\newif \ifELS

\documentclass[10pt,a4paper]{article}
\usepackage[a4paper,left=2.5cm,right=2.5cm,top=2.5cm,bottom=2.5cm]{geometry} 

\usepackage{graphicx}
\usepackage{amssymb}
\usepackage{amsmath}
\usepackage{amsthm}
\usepackage{bm}
\usepackage{mathtools}
\usepackage{booktabs}
\usepackage{paralist}
\usepackage[usenames]{color}
\usepackage{empheq}
\usepackage{colortbl}
\usepackage{tcolorbox}
\tcbuselibrary{breakable}
\usepackage{enumitem}
\usepackage{pgfplots}
\usepackage{tikz}
\usepackage{subfig}
\usepackage{algorithm}
\usepackage{algorithmic}
\usepackage[normalem]{ulem} 
\usepackage[textsize=footnotesize,color=blue!30!white]{todonotes} 
\usepackage{cancel}
\usepackage[breaklinks,bookmarks=false]{hyperref}

\hypersetup{colorlinks, linkcolor=blue, citecolor=blue,
urlcolor=blue, plainpages=false, pdfwindowui=false, pdfstartview={FitH}, pdftitle={}, pdfauthor={Martin Vohral\'ik} }
\pdfcatalog{/PageLayout /OneColumn}

\makeatletter
\newcommand{\MS[1]}{\bBigg@{#1}}
\makeatother

\numberwithin{equation}{section}

\newcommand*\widefbox[1]{\fbox{\hspace{0.2em}#1\hspace{0.2em}}}

\ifELS
\newtheorem{remark}[theorem]{Remark}
\newtheorem{assumption}[theorem]{Assumption}
\newtheorem{construction}[theorem]{Construction}

\else

\newtheorem{lemma}{Lemma}[section]
\newtheorem{corollary}[lemma]{Corollary}
\newtheorem{assumption}[lemma]{Assumption}
\newtheorem{theorem}[lemma]{Theorem}
\newtheorem{proposition}[lemma]{Proposition}
\newtheorem{definition}[lemma]{Definition}
\newtheorem{remark}[lemma]{Remark}
\newtheorem{example}[lemma]{Example}
\newtheorem{motivation}[lemma]{Motivation}
\newtheorem{construction}[lemma]{Construction}
\newtheorem{algo}[lemma]{Algorithm}

\fi

\newtcolorbox[number within=section, use counter=lemma]{ctheorem}[3][]{
left=0mm, right=0mm, colback=black!5, colframe=black!30, 
colbacktitle=black!30, coltitle=black, label=#3, breakable,
title={{\bf Theorem~\thelemma}~(#2){\bf .}},#1}

\newcommand{\be}{\begin{equation}}
\newcommand{\ee}{\end{equation}}
\newcommand{\bea}{\begin{eqnarray}}
\newcommand{\eea}{\end{eqnarray}}
\newcommand{\bean}{\begin{eqnarray*}}
\newcommand{\eean}{\end{eqnarray*}}
\def\ba#1\ea{\begin{align}#1\end{align}}
\def\ban#1\ean{\begin{align*}#1\end{align*}}
\def\bat#1\eat{\begin{alignat}#1\end{alignat}}
\def\batn#1\eatn{\begin{alignat*}#1\end{alignat*}}
\def\bs#1\es{\begin{split}#1\end{split}}
\newcommand{\bse}{\begin{subequations}}
\newcommand{\ese}{\end{subequations}}
\newcommand{\bt}{\begin{theorem}}
\newcommand{\et}{\end{theorem}}
\newcommand{\bpr}{\begin{proposition}}
\newcommand{\epr}{\end{proposition}}
\newcommand{\bl}{\begin{lemma}}
\newcommand{\el}{\end{lemma}}
\newcommand{\bc}{\begin{corollary}}
\newcommand{\ec}{\end{corollary}}
\newcommand{\bp}{\begin{proof}}
\newcommand{\ep}{\end{proof}}
\newcommand{\bd}{\begin{definition}}
\newcommand{\ed}{\end{definition}}
\newcommand{\br}{\begin{remark}}
\newcommand{\er}{\end{remark}}
\newcommand{\bmo}{\begin{motivation}}
\newcommand{\emo}{\end{motivation}}
\newcommand{\bco}{\begin{construction}}
\newcommand{\eco}{\end{construction}}
\newcommand{\bas}{\begin{assumption}}
\newcommand{\eas}{\end{assumption}}
\newcommand{\bex}{\begin{example}}
\newcommand{\eex}{\end{example}}
\newcommand{\bqo}{\begin{quote}}
\newcommand{\eqo}{\end{quote}}
\newcommand{\bdc}{\begin{description}}
\newcommand{\edc}{\end{description}}
\newcommand{\bi}{\begin{itemize}}
\newcommand{\ei}{\end{itemize}}
\newcommand{\ben}{\begin{enumerate}}
\newcommand{\een}{\end{enumerate}}

\newcommand\ie{i.e.}
\newcommand\cf{cf.}
\newcommand\eg{e.g.}
\newcommand\eal{{\em et al.}}
\newcommand\eq{:=}
\newcommand\qe{=:}

\newcommand\ds{\displaystyle}
\newcommand\nn{\nonumber}
\newcommand\pt{\partial}
\newcommand{\<}{\langle}
\renewcommand{\>}{\rangle}
\newcommand\ra{\rightarrow}

\newcommand\reff[2]{\stackrel{\eqref{#1}}{#2}}

\newcommand\Om{\Omega}
\newcommand\om{\omega}
\newcommand\omK{{\omega_\elm}}

\newcommand{\elm}{{K}}
\newcommand{\elmt}{{L}}
\newcommand{\sd}{{\sigma}}
\newcommand{\sdt}{{\sigma'}}
\newcommand{\ver}{{\ta}}
\newcommand{\vertt}{{\ta'}}
\newcommand{\htt}{{h\tau}}

\newcommand\Fh{\mathcal{F}_h}
\newcommand\Fhint{\mathcal{F}_h^\mathrm{int}}
\newcommand\Fhext{\mathcal{F}_h^\mathrm{ext}}
\newcommand\FK{\mathcal{F}_{\elm}}
\newcommand\FKh{\mathcal{F}_{\elm,h}}
\newcommand\FKt{\mathcal{F}_\elmt}
\newcommand\FKhext{\FKh^{\rm ext}}
\newcommand\FKhint{\FKh^{\rm int}}
\newcommand\FKhextn{\FKh^{n,{\rm ext}}}
\newcommand\FKhintn{\FKh^{n,{\rm int}}}

\newcommand\Th{\mathcal{T}_h}
\newcommand\TK{\mathcal{T}_\elm}

\newcommand\Vh{\mathcal{V}_h}
\newcommand\VK{\mathcal{V}_\elm}
\newcommand\Ve{\mathcal{V}_\sd}
\newcommand\VKh{\mathcal{V}_{\elm,h}}

\newcommand\Set{\mathcal{S}}

\newcommand\Gr{\nabla}

\newcommand\Dv{\nabla {\cdot}}

\newcommand\dv{\mathrm{div}}
\newcommand\scp{{\cdot}}
\newcommand\trp{\mathrm{t}}
\newcommand\Lap{\Delta}

\newcommand{\matr}[1]{\mathbb{#1}}
\newcommand\il{|\hspace{-0.02cm}|\hspace{-0.02cm}|}

\DeclarePairedDelimiter\norm{\|}{\|} 
\DeclarePairedDelimiter\enorm{\il}{\il} 

\newcommand\Ho{H^1(\Om)}
\newcommand\Hoi[1]{H^1(#1)}
\newcommand\Hoo{H^1_0(\Om)}

\newcommand\Hmo{H^{-1}(\Om)}
\newcommand\HTh{H^1(\Th)}
\newcommand\HdTh{\bm{H}(\dv,\Th)}

\newcommand\Cz{C^0(\overline\Om)}
\newcommand\Co{C^1(\overline\Om)}
\newcommand\Cinf{C^\infty(\Om)}

\newcommand\Lt{L^2(\Om)}
\newcommand\tLt{\bm{L}^2(\Om)}
\newcommand\Lti[1]{L^2(#1)}
\newcommand\tLti[1]{\bm{L}^2(#1)}

\newcommand\Do{{\mathcal D}(\Om)}

\newcommand\Hdv{\bm{H}(\dv,\Om)}
\newcommand\Hdvi[1]{\bm{H}(\dv,#1)}

\newcommand\RT{\bm{\mathcal{R\hspace{-0.1em}T}}\hspace{-0.25em}}

\newcommand\pth{\tilde p_h}

\newcommand\pr{\zeta}
\newcommand\prh{\zeta_h}
\newcommand\prhki{\zeta_h^{k,i}}
\newcommand\prht{\zeta_{p,\htt}}

\newcommand\tuh{\tu_h}

\newcommand\tuhki{\tu_h^{k,i}}
\newcommand\tuhtki{\tu_{h\tau}^{k,i}}

\newcommand\frh{\bsig_h} 

\newcommand\vf{\varphi}

\newcommand\ta{\bm{a}}

\newcommand\tn{\bm{n}}

\newcommand\tu{\bm{u}}
\newcommand\tv{\bm{v}}
\newcommand\tw{\bm{w}}
\newcommand\tx{\bm{x}}

\newcommand\Km{{\bm{\underline K}}}
\newcommand\Kmt{{\tilde {\bm{\underline K}}}}

\newcommand\tV{\bm{V}}

\newcommand{\bsig}{\bm{\sigma}}

\newcommand{\bxi}{{\bm\xi}}

\newcommand\A{\mathcal{A}}

\newcommand\C{\mathcal{C}}

\newcommand\F{\mathcal{F}}

\newcommand\Pp{\mathcal{P}}

\newcommand\T{\mathcal{T}}

\newcommand\Idd{{\bm{\underline I}}}

\newcommand\RR{{\mathbb R}}

\newcommand\PP{{\mathcal{P}}}
\newcommand\QQ{{\mathcal{Q}}}
\def\sol{{\bm{X}}}
\def\solM{{\bm{X}}^n_{\elm}}
\def\solMp{{\bm{X}}^{n-1}_{\elm}}
\def\solMpp{{\bm{X}
}^{n,k-1}_{\elm}}
\def\solMppp{{\bm{X}}^{n,k,i}_{\elm}}
\def\solMK{{\bm{X}}^n_{\elm'}}
\newcommand{\solcM}[1][n]{{\bm{X}}_{\T_H}^{#1}}
\newcommand{\ffd}{\VEC{\theta}_{\mathrm{\cor{upw}},c,h}}
\newcommand{\ffl}{\VEC{\theta}_{\mathrm{lin},c,h}}
\newcommand{\ffr}{\VEC{\theta}_{\mathrm{alg},c,h}}

\def\d{\mathrm{d}}

\newcommand{\perm}{\Km}

\newcommand{\VEC}[1]{\boldsymbol{#1}}
\newcommand{\DIV}{{\nabla{\cdot}}}
\newcommand{\GRAD}{{\nabla}}
\newcommand{\SCAL}{{\cdot}}
\newcommand{\Nnn}[1][]{\mathcal{N}^{#1}}

\newcommand{\Phicht}[1][]{\VEC{\theta}_{c,h\tau}^{#1}}
\newcommand{\Phich}[1][]{\VEC{\theta}_{c,h}^{#1}}

\newcommand{\normall}[1][]{\tn_{#1}}

\newcommand{\Nn}[1]{\mathcal{N}^{#1}}

\newcommand{\Ll}{{l}}
\newcommand{\Fl}{\theta}
\newcommand{\Rl}{{R}}

\newcommand{\FpMs}{v_{p,\elm,\sigma}}
\newcommand{\FcMs}{\Fl_{c,\elm,\sigma,\uparrow}}
\newcommand{\FcMsNL}{\Fl_{c,\elm,\sigma}}

\newcommand\ft{{1/2}}
\newcommand\mft{{-1/2}}

\newcommand\dx{\, \mathrm{d} \tx}
\newcommand\dt{\, \mathrm{d} t}

\newcommand{\err}[3]{\ensuremath{{\rm{err}}^{#3}_{{\rm #1}}}}
\newcommand{\est}[3]{\ensuremath{\eta_{{\rm #1},#2}^{#3}}}
\newcommand{\Est}[2]{\ensuremath{\eta_{{\rm #1}}^{#2}}}
\newcommand{\Esti}[2]{\ensuremath{\eta_{{\rm #1}}^{#2}}}

\newcommand{\Param}[1]{\Gamma_{\rm #1}}
\newcommand{\param}[1]{\gamma_{\rm #1}}

\newcommand{\alg}[1]{\mathsf{#1}}
\newcommand{\algUi}[1]{\mathsf{U}_{{\rm #1},\elm}^{ {k,i}}}
\newcommand{\algP}[1]{\mathsf{\Theta}_{{\rm #1},\elm,c}^{n,k,i}}
\newcommand{\algPt}[1]{\mathsf{\Theta}_{{\rm #1},\elm,c}^{t,n,k,i}}
\newcommand{\Usd}{\alg{U}_\sd}
\newcommand{\Usdt}{\alg{U}_\sdt}
\newcommand{\Peo}{\alg{P}^0_\elm}
\newcommand{\Pelm}{\alg{P}_\elm}
\newcommand{\Pelmt}{\alg{P}_\elmt}
\newcommand{\Sver}{\alg{Z}_{\ver}}
\newcommand{\Ssd}{\alg{Z}_{\sd}}

\newcommand{\SeExt}{\alg{Z}^{\mathrm{ext}}_{\elm}}

\newcommand{\Selm}{\alg{Z}_\elm}
\newcommand{\Sel}{\alg{Z}}

\newcommand\bzeta{\bm{\zeta}}

\newcommand{\cor}[1]{{#1}}

\begin{document}

\title{A posteriori error estimates and adaptivity for locally conservative methods\\ {\large Inexpensive implementation and evaluation, polytopal meshes, iterative linearization and algebraic solvers, and applications to complex porous media flows}}

\author{Martin Vohral\'ik\footnotemark[2] \and Soleiman Yousef\footnotemark[3]}

\maketitle

\renewcommand{\thefootnote}{\fnsymbol{footnote}}

\footnotetext[2]{Project-team SERENA, Inria Paris, 48 rue Barrault, 75647 Paris, France \&
CERMICS, Ecole nationale des ponts et chauss\'ees, IP Paris, 77455 Marne la Vall\'ee, France (\href{mailto:martin.vohralik@inria.fr}{\texttt{martin.vohralik@inria.fr}}).}

\footnotetext[3]{IFP Energies nouvelles, 1 \& 4 av. Bois Pr\'eau, 92852 Rueil-Malmaison, France
(\href{soleiman.yousef@ifpen.fr}{\texttt{soleiman.yousef@ifpen.fr}}).}

\renewcommand{\thefootnote}{\arabic{footnote}}

\begin{abstract} A posteriori estimates give bounds on the error between the unknown solution of a partial differential equation and its numerical approximation. We present here the methodology based on $H^1$-conforming potential and $\bm{H}(\dv)$-conforming equilibrated flux reconstructions, where the error bounds are guaranteed and fully computable. 
We consider any lowest-order locally conservative method of the finite volume type and treat general polytopal meshes. We start by a pure diffusion problem and first address the discretization error. 
We then progressively pass to more complicated model problems and also take into account the errors arising in iterative linearization of nonlinear problems and in algebraic resolution of systems of linear algebraic equations. 
We focus on the ease of implementation and evaluation of the estimates. In particular, the evaluation of our estimates is explicit and inexpensive, since it merely consists in some local matrix-vector multiplications. 
Here, on each mesh element, the matrices are either directly inherited from the given numerical method, or easily constructed from the element geometry, while the vectors are the \cor{algebraic unknowns of the} flux and potential approximations on the given element. 
Our last problem is a real-life unsteady nonlinear coupled degenerate advection--diffusion--reaction system describing a complex multiphase multicomponent flow in porous media. 
Here, on each step of the time-marching scheme, on each step of the iterative linearization procedure, and on each step of the linear algebraic solver, the estimate gives a guaranteed upper bound on the total intrinsic error, still takes the simple matrix-vector multiplication form, and distinguishes the different error components. 
It leads to an easy-to-implement and fast-to-run adaptive algorithm with guaranteed overall precision, adaptive stopping criteria for nonlinear and linear solvers, and adaptive space and time mesh refinements and derefinements. 
Progressively along the theoretical exposition, numerical experiments on academic benchmarks as well as on real-life problems in two and three space dimensions illustrate the performance of the derived methodology. The presentation is largely self-standing, developing all the details and recalling all necessary basic notions. 

\end{abstract}



\tableofcontents

\section{Introduction} \label{sec_intr}

The purpose of this contribution is to present the theory of {\em a posteriori
error estimates} and {\em adaptivity} in numerical approximation of partial differential equations (PDEs) by locally conservative methods of the {\em finite volume} type. We focus on estimates that {\em certify} the {\em error} and can be {\em easily coded}, {\em cheaply evaluated}, and {\em efficiently used}, even in complex {\em practical simulations}. We consider {\em polytopal meshes}, \ie, meshes formed by general polygonal or polyhedral elements, appealing in various applications. We cover the overall chain of computational practice, including {\em iterative linearization} of nonlinear problems and {\em algebraic resolution} of systems of linear algebraic equations. 

An {\em optimal} a posteriori estimate can be described by a set of requested properties, including namely the bound on the error between the unknown solution of the partial differential equation and the available numerical approximation that is {\em guaranteed} and {\em fully computable} from the approximate solution. 
We introduce this concept in Section~\ref{sec_a_post} at an abstract level. 
We start from the notion of a PDE, discuss some intrinsic physical and mathematical properties of the weak solutions of PDEs, the notion of a total error and its components, and of adaptivity which aims at balancing these error components. One particular property of the methodology we develop here is that it leads to (exact or approximate) mass balance recovery at any step of numerical resolution, namely during the iteration of the nonlinear and linear solvers.

This contribution is made as self-standing as possible. For this purpose, in Section~\ref{sec_not}, we set up the notation, recall the intrinsic function spaces related to the considered partial differential equations, namely the infinite-dimensional {\em Sobolev spaces} $\Hoo$ and $\Hdv$, and recall their finite-dimensional subspaces composed of {\em piecewise polynomials} with respect to a simplicial mesh. 
We also describe there the polytopal meshes we consider and introduce the basic principle of the finite volume and related lowest-order locally conservative methods.

The heart of our exposition starts in Section~\ref{sec_Pois}, where we treat in detail the model steady linear diffusion problem with homogeneous and isotropic diffusion tensor (the {\em Poisson equation}). Only the basic cell-centered finite volume discretization on a simplicial mesh and only the {\em spatial discretization} error are investigated here, to expose the construction principles as clearly as possible. 
We namely motivate and present the methodology based on $\Ho$-conforming potential and $\Hdv$-conforming equilibrated flux reconstructions, which yields \cor{a} guaranteed and fully computable error upper bound. 
These reconstructions are practically obtained here in piecewise polynomial spaces. 
This is perfectly fine for academic purposes and has an asymptotically optimal evaluation cost, since it is in particular explicit and no local problems need to be solved. \cor{A} cheaper implementation and evaluation of a posteriori error estimates is, however, possible and suitable for more complicated model or real-life problems, which we treat later. 
Numerical experiments illustrate these basic developments. 

In Section~\ref{sec_Darcy}, we extend the above analysis to tensor-valued steady linear diffusion problems (the {\em singlephase steady linear Darcy flow}), general polytopal meshes, and any lowest-order locally conservative method such as mixed finite elements, mixed and hybrid finite volumes, mimetic finite differences, mixed virtual finite elements or hybrid high-order methods. In particular, though we still rely on the notions of potential and equilibrated flux reconstructions, there become {\em virtual} here, since they are not constructed in practice anymore. 
Similarly, though we rely on a notion of a simplicial submesh of a general polytopal mesh, this is also virtual, where the physical construction and computer implementation of a simplicial submesh is avoided. 
In particular, the evaluation of our estimates merely consists in some local matrix-vector multiplications, where, on each mesh element, the matrices are either directly inherited from the given numerical method, or easily constructed from the element geometry, while the vectors are the \cor{algebraic unknowns of the} flux and potential \cor{approximations} on the given element. This gives an easy and practically accessible application to polytopal meshes of the general methodology of $H^1$-conforming potential reconstruction and $\Hdv$-conforming flux reconstruction. Numerical experiments focus on general polygonal meshes \cor{and} illustrate the performance of this methodology.

In Section~\ref{sec_Darcy_NL}, we extend the previous developments to encompass two cornerstones of numerical approximation in applications: the {\em nonlinear} and {\em linear solvers}. Our problem here is the {\em singlephase steady nonlinear Darcy flow}, where an {\em iterative linearization} is applied and where the {\em algebraic resolution} of the system of arising linear algebraic equations is addressed. 
We still derive a guaranteed upper bound on the total error. 
Moreover, we identify the three arising {\em error components}, related to discretization, linearization, and algebraic solver, and design {\em adaptive stopping criteria} for the iterative solvers. Remarkably, the estimators still take the form of a {\em simple matrix-vector multiplication}, with the same local matrices as in the linear case. Numerical illustrations confirm the theoretical findings.

In the last part of this contribution, Section~\ref{sec_MP_MC}, we consider an unsteady nonlinear coupled degenerate advection--diffusion--reaction problem, the {\em multiphase compositional unsteady Darcy flow}. 
We apply the entire methodology developed above to obtain estimates that are still guaranteed and fully computable for the intrinsic error measure, apply to any lowest-order locally conservative method on a polytopal mesh, and take a simple matrix-vector multiplication form. 
Moreover, they are valid on each stage of the overall solution algorithm: on {\em each time} step $n$, each {\em linearization} step $k$, and each {\em linear solver step} $i$. 
They also allow to
distinguish different error components and design {\em adaptive stopping criteria} for the solvers, as well as {\em adaptive choice} of {\em space} and {\em time meshes}. A comprehensive numerical study on realistic porous media problems is included.

The present a posteriori error estimates can be readily implemented into numerical codes on general polytopal meshes, with a minimal overhead. These estimates allow for a very fast evaluation and, according to the main results summarized in Theorems~\ref{thm_est_Pois}, \ref{thm_est_Darcy}, \ref{thm_estim_Darcy_NL}, and~\ref{thm_estim_MP_MC}, give a guaranteed control over the intrinsic error committed in the numerical approximation. 
Additionally, all the different error components
(time and space discretizations, linearization, algebraic) are identified, leading to fully adaptive algorithms with all adaptive stopping criteria for linear and nonlinear solvers, adaptive time step management, and adaptive mesh refinement. Numerical experiments on real-life problems confirm important {\em computational speed-ups} that can be achieved with our methodology, in addition to the {\em certification of the computed output}.

Our presentation proceeds along the lines of Vohral{\'{\i}}k~\cite{Voh_apost_FV_08}, Ern and Vohral{\'{\i}}k~\cite{Ern_Voh_adpt_IN_13}, and Vohral{\'{\i}}k and Yousef~\cite{Voh_Yous_polyt_18}, using some central concepts for general polytopal meshes from Vohral{\'{\i}}k and Wohlmuth~\cite{Voh_Wohl_MFE_1_unkn_el_rel_13}. Ample bibliographic references age given in the dedicated sections below.

\section{A posteriori error estimation and adaptivity -- the principles} \label{sec_a_post}

In this section, we set up the main ideas. We proceed informally; notation is in detail set in Section~\ref{sec_not}.

\subsection{Partial differential equations}

Let $\Om \subset \RR^d$, $d \geq 1$, be a polytopal domain where our problem is set. 
Consider an abstract (steady \cor{but possibly non}linear) partial differential equation in the form: for a source term $f$, find a real-valued function $p: \Om \ra \RR$ such that
\bse \label{eq_PDE} \ba
    \A(p) & = f \qquad \mbox{ in } \, \Om, \label{eq_PDE_eq} \\
    p & = 0 \qquad \mbox{ on } \, \pt \Om. \label{eq_PDE_BC}
\ea \ese
The prominent example that we will consider below is the Laplace operator 
\be \label{eq_Lapl_op}
    \A(p) = - \Dv(\Gr p) = -\Lap p = - \sum_{i=1}^d \pt_{\tx_i}^2 p.
\ee
In general, $\A(\cdot)$ is a (\cor{non}linear) operator including {\em partial derivatives} of the unknown solution $p$. 

More generally, for a final time $t_{\rm{F}} > 0$, a source term $f$, and an initial condition $p_0$, consider an abstract (unsteady nonlinear) partial differential equation in the form: find a real-valued function $p: \Om \times (0,t_{\rm{F}}) \ra \RR$ such that
\bse \label{eq_PDE_unst} \ba
    \A(p) & = f \qquad \mbox{ in } \, \Om \times (0,t_{\rm{F}}), \label{eq_PDE_eq_unst} \\
    p & = 0 \qquad \mbox{ on } \, \pt \Om \times (0,t_{\rm{F}}), \label{eq_PDE_BC_unst} \\
    p(\cdot,0) & = p_0 \qquad \mbox{ in } \, \Om. \label{eq_PDE_IC_unst}
\ea \ese
Here, $\A(\cdot)$ is a (nonlinear) operator including partial derivatives of $p$ and namely a partial derivative with respect to the time, for example 
\be \label{eq_nonl_heat_op}
    \A(p) = \pt_t p - \Dv(\Km(|\Gr p|) \Gr p),
\ee
where $|\Gr p|$ is the Euclidean size of the vector $\Gr p$ and $\Km(\cdot)$ is a scalar- or matrix-valued {\em nonlinear function} of a real variable.

\subsection{Physical and mathematical properties of the weak solutions of PDEs} \label{sec_props}

In order to properly describe the solution $p$ of~\eqref{eq_PDE} or~\eqref{eq_PDE_unst}, adequate mathematical {\em function spaces} are necessary. We develop this task in details in Section~\ref{sec_spaces} below, introducing namely the Sobolev spaces $\Hoo$ and $\Hdv$. These spaces also reflect the {\em physical principles} behind problems~\eqref{eq_PDE} and~\eqref{eq_PDE_unst}. For instance, for the Laplace operator~\eqref{eq_Lapl_op} in~\eqref{eq_PDE}, it turns out that the solution $p$ satisfies
\be \label{eq_Lapl_prop}
    p \in \Hoo, \quad \tu \eq - \Gr p \in \Hdv, \quad \Dv \tu = f.
\ee
This reflects that the scalar-valued field $p$ (pressure, potential, primal variable) is {\em continuous} in the appropriate trace sense (cannot jump across interfaces), the normal trace of the vector-valued field $\tu \eq - \Gr p$ (Darcy velocity, flux or flow field, dual variable) is {\em normal-component continuous} in the appropriate trace sense (what flows out from a subdomain through an interface flows in through this interface in the neighboring subdomain), the Darcy velocity is related to the pressure by the {\em constitutive law} $\tu \eq - \Gr p$, and the Darcy velocity is in {\em equilibrium} with the load $f$ as \cor{expressed by} $\Dv \tu = f$.

\subsection{Numerical approximation: spatial discretization, temporal discretization, iterative linearization, and iterative linear algebraic resolution}

It is in general not possible to find analytically the exact potential $p$ and \cor{the} exact flux $\tu$ of problems like~\eqref{eq_PDE} or~\eqref{eq_PDE_unst}. Then, we want to find an approximation obtained by a mathematical algorithm implemented on a computer. 
Typically, we will take a subdivision of $\Om$ into computational cells, a spatial mesh $\Th$ of $\Om$, with a characteristic size $h$, and use some {\em spatial discretization} scheme. For~\eqref{eq_PDE_unst}, we additionally consider a subdivision of the time interval $(0,t_{\rm{F}})$ into time steps, a temporal mesh of $(0,t_{\rm{F}})$, with a characteristic size $\tau$, and proceed to a {\em temporal discretization}.
When a nonlinear function as $\Km(\cdot)$ in~\eqref{eq_nonl_heat_op} is present, we typically need an {\em iterative linearization}, to, on each iteration step $k$, approximate a nonlinear operator by a linear one. 
Finally, for steady problems or implicit discretizations of unsteady problems, we need an {\em iterative linear algebraic solver} for an approximate solution of a (large sparse) system of linear algebraic equations, with an iteration index $i$. In this text, we focus on discretizations yielding {\em numerical approximations} to the exact flux $\tu$ \cor{(rather than to the exact potential $p$)} that we denote by $\tuhki$ or $\tuhtki$.

\subsection{Total error and its components} \label{sec_err_comps_all}

In principle, the numerical approximation is not exact, so that $\tuhki \neq \tu$ for~\eqref{eq_PDE} and $\tuhtki \neq \tu$ for~\eqref{eq_PDE_unst}. Evaluating the difference $\tu - \tuhki$ or $\tu-\tuhtki$ in some norm then gives the notion of the {\em total error}. This is in our exposition composed of the spatial discretization error, temporal discretization error, iterative linearization error, and iterative linear algebraic solver error.\footnote{In this text, we suppose exact computer implementation, \ie, we neglect rounding errors coming from approximate computer work with real numbers via floating point arithmetic. \cor{Such errors can also be taken into account, as well as additional errors such as those from (iterative) regularization or model selection.}}

\subsubsection{Intrinsic error measure} \label{sec_intrin}

In the choice of the norm to measure the distances $\tu - \tuhki$ or $\tu-\tuhtki$, we in this contribution stick to problem-induced norms, written as
\be \label{eq_intr_norm}
    \enorm{\tu - \tuhki} \quad \text { or } \quad \enorm{\tu-\tuhtki}.
\ee
These {\em intrinsic error measures} arise from the problem considered, the mathematical function spaces adopted, and are defined so as to yield the most direct link to the sum of the {\em residual} and {\em nonconformity} in $\tuhki$ or $\tuhtki$, \cf\ Theorem~\ref{thm_err_char} below. In the Laplace case~\eqref{eq_Lapl_op}, our choice of the intrinsic error measure in particular gives
\be \label{eq_norm_Lapl}
    \enorm{\tu-\tuhki} \cor{=} \norm{\tu-\tu_h^{k,i}},
\ee
where $\norm{\cdot}$ is the $\tLt$ norm on $\Om$. Th\cor{is $\Lt$} flux error $\norm{\tu-\tu_h^{k,i}}$ is a physically crucial\cor{ly important} quantity in the applications. 

\subsubsection{Spatial discretization error} \label{sec_sp_err}

The {\em spatial discretization error} is the error arising from the discretization on the spatial mesh $\Th$ of $\Om$. 
For steady problems~\eqref{eq_PDE} and supposing that we have obtained the approximate \cor{flux} $\tu_h$ as specified ``on paper'', \ie, in particular solving ``exactly'' the associated \cor{(non)}linear systems (so that there is no \cor{linearization index $k$ and no} algebraic solver index $i$), this is the only error component.
In the Laplace case~\eqref{eq_Lapl_op}, the spatial discretization error writes 
\be \label{eq_norm_Lapl_sp}
    \norm{\tu-\tu_h}.
\ee
In the setting of locally conservative methods that we consider, the \cor{approximate} flux $\tu_h$ will satisfy
\be \label{eq_loc_cons}
    \tu_h \in \Hdv \quad \text{with} \quad \Dv \tu_h = f
\ee
for a piecewise constant source term $f$. Property~\eqref{eq_loc_cons} is precisely the ``local conservation'' from ``locally conservative'' methods. \cor{Additionally, we will be able to find an approximate potential $\pth$ such that 
\[
- \Gr \pth|_\elm = \tu_h|_\elm
\]
for all elements $\elm$ of the mesh $\Th$, so that the constitutive law is satisfied.} In comparison with the continuous level~\eqref{eq_Lapl_prop}, we see that \cor{the only} deficiency of the discrete level is that the \cor{approximate potential $\pth$} is {\em nonconforming}, \cor{does not belong to the correct space},
\be \label{eq_pth_disc}
    \pth \not \in \Hoo,
\ee
in contrast to the situation at the continuous level, where \cor{$p \in \Hoo$} as per~\eqref{eq_Lapl_prop}.

\subsubsection{Iterative linear algebraic solver error} \label{sec_alg_err}

In order to obtain the approximate solution $\tu_h$ as specified by some spatial discretization scheme, one typically needs to solve a (large, sparse) system of linear algebraic equations. In practice, we typically proceed iteratively, yielding an approximation $\tu_h^i$ to $\tu_h$ on each algebraic solver step $i \geq 0$. Then the {\em iterative linear algebraic solver error} is
\be \label{eq_norm_Lapl_alg}
    \norm{\tu_h-\tu_h^i}.
\ee

\subsubsection{Iterative linearization error} \label{sec_lin_err}

\cor{In presence of nonlinearities, one typically employs an iterative approximation} of a nonlinear operator $\A(\cdot)$ by linear operators $\A^k(\cdot)$ \cor{(on the continuous or on the discrete level). For steady problems, this yields an approximation $\tu_h^k$ to $\tu_h$ on each iterative linearization step $k \geq 0$.} The corresponding {\em iterative linearization error} can \cor{then} be structurally expressed as
\be \label{eq_norm_nonl_Lapl_lin}
    \norm{\tu_h-\tu_h^k}.
\ee

\subsubsection{Temporal discretization error} \label{sec_temp_err}

The {\em temporal discretization error} is the error arising from the discretization on the temporal mesh of the time interval $(0,t_{\rm{F}})$. Naturally, it only appears for unsteady problems of the form~\eqref{eq_PDE_unst}, when some temporal discretization has been applied.

\subsubsection{\cor{Error components}} \label{sec_err_comps}

\cor{Wrapping up the previous sections, the total error is typically composed of several error components. For problems of the form~\eqref{eq_PDE_unst}--\eqref{eq_nonl_heat_op} namely,} we still have the potential and the Darcy velocity
\be  \label{eq_Lapl_prop_unst}
    \cor{p \in \Hoo,} \quad \tu \eq - \Km(|\Gr p|) \Gr p \in \Hdv, \quad \cor{\pt_t p +} \Dv \tu = f\cor{,}
\ee
as in~\eqref{eq_Lapl_prop}. Then, for a numerical approximation \cor{$\tuhtki$}, we will have 
\bse \label{eq_err_comps} \begin{alignat}{2}
& \text{total error} && \enorm{\tu-\tu_\htt^{k,i}},\\
& \text{\cor{temporal discretization error}} && \enorm{\tu-\tu_\tau},\\ 
& \text{spatial discretization error} && \enorm{\tu_\tau-\tu_\htt},\\ 
& \text{iterative linearization error} && \enorm{\tu_\htt-\tu_\htt^k},\\ 
& \text{iterative linear algebraic solver error} \qquad && \enorm{\tu_\htt^k-\tu_\htt^{k,i}},
\end{alignat}\ese
in extension of~\eqref{eq_norm_Lapl_sp}, \eqref{eq_norm_nonl_Lapl_lin}, and~\eqref{eq_norm_Lapl_alg}. Details are developed in Sections~\ref{sec_Darcy_NL} and~\ref{sec_MP_MC}.

\subsection{A posteriori error estimate} \label{sec_a_post_props}

An a posteriori error estimate is \cor{real number} $\eta(\tuhki)$ or $\eta(\tuhtki)$ \cor{that bounds} the error between the unknown \cor{exact} flux $\tu$ of the partial differential equation like~\eqref{eq_PDE} or~\eqref{eq_PDE_unst} and the available numerical approximation $\tuhki$ or $\tuhtki$ \cor{of $\tu$. Importantly,} $\eta(\tuhki)$ or $\eta(\tuhtki)$ \cor{has to be} {\em fully computable} from the approximate solution $\tuhki$ or $\tuhtki$. The following are some highly desirable properties.

\subsubsection{Guaranteed error upper bound (reliability)} \label{sec_rel}

We say that an a posteriori error estimate gives a {\em guaranteed error upper bound} when the intrinsic error measure~\eqref{eq_intr_norm} is bounded by it as
\be \label{eq_est_rel}
    \enorm{\tu-\tuhki} \leq \eta(\tuhki) \quad \text { or } \quad \enorm{\tu-\tuhtki} \leq \eta(\tuhtki).
\ee
This allows to {\em certify} the {\em error}.
This may appear as a miracle, since the left-hand side of~\eqref{eq_est_rel} is unknown, whereas the right-hand side of~\eqref{eq_est_rel} is known, but is indeed possible, as we will see in the developments below. 
The key will be to follow at the discrete level the physical and mathematical properties of the weak solutions of PDEs as \cor{expressed by~\eqref{eq_Lapl_prop} and~\eqref{eq_Lapl_prop_unst}}. In the setting of Section~\ref{sec_sp_err}, namely, from~\eqref{eq_Lapl_prop} and since $\tu_h \in \Hdv$ with $\Dv \tu_h = f$, we easily show that
\be \label{eq_a_post_Lapl_sp}
    \norm{\tu-\tu_h} \leq \eta(\tuh) \eq \norm{\tu_h + \Gr \prh}
\ee
for an arbitrary
\be \label{eq_pot_rec}
    \prh \in \Hoo;
\ee
we will see the details in Theorem~\ref{thm_est_Pois} below. Here, $\prh$ is discrete, piecewise polynomial, and obtained by some {\em local} modification from \cor{$\pth$} (postprocessing and averaging); we call it a {\em potential reconstruction}. 
Similarly, when $\tuhki \not \in \Hdv$ or $\Dv \tuhki \neq f$ (namely in presence of linear and nonlinear solvers), an {\em equilibrated flux reconstruction} $\bsig_h^{k,i}$ such that 
\be \label{eq_eq_fl_rec}
    \bsig_h^{k,i} \in \Hdv \quad \text{and} \quad \Dv \bsig_h^{k,i} = f
\ee
(we will also consider merely $\Dv \bsig_h^{k,i} = f + \rho_h^{k,i}$ where $\rho_h^{k,i}$ is a remainder made as small as necessary) will be the key together with the potential reconstruction~\eqref{eq_pot_rec}.

\subsubsection{Error lower bound (efficiency)} \label{sec_eff}

Relation~\eqref{eq_est_rel} is only one-sided, so that it might potentially happen that the left-hand side of~\eqref{eq_est_rel} is zero (small), whereas the right-hand side is nonzero (big). A mathematical equivalence between the error and the a posteriori error estimator is ensured when, in addition to~\eqref{eq_est_rel}, there also holds 
\be \label{eq_est_eff}
    \eta(\tuhki) \leq C \enorm{\tu-\tuhki} \quad \text { or } \quad \eta(\tuhtki) \leq C \enorm{\tu-\tuhtki},
\ee
which is called (global) {\em efficiency}. Here, $C$ is a generic constant. Actually, \eqref{eq_est_eff} may only hold when on the right-hand side, there are some additional terms without structural importance, \ie, of vanishing importance with decreasing the mesh size $h$ and the time step $\tau$ (so-called ``data oscillation'' terms). 

Strengthening~\eqref{eq_est_eff} to hold locally in space (and in time) is actually often possible, and highly desirable. This then means that when we {\em predict} the {\em presence of error} in some {\em mesh element} or {\em time step}, then it is {\em indeed localized there}, or in a small neighborhood. 
In the example~\eqref{eq_a_post_Lapl_sp}, in particular, we easily have
\be \label{eq_est_loc}
    \underbrace{\norm{\tu_h + \Gr \prh}^2}_{[\eta(\tuh)]^2} = \sum_{\elm \in \Th} \underbrace{\norm{\tu_h + \Gr \prh}_\elm^2}_{[\eta_\elm(\tuh)]^2},
\ee
and it is possible to strengthen~\eqref{eq_est_eff} to
\be \label{eq_est_eff_loc}
    \eta_\elm(\tuh) \leq C \norm{\tu-\tuh}_\omK,
\ee
where $\omK$ are the simplices of the mesh $\Th$ sharing a vertex with the simplex $\elm$. This is called {\em local efficiency}.

\subsubsection{Robustness with respect to the spatial and temporal domains as well as physical and numerical parameters} \label{sec_rob}

There may be various ways how to simultaneously achieve~\eqref{eq_est_rel} and~\eqref{eq_est_eff} \cor{(and also~\eqref{eq_est_eff_loc})}. Then, a crucial indicator of the quality of the derived estimate $\eta(\tuhki)$ or $\eta(\tuhtki)$ is the nature and behavior of the ``generic'' constant $C$ from~\eqref{eq_est_eff} (and~\eqref{eq_est_eff_loc}) (recall that there is no unknown constant in~\eqref{eq_est_rel}). Ideally, $C$ is truly generic, only depending on the space dimension $d$ and the shape regularity (smallest angle) of the spatial meshes $\Th$ (say taking a value $2 - 10$ in practical applications when $1 \leq d \leq 3$). In any case, it is requested that $C$ is independent of the spatial domain $\Om$ (namely its size) and the temporal domain $(0,t_{\rm{F}})$ (the final simulation time $t_{\rm{F}}$) and of any numerical parameters, namely the mesh size $h$ and the time step $\tau$. It is also desirable to have $C$ independent of the physical parameters such us medium properties (the smallest and highest (eigen)values of diffusion tensor $\Km$ or specifications of the nonlinearity in $\Km(\cdot)$). We then speak of {\em robustness}.

\subsubsection{Asymptotic exactness} \label{sec_ass_ex}

One characterization of the quality of an a posteriori error estimate is the so-called effectivity index given by
\be \label{eq_I_eff}
    I_{\mathrm{eff}} \eq \frac{\eta(\tuhki)}{\enorm{\tu-\tuhki}} \quad \text { or } \quad I_{\mathrm{eff}} \eq \frac{\eta(\tuhtki)}{\enorm{\tu-\tuhtki}}.
\ee
Recall that in our setting, from~\eqref{eq_est_rel}, $I_{\mathrm{eff}}$ is greater than or equal to $1$. If $I_{\mathrm{eff}}$ goes to the optimal value of $1$ with decreasing mesh size $h$ and time step $\tau$, we speak of {\em asymptotic exactness}. We will see such cases below.

\subsubsection{Inexpensive evaluation} \label{sec_inexp}

There exist very sharp estimates $\eta(\tuhki)$ or $\eta(\tuhtki)$ in~\eqref{eq_est_rel}\cor{, which, however,} are expensive to compute and evaluate, comparatively as expensive as to compute the numerical approximations $\tuhki$ or $\tuhtki$ themselves. 
\cor{Indeed, from~\eqref{eq_a_post_Lapl_sp} and~\eqref{eq_Lapl_prop}, as sharp as requested upper bound can be obtained upon a solution of some global problem with increasing size, see Theorem~\ref{thm_err_char} below for more details.}
Such classes of estimates are\cor{, however,} excluded from our considerations, as not suitable in practice (they may serve to certify the error but it is hard to think of their use in adaptivity \cor{because of their cost}). 
Some other classes of estimates request a solution of some small\cor{,} local\cor{,} mutually independent problems in patches of mesh elements. These are termed implicit and can be labelled ``inexpensive''. 
In this contribution, however, we do not consider such estimates either, and we only focus on estimates whose evaluation from the approximate solution $\tuhki$ or $\tuhtki$ is yet cheaper, {\em explicit}. 
For example, the potential reconstruction $\prh$ from~\eqref{eq_pot_rec} is obtained from $\tuh$ \cor{(passing through $\pth$)} by some local postprocessing and averaging. 
On general polytopal meshes $\T_H$, we will then develop estimates whose evaluation merely consists in {\em multiplications} of {\em local vectors} by {\em local matrices}, where the local matrices are pre-processed from the scheme or cell geometry at hand and the local vectors are the degrees of freedom representing the available values of the local fluxes and potentials. These estimate are then as {\em inexpensive} as possible.

\subsubsection{Estimating the error components} \label{sec_comps}

We have seen in Section~\ref{sec_err_comps} that the total error in a numerical approximation of, say~\eqref{eq_PDE_unst}--\eqref{eq_nonl_heat_op}, has several components, namely due to spatial discretization, temporal discretization, iterative linearization, and iterative linear algebraic resolution. Congruently, an a posteriori error estimator should {\em identify} and estimate these {\em error components}, developing~\eqref{eq_est_rel} structurally into
\bse \label{eq_est_rel_comps} \ba 
    \enorm{\tu-\tuhki} & \leq \eta(\tuhki) \leq \eta_{\mathrm{sp}}(\tuhki) + \eta_{\mathrm{lin}}(\tuhki) + \eta_{\mathrm{alg}}(\tuhki) \label{eq_est_rel_comps_st}  \\
    \text{or} & \nn \\
    \enorm{\tu-\tuhtki} & \leq \eta(\tuhtki) \leq \eta_{\mathrm{sp}}(\tuhtki) + \eta_{\mathrm{tm}}(\tuhtki) + \eta_{\mathrm{lin}}(\tuhtki) + \eta_{\mathrm{alg}}(\tuhtki), \label{eq_est_rel_comps_unst} 
\ea \ese
in congruence with the example of error components in~\eqref{eq_err_comps}.

\subsection{Adaptivity: balancing the error component estimates} \label{sec_adapt}

Having approximately identified the total error and its components as per~\eqref{eq_est_rel_comps}, one can think of how to use this information. We are in particular going to develop {\em adaptivity} consisting in {\em balancing} of the estimated {\em error components}. For~\eqref{eq_est_rel_comps_unst}, for example, this consists in making all the estimators $\eta_{\mathrm{sp}}(\tuhtki)$, $\eta_{\mathrm{tm}}(\tuhtki)$, $\eta_{\mathrm{lin}}(\tuhtki)$, and $\eta_{\mathrm{alg}}(\tuhtki)$ of comparable size, uniformly throughout the entire simulation. Indeed, shall for example the iterative linear algebraic solver error estimate $\eta_{\mathrm{alg}}(\tuhtki)$ be orders of magnitude smaller than the spatial or temporal discretization error estimates, we may be wasting our computational resources (unfortunately, this is common in practice). On a more common note, say for a steady problem, we will develop
\[
    \big[\eta_{\mathrm{sp}}(\tuhki)\big]^2 = \sum_{\elm \in \Th} \big[\eta_{\elm, \mathrm{sp}}(\tuhki)\big]^2,
\]
in extension of~\eqref{eq_est_loc}. We then design algorithms ensuring that the elementwise a posteriori error estimators $\eta_{\elm, \mathrm{sp}}(\tuhtki)$ have all comparable values. This is termed {\em mesh adaptivity} (adaptive mesh refinement). 

\subsection{Mass balance recovery at any step of numerical resolution} \label{sec_mass_bal}

A locally conservative method produces an approximation to the Darcy velocity $\tu$ from~\eqref{eq_Lapl_prop} in a way that 
\[
    \tu_h \in \Hdv \quad \text{with} \quad \Dv \tu_h = f
\]
for a piecewise constant source term $f$, see~\eqref{eq_loc_cons}. Unfortunately, during an iterative linearization or during iterative linear algebraic resolution, this may be lost in that 
\[
    \tu_h^{k,i} \not \in \Hdv \quad \text{\cor{and/or}} \quad \Dv \tu_h^{k,i} \neq f.
\]
One of the specificities of our approach is that in the construction of our a posteriori error estimates giving a guaranteed error upper bound and distinguishing the total error components, we actually recover a locally conservative flux field (Darcy velocity) $\bsig_h^{k,i}$ such that 
\[
    \bsig_h^{k,i} \in \Hdv \quad \text{and} \quad \Dv \bsig_h^{k,i} = f,
\]
see~\eqref{eq_eq_fl_rec}. In practice, we may go for a cheaper variant where $\Dv \bsig_h^{k,i} = f + \rho_h^{k,i}$ with $\rho_h^{k,i}$ is a remainder made as small as necessary.

\subsection{Bibliographic resources}\label{sec_biblio_a_post}

Focusing on the spatial discretization error, it has been rather soon understood how to structurally obtain a guaranteed error upper bound as in~\eqref{eq_a_post_Lapl_sp}: this follows from the Prager--Synge equality~\cite{Prag_Syng_47}, relying on the physical and mathematical properties of the weak solutions of PDEs as discussed in Section~\ref{sec_props}. In particular, in finite element methods, where one obtains an approximate solution $p_h \in \Hoo$, one was in quest of an equilibrated flux reconstruction $\bsig_h \in \Hdv$ with $\Dv \bsig_h = f$ as in~\eqref{eq_eq_fl_rec}. A {\em local} construction, is however, a bit involved in finite elements, so that its establishment was rather labyrinthine, see
Ladev{\`e}ze~\cite{Lad_these_75}, Ladev{\`e}ze and Leguillon~\cite{Lad_Leg_83}, Bossavit~\cite{Boss_hypercirc_a_post_98}, Destuynder and M{\'e}tivet~\cite{Dest_Met_expl_err_CFE_99}, Larson and Niklasson~\cite{Lars_Nikl_cons_flux_FEs_04}, Luce and Wohlmuth~\cite{Luce_Wohl_local_a_post_fluxes_04}, Vohral{\'{\i}}k~\cite{Voh_a_post_FE_loc_cons_min_CRAS_08}, Braess and Sch{\"o}berl~\cite{Braess_Scho_a_post_edge_08}, Cottereau~\eal\ \cite{Cott_Diez_Huer_strict_lin_09}, Ern and Vohral{\'{\i}}k~\cite{Ern_Voh_adpt_IN_13, Ern_Voh_p_rob_15}, and the references therein.

In contrast, in methods locally conservative by construction, we already have $\tu_h \in \Hdv$ with $\quad \Dv \tu_h = f$: \cor{directly from mixed finite elements and by local postprocessing for finite volume-type methods, as detailed below.} And since a potential reconstruction $\prh \in \Hoo$ as in~\eqref{eq_pot_rec} is rather easy to design, a posteriori error estimates of the form~\eqref{eq_a_post_Lapl_sp} were put in place in Achdou~\eal\ \cite{Ach_Ber_Coq_FV_Darcy_03}, Ainsworth~\cite{Ains_rob_a_post_NCFE_05, Ains_a_post_MFE_07}, Kim~\cite{Kim_a_post_MFE_07}, and Vohral{\'{\i}}k~\cite{Voh_apost_MFE_07, Voh_apost_FV_08}, see also the references therein. 
In discontinuous Galerkin and related methods, $\tu_h \in \Hdv$ with $\Dv \tu_h = f$ is not directly available but can be readily reconstructed, see Bastian and Rivi{\`e}re~\cite{Bas_Riv_DG_H_div_03}, Kim~\cite{Kim_a_post_MFE_07}, and Ern~\eal\ \cite{Ern_Nic_Voh_DG_flux_rec_CRAS_07}. More references specifically for lowest-order locally conservative methods/finite volumes are discussed in Sections~\ref{sec_biblio_Pois} and~\ref{sec_biblio_Darcy} below.

Unsteady problems were analyzed in Picasso~\cite{Pic_adpt_par_98}, Verf{\"u}rth~\cite{Ver_a_post_heat_03}, Bergam~\eal\ \cite{Ber_Ber_Mgh_a_post_par_04}, Ern and Vohral{\'{\i}}k~\cite{Ern_Voh_a_post_par_10}, and Ern~\eal\ \cite{Ern_Sme_Voh_heat_HO_Y_17}, see also the references therein. Iterative linear algebraic solver error and iterative linearization error were in particular addressed in Becker~\eal\ \cite{Beck_John_Ran_95}, Jir{\'a}nek~\eal\ \cite{Jir_Strak_Voh_a_post_it_solv_10}, Arioli~\eal\ \cite{Ar_Geor_Log_st_crit_cvg_13}, 
Ern and Vohral{\'{\i}}k~\cite{Ern_Voh_adpt_IN_13}, Heid and Wihler~\cite{Heid_Wih_it_lin_NL_20}, 
Gantner~\eal\ \cite{Gant_Hab_Praet_Schi_opt_cost_NL_21}, Haberl~\eal\ \cite{Hab_Praet_Schim_Voh_inex_Newt_opt_cost_21}, Mitra and Vohral{\'{\i}}k~\cite{Mitra_Voh_Richards_24}, and the references therein. In-depth presentation of a posteriori error estimates and adaptivity can then be found in the monographs by Synge~\cite{Synge_hypercirc_57}, Hlav{\'a}{\v{c}}ek~\eal\ \cite{Hlav_Has_Nec_Lov_VI_88}, Zeidler~\cite{Zeid_FA_90}, Ainsworth and Oden~\cite{Ainsw_Oden_a_post_FE_00}, Babu{\v{s}}ka and Strouboulis~\cite{Bab_Stroub_FE_rel_01}, Bangerth and Rannacher~\cite{Bang_Ran_AFEM_03}, Han~\cite{Han_a_post_dual_05}, Ladev\`eze and Pelle~\cite{Lad_Pel_mech_05}, 
Repin~\cite{Repin_book_08}, Nochetto~\eal\ \cite{Noch_Sieb_Vees_09}, Deuflhard and Weiser~\cite{Deuf_Weiss_adpt_book_12}, 
and Verf{\"u}rth~\cite{Verf_13}, see also Chamoin and Legoll~\cite{Cham_Legoll_a_post_rev_23}, Bartels and Kaltenbach~\cite{Bart_Kalt_a_post_conv_dual_24}, Mghazli~\cite{Mghaz_a_post_lin_25, Mghaz_a_post_nonlin_25}, and Smears~\cite{Smears_a_post_par_25}.

\section{Sobolev spaces and their piecewise polynomial subspaces, me\-shes, and finite volume methods} \label{sec_not}

In this section, we describe the infinite-dimensional Sobolev spaces $\Hoo$ and $\Hdv$ and their finite-dimensional subspaces formed by piecewise polynomials. We also set up the notation, define a simplicial mesh, and describe the considered polytopal meshes. We finally introduce the basic principle of the finite volume and related lowest-order locally conservative methods.

\subsection{Basic notation}

We introduce here the basic notation used throughout the manuscript.

\subsubsection{Domain $\Om$} \label{sec_dom_Om}

We let $\Om \subset \RR^d$ be an open bounded connected set in a form of an interval or polygon or polyhedron for respectively $d=1,2,3$ space dimensions, or a polytope in general for $d \geq 1$, with a Lipschitz-continuous boundary $\pt \Om$. The notation $\om$ is reserved for a subdomain $\om \subset \Om$.

\subsubsection{Lebesgue spaces}

For $\om \subset \Om$, we let $\Lti{\om}$ stand for square-integrable scalar-valued functions $v,w: \om \ra \RR$ with the scalar product 
\bse \label{eq_not} \be 
    (v,w)_\om \eq \int_\om v(\tx) w(\tx) \dx. \label{eq_not_scal} 
\ee
Square-integrable vector-valued functions $\tv,\tw: \om \ra \RR^d$ are then collected in $\tLti{\om} \eq [\Lti{\om}]^d$, with the scalar product 
\be 
    (\tv,\tw)_\om \eq \int_\om \tv(\tx) \scp \tw(\tx) \dx. \label{eq_not_vect} 
\ee \ese
The associated norms are denoted by $\norm{\cdot}_\om$. When $\om = \Om$, we drop the subscripts.

\subsubsection{Measures, sizes, and cardinalities}

We write $|\om|$ for the Lebesgue measure of a domain $\om$, $|\tv|$ for the Euclidean size of a vector $\tv \in \RR^d$, and $|\Set|$ for the cardinality of a set $\Set$. 
We denote by $\tn_\om$ the exterior unit normal of a domain $\om$.

\subsection{Infinite-dimensional Sobolev spaces $\Hoo$ and $\Hdv$} \label{sec_spaces}

In the mathematical sense, the solutions of partial differential equations live in specific infinite-dimensional spaces that admit generalized (weak) partial derivatives and divergences. This is in particular the Sobolev space $\Hoo$ for scalar-valued functions and the Sobolev space $\Hdv$ for vector-valued functions. We recapitulate some basics about these spaces here. 
We follow Allaire~\cite{Allaire_05}, Ern and Guermond~\cite{Ern_Guermond_FEs_I_21}, and Vohral\'ik~\cite{Voh_a_post_LN_24, Voh_FEM_LN_25}. 

\subsubsection{Sobolev space $\Hoo$}\label{sec_H1}

We start by the notion of a weak partial derivative. Let $\Do$ be the space of
functions from $\Cinf$ with a compact support in $\Om$.

\bd[Weak partial derivative] \label{def_WPD} Let a scalar-valued function $v
: \Om \ra \RR$ be given. We say that $v$ admits a weak $i$-th partial
derivative, $1 \leq i \leq d$, if
\ben[parsep=1pt, itemsep=1pt, topsep=5pt, partopsep=5pt]

\item \label{prop_L2} $v \in \Lt$;

\item \label{prop_wi} there exists a function $w_i: \Om \ra \RR$
such that

\ben[parsep=1pt, itemsep=1pt, topsep=1pt, partopsep=1pt]

\item \label{prop_L2_der} $w_i \in \Lt$;

\item \label{prop_Green} $(v, \pt_{\tx_i} \vf) = - (w_i, \vf)
\qquad \forall \vf \in \Do$.

\een

\een
The function $w_i$ is called the {\em weak $i$-th partial
derivative} of $v$. We use the
notation $\pt_{\tx_i} v = w_i$. \ed

We then define a weak gradient:

\bd[Weak gradient] Let a scalar-valued function $v : \Om \ra \RR$ be given.
We say that $v$ admits a {\em weak gradient} if
$v$ admits the weak $i$-th partial derivative for all $1 \leq i \leq
d$. We set
\be \label{eq_weak_grad}
    \Gr v \eq (\pt_{\tx_1} v, \ldots, \pt_{\tx_d} v)^t.
\ee
\ed

\bd[The space $\Ho$] \label{def_Ho} The space $\Ho$ is the space of all the functions that admit the weak gradient. \ed

Let us recall from~\cite{Allaire_05, Ern_Guermond_FEs_I_21} that
$\Ho$ is a Hilbert space for the scalar product given by $(u,v) + (\Gr u, \Gr v)$.

\bd[The space $\Hoo$] \label{def_Hoo} The space
$\Hoo$ is the space of functions $v \in
\Ho$ such that $v|_{\pt \Om} = 0$. \ed

Recalling again from~\cite{Allaire_05, Ern_Guermond_FEs_I_21}, $\Hoo$ is a Hilbert space for the scalar product given by $(\Gr u, \Gr v)$.

\begin{figure}
\centerline{\includegraphics[width=0.48\textwidth]{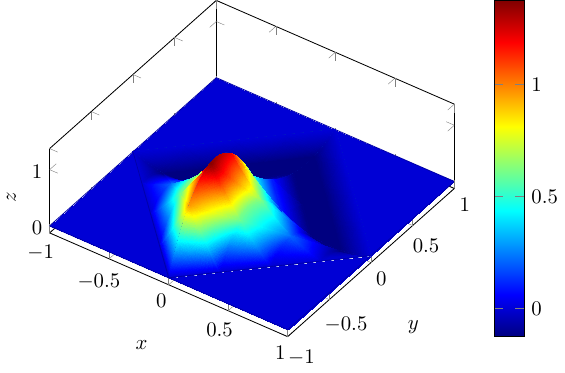} \quad \includegraphics[width=0.48\textwidth]{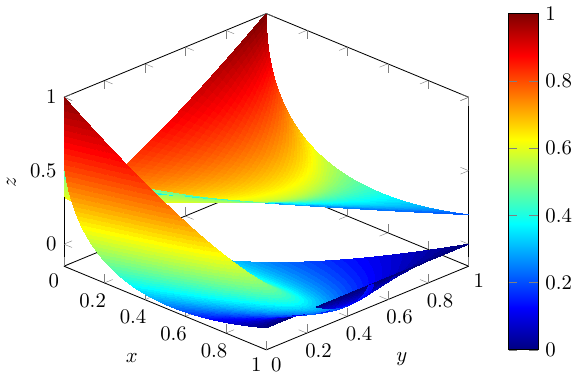}}
\caption{Example of a function belonging to the Sobolev space $\Hoo$ (left; the function actually lies in $\Cz$, but not in $\Co$). Example of a function not belonging to the Sobolev space $\Hoo$ but belonging to the broken Sobolev space $\HTh$ for a mesh $\Th$ composed of two triangles (right; note that there is no trace-continuity on the common face of the two mesh elements)}\label{fig_Ho}
\end{figure}

\br[Trace and trace continuity] \label{rem_trace} In the above definition, we have used
the notation $v|_{\pt \Om}$. As the functions from $\Ho$ are taken from the Lebesgue space $\Lt$, $v|_{\pt \Om}$ in the sense of a restriction is not defined since $\pt \Om$ is a set of
measure zero. In our notation, by $v|_{\pt \Om}$, we mean the {\em trace} of $v \in \Ho$ on the boundary of $\Om$, a key notion for the Sobolev space $\Ho$. Crucially, $v|_{\pt \Om}$ in the sense of traces and $v|_{\pt \Om}$ as a restriction coincide for any function $v \in \Ho$ that is continuous, $v \in \Ho \cap \Cz$. Let us also recall that $\Ho
\subset \Cz$ in one space dimension but $\Ho \not \subset \Cz$ in multiple space dimensions. Crucially, functions $v$ in $\Ho$, though not necessarily continuous, are continuous in the trace sense.
An illustration of a function $v$ belonging to $\Ho$ is given in Figure~\ref{fig_Ho}, left, and of a function $v$ not belonging to $\Ho$ in Figure~\ref{fig_Ho}, right. We refer for 
details to~\cite{Allaire_05, Ern_Guermond_FEs_I_21, Voh_a_post_LN_24}. \er

\subsubsection{Sobolev space $\Hdv$}\label{sec_Hdv}

We now turn to weak a divergence. 

\bd[Weak divergence] \label{def_WDG} Let a vector-valued function $\tv: \Om
\ra \RR^d$ be given. We say that $\tv$ admits a weak divergence if
\ben[parsep=1pt, itemsep=1pt, topsep=5pt, partopsep=5pt]

\item \label{prop_L2_div} $\tv \in \tLt$;

\item there exists a function $w: \Om \ra \RR$ such that

\ben[parsep=1pt, itemsep=1pt, topsep=1pt, partopsep=1pt]

\item \label{prop_L2_dv} $w \in \Lt$;

\item \label{prop_Green_dv} $(\tv, \Gr \vf) = - (w, \vf) \qquad
\forall \vf \in \Do$.

\een

\een
The function $w$ is called the {\em weak
divergence} of $\tv$. We use the notation $\Dv
\tv = w$. \ed

\bd[The space $\Hdv$] \label{def_Hdv} The space
$\Hdv$ is the space of all the functions that admit the weak divergence. \ed

Let us recall from~\cite{Allaire_05, Ern_Guermond_FEs_I_21} that
$\Hdv$ is a Hilbert space for the scalar product given by $(\tu,\tv) + (\Dv \tu, \Dv \tv)$. Similarly to Remark~\ref{rem_trace}, there holds:

\begin{figure}
\centerline{\includegraphics[width=0.45\textwidth]{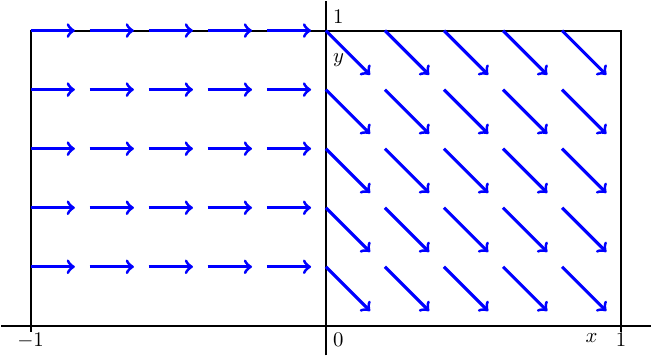} \quad \includegraphics[width=0.45\textwidth]{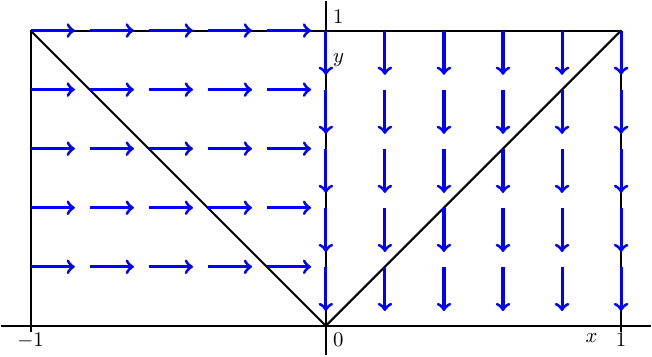}}
\caption{Example of a function belonging to the Sobolev space $\Hdv$ (left; on the interface $x=0$ between the two subdomains $\Om_1 = (-1,0)\times(0,1)$ and $\Om_2 = (0,1)\times(0,1)$, $\tv|_{\Om_1} \scp \tn = \tv|_{\Om_2} \scp \tn$, with $\tn=(1,0)^\trp$; $\tv$ is not continuous for each component (the $y$ component of $\tv$ is discontinuous, as it passes from value $0$ in $\Om_1$ to a nonzero value in $\Om_2$) but $\tv$ is normal-trace continuous). Example of a function not belonging to the Sobolev space $\Hdv$ but belonging to the broken Sobolev space $\HdTh$ (right, there is no normal trace-continuity across the mesh face at $x=0$).}\label{fig_H_div}
\end{figure}

\br[Normal trace and normal-trace continuity] \label{rem_norm_trace} On $\pt \Om$, functions $\tv \in \Hdv$ admit a normal component in appropriate sense, called normal trace and denoted as $(\tv \scp \tn)|_{\pt \Om}$. Namely, if $\tv \in \Hdv \cap [\Cz^d]$, then $\tv \scp \tn$ is the usual normal component of $\tv$ on the boundary $\pt \Om$. Functions $\tv$ in $\Hdv$ are not necessarily continuous in each component, not necessarily normal-component continuous, but are normal-trace continuous. An illustration of a function $\tv$ belonging to $\Hdv$ is given in Figure~\ref{fig_H_div}, left, and of a function $\tv$ not belonging to $\Hdv$ in Figure~\ref{fig_H_div}, right. We refer for details to~\cite{Allaire_05, Ern_Guermond_FEs_I_21}. \er

The Sobolev spaces $\Hoo$ from Definition~\ref{def_Hoo} and $\Hdv$ from Definition~\ref{def_Hdv} form exactly the right setting for the following Green theorem: 

\bt[Green theorem] \label{thm_green} Let $v \in \Hoo$ and let $\tw \in \Hdv$. Then
\be \label{eq_Green_Hdv}
    (\tw, \Gr v) + (\Dv \tw, v) = 0.
\ee
\et

\cor{In the following, we will also employ the above notation and developments for subdomains $\om \subset \Om$.}

\subsection{Simplicial meshes}\label{sec_Th}

\begin{figure}
\centerline{\includegraphics[width=0.43\textwidth]{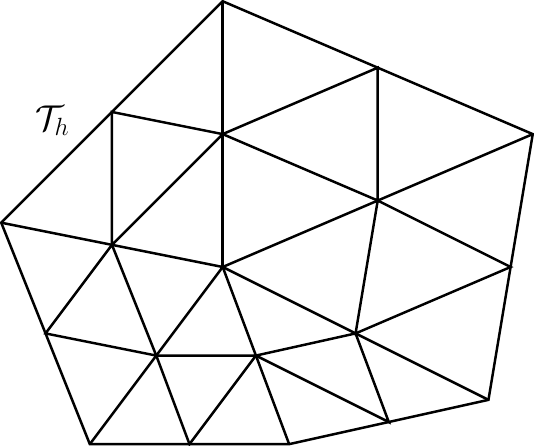} \qquad \includegraphics[width=0.43\textwidth]{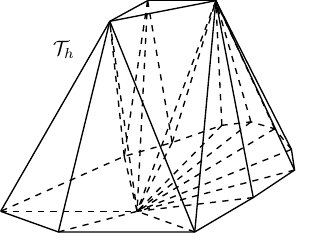}}
\caption{Simplicial mesh $\Th$ for $d=2$ (left) and $d=3$ (right)}\label{fig_Th}
\end{figure}

In order to form finite-dimensional subspaces of $\Hoo$ and $\Hdv$, we will need a simplicial mesh $\Th$ of $\Om$. This a finite collection of closed $d$-simplices (intervals for $d=1$, triangles for $d=2$, tetrahedra for $d=3$) such that $\cup_{\elm \in \Th} \elm = \overline \Om$ and such that the intersection of two different simplices is either empty or their entire $l$-dimensional face, $0 \leq l \leq (d-1)$. Here, for example for space dimension $d=3$, a vertex is a $0$-dimensional face, an edge is a $1$-dimensional face, and a face is a $2$-dimensional face. Henceforth, by face, we mean $(d-1)$-dimensional face. Figure~\ref{fig_Th} gives illustrations for $d=2$ and $d=3$. 
We denote by
\be \label{eq_shape_reg}
    \theta_{\Th} \eq \max_{\elm \in \Th} \frac{h_\elm}{\iota_\elm}
\ee
the shape-regularity parameter of the mesh $\Th$, where $h_\elm$ is the diameter of the simplex $\elm \in \Th$ and $\iota_\elm$ is the diameter of largest ball inscribed in $\elm$. 
\cor{This parameter expresses how much the mesh is ``distorted'': it is close to $1$ for meshes of equilateral simplices and increases in presence of small angles. For each mesh $\Th$, $\theta_{\Th}$ is a bounded number. Our theoretical constants such as $C$ from~\eqref{eq_est_eff} and~\eqref{eq_est_eff_loc} will depend on $\theta_{\Th}$. If we consider a sequence of meshes, then we need $\theta_{\Th}$ uniformly bounded.}

The set of the ($(d-1)$-dimensional) faces of $\Th$ is denoted as $\Fh$. It is divided into interior faces $\Fhint$ and boundary faces $\Fhext$: $\sd \in \Fhint$ if there exist $\elm, \elmt \in \Th$, $\elm \neq \elmt$, such that $\sd = \elm \cap \elmt$, and $\sd \in \Fhext$ if there exist $\elm \in \Th$ such that $\sd \cor{\subset} \elm \cap \pt \Om$. We will also need the notation $\Vh$ for the vertices of the mesh $\Th$ and $\Ve$ for the vertices of the face $\sd \in \Fh$.
For an interior face $\sd \in \Fhint$, we fix an arbitrary orientation and denote the corresponding unit normal vector by $\tn_\sd$. 
For a boundary face $\sd \in \Fhext$, we make $\tn_\sd$ coincide with the exterior unit normal $\tn_\Omega$ of $\Omega$. 
Finally, for all $\elm \in \Th$, we denote by $\tn_\elm$ the unit normal vector to $\sd$ pointing out of $\elm$ and by $\FK$ the set of all faces $\sd$ of $\elm$.

\subsection{Broken Sobolev spaces $\HTh$ and $\HdTh$}\label{sec_HTh}

Before we pass to finite-dimensional subspaces of $\Hoo$ and $\Hdv$, it is useful to define the broken Sobolev space $\HTh$ as
\be \label{eq_H1Th}
    \HTh \eq \{v \in \Lt; \, v|_\elm \in \Hoi{\elm} \quad \forall \elm  \in \Th \},
\ee
and the broken Sobolev space $\HdTh$ as
\be \label{eq_HdTh}
    \HdTh \eq \{\tv \in \tLt; \, \tv|_\elm \in \Hdvi{\elm} \quad \forall \elm  \in \Th \}.
\ee
An illustration of a function $v \in \HTh$ is given in Figure~\ref{fig_Ho}, right, and an illustration of a function $\tv \in \HdTh$ is given in Figure~\ref{fig_H_div}, right. In this setting, Remarks~\ref{rem_trace} and~\ref{rem_norm_trace} can \cor{respectively} be developed into, \cor{see~\cite[Lemma~1.23]{Di_Pietr_Ern_book_12} or~\cite[Theorem~18.8]{Ern_Guermond_FEs_I_21}:}

\bl[Trace continuity over mesh faces] \label{lem_trace_cont} 
There holds
\bse \label{eq_trace_cont} \ba 
    v \in \HTh \text{ and } (v|_\elm)|_\sd = (v|_\elmt)|_\sd \quad \forall \sd = \elm \cap \elmt \in \Fhint & \Longleftrightarrow v \in \Ho, \\
    v \in \HTh \text{ and } (v|_\elm)|_\sd = (v|_\elmt)|_\sd \quad \forall \sd = \elm \cap \elmt \in \Fhint & \nn \\
    \text{ and } (v|_\elm)|_\sd = 0 \quad \forall \sd = \elm \cap \pt \Om \in \Fhext & \Longleftrightarrow v \in \Hoo.
\ea \ese
In words, functions from the broken Sobolev space $\HTh$ belong to the Sobolev space $\Ho$ if and only if their traces are uniquely defined on all mesh interior faces, with the same trace from both elements sharing the given interior face. We say that the functions from $\Ho$ are trace continuous. Functions from $\Hoo$ then additionally have the traces zero on boundary faces. 
\el

\cor{Similarly, following~\cite[Lemma~1.24]{Di_Pietr_Ern_book_12} or~\cite[Theorem~18.10]{Ern_Guermond_FEs_I_21}, we have:}

\bl[Normal trace continuity over mesh faces] \label{rem_normal_trace_cont} Functions from $\HdTh$ belong to $\Hdv$ if and only if their normal traces are continuous over mesh interior faces in appropriate normal-trace sense. \el

\subsection{Finite-dimensional subspaces of $\Hoo$ and $\Hdv$} \label{sec_discr_spaces}

We now introduce standard finite-dimensional subspaces of $\Hoo$ and $\Hdv$, again following Allaire~\cite{Allaire_05} and Ern and Guermond~\cite{Ern_Guermond_FEs_I_21}. 

\subsubsection{Piecewise polynomial subspaces of $\Lt$} \label{sec_pw_pol}

\begin{figure}
\centerline{\includegraphics[height=0.2\textwidth]{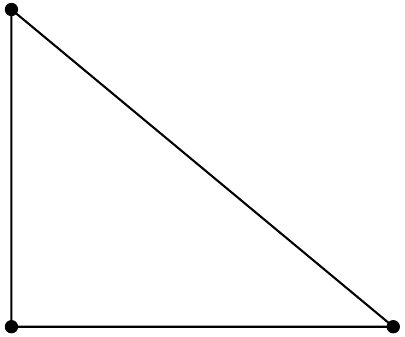} \quad \includegraphics[height=0.2\textwidth]{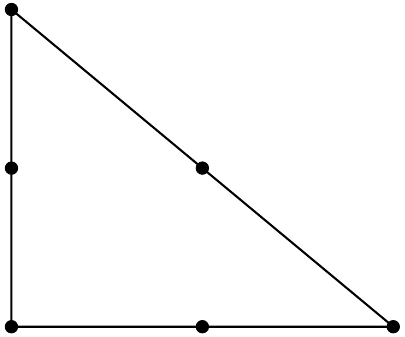} \quad \includegraphics[height=0.22\textwidth]{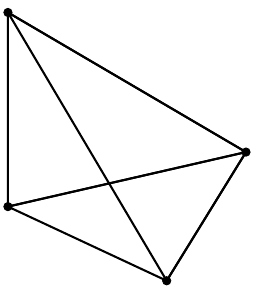} \quad \includegraphics[height=0.22\textwidth]{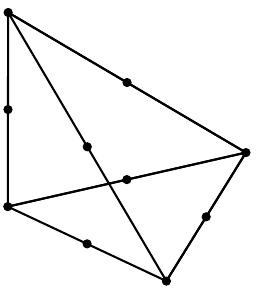}}
\caption{Spaces $\PP_k(\elm)$ and point values uniquely fixing a function $v_h \in \PP_k(\elm)$ (the so-called Lagrange nodes) for a mesh element $\elm \in \Th$. Polynomial degrees $k=1$ and $k=2$, space dimensions $d=2$ (left) and $d=3$ (right)}\label{fig_Lagr_K}
\end{figure}

Let $\elm \in \Th$ be a simplex \cor{from the mesh $\Th$} and \cor{let} $k \geq 0$ \cor{be} an integer. 
We denote by $\PP_k(\elm)$ the space of scalar-valued polynomials on $\elm$ of total degree at most $k$. In particular, $\PP_0(\elm)$ stands for constants on $\elm$, $\PP_1(\elm)$ for affine functions on $\elm$, and $\PP_2(\elm)$ for quadratic functions on $\elm$. The spaces for $k=1,2$ and $d=2,3$ are schematically visualized in Figure~\ref{fig_Lagr_K}. In particular, fixing the point values as \cor{per} Figure~\ref{fig_Lagr_K} uniquely determines a function in $\PP_k(\elm)$.
We will often use the space of piecewise $k$-degree polynomials on the simplicial mesh $\Th$, 
\be \label{eq_pw_pols}
    \PP_k(\Th) \eq \{v_h \in \Lt; \, v_h|_\elm \in \PP_k(\elm) \quad \forall \elm  \in \Th \}.
\ee

\subsubsection{Lagrange piecewise polynomial subspaces of $\Hoo$} \label{sec_Lagr}

The space $\PP_k(\Th)$ from~\eqref{eq_pw_pols} is a finite-dimensional subspace of the space $\HTh$ from~\eqref{eq_H1Th} but not of the Sobolev space $\Hoo$. To create a subspace of $\Hoo$ of the so-called {\em Lagrange} piecewise polynomials, we write\cor{, in several equivalent ways,}
\be \label{eq_Lagr} \begin{split}
    \PP_k(\Th) \cap \Hoo = {} & \{v_h \in \Hoo; \, v_h|_\elm \in \PP_k(\elm) \quad \forall \elm  \in \Th \} \\
     = {} & \{v_h \in \PP_k(\Th); \, (v_h|_\elm)|_\sd = (v_h|_\elmt)|_\sd \quad \forall \sd = \elm \cap \elmt \in \Fhint \\
    {} & \text{ and } (v_h|_\elm)|_\sd = 0 \quad \forall \sd \cor{\subset} \elm \cap \pt \Om \in \Fhext\} \\
    = {} & \cor{\{v_h \in }\PP_k(\Th) \cap \Cz; \, \cor{v_h = 0} \text{ on } \pt \Om \cor{ \}}.
\end{split} \ee
Indeed, Lemma~\ref{lem_trace_cont} holds with $\HTh$ replaced by $\PP_k(\Th)$ and $\Hoo$ replaced by $\PP_k(\Th) \cap \Hoo$, where the traces on faces $\sd$ now equal the \cor{usual} restrictions to $\sd$. Congruently, fixing the point values in each simplex $\elm$ \cor{of the mesh} $\Th$ as schematically visualized in Figure~\ref{fig_Lagr_K} by \cor{a unique} bullet \cor{value} yields the $\Hoo$- or $\Cz$-conformity for piecewise $k$-degree polynomials. 
We again refer for details to~\cite{Allaire_05, Ern_Guermond_FEs_I_21}. 

\subsubsection{Raviart--Thomas piecewise polynomial subspaces of $\Hdv$} \label{sec_RT}

\begin{figure}
\centerline{\includegraphics[height=0.25\textwidth]{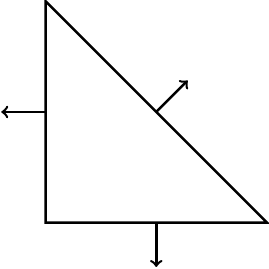} \qquad \includegraphics[height=0.3\textwidth]{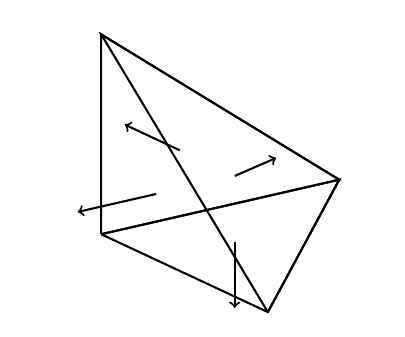}}
\caption{Spaces $\RT_0(\elm)$ and normal fluxes uniquely fixing a function $v_h \in \RT_0(\elm)$ for a \cor{simplex} $\elm \in \Th$. Space dimensions $d=2$ (left) and $d=3$ (right)}\label{fig_RTN_K}
\end{figure}

\begin{figure}
\centerline{\includegraphics[height=0.3\textwidth]{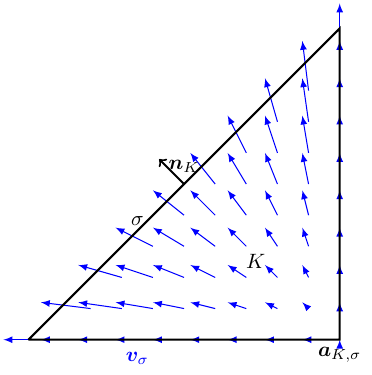} \qquad \qquad \includegraphics[height=0.3\textwidth]{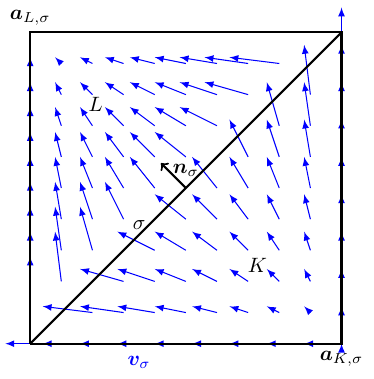}}
\caption{A basis function $\tv_\sd$ of $\RT_0(\elm)$ \cor{(supported on the simplex $\elm$)} with nonzero normal flux only through the face $\sd \in \FK$ (left). A basis function $\tv_\sd$ of $\RT_0(\Th) \cap \Hdv$ \cor{(supported on all (two here) simplices sharing the face $\sd$)} with nonzero normal flux only through the face $\sd \in \Fh$ (right).}\label{fig_RTN_bas}
\end{figure}

Recall that $\PP_0(\elm)$ stands for constants on the simplex $\elm \in \Th$. Following Whitney~\cite{Whit_geom_int_57}, Raviart and Thomas~\cite{Ra_Tho_MFE_77}, and N\'ed\'elec~\cite{Ned_mix_R_3_80}, see alternatively Boffi~\eal\ \cite{Bof_Brez_For_MFEs_13} or Ern and Guermond~\cite{Ern_Guermond_FEs_I_21}, we define the space of vector-valued polynomials
\be \label{eq_RT0} \bs
    \RT_0(\elm) & \eq [\PP_0(\elm)]^d + \tx \PP_0(\elm) \\
     & = \{\tv_h(\tx) + w_h(\tx) \tx; \, \tv_h \in [\PP_0(\elm)]^d, \, w_h \in \PP_0(\elm)\} \\
     & \cor{= \MS[2]\{\begin{pmatrix} a \\ b \end{pmatrix} + c \begin{pmatrix} x \\ y \end{pmatrix} \quad a,b,c \in \RR, \, \tx=(x,y)^{\mathrm{t}}\MS[2]\} \quad \text{ for } d=2,}\\
     & \cor{= \MS[2.8]\{\begin{pmatrix} a \\ b \\ c \end{pmatrix} + d \begin{pmatrix} x \\ y \\ z \end{pmatrix} \quad a,b,c,d \in \RR, \, \tx=(x,y,z)^{\mathrm{t}}\MS[2.8]\} \quad \text{ for } d=3.}
\es \ee
This space is slightly bigger than just constant vector-valued polynomials $[\PP_0(\elm)]^d$ but slightly smaller than affine vector-valued polynomials $[\PP_1(\elm)]^d$. \cor{Its dimension equals $d+1$, \ie, the number of faces of the simplex $\elm$, which means $3$ in two space dimensions and $4$ in three space dimensions.} Its distinctive properties are
\be \label{eq_RT0_prop}
    \Dv \tv_h \in \PP_0(\elm) \quad \text{ and } \quad (\tv_h \scp \tn_\elm)|_\sd \in \PP_0(\sd) \quad \forall \sd \in \FK \qquad \text{ for } \tv_h \in \RT_0(\elm),
\ee
\ie, the divergence is constant and the normal trace is constant on each face \cor{of the simplex $\elm$} for any function $\tv_h \in \RT_0(\elm)$. Let us recall here that for a polynomial $\tv_h$, the normal trace on a face $\sd$ equals the restriction to $\sd$ of the normal component $\tv_h \scp \tn_\elm$. \cor{Moreover, a} function $\tv_h$ in $\RT_0(\elm)$ is uniquely defined by prescribing the face normal fluxes, \ie,
\be \label{eq_RT0_DoF}
    \<\tv_h \scp \tn_\elm, 1\>_\sd \cor{\eq \int_\sd \tv_h \scp \tn_\elm} \qquad \forall \sd \in \FK,
\ee
which is a part of the schematic visualization of the spaces $\RT_0(\elm)$ for $d=2,3$ in Figure~\ref{fig_RTN_K}. We choose the $d+1$ basis functions $\tv_\sd$ as associated with the faces $\sd \in \FK$. Setting $\<\tv_\sd \scp \tn_\elm, 1\>_\sd = 1$ and $\<\tv_\sd \scp \tn_\elm, 1\>_\sdt = 0$ for $\sdt \neq \sd$, we obtain the explicit expression
\be \label{eq_RT0_bas}
    \tv_\sd (\tx) = \frac{1}{d |\elm|}(\tx - \ver_{\elm, \sd}), \qquad \tx \in \elm,
\ee
\cf\ the illustration in Figure~\ref{fig_RTN_bas}, left, where $\ver_{\elm, \sd}$ is the vertex of $\elm$ opposite to $\sd$.

On the simplicial mesh $\Th$, we now define
\be \label{eq_pw_RT+pols}
    \RT_0(\Th) \eq \{\tv_h \in \tLt; \, \tv_h|_\elm \in \RT_0(\elm) \quad \forall \elm  \in \Th \}.
\ee
This is a finite-dimensional subspace of the space $\HdTh$ from~\eqref{eq_HdTh} but not of the Sobolev space $\Hdv$. To create a subspace of $\Hdv$, we set
\be \label{eq_RTN} \begin{split}
    \RT_0(\Th) \cap \Hdv = {} & \{\tv_h \in \Hdv; \, \tv_h|_\elm \in \RT_0(\elm) \quad \forall \elm  \in \Th \}\\
     = {} & \{\tv_h \in \RT_0(\Th); \, (\tv_h|_\elm \scp \tn_\elm)|_\sd = - (\tv_h|_\elmt \scp \tn_\elmt)|_\sd \quad \forall \sd = \elm \cap \elmt \in \Fhint\};
\end{split} \ee
it follows by Lemma~\ref{rem_normal_trace_cont} that the above characterizations coincide.
Congruently, fixing the face normal values/integrals in each simplex $\elm \in \Th$ as schematically visualized in Figure~\ref{fig_RTN_K} by the arrows in particular yields the $\Hdv$-conformity for piecewise $\RT_0(\elm)$-polynomials. 
The space~\eqref{eq_RTN} of normal-component-continuous piecewise vector-valued polynomials is called the {\em Raviart--Thomas}(--N\'ed\'elec) space. Its dimension is the number of mesh faces in the mesh $\Th$. There is one basis function $\tv_\sd$ for each face $\sd \in \Fh$, prescribed such that $\<\tv_\sd \scp \tn_\sd, 1\>_\sd = 1$ and $\<\tv_\sd \scp \tn_\sdt, 1\>_\sdt = 0$ for all $\sdt \in \Fh$ with $\sdt \neq \sd$. An illustration is given in Figure~\ref{fig_RTN_bas}, right. We again refer for details to~\cite{Allaire_05, Ern_Guermond_FEs_I_21}. 

\cor{\br[Meshes consisting of rectangular parallelepipeds] \label{rem_mesh_rect} Above, we only consider simplicial meshes since these are our key to treat general polytopal meshes below. In practice, meshes consisting of rectangles, blocks, or in general rectangular parallelepipeds are often used. On rectangles or blocks, in place of~\eqref{eq_pw_pols}, piecewise polynomial subspaces of $\Lt$ are typically chosen as 
\be \label{eq_pw_pols_Q}
    \QQ_k(\Th) \eq \{v_h \in \Lt; \, v_h|_\elm \in \QQ_k(\elm) \quad \forall \elm  \in \Th \},
\ee
where $\QQ_k(\elm)$ is the space of scalar-valued polynomials on $\elm$ of degree at most $k$ in each space coordinate. Note that $\PP_k(\Th)$ and $\QQ_k(\Th)$ coincide for the lowest polynomial degree $k=0$, being simply piecewise constants with respect to the mesh $\Th$. The Raviart--Thomas spaces $\RT_0(\elm)$ then take the form
\bse \label{eq_RT0_Q}\be 
    \RT_0(\elm) \eq \MS[2]\{\begin{pmatrix} a + c x \\ b + d y \end{pmatrix} \quad a,b,c,d \in \RR, \, \tx=(x,y)^{\mathrm{t}}\MS[2]\} \quad \text{ for } d=2,
\ee
and
\be 
    \RT_0(\elm) \eq \MS[2.8]\{\begin{pmatrix} a + d x \\ b + e y \\ c + f z \end{pmatrix} \qquad a,b,c,d,e,f \in \RR, \, \tx=(x,y,z)^{\mathrm{t}}\MS[2.8]\} \quad \text{ for } d=3,
\ee \ese
in place of~\eqref{eq_RT0}. Subspaces of $\Hdv$ are then still created by~\eqref{eq_pw_RT+pols}--\eqref{eq_RTN}. \er}

\cor{In the following, we will sometimes employ the above notation and developments for subdomains $\om \subset \Om$ and the corresponding submeshes $\T_\om$ of $\om$.}

\subsection{Polytopal meshes with virtual simplicial submeshes} \label{sec_meshes}

Recall that $\Om \subset \RR^d$, $d \geq 1$, is an open bounded connected polytope with a
Lipschitz-continuous boundary $\pt \Om$. Meshes with mesh elements formed by polygons (in two space dimensions), polyhedrons (in three space dimensions), or in general polytopes (in $d$ space dimensions) are useful in practice. \cor{An example of a general polygonal element ($d=2$) is given in Figure~\ref{fig:notation_types} (left); some examples of} polygonal meshes can be found in Figures~\ref{fig:mesh} (left) and~\ref{fig:spe10.front} and of polyhedral meshes in Figures~\ref{fig:BO.perm} and~\ref{fig:BO.amr}. These meshes stand in contrast to regular meshes of rectangular parallelepipeds by their generality and, in comparison with meshes formed by simplicial elements, they may contain (much) fewer elements, \cf\ Figure~\ref{fig:mesh} (right).

\subsubsection{Polytopal mesh with a virtual simplicial submesh} \label{sec_meshes_polyt}

\begin{figure}
   \centering
    \includegraphics[width=0.3\linewidth]{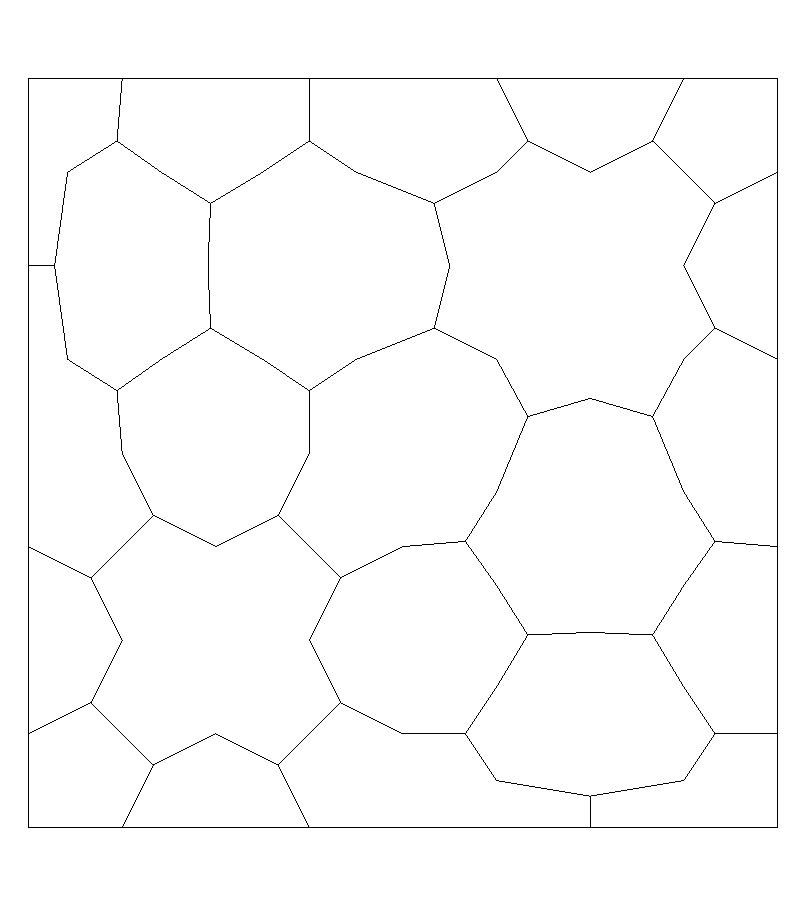}
\hspace{1.5cm}
    \includegraphics[width=0.3\linewidth]{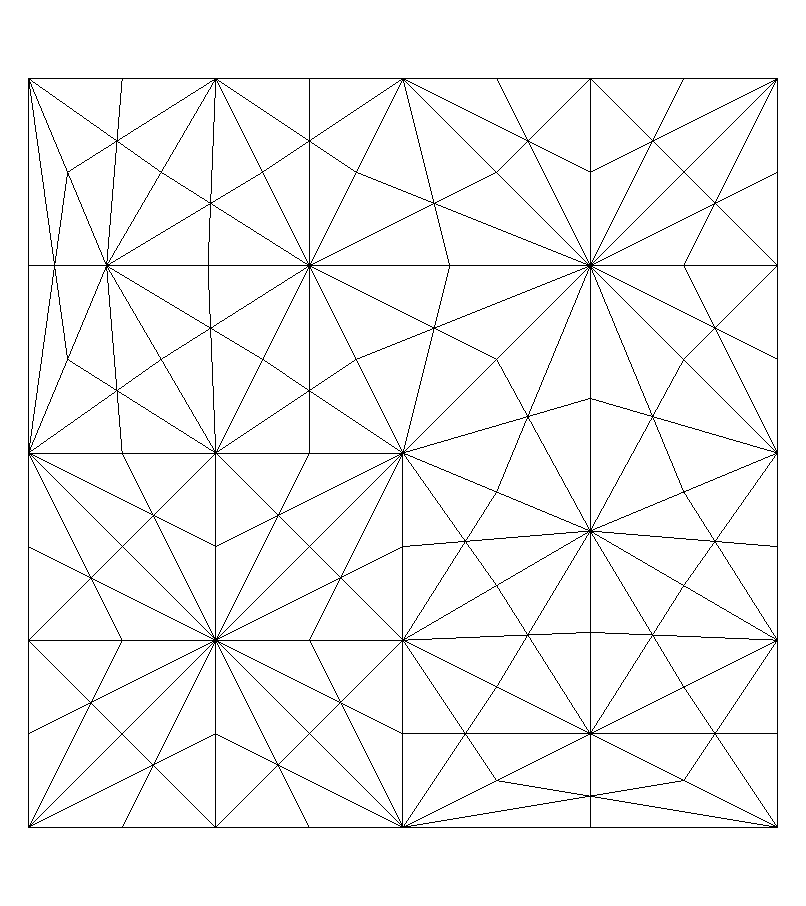}
    \caption{A polygonal mesh $\T_H$ and the corresponding
    triangular submesh $\Th$}
  \label{fig:mesh}
\end{figure}

We call a polytopal mesh a partition $\T_H$ of the domain $\Om$ into a finite collection of closed $d$-dimensional polytopes $\elm$ homotopic to a ball and with Lipschitz-continuous boundary $\pt \elm$ such that $\cup_{\elm \in \T_H} \elm = \overline \Om$. We suppose that the intersection of two different polytopes is either empty or their entire $l$-dimensional face, $0 \leq l \leq (d-1)$, or a collection of their entire $l$-dimensional faces, $0 \leq l \leq (d-1)$.
The elements $\elm$ can in particular be nonconvex and non star-shaped.
The structural assumption we impose is that there exists a virtual, not to be constructed in practice, simplicial mesh $\Th$ of $\Om$ in the sense of Section~\ref{sec_Th} such that $\TK \eq \Th|_\elm$ is a simplicial mesh of the polytope $\elm$. 
An illustration is provided in Figures~\ref{fig:notation_types} and~\ref{fig:mesh}.
Altogether, each element $\kappa \in \Th$ is a $d$-simplex, $\cup_{\kappa \in \Th} \kappa = \overline \Om$, $\cup_{\kappa \in \TK} \kappa = \elm$, and $\cup_{\elm \in \T_H} \TK = \Th$.

\def\sz{1}

\begin{figure}
\begin{tikzpicture}
       \draw[line width=\sz pt] (5,-0.7) -- (5,1)  ;
        \draw[line width=\sz pt] (5,1) -- (2.5,2)  ;
        \draw[line width=\sz pt] (2.5,2) -- (0.5,1.5)  ;
        \draw[line width=\sz pt] (0.5,1.5) -- (0.,0)  ;
        \draw[line width=\sz pt] (0,0) -- (2,-0.5) ;
        \draw[line width=\sz pt] (2,-0.5) -- (3.5,-1.5) ;
        \draw[line width=\sz pt] (3.5,-1.5) -- (5,-0.7) ;
\end{tikzpicture}
\hspace{0.2cm}
\begin{tikzpicture}
       \draw[line width=\sz pt] (5,-0.7) -- (5,1)  ;
         \draw[dashed, line width=\sz pt] (2.6,0.5) -- (5,-0.7) ;
         \draw[line width=\sz pt] (4.3,0.) ;
        \draw[line width=\sz pt] (5,1) -- (2.5,2)  ;
         \draw[dashed, line width=\sz pt] (2.6,0.5) -- (5,1) ;
        \draw[line width=\sz pt] (2.5,2) -- (0.5,1.5)  ;
         \draw[dashed, line width=\sz pt] (2.6,0.5) -- (2.5,2) ;
       \draw[line width=\sz pt] (0.5,1.5) -- (0.,0)  ;
         \draw [dashed, line width=\sz pt](2.6,0.5) -- (0.5,1.5) ;
        \draw[line width=\sz pt] (0,0) -- (2,-0.5) ;
          \draw [dashed, line width=\sz pt](2.6,0.5) -- (0,0) ;
        \draw[line width=\sz pt] (2,-0.5) -- (3.5,-1.5) ;
         \draw [dashed, line width=\sz pt](2.6,0.5) -- (2,-0.5) ;
        \draw[line width=\sz pt] (3.5,-1.5) -- (5,-0.7)  ;
         \draw[dashed, line width=\sz pt] (2.6,0.5) -- (3.5,-1.5) ;
\end{tikzpicture}
\hspace{0.2cm}
\begin{tikzpicture}
       \draw[line width=\sz pt] (5,-0.7) -- (5,1)  ;
         \draw[dashed, line width=\sz pt] (3,0.5) -- (5,-0.7) ;
         \draw[dashed, line width=\sz pt] (3,0.5) -- (5,0.15) ;
         \draw[line width=\sz pt] (4.3,0.) ;
        \draw[line width=\sz pt] (5,1) -- (2.5,2)  ;
         \draw[dashed, line width=\sz pt] (3,0.5) -- (5,1) ;
        \draw[line width=\sz pt] (2.5,2) -- (0.5,1.5)  ;
         \draw[dashed, line width=\sz pt] (1.6,0.5) -- (2.5,2) ;
         \draw[dashed, line width=\sz pt] (3,0.5) -- (2.5,2) ;
       \draw[line width=\sz pt] (0.5,1.5) -- (0.,0)  ;
         \draw [dashed, line width=\sz pt](1.6,0.5) -- (0.5,1.5) ;
        \draw[line width=\sz pt] (0,0) -- (2,-0.5) ;
          \draw [dashed, line width=\sz pt](1.6,0.5) -- (0,0) ;
         \draw [dashed, line width=\sz pt](1.6,0.5) -- (0.25,0.75) ;
        \draw[line width=\sz pt] (2,-0.5) -- (3.5,-1.5) ;
         \draw [dashed, line width=\sz pt](1.6,0.5) -- (2,-0.5) ;
        \draw[line width=\sz pt] (3.5,-1.5) -- (5,-0.7)  ;
         \draw[dashed, line width=\sz pt] (2,-0.5) -- (3,0.5) ;
         \draw[dashed, line width=\sz pt] (3,0.5) -- (3.5,-1.5) ;
         \draw[dashed, line width=\sz pt] (3,0.5) -- (1.6,0.5) ;
\end{tikzpicture}
\caption{\cor{Polygonal element $\elm \in \T_H$ (left). Virtual simplicial submesh $\TK$ of $\elm$ (middle and right). One point inside $\elm \in \T_H$ shared by all the simplices in $\TK$ and faces of $\Th$ do not subdividing the faces of $\T_H$ (middle). General situation (right).}}
\label{fig:notation_types}
\end{figure}
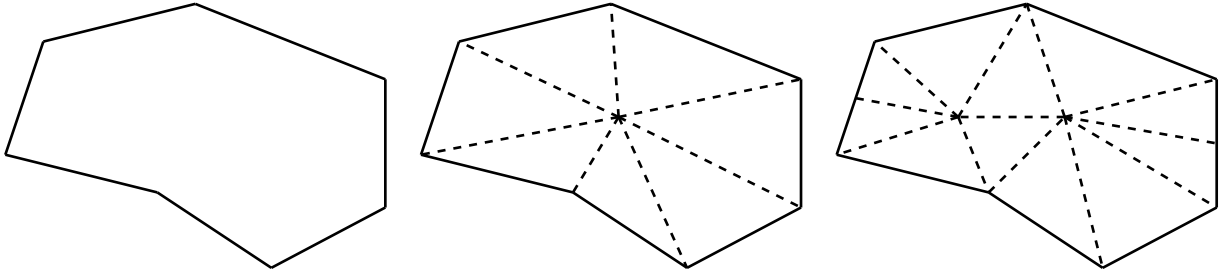

\def\sz{1}

\begin{figure}
\begin{tikzpicture}
          \draw (2.8,0.3) node{\small$(\mathsf{P})_K$};
       \draw[line width=\sz pt] (5,-0.7) -- (5,1)  ;
        \draw[->, line width=\sz pt]  (5,0.2) -- (5.6,0.2) ;
        \draw  (5.6,0.2) node[below]{\footnotesize$ (\mathsf{U}_K)_{\sigma_1}$};
        \draw[line width=\sz pt] (5,1) -- (2.5,2)  ;
        \draw[->, line width=\sz pt]   (3.75,1.5) -- (4,2.2) ;
        \draw  (4,2.2) node[right]{\footnotesize$ (\mathsf{U}_K)_{\sigma_2}$};
        \draw[line width=\sz pt] (2.5,2) -- (0.5,1.5)  ;
        \draw[->, line width=\sz pt]   (1.5,1.75) -- (1.3,2.3) ;
        \draw  (1.3,2.3) node[right]{\footnotesize$ (\mathsf{U}_K)_{\sigma_3}$};
        \draw[line width=\sz pt] (0.5,1.5) -- (0.,0)  ;
        \draw[->, line width=\sz pt]   (0.25,0.75) -- (-.5,1) ;
        \draw  (0.5,1.2) node[left]{\footnotesize$ (\mathsf{U}_K)_{\sigma_4}$};
        \draw[line width=\sz pt] (0,0) -- (2,-0.5) ;
        \draw[->, line width=\sz pt]   (1,-0.25) -- (0.75,-1) ;
        \draw  (0.8,-0.5) node[left]{\footnotesize$ (\mathsf{U}_K)_{\sigma_5}$};
        \draw[line width=\sz pt] (2,-0.5) -- (3.5,-1.5) ;
        \draw[->, line width=\sz pt]   (2.75,-1) -- (2.25,-1.5) ;
        \draw  (1.2,-1.2) node[right]{\footnotesize$ (\mathsf{U}_K)_{\sigma_6}$};
        \draw[line width=\sz pt] (3.5,-1.5) -- (5,-0.7) ;
        \draw[->, line width=\sz pt]   (4.25,-1.1) -- (4.6,-1.6) ;
        \draw  (4,-2) node[right]{\footnotesize$ (\mathsf{U}_K)_{\sigma_7}$};
         \draw (7.6,1.6) node{$\FK =  \{\sigma_i\}^7_{i=1} $} ;
         \draw (7.6,1.) node{$\mathsf{U}^{\mathrm{ext}}_K = \{(\mathsf{U}_K)_{\sigma_i}\}^7_{i=1}   $} ;
         \draw (7.5,0.4) node{$\TK =  \{\kappa_i\}^7_{i=1}$} ;
         \draw (7.5,-0.2) node{$\FKhext  =  \{\sigma_i\}^7_{i=1} $};
          \draw (7.8,-0.8) node{$\FKhint =  \{\sigma_i\}^{14}_{i=8}$};
        \draw (7.9,-1.4) node{$\FKh = \FKhext  \cup \FKhint $};
\end{tikzpicture}
\hspace{0.1cm}
\begin{tikzpicture}
       \draw[line width=\sz pt] (5,-0.7) -- (5,1)  ;
         \draw[dashed, line width=\sz pt] (2.6,0.5) -- (5,-0.7) ;
         \draw[line width=\sz pt] (4.3,0.) node[left]{$\sigma_8$};
         \draw (4.7,0.2) node[left]{$\kappa_1$};
        \draw[line width=\sz pt] (5,1) -- (2.5,2)  ;
         \draw[dashed, line width=\sz pt] (2.6,0.5) -- (5,1) ;
         \draw (3.5,1.1) node[below]{$\sigma_9$};
         \draw (3.4,1.4) node[below]{$\kappa_2$};
        \draw[line width=\sz pt] (2.5,2) -- (0.5,1.5)  ;
         \draw[dashed, line width=\sz pt] (2.6,0.5) -- (2.5,2) ;
         \draw (2.3,1.4) node[below]{$\sigma_{10}$};
          \draw (1.8,1.7) node[below]{$\kappa_3$};
       \draw[line width=\sz pt] (0.5,1.5) -- (0.,0)  ;
         \draw [dashed, line width=\sz pt](2.6,0.5) -- (0.5,1.5) ;
         \draw (1.4,1.1) node[below]{$\sigma_{11}$};
         \draw (0.4,0.6) node[right]{$\kappa_4$};
        \draw[line width=\sz pt] (0,0) -- (2,-0.5) ;
          \draw [dashed, line width=\sz pt](2.6,0.5) -- (0,0) ;
          \draw (1.6,0.3) node[below]{$\sigma_{12}$};
        \draw (1.5,-0.4) node[above]{$\kappa_5$};
        \draw[line width=\sz pt] (2,-0.5) -- (3.5,-1.5) ;
         \draw [dashed, line width=\sz pt](2.6,0.5) -- (2,-0.5) ;
          \draw (2.5,0.) node[below]{$\sigma_{13}$};
         \draw (2.3,-0.7) node[right]{$\kappa_6$};
        \draw[line width=\sz pt] (3.5,-1.5) -- (5,-0.7)  ;
         \draw[dashed, line width=\sz pt] (2.6,0.5) -- (3.5,-1.5) ;
           \draw (3.3,-0.2) node[below]{$\sigma_{14}$};
        \draw (3.8,-1.) node[above]{$\kappa_7$};
\end{tikzpicture}
\caption{Example of a polygonal element $\elm$ with its faces $\FK$,
corresponding face fluxes $\mathsf{U}^{\mathrm{ext}}_K$, and pressure head $(\mathsf{P})_K$ (left);
virtual simplicial submesh $\TK$ of $\elm$ (right)}
\label{fig:notations}
\end{figure}
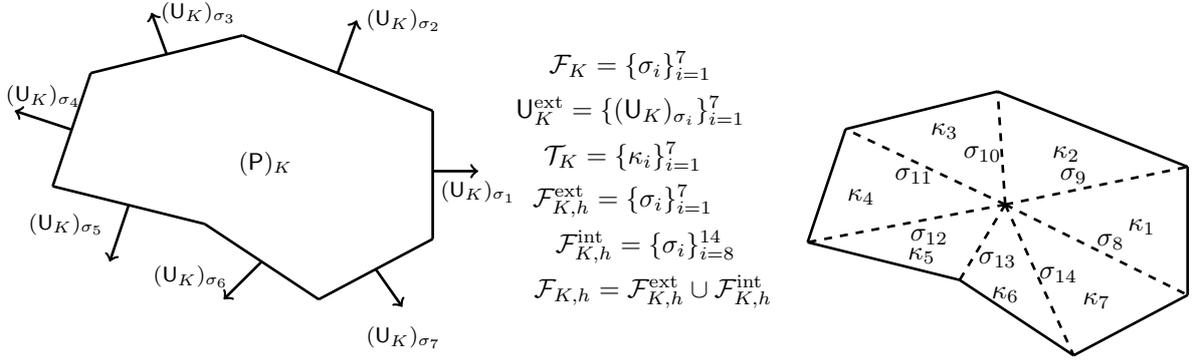

\subsubsection{Mesh faces and mesh vertices}

Let $\F_H$ be the set of the $(d-1)$-dimensional faces of $\T_H$. We divide
it into interior faces $\F_H^{\rm int}$ and boundary faces $\F_H^{\rm ext}$:
$\sd \in \F_H^{\rm int}$ if there exist $\elm, \elmt \in \T_H$, $\elm \neq
\elmt$, such that $\sd \subset \elm \cap \elmt$, and $\sd \in \F_H^{\rm ext}$ if
there exist $\elm \in \T_H$ such that $\sd \subset \elm \cap \pt \Om$. 
Recall from Section~\ref{sec_Th} that $\Fh$ is the set of the ($d-1$)-dimensional faces of the simplicial mesh $\Th$. We denote by $\F_{H,h}$ the set of such faces from $\Fh$ that lie in
some polytopal face of $\F_H$.  
Let $\FK {\subset} \F_H$ be the set
of the $(d-1)$-dimensional faces of the polytope $\elm \in \T_H$, and let
$\FKhext$ collect those faces of $\F_{H,h}$ that lie on the boundary of the
element $\elm \in \T_H$, whereas $\FKhint$ those faces of $\Fh$ that lie
inside the element $\elm \in \T_H$. \cor{We also set $\FKh \eq \FKhext \cup \FKhint$} see Figure~\ref{fig:notations}. 
If the faces of $\Th$ do not subdivide the faces of $\T_H$, then $\F_{H,h} = \F_H$ \cor{and also
$\FKhext = \FK$}.
For an interior face $\sd \in \F_H^{\rm int}$, we fix an arbitrary
orientation and denote the corresponding unit normal vector by $\tn_\sd$. For
a boundary face $\sd \in \F_H^{\rm ext}$, $\tn_\sd$ coincides with the
exterior unit normal $\tn_\Omega$ of $\Omega$. Finally, for all $\elm \in
\T_H$, we denote by $\tn_\elm$ the unit normal vector to $\sd$ pointing out of $\elm$.
To complete the notation, we let 
$\VKh$ to be the set of vertices of all the elements $\kappa \in \T_\elm$\cor{, and, as for simplicial meshes, $\Ve$ the set of vertices of the given face $\sd$}.

\subsection{Basic principle of finite volume and related lowest-order locally conservative methods} \label{sec_FV}

We now introduce the basic principle of the finite volume and related lowest-order locally conservative methods considered in this contribution. 
In Sections~\ref{sec_Pois}--\ref{sec_MP_MC} below, for each considered partial differential equation, a general lowest-order locally conservative discretization is introduced in detail.

Consider a steady diffusion partial differential equation of the form~\eqref{eq_PDE}--\eqref{eq_Lapl_op} and recall~\eqref{eq_Lapl_prop}. A general lowest-order locally conservative method looks for real values $U_{\elm,\sd}$ which for each polytopal mesh element $\elm \in \T_H$ and each of its faces $\sd \in \FK$ approximate the normal outflux $\<\tu \scp \tn_\elm, 1\>_\sd \cor{= \int_\sd \tu \scp \tn_\elm}$ from $\elm$ through the face $\sd$ such that
\bse \label{eq_FV} \begin{equation} \label{eq_FV_scheme_princ}
    \sum_{\sd \in \FK} U_{\elm,\sd} = (f,1)_\elm \qquad \forall \elm \in \T_H.
\end{equation}
Here the sum runs through all the faces of the polytope $\elm$ and, recalling \cor{the notation}~\eqref{eq_not_scal}, $(f,1)_\elm = \int_\elm f(\tx) \dx$ is the integral of the source term. From a neighboring element $\elmt \in \T_H$ sharing the same face $\sd$, we then impose, 
\begin{equation} \label{eq_FV_normal_cont}
    U_{\elmt,\sd} = - U_{\elm,\sd}.
\end{equation}\ese
Requirement~\eqref{eq_FV_normal_cont} \cor{ensures the} normal-trace continuity, as in the Raviart--Thomas space $\RT_0(\Th) \cap \Hdv$ of~\eqref{eq_RTN}. 
Then~\eqref{eq_FV_normal_cont} expresses the conservation of mass across mesh faces and~\eqref{eq_FV_scheme_princ} the equilibrium with the load, respectively the discrete versions of $\tu \in \Hdv$ and $\Dv \tu = f$. \cor{Mimicking the} remaining items from~\eqref{eq_Lapl_prop}, namely $p \in \Hoo$ and $\tu = - \Gr p$, then have to be built in the construction of the discrete normal fluxes $U_{\elm,\sd}$ from elementwise approximations of the primal variable $p_\elm$.

For an unsteady problem of the form~\eqref{eq_PDE_unst}--\eqref{eq_nonl_heat_op}, we consider an implicit Euler time discretization on discrete times $0=t^0 < \ldots < t^n < \ldots t^N = t_{\rm{F}}$ with time steps $\tau^n \eq t^n - t^{n-1}$. We then get, for each discrete time $t^n$, $1\le n\le N$,
\bse \label{eq_FV_unst} \begin{equation} \label{eq_FV_scheme_princ_unst}
    \frac{p_\elm^n - p_\elm^{n-1}}{\tau^n} + \sum_{\sd \in \FK} U_{\elm,\sd}^n = (f^n,1)_\elm \qquad \forall \elm \in \T_H
\end{equation}
and
\begin{equation} \label{eq_FV_normal_cont_unst}
    U_{\elmt,\sd}^n = - U_{\elm,\sd}^n,
\end{equation}\ese
where $f^n \eq \int_{t^{n-1}}^{t^n} f / \tau^n$. Again, the heart of each method is the construction of the discrete normal fluxes $U_{\elm,\sd}^n$ for each \cor{polytopal mesh element} $\elm \in \T_H$ and \cor{each face} $\sd \in \FK$ from the elementwise approximations of the primal variable $p_\elm^n$.

\subsection{Bibliographic resources}\label{sec_biblio_setting}

Details on the material rapidly exposed in Sections~\ref{sec_spaces}--\ref{sec_meshes} 
can be found in Adams~\cite{Adams_75}, Thomas~\cite{Thomas_dis_77}, Ciarlet~\cite{Ciar_78}, Allaire~\cite{Allaire_05}, Brenner and Scott~\cite{Bren_Scott_FEs_08}, Boffi~\eal\ \cite{Bof_Brez_For_MFEs_13}, Ern and Guermond~\cite{Ern_Guermond_FEs_I_21}, and Vohral\'ik~\cite{Voh_a_post_LN_24, Voh_FEM_LN_25}.

Finite volume schemes like~\eqref{eq_FV} and~\eqref{eq_FV_unst} are very natural and used traditionally in the engineering practice. Amongst first mathematical analyses for diffusion problems, one can cite Manteuffel and White~\cite{Mant_White_FVs_86}, Rose~\cite{Rose_FV_diff_89}, Morton and S\"{u}li~\cite{Mort_Suli_FVs_91}, Faille~\cite{Faille_FV9_92}, Herbin~\cite{Herb_FV_ADR_95}, and the overview in Eymard~\eal\ \cite{Eym_Gal_Her_00}.
For the closely related mixed finite elements, the construction of the spaces in the lowest-order case can be traced back to Whitney~\cite{Whit_geom_int_57}, whereas higher-order extensions and analysis have been undertaken in Raviart and Thomas~\cite{Ra_Tho_MFE_77} and N\'ed\'elec~\cite{Ned_mix_R_3_80}, see Brezzi and Fortin~\cite{Brez_For_91} for an early overview. The ample connections between finite volumes and mixed finite elements are discussed in, \eg, Russell and Wheeler~\cite{Rus_Whe_MFE_FD_83}, Agouzal~\eal\ \cite{Agou_Bar_Mai_Oud_FV_MFE_95}, Baranger~\eal\ \cite{Bar_Mai_Oud_FV_MFE_96}, Youn\`es~\eal\ \cite{You_Mos_Ack_Chav_MFE_FV_99, You_Ack_Chav_MFE_FV_3D_04}, Chavent~\eal\ \cite{Chav_You_Ack_MFE_FV_03}, Klausen and Russell~\cite{Klaus_Rus_rel_loc_cons_04}, and Vohral{\'{\i}}k~\cite{Voh_eq_MFE_FV_06}.

Numerous extensions of the initial approaches together with generalizations to polytopal meshes, including a priori analysis (existence, uniqueness, convergence, a priori error estimates), have been undertaken recently. 
Let us cite multi-point finite volumes by Aavatsmark~\eal\ \cite{Aava_98_I, Aava_Eig_Klau_Whee_Yot_cvg_MPFA_07}, Edwards~\cite{Edw_CVFE_02}, and Breil and Maire~\cite{Bre_Mai_CCFV_sym_06}, mimetic finite differences by Brezzi~\eal\ \cite{Brez_Lip_Shash_MFD_cvg_05,
Brez_Lip_Sim_MFD_05}, mixed and hybrid finite volumes by Eymard~\eal\ \cite{Dro_Ey_MFV_06, Eym_Gal_Her_SUSHI_10}, multipoint flux mixed finite elements by Wheeler and Yotov~\eal\ \cite{Wh_Yot_MPFMFE_06, Whee_Xue_Yot_MPMFE_3D_12}, compatible discrete operator schemes by Bonelle and Ern~\cite{Bon_Ern_op_pol_14}, mixed virtual elements by Brezzi~\eal\ \cite{Brez_Falk_Mar_MVEM_14} and Beir\~ao da Veiga~\eal\ \cite{Bei_Brez_Mar_Rus_MVEM_16}, hybrid high-order methods by Di Pietro, Ern, and Lemaire~\cite{DiPiet_Ern_Lem_HHO_14, Di_Pietr_Ern_HHO_el_15}, or mixed finite elements in~\cite{Kuzn_Rep_MFE_pol_03, Kuzn_MFE_pol_08, Sbou_Jaff_Rob_comp_MFE_09, Christ_MFE_pol_10, Voh_Wohl_MFE_1_unkn_el_rel_13}, see also~\cite{Boch_Hym_MFD_CDO_06, Her_Hub_FVCA_benchm_08, Christ_CDO_08, Dro_rev_14}; unifying frameworks can be found in Droniou~\eal\ \cite{Dro_Ey_Gal_Her_un_10, Dro_Eym_Her_grad_sch_16}, Vohral{\'{\i}}k and Wohlmuth~\cite{Voh_Wohl_MFE_1_unkn_el_rel_13}, Cockburn~\eal\ \cite{Cock_DiPi_ern_HHO_HDG_16}, Boffi and Di Pietro~\cite{Boff_DiPietr_unif_mixed_hybrid_polyt_18}, and the references therein. Amongst recent books, let us cite Jovanovi\'{c} and S\"{u}li~\cite{Jov_Suli_FDs_14} for finite differences, Beir\~{a}o da Veiga~\eal\ \cite{Bei_Lip_Manz_MFD_14} for mimetic finite differences, Droniou~\eal\ \cite{Dro_Eym_Gal_Guich_Her_grad_sch_18} for a unifying approach, and Nordbotten and Keilegavlen~\cite{Nord_Keil_MPFA_PM_21} for multi-point finite volumes.

\section{The Poisson equation. Simplicial meshes and spatial discretization error} \label{sec_Pois}

In this section, we start with the simplest partial differential equation for diffusion processes: the Poisson equation. We first recall the definition of a weak solution based on the Sobolev space $\Hoo$ of Section~\ref{sec_H1} and show that the associated flux lies in the Sobolev space $\Hdv$ of Section~\ref{sec_Hdv}. This is very instrumental and will be repeated in all generalizations of Sections~\ref{sec_Darcy}--\ref{sec_MP_MC} below. We then define the energy error and explore its structure for a completely generic approximate solution. We subsequently consider a basic cell-centered finite volume discretization incarnating the principles of Section~\ref{sec_FV}. We only employ in this section a simplicial mesh as per Section~\ref{sec_Th} and define the approximate solution by lifting the face fluxes into a Raviart--Thomas piecewise polynomial space of Section~\ref{sec_RT}. Next, we introduce and discuss two central notions of this work: the flux and potential reconstructions. These in particular lead to a posteriori error estimates fulfilling all the building principles of Section~\ref{sec_a_post}. Illustrative numerical experiments close this section.

\subsection{The Poisson equation}

In this section, for the domain $\Om$ as specified in Section~\ref{sec_dom_Om}, we consider the problem of finding $p: \Om \ra \RR$ such that
\bse \label{eq_Pois} \begin{empheq}[box=\widefbox]{align}
    \ds - \Dv(\Gr p) & = f \qquad \mbox{ in } \, \Om, \label{eq_Pois_eq} \\
    \ds p & = 0 \qquad \mbox{ on } \, \pt \Om. \label{eq_Pois_BC}
\end{empheq} \ese
Here $f \in \Lt$ is a source term. 

\subsection{Weak solution, potential, and flux} \label{sec_Pois_WF}

The Sobolev space $\Hoo$ from Definition~\ref{def_Hoo} gives a proper mathematical setting to define a unique solution $p$ of~\eqref{eq_Pois}. Indeed, there is one and only one $p \in \Hoo$ such that
\be \label{eq_Pois_WF}
    (\Gr p, \Gr v) = (f, v) \qquad \forall v \in \Hoo;
\ee
note that this follows trivially from the Riesz representation theorem, since $\Hoo$ is a Hilbert space for the scalar product given by $(\Gr p, \Gr v)$ and since $(f, v)$ is bounded linear form over $\Hoo$. We call~\eqref{eq_Pois_WF} the weak formulation of~\eqref{eq_Pois}. From the weak solution $p$, we can define the dual variable
\be \label{eq_velocity_Pois}
    \tu \eq - \Gr p.
\ee
In relation to the applications below, $p$ may be a pressure and $\tu$ a Darcy velocity, but we rather generically refer to $p$ as {\em potential} and to $\tu$ as {\em flux}.

\subsection{Properties of the exact potential and flux} \label{sec_prop_WS}

Using Definitions~\ref{def_Hoo} and~\ref{def_Hdv} of the spaces $\Hoo$ and $\Hdv$, we have the following important result, detailing what we have discovered in~\eqref{eq_Lapl_prop}:

\bpr[Properties of the weak solution~\eqref{eq_Pois_WF}] \label{pr_prop_WS} Let
$p$ be the solution of~\eqref{eq_Pois_WF}. Let $\tu$ be given by~\eqref{eq_velocity_Pois}. Then
\be \label{eq_prop_WS}
    p \in \Hoo, \quad \tu \in \Hdv, \quad \Dv \tu = f.
\ee
\epr

\bp The weak solution $p$ belongs to $\Hoo$ by definition. In order
to verify that $\tu \in \Hdv$ and $\Dv \tu = f$, we need to check the three conditions
of Definition~\ref{def_WDG}. Condition~\ref{prop_L2_div} is obvious,
as $p \in \Hoo$ and thus $- \Gr p = \tu$ is square-integrable \cor{in view of Definition~\ref{def_WPD}, condition~\ref{prop_L2_der}}.
For the function $w$ of condition~\ref{prop_L2_dv} in Definition~\ref{def_WDG}, choose $w \eq f$
and note that $f \in \Lt$ by assumption. Then
condition~\ref{prop_Green_dv} follows immediately
from~\eqref{eq_Pois_WF} and the fact that $\Do \subset
\Hoo$. \ep

\br[Mathematical properties of the weak solution] It is remarkable that the flux $\tu$ prescribed from the potential $p$ by~\eqref{eq_velocity_Pois} belongs to $\Hdv$ and has a divergence in equilibrium with the load $f$, since this is not explicitly stipulated in the weak formulation~\eqref{eq_Pois_WF}. Congruently, it immediately follows \cor{(\cf\ Proposition~\ref{thm_Prag_Syng} below)} that the three properties~\eqref{eq_prop_WS} together with the constitutive relation~\eqref{eq_velocity_Pois} actually define the weak solution $p$ of~\eqref{eq_Pois_WF} in a unique way. This mathematically describes the intrinsic properties of the exact potential $p$ and flux $\tu$: $p$ has to be trace continuous and $\tu$ has to be related by the constitutive relation~\eqref{eq_velocity_Pois} to it, normal-trace continuous, and in equilibrium with the load. \er

\br[Physical properties of the weak solution] \label{rem_p_u} Scalar-valued physical variables such as the potential (pressure) $p$ are naturally trace continuous, whereas vector-valued physical variables such as the flux (Darcy velocity) $\tu = - \Gr p$ are naturally normal-trace continuous, \cf\ Figures~\ref{fig_Ho} and~\ref{fig_H_div}. Moreover, $\tu$ has to be in equilibrium with the load (locally mass conservative) which amounts to
$\Dv \tu = f$. This restates that the three properties~\eqref{eq_prop_WS} together with the constitutive relation~\eqref{eq_velocity_Pois} characterize the exact solution. \er

\subsection{Approximate solution and spatial discretization error}\label{sec_err_Poiss}

Let $\tu_h \in \tLt$ be an arbitrary function that we think of a (numerical) approximation of the exact flux $\tu$ from~\eqref{eq_Pois_WF}--\eqref{eq_velocity_Pois}. We will give a specific example for the cell-centered finite volume discretizations on a simplicial mesh in Section~\ref{sec_fl_rec_Pois} below, but we want to proceed abstractly for the moment, since there is no link to a mesh or to the discrete world of piecewise polynomials at this stage. We will call the $\tLt$-distance of $\tu_h$ to $\tu$,
\be \label{eq_sp_dis_err}
    \norm{\tu - \tu_h},
\ee
the {\em spatial discretization error}.

\subsection{Prager--Synge equality}

We now state the following important result dating back to Prager and Synge~\cite{Prag_Syng_47}:

\bpr[Prager--Synge equality] \label{thm_Prag_Syng} 
Let $p \in \Hoo$ be the weak solution of~\eqref{eq_Pois_WF} and let $\prh \in \Hoo$ and $\frh \in \Hdv$ with $\Dv \frh = f$ be arbitrary. Then
\be \label{eq_Prag_Syng}
    \norm{\Gr(p - \prh)}^2 + \norm{\Gr p + \frh}^2 = \norm{\Gr \prh +
    \frh}^2.
\ee
\epr

\bp Adding and subtracting $\Gr p$, we develop
\ban
    \norm{\Gr \prh + \frh}^2 {} & = \norm{\Gr (\prh - p) + \Gr p  + \frh}^2 \\
    & {} = \norm{\Gr (\prh - p)}^2 + \norm{\Gr p + \frh}^2
    + 2(\Gr (\prh - p), \Gr p  + \frh).
\ean
Note from Proposition~\ref{pr_prop_WS} that $\Gr p \in \Hdv$ with $\Dv
(\Gr p ) = - f$. Thus $(\Gr p  + \frh) \in \Hdv$ and in
particular $\Dv (\Gr p  + \frh) = 0$. Thus, using that $\prh - p
\in \Hoo$, the Green theorem~\eqref{eq_Green_Hdv} gives
\[
    (\underbrace{\Gr p + \frh}_{\in \Hdv}, \Gr (\underbrace{\prh - p}_{\in \Hoo})) \reff{eq_Green_Hdv}= - (\underbrace{\Dv (\Gr p  + \frh)}_{=0}, \prh - p) = 0,
\]
whence the assertion follows. \ep

\br[Prager--Synge equality~\eqref{eq_Prag_Syng}] We stress that~\eqref{eq_Prag_Syng} is an equality, where on the right-hand side, there are, in practice, two known, discrete objects from the piecewise polynomial subspaces of $\Hoo$ and $\Hdv$, whereas the left-hand side are their spatial discretization errors with respect to the weak solution of~\eqref{eq_Pois_WF}. In other words, if we want to trace a simultaneous potential and flux error, then we can {\em compute it}, there is no (a posteriori) estimate. \er

From~\eqref{eq_Prag_Syng}, we immediately observe (\cf\ Theorem~6.1 and Corollary~6.6 in~\cite{Voh_un_apr_apost_MFE_10} and the references therein):

\bc[Prager--Synge (in)equality] \label{cor_Prag_Syng} Let $p \in \Hoo$ be the solution of~\eqref{eq_Pois_WF}, let $\tu = - \Gr p$ from~\eqref{eq_velocity_Pois}, and let $\tu_h \in \Hdv$ with $\Dv \tu_h = f$ be arbitrary. Then
\be \label{eq_Prag_Syng_cons_1}
    \norm{\tu - \tu_h} = \min_{v \in \Hoo} \norm{\tu_h + \Gr v}\cor{,}
\ee
\cor{so that in particular}
\be \label{eq_Prag_Syng_cons_2}
    \norm{\tu - \tu_h} \leq \norm{\tu_h + \Gr \prh}
\ee
for an arbitrary $\prh \in \Hoo$. \ec

\bp The second assertion~\eqref{eq_Prag_Syng_cons_2} follows from~\eqref{eq_Prag_Syng} upon taking $\frh = \tu_h$ and noticing that
\[
    \norm{\Gr p + \tu_h}^2 = \norm{\Gr \prh +
    \tu_h}^2 - \norm{\Gr(p - \prh)}^2 \leq \norm{\Gr \prh +
    \tu_h}^2
\]
for any $\prh \in \Hoo$, using $\tu = - \Gr p$ from~\eqref{eq_velocity_Pois}. This immediately gives~\eqref{eq_Prag_Syng_cons_1} with infimum in place of minimum and the $\leq$ sign in place of the equality. Since we can take $\prh = p$ and since $\tu = - \Gr p$, the minimum in~\eqref{eq_Prag_Syng_cons_1} is actually attained and the right-hand side equals the left-hand side. \ep

\br[Prager--Synge equality~\eqref{eq_Prag_Syng_cons_1}] \label{rem_PS} According to~\eqref{eq_velocity_Pois} and~\eqref{eq_prop_WS}, to be equal to \cor{the exact flux} $-\Gr p$, \cor{any} $\tu_h \in \Hdv$ with $\Dv \tu_h = f$ only misses being a minus weak gradient of some function from $\Hoo$. This is the meaning of the characterization~\eqref{eq_Prag_Syng_cons_1}. Notice that~\eqref{eq_Prag_Syng_cons_1} expresses the distance of $\tu_h$ to $\Gr \Hoo$. \er 

\br[Prager--Synge inequality~\eqref{eq_Prag_Syng_cons_2}] Having $\tu_h \in \Hdv$ with $\Dv \tu_h = f$ at disposal, which we will see in Section~\ref{sec_fl_rec_Pois} below is immediate in finite volume methods for the Poisson equation~\eqref{eq_Pois}, when $f$ is piecewise constant, the Prager--Synge inequality~\eqref{eq_Prag_Syng_cons_2} can be readily used to produce a guaranteed error upper bound also satisfying the other requirements formulated in Section~\ref{sec_a_post_props}. 
This will still be the situation for the singlephase steady linear Darcy flow in Section~\ref{sec_Darcy} below, but when we want to take into account iterative algebraic solvers in Section~\ref{sec_Darcy_NL} and component fluxes together with nonlinear zero-order terms in Section~\ref{sec_MP_MC} for nonlinear compositional Darcy flows, it may not be feasible to locally construct $\tu_h \in \Hdv$ with $\Dv \tu_h = f$ as an approximate solution. There, we will in general only have $\tu_h \in \tLt$. \cor{In this case, we will reconstruct $\bsig_h^{k,i} \in \Hdv$ with $\Dv \bsig_h^{k,i} = f$ or $\Dv \bsig_h^{k,i} = f + \rho_h^{k,i}$ with $\rho_h^{k,i}$ a remainder made as small as necessary.}\er

\subsection{Equivalence of the spatial discretization error to the flux nonconformity (dual norm of the residual) plus potential nonconformity}

Following~\cite[Theorem~3.1]{Pench_Voh_Whee_Wild_a_post_MS_MN_M_13}, \cite[Theorem 3.3]{Ern_Voh_p_rob_15}, and the references therein, we now extend Corollary~\ref{cor_Prag_Syng} to the general case $\tu_h \in \tLt$.

\bt[Error characterization for $\tu_h \in \tLt$] \label{thm_err_char} Let $p \in \Hoo$ be the weak solution of~\eqref{eq_Pois_WF}, let $\tu = - \Gr p$ from~\eqref{eq_velocity_Pois}, and let $\tu_h \in \tLt$ be arbitrary. Then
\be \label{eq_err_char}
    \norm{\tu - \tu_h}^2 = \underbrace{\min_{\substack{\tv \in \Hdv\\\Dv \tv = f}} \norm{\tu_h - \tv}^2}_{\substack{\max_{\vf \in \Hoo; \, \norm{\Gr \vf} = 1} \{(f,\vf) + (\tu_h, \Gr \vf)\}^2,\\ \text{dual norm of the residual of } \tu_h}} + \underbrace{\min_{v \in \Hoo} \norm{\tu_h + \Gr v}^2}_{\substack{\text{nonconformity}\\ (\text{distance of } \tu_h \text{ to } \Gr \Hoo)}}.
\ee
\et

\br[Error characterization~\eqref{eq_err_char}] Relation~\eqref{eq_err_char} extends~\eqref{eq_Prag_Syng_cons_1} from $\tu_h \in \Hdv$ with $\Dv \tu_h = f$ to any $\tu_h \in \tLt$. The supplementary arising term is naturally the distance of $\tu_h$ to $\Hdv$ under the constraint of the divergence being equal to $f$. \er

\br[Dual norm of the residual] The first term in~\eqref{eq_err_char} has an equivalent expression investigating ``how much, for a test function $\vf \in \Hoo$, $\tu_h$ misses to fulfill the weak formulation~\eqref{eq_Pois_WF} \cor{when put in place of $- \Gr p$}'', which is related to the {\em residual} of~\eqref{eq_Pois_WF} for $\vf$. The overall maximum \cor{over all $\vf$ from $\Hoo$ with $\norm{\Gr \vf} = 1$} is then the {\em dual norm} of the residual. \er

\br[Potential nonconformity]The second term in~\eqref{eq_err_char} expresses the distance of $\tu_h$ to $\Gr \Hoo$, \cf\ Remark~\ref{rem_PS}.\er

\bp[Proof of Theorem~\ref{thm_err_char}] Let us define a function $\pr \in \Hoo$
by
\be \label{eq_s}
    (\Gr \pr, \Gr v) = - (\tu_h, \Gr v) \qquad \forall v \in \Hoo.
\ee
There exists one and only one $\pr$ by the Riesz representation theorem; indeed, recall from Section~\ref{sec_H1} that the left-hand side of~\eqref{eq_s} is a scalar product on $\Hoo$ and notice that from $|- (\tu_h, \Gr v)| \leq \norm{\tu_h}\norm{\Gr v}$, the right-hand side of~\eqref{eq_s} is a bounded linear form on $\Hoo$ \cor{($\tu_h$ is here a fixed datum)}. The function $\pr$ can be seen as the orthogonal projection of the approximate solution $\tu_h$ onto $\Gr \Hoo$. With the aid of $\pr$, we can thus write the Pythagorean equality
\be \label{eq_Pyth}
    \norm{\Gr p + \tu_h}^2 = \norm{\Gr(p - \pr)}^2 + \norm{\Gr \pr + \tu_h}^2.
\ee
Indeed,
\[
    \norm{\Gr p + \tu_h}^2 = \norm{\Gr p - \Gr \pr + \Gr \pr + \tu_h}^2 =
    \norm{\Gr(p - \pr)}^2 + \norm{\Gr \pr + \tu_h}^2 + 2(\Gr(p - \pr), \Gr \pr + \tu_h),
\]
and the last term in the above expression vanishes in view of the orthogonality~\eqref{eq_s}, since $p - \pr$ can be taken as a test function $v \in \Hoo$ in~\eqref{eq_s}. We continue in three steps.

1) Since $\pr$ can be seen as a projection of $\tu_h$,
\be \label{eq_NC_min}
    \norm{\Gr \pr + \tu_h}^2 = \min_{v \in \Hoo} \norm{\Gr v + \tu_h}^2.
\ee
Indeed, from~\eqref{eq_Pyth} used for \cor{an arbitrary} function $v \in \Hoo$ in place of $p \in \Hoo$, we see
\[
    \norm{\Gr v + \tu_h}^2 = \norm{\Gr(v - \pr)}^2 + \norm{\Gr \pr + \tu_h}^2.
\]
Therefrom, we get
\[
    \norm{\Gr \pr + \tu_h}^2 = \norm{\Gr v + \tu_h}^2 - \norm{\Gr(v - \pr)}^2 \leq \norm{\Gr v + \tu_h}^2 \qquad \forall v \in \Hoo\cor{,}
\]
\cor{and there is an equality for $v = \pr$.} 
This \cor{shows that} the second term in~\eqref{eq_Pyth} \cor{indeed takes} the form needed in~\eqref{eq_err_char}.

2) For the first term in~\eqref{eq_Pyth}, we first use that $p - \pr \in \Hoo$. Thus, the dual norm characterization and~\eqref{eq_s} give
\be \label{eq_1}
    \norm{\Gr(p - \pr)} = \max_{\substack{\vf \in \Hoo\\ \norm{\Gr \vf} = 1}} (\Gr (p - \pr), \Gr \vf) \reff{eq_s}= \max_{\substack{\vf \in \Hoo\\ \norm{\Gr \vf} = 1}} (\Gr p + \tu_h, \Gr \vf).
\ee
Let now $\vf \in \Hoo$ with $\norm{\Gr \vf} = 1$ be fixed. Using the
characterization~\eqref{eq_Pois_WF} of the weak solution, we have
\be \label{eq_2}
    (\Gr p + \tu_h,\Gr\vf) \reff{eq_Pois_WF}= (f, \vf) + (\tu_h,\Gr\vf).
\ee
Thus, we have~\eqref{eq_err_char} with the first term in the form expressed in the underbrace. 

3) We are left to show the equivalence
\be \label{eq_dist_H_div_res}
    \min_{\substack{\tv \in \Hdv\\\Dv \tv = f}} \norm{\tu_h - \tv} = \max_{\vf \in \Hoo; \, \norm{\Gr \vf} = 1} \{(f,\vf) + (\tu_h, \Gr \vf)\}.
\ee

3a) First, to motivate, we show that the right-hand side of~\eqref{eq_dist_H_div_res} is bounded by the left-hand side one. Let $\tv \in \Hdv$ such that $\Dv \tv = f$ be arbitrary. Then the Green theorem~\eqref{eq_Green_Hdv} gives
\[
    (f, \vf) + (\tu_h,\Gr\vf) = (\Dv \tv, \vf) + (\tu_h,\Gr\vf) = (\tu_h - \tv, \Gr\vf).
\]
Consequently, by the Cauchy--Schwarz inequality,
\be \label{eq_flux_min_est}
    \max_{\vf \in \Hoo; \, \norm{\Gr \vf} = 1} \{(f,\vf) + (\tu_h, \Gr \vf)\} \leq \min_{\substack{\tv \in \Hdv\\\Dv \tv = f}} \norm{\tu_h - \tv}.
\ee

3b) We now show the equality~\eqref{eq_dist_H_div_res}. The argument of the \cor{constrained minimization} in~\eqref{eq_dist_H_div_res} is
\[
    \bsig \eq \arg\min_{\substack{\tv \in \Hdv\\\Dv \tv = f}} \norm{\tu_h - \tv}
\]
and is characterized by the Euler--Lagrange conditions as a function $\bsig \in \Hdv$ with $\Dv \bsig = f$ such that
\[
    (\bsig, \tv) = (\tu_h, \tv) \qquad \forall \tv \in \Hdv \text{ with } \Dv \tv = 0.
\]
Imposing the constraint $\Dv \bsig = f$ equivalently with test functions $q \in \Lt$ and introducing the Lagrange multiplier $r \in \Lt$, using that $\Dv \Hdv = \Lt$, the above problem is equivalent to finding $\bsig \in \Hdv$ and $r \in \Lt$ such that
\bse \bat{2}
    (\bsig, \tv) - (r, \Dv \tv) & = (\tu_h, \tv) \qquad & & \forall \tv \in \Hdv, \label{eq_Dual_1}\\
    (\Dv \bsig, q) & = (f,q) && \forall q \in \Lt. \label{eq_Dual_2}
\eat \ese
Now, \eqref{eq_Dual_1} implies by Definition~\ref{def_WPD} of the weak partial derivative that $r \in \Hoo$ with $\Gr r = - (\bsig - \tu_h)$. Consequently, by the Green theorem~\eqref{eq_Green_Hdv},
\be \label{eq_char} \bs
        \min_{\substack{\tv \in \Hdv\\\Dv \tv = f}} \norm{\tu_h - \tv} & = \norm{\tu_h - \bsig} = \norm{\Gr r} = \max_{\substack{\vf \in \Hoo\\ \norm{\Gr \vf} = 1}} (\Gr r, \Gr \vf)\\
        & = \max_{\substack{\vf \in \Hoo\\ \norm{\Gr \vf} = 1}} (- \bsig + \tu_h, \Gr \vf) = \max_{\substack{\vf \in \Hoo\\ \norm{\Gr \vf} = 1}} \{(f, \vf) + (\tu_h,\Gr\vf)\},
\es \ee
which is~\eqref{eq_dist_H_div_res}. \ep

\subsection{Cell-centered finite volume discretizations on a simplicial mesh} \label{sec_FV_Pois}

\begin{figure}
\centerline{\includegraphics[width=0.5\textwidth]{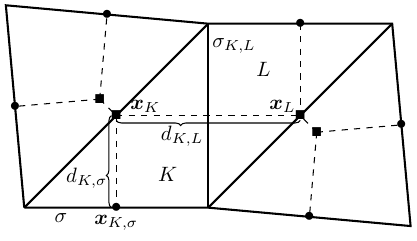}}
\caption{Admissible triangular mesh $\Th$ and notation for the cell-centered finite volume scheme}\label{fig_FV_mesh}
\end{figure}

Consider now a simplicial mesh $\Th$ of the domain $\Om$ as defined in Section~\ref{sec_Th}. We additionally assume that 1) there exists a point $\tx_\elm$ associated with each simplex $\elm \in \Th$; 2) there exists a point $\tx_{\elm, \sd}$ associated with each face $\sd$ lying in the boundary $\pt \Om$, $\sd \in \FK \cap \Fhext$, $\elm \in \Th$, lying in the interior of the face $\sd$; 3) all the points $\tx_\elm$ and $\tx_{\elm, \sd}$ are distinct; 4) the straight line connecting $\tx_\elm$ and $\tx_\elmt$ for two neighboring elements $\elm$ and $\elmt$ is orthogonal to their common face $\sd_{\elm,\elmt} = \pt \elm \cap \pt \elmt \in \Fhint$; 5) the straight line connecting $\tx_\elm$ and $\tx_{\elm, \sd}$ is orthogonal to $\sd$ for $\sd \in \FK \cap \Fhext$. The points $\tx_\elm$ typically, but not necessarily, lie in the interior of the element $\elm \in \Th$. In particular, Delaunay meshes satisfy this requirement of ``admissibility'' in the sense of~\cite[Definition~9.1]{Eym_Gal_Her_00}. An illustration is provided in Figure~\ref{fig_FV_mesh}.

The cell-centered finite volume scheme for the Poisson equation~\eqref{eq_Pois} on the simplicial mesh $\Th$ reads: find the real values $p_\elm$, $\elm \in \Th$, the approximations to the mean values of $p$ in the mesh elements $\elm$, such that
\bse\label{eq_FVs} \begin{equation} \label{eq_FV_scheme}
    \boxed{\sum_{\sd \in \FK} U_{\elm,\sd} = (f,1)_\elm \qquad \forall \elm \in \Th,}
\end{equation}
where $U_{\elm,\sd} \in \RR$ for each face $\sd \in \FK$ approximates the normal (out)flux $\<\tu \scp \tn_\elm, 1\>_\sd$ from $\elm$ over the face $\sd$ \cor{(recall the notation~\eqref{eq_RT0_DoF})} by
\begin{alignat}{2}
    \label{eq_dif_fl_int}
        U_{\elm,\sd} & \eq - \frac{|\sd_{\elm,\elmt}|}{d_{\elm,\elmt}}(p_\elmt - p_\elm) \qquad & & 
        \sd = \sd_{\elm,\elmt} \in \Fhint, \\
    \label{eq_dif_fl_Dir}
        U_{\elm,\sd} & \eq - \frac{|\sd|}{d_{\elm,\sd}}(0-p_\elm) & & 
        \sd \in \FK \cap \Fhext.
\end{alignat}\ese
Recall that $|\sd|$ is the $(d-1)$-dimensional measure of the face $\sd$. Moreover, $d_{\elm,\elmt} = |\tx_\elm - \tx_\elmt|$ is the Euclidean distance of
the points $\tx_\elm$ and $\tx_\elmt$ for an interior face $\sd = \sd_{\elm,\elmt} \in \Fhint$, unless the positions of $\tx_\elm$ and $\tx_\elmt$ are interchanged with respect to the positions of $\elm$ and $\elmt$, where we take $d_{\elm,\elmt} = - |\tx_\elm - \tx_\elmt|$. \cor{Since we suppose $\tx_\elm$ and $\tx_\elmt$ distinct, $d_{\elm,\elmt} \neq 0$.} Similarly, $d_{\elm,\sd} = |\tx_\elm - \tx_{\elm, \sd}|$ is the Euclidean distance of the points $\tx_\elm$ and $\tx_{\elm, \sd}$ for a boundary face $\sd \in \FK \cap \Fhext$, unless $\tx_\elm$ lies outside of $\elm$ past $\sd$, where we take $d_{\elm,\sd} = - |\tx_\elm - \tx_{\elm, \sd}|$. Figure~\ref{fig_FV_mesh} provides an illustration.

\br[Flux conservation and equilibrium in finite volumes \cor{on simplicial meshes}] From~\eqref{eq_dif_fl_int}, there holds 
\be \label{eq_FV_cons}
    U_{\elm,\sd} = - U_{\elmt,\sd} \qquad \forall \sd = \sd_{\elm,\elmt} \in \Fhint.
\ee
This is the flux conservation in finite volumes: what flows from the element $\elm$ out through the face $\sd_{\elm,\elmt}$ flows in the neighbor element $\elmt$ through the same face $\sd_{\elm,\elmt}$, \cf\ \cor{Lemma}~\ref{rem_normal_trace_cont} and definition~\eqref{eq_RTN}. The requirement~\eqref{eq_FV_scheme} then stipulates the equilibrium of these face fluxes with the load $f$, mimicking $\Dv \tu = f$ from Proposition~\ref{pr_prop_WS}. \er

\subsection{Flux reconstruction by lifting the face normal fluxes} \label{sec_fl_rec_Pois}

\begin{figure}
\centerline{\includegraphics[width=0.5\textwidth]{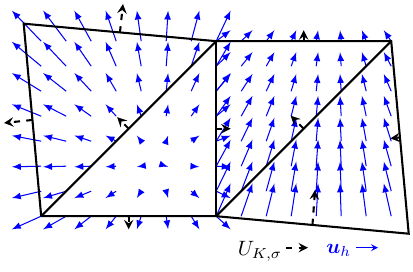}}
\caption{Finite volume face fluxes $U_{\elm,\sd}$ (dashed black arrows) and flux reconstruction $\tu_h$ \cor{of Definition~\ref{def_fr_Pois}} (blue arrows)}\label{fig_flux_rec}
\end{figure}

Following~\cite{Eym_Gal_Her_grad_01}, we now lift the normal fluxes $U_{\elm,\sd}$ from~\eqref{eq_dif_fl_int}--\eqref{eq_dif_fl_Dir} to form a vector-valued approximate flux field $\tu_h$ that is $\Hdv$-conforming and \cor{in equilibrium with the load $f$}. More precisely, we construct $\tu_h$ in the piecewise polynomial Raviart--Thomas space $\RT_0(\Th) \cap \Hdv$ from~\eqref{eq_RTN}.

Recalling the definition~\eqref{eq_RT0} of the space $\RT_0(\elm)$ for a simplex $\elm$ and the notation from Figure~\ref{fig_RTN_bas}, we start with:

\bd[Flux reconstruction] \label{def_fr_Pois} Let the cell-centered finite volume discretization on a simplicial mesh be given by~\eqref{eq_FVs}. Then, on each simplex $\elm \in \Th$, define 
\be \label{eq_fr_Pois}
    \tu_h|_\elm \in \RT_0(\elm), \quad \<\tu_h \scp \tn_\elm, 1\>_\sd = U_{\elm,\sd} \qquad \forall \sd \in \FK,
\ee
or, equivalently, \cor{using the Raviart--Thomas basis functions $\tv_\sd$ from~\eqref{eq_RT0_bas},}
\be \label{eq_fr_Pois_equiv}
    (\tu_h|_\elm) (\tx) \eq \sum_{\sd \in \FK} U_{\elm,\sd} \cor{\tv_\sd (\tx) = \sum_{\sd \in \FK} U_{\elm,\sd}} \frac{1}{d |\elm|}(\tx - \ver_{\elm, \sd}), \qquad \tx \in \elm.
\ee
\ed 

An illustration of the flux reconstruction $\tu_h$ from the finite volume face fluxes $U_{\elm,\sd}$ is provided in Figure~\ref{fig_flux_rec}. 
Crucially, we have:

\bl[Flux reconstruction] \label{lem_fr_Pois_prop} Let $\tu_h$ be given by Definition~\ref{def_fr_Pois}. Then
\be \label{eq_fr_Pois_prop}
    \tu_h \in \RT_0(\Th) \cap \Hdv, \qquad \cor{(}\Dv \tu_h\cor{)|_\elm} = \frac{(f,1)_\elm}{|\elm|} \qquad \forall \elm \in \Th.
\ee
\el

\bp From the finite volume conservation~\eqref{eq_FV_cons} and definition~\eqref{eq_fr_Pois}, we have the normal trace continuity needed in~\eqref{eq_RTN}; this yields the first property in~\eqref{eq_fr_Pois_prop}. As for the second property in~\eqref{eq_fr_Pois_prop}, we write, using the Green theorem, definition~\eqref{eq_fr_Pois}, and the finite volume scheme~\eqref{eq_FV_scheme}
\[
    (\Dv \tu_h,1)_\elm = \sum_{\sd \in \FK} \<\tu_h \scp \tn_\elm, 1\>_\sd \reff{eq_fr_Pois}= \sum_{\sd \in \FK} U_{\elm,\sd} \reff{eq_FV_scheme}= (f,1)_\elm.
\]
The conclusion follows since $\Dv \tu_h$ is constant by virtue of~\eqref{eq_RT0_prop}. \ep

\subsection{Potential reconstruction by elementwise postprocessing and averaging} \label{sec_pot_rec_Pois}

Supposing for the moment that $f$ is piecewise constant, so that $\tu_h \in \Hdv$ with $\Dv \tu_h = f$ from~\eqref{eq_fr_Pois_prop}, we will be in position to use Corollary~\ref{cor_Prag_Syng} to estimate the error between the finite volume flux reconstruction $\tu_h$ of Definition~\ref{def_fr_Pois} and the exact flux $\tu = - \Gr p$. Inequality~\eqref{eq_Prag_Syng_cons_2} tells us that we will need a suitable scalar-valued object in the primal space $\Hoo$ whose negative gradient is as close as possible to $\tu_h$. Since the piecewise constant function given by the finite volume values $p_\elm$ on each $\elm \in \Th$ cannot be used directly, as it does not lie in the space $\Hoo$, we proceed in two steps. 

Following~\cite[Section~4.1]{Voh_apost_MFE_07} and~\cite[Section~3.2]{Voh_apost_FV_08}, we first create a scalar-valued piecewise polynomial whose negative gradient equals $\tu_h$ and whose mean values are given by $p_\elm$ on each mesh element $\elm \in \Th$. Recalling the piecewise polynomial space $\PP_k(\Th)$ from~\eqref{eq_pw_pols}, we define:

\bd[Potential postprocessing] \label{def_postpr} Let $p_\elm$, $\elm \in \Th$, be the finite volume approximate solution of~\eqref{eq_FVs}. Let $\tu_h$ be given by Definition~\ref{def_fr_Pois}. Define the potential postprocessing $\pth$ as a piecewise quadratic polynomial on the simplicial mesh $\Th$, $\pth \in \PP_2(\Th)$, given by
\be \label{eq_postpr}
    - \Gr \pth|_\elm = \tu_h|_\elm, \qquad \frac{(\pth,1)_\elm}{|\elm|} = p_\elm
    \qquad \forall \elm \in \Th.
\ee
\ed

\begin{figure}
\centerline{\includegraphics[width=0.4\textwidth]{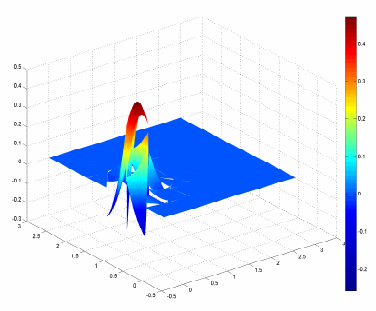} \quad \includegraphics[width=0.4\textwidth]{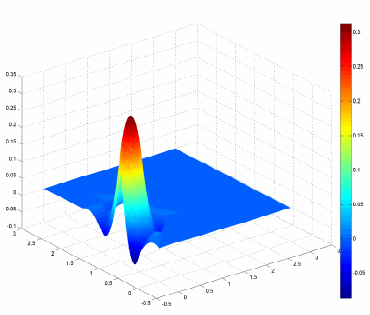}}
\caption{Potential postprocessing $\pth \in \PP_2(\Th)$ (left) and potential reconstruction by averaging $\prh \in \PP_2(\Th) \cap \Ho$ (right)}\label{fig_pospr_rec}
\end{figure}

It is easy to verify that any object from $\RT_0(\elm)$ for a simplex $\elm$ is a gradient of some polynomial \cor{which is at most second-order, since fields in $\RT_0(\elm)$ are curl-free and at most affine. Thus}~\eqref{eq_postpr} is well-posed. 

Unfortunately, in general $\pth$ \cor{from Definition~\ref{def_postpr}} does not lie in $\Hoo$\cor{; if this was the case, then, actually, by Proposition~\ref{thm_Prag_Syng} or Theorem~\ref{thm_err_char}, $\pth$ would be the weak solution $p$. An} illustration \cor{of a prototypical situation is provided} in Figure~\ref{fig_pospr_rec}, left. For this reason, we continue with a second step, in the spirit of averaging operators in Dur{\'a}n and Padra~\cite{Dur_Padr_err_est_NC_NL_94}, Achdou~\eal\ \cite{Ach_Ber_Coq_FV_Darcy_03}, Karakashian and Pascal~\cite{Karak_Pasc_DG_a_post_03}, Ainsworth~\cite{Ains_rob_a_post_NCFE_05}, or Burman and Ern~\cite{Bur_Ern_int_pen_DG_07}. Recall the Lagrange spaces and nodes of Section~\ref{sec_Lagr} and Figure~\ref{fig_Lagr_K}.

\bd[Potential reconstruction by averaging] \label{def_pr_Pois} Let $\pth \in \PP_2(\Th)$ be given by Definition~\ref{def_postpr}. We call a potential reconstruction by averaging a piecewise quadratic and continuous polynomial $\prh \in \PP_2(\Th) \cap \Hoo$ that is prescribed by its values in the Lagrange nodes of the space $\PP_2(\Th) \cap \Hoo$ by
\bse \label{eq_sh} \bat{2}
    \prh(\tx) & \eq \frac{1}{|\T_{\tx}|} \sum_{\elm \in \T_{\tx}}
        \pth|_\elm(\tx) \qquad & & \tx \text{ is a Lagrange node of $\PP_2(\Th) \cap \Hoo$ included in } \Om, \label{eq_sh_int} \\
    \prh(\tx) & \eq 0 & & \tx \text{ is a Lagrange node of $\PP_2(\Th) \cap \Hoo$ included in } \pt \Om, \label{eq_sh_BC}
\eat \ese
where $\T_{\tx}$ denotes the set of elements of the mesh $\Th$ that contain the point $\tx$ and $|\T_{\tx}|$ is the cardinality (number of elements) of this set.
\ed 

We simply average all the (typically) different values that $\prh$ takes in the Lagrange nodes. An illustration is given in Figure~\ref{fig_pospr_rec}, right (where the $0$ values at the boundary $\pt \Om$ are not imposed).
 
\subsection{A guaranteed a posteriori error estimate}
\label{sec_a_post_Pois}

With the above developments, we are now in a position to present a guaranteed a posteriori error estimate for the cell-centered finite volume discretization~\eqref{eq_FVs} of the Poisson equation~\eqref{eq_Pois} on a simplicial mesh $\Th$. We follow~\cite[Theorem~4.1]{Voh_apost_FV_08}.

\begin{ctheorem}{A guaranteed a posteriori error estimate}{thm_est_Pois}
Let $p \in \Hoo$ be the exact potential given by~\eqref{eq_Pois_WF} and let $\tu$ be the exact flux given by~\eqref{eq_velocity_Pois}. Let the cell-centered finite volume discretization be given by~\eqref{eq_FVs} and let the approximate flux $\tu_h \in \RT_0(\Th) \cap \Hdv$ be constructed following Definition~\ref{def_fr_Pois}. Let the potential postprocessing $\pth \in \PP_2(\Th)$ be given by Definition~\ref{def_postpr} and the potential reconstruction by averaging $\prh \in \PP_2(\Th) \cap \Hoo$ by Definition~\ref{def_pr_Pois}. Then there holds
\bse\label{eq_est_Pois} \begin{equation} \label{eq_est_Pois_1}
    \norm{\tu - \tu_h} \leq \eta \eq\Biggl\{\sum_{\elm \in \Th} \eta_\elm^2\Biggr\}^\ft,
\end{equation}
where
\be \label{eq_est_Pois_2}
    \eta_\elm^2 \eq \norm{\tu_h + \Gr \prh}_\elm^2 + \frac{h_\elm^2}{\pi^2}\norm[\bigg]{f-\frac{(f,1)_\elm}{|\elm|}}_\elm^2.
\ee\ese
\end{ctheorem}

\bp If $f$ is piecewise constant, then $\tu_h \in \Hdv$ with $\Dv \tu_h = f$ from~\eqref{eq_fr_Pois_prop}. In this case, we merely employ inequality~\eqref{eq_Prag_Syng_cons_2} from Corollary~\ref{cor_Prag_Syng}.

In the general case $f \in \Lt$, we rather use equality~\eqref{eq_err_char} from Theorem~\ref{thm_err_char}. The second term on the right-hand side of~\eqref{eq_err_char} gives rise to $\norm{\tu_h + \Gr \prh}^2$. As for the first one, let $\vf \in \Hoo$ with $\norm{\Gr \vf} = 1$ be fixed. By the Green theorem~\eqref{eq_Green_Hdv}, the flux reconstruction property~\eqref{eq_fr_Pois_prop}, the Cauchy--Schwarz inequality, and the Poincar\'e inequality
\be \label{eq_Poinc}
    \norm[\bigg]{\vf - \frac{(\vf,1)_\elm}{|\elm|}}_\elm \leq \frac{h_\elm}{\pi} \norm{\Gr \vf}_\elm \qquad \vf \in \Hoi{\elm},
\ee 
\cf~\cite{Pay_Wei_Poin_conv_60, Beben_Poin_conv_03}, we see
\ban
    (f,\vf) + (\tu_h, \Gr \vf) & \reff{eq_Green_Hdv}= (f - \Dv \tu_h,\vf) = \sum_{\elm \in \Th}(f - \Dv \tu_h,\vf)_\elm \reff{eq_fr_Pois_prop}= \sum_{\elm \in \Th}\bigg(f - \frac{(f,1)_\elm}{|\elm|},\vf\bigg)_\elm \\
    & = \sum_{\elm \in \Th}\bigg(f - \frac{(f,1)_\elm}{|\elm|},\vf - \frac{(\vf,1)_\elm}{|\elm|}\bigg)_\elm \reff{eq_Poinc}\leq \sum_{\elm \in \Th} \bigg(\norm[\bigg]{f - \frac{(f,1)_\elm}{|\elm|}}_\elm \frac{h_\elm}{\pi} \norm{\Gr \vf}_\elm \bigg) \\
    & \leq \Bigg\{\sum_{\elm \in \Th} \bigg(\frac{h_\elm^2}{\pi^2}\norm[\bigg]{f-\frac{(f,1)_\elm}{|\elm|}}_\elm^2\bigg)\Bigg\}^\ft \underbrace{\norm{\Gr \vf}}_{=1} = \Bigg\{\sum_{\elm \in \Th} \bigg(\frac{h_\elm^2}{\pi^2}\norm[\bigg]{f-\frac{(f,1)_\elm}{|\elm|}}_\elm^2\bigg)\Bigg\}^\ft.
\ean
Thus~\eqref{eq_est_Pois} follows. \ep

\br[Local efficiency] Local efficiency of the estimate of~\eqref{eq_est_Pois} as per Section~\ref{sec_eff} holds true, see~\cite[Theorem~4.2]{Voh_apost_FV_08}.\er

\cor{\br[Rectangular meshes] \label{rem_a_post_rect} All the above results are in~\cite{Voh_apost_FV_08} also presented on rectangular meshes.\er}

\subsection{Numerical experiments} \label{sec_num_Pois}

In this Section, we \cor{illustrate} the application of the a posteriori error estimates of \cor{the form of} Theorem~\ref{thm_est_Pois} to the finite volume method on two test cases in two space dimensions. 

\subsubsection{Regular solution}\label{sec_reg_sol}

\begin{figure}
   \centering
    \includegraphics[width=0.31\linewidth]{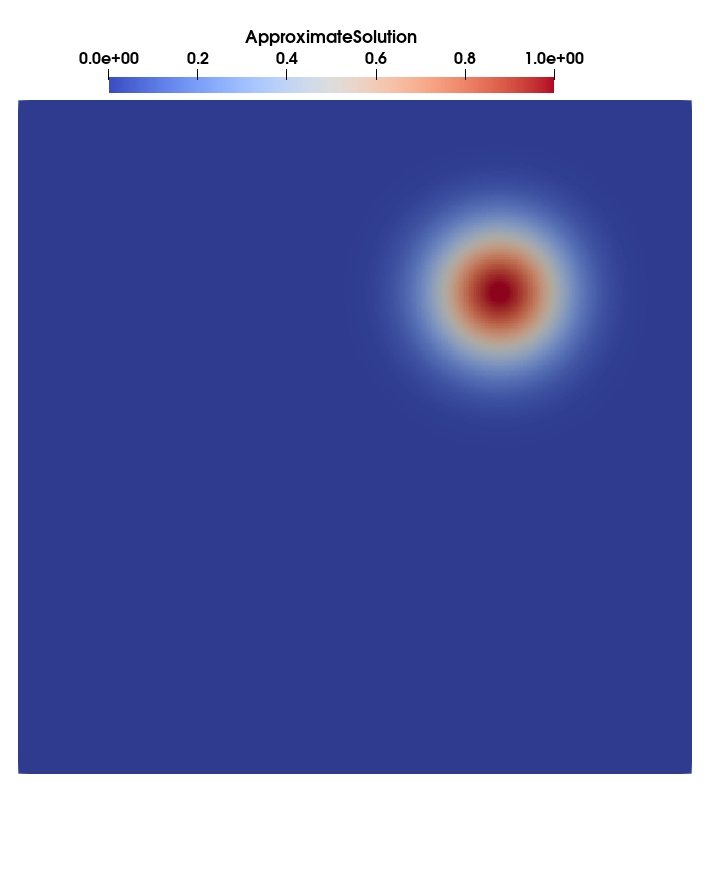}
\hspace{0.1cm}
    \includegraphics[width=0.31\linewidth]{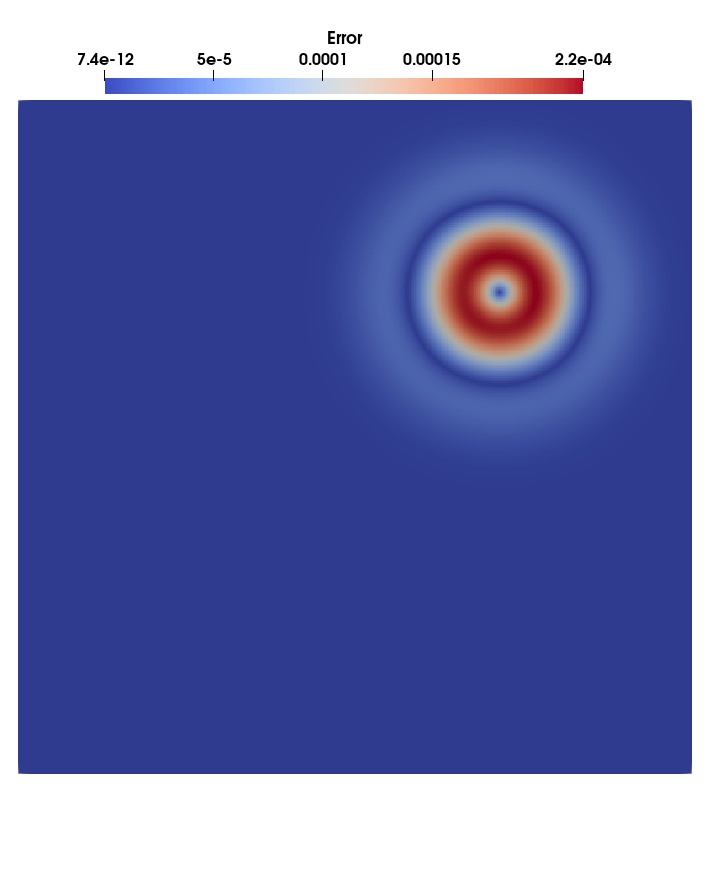}
\hspace{0.1cm}
    \includegraphics[width=0.31\linewidth]{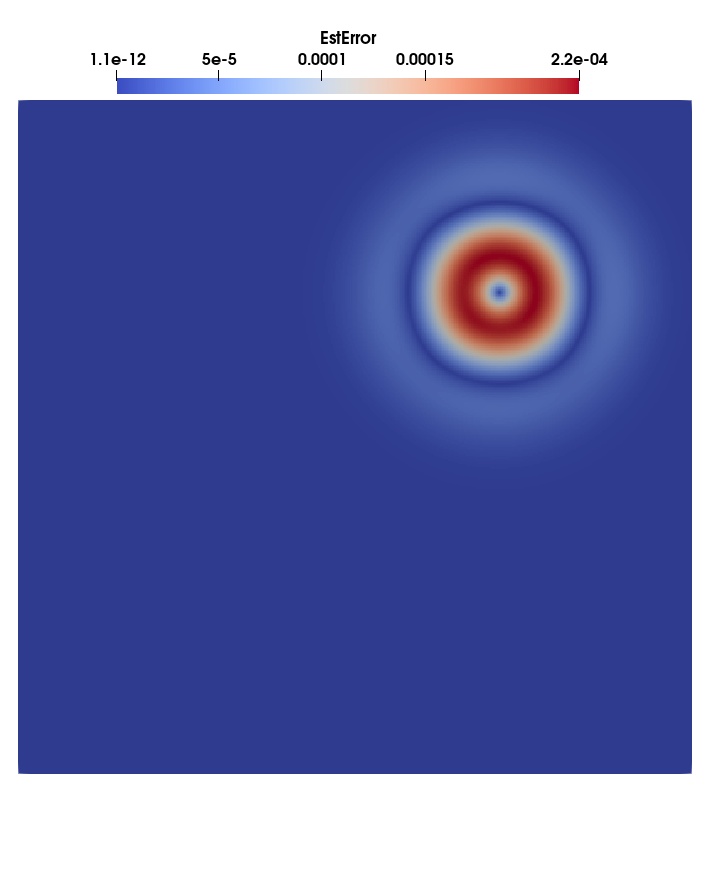}
    \caption{[Section~\ref{sec_reg_sol}, mesh $\T_3$] The approximate solution given by the values $p_\elm$, $\elm \in \Th$, ({\em left}), the exact errors $\norm{\tu - \tu_h}_\elm$ ({\em middle}), and the error estimators $\eta_\elm$ from Theorem~\ref{thm_est_Pois} ({\em right}).}
  \label{fig:reg.sol}
\end{figure}

In this test, we consider a model problem from~\cite{Mozo_Prudh_goal_est_15} of form~\eqref{eq_Pois} with $\Omega = (0,1) \times (0,1)$ and the load term $f$ given such that the (regular exact solution reads
\begin{equation}
p(x,y)=\frac{10^2}{4}x(1-x)y(1-y)\exp\left(-100\left((x-0.75)^2+(y-0.75)^2\right)\right).
\end{equation}
\cor{This solution is regular but contains a localized exponential peak.} 
The numerical tests are performed on a sequence of uniformly refined \cor{rectangular} meshes $\T_0,\T_1,\ldots,\T_{J}$, $J=3$\cor{, \cf\ Remarks~\ref{rem_mesh_rect} and~\ref{rem_a_post_rect}}. Figure~\ref{fig:reg.sol} illustrates the approximate solution given by the values $p_\elm$, $\elm \in \Th$, the elementwise errors $\norm{\tu - \tu_h}_\elm$, and the corresponding a posteriori error estimators $\eta_\elm$ \cor{arising} from Theorem~\ref{thm_est_Pois} \cor{on mesh $\T_3$}. We observe that the error is localized in a circular zone around the peak and that the actual and predicted error distributions match very closely.

\begin{figure}
\centering
     \begin{tikzpicture}[scale=0.8]
      \begin{semilogxaxis}[
         max space between ticks=30,
         yticklabel style={/pgf/number format/fixed},
          xlabel = {Number of unknowns},
          ylabel = {\cor{Relative e}rrors and estimates},
          legend style ={at = { (0.99,0.99)}} 
        ]
        \addplot +[green!50!black,  mark options={solid}, mark=diamond*, mark size=3]
        table[x=N,y=relativeErr]{Figs/res-convergence.tex};
        \addplot +[mark=square*, mark size=1.5, mark options={solid}, black]table[x=N,y=relativeEst]{Figs/res-convergence.tex};
        \legend{Error, Estimate};
      \end{semilogxaxis}
    \end{tikzpicture}
\hspace{0.5cm}
     \begin{tikzpicture}[scale=0.8]
      \begin{semilogxaxis}[
         max space between ticks=30,
          xlabel = {Number of unknowns},
          ylabel = {Effectivity indices},
          legend style ={at = { (0.92,0.92)}} 
        ]
        \addplot +[line width=0.2mm] table[x=N,y=effId]{Figs/res-convergence.tex};
        \legend{Effectivity index};
      \end{semilogxaxis}
    \end{tikzpicture}
  \caption{[Section~\ref{sec_reg_sol}, uniformly refined meshes $\T_0$--$\T_3$] Relative errors $\norm{\tu - \tu_h}/\norm{\tu_h}$ and relative estimates $\eta/\norm{\tu_h}$ from Theorem~\ref{thm_est_Pois} ({\em left}), effectivity indices $I_{\mathrm{eff}}$ from~\eqref{eq_eff_ind} ({\em right})}
\label{fig_estim_reg_sol}
\end{figure}
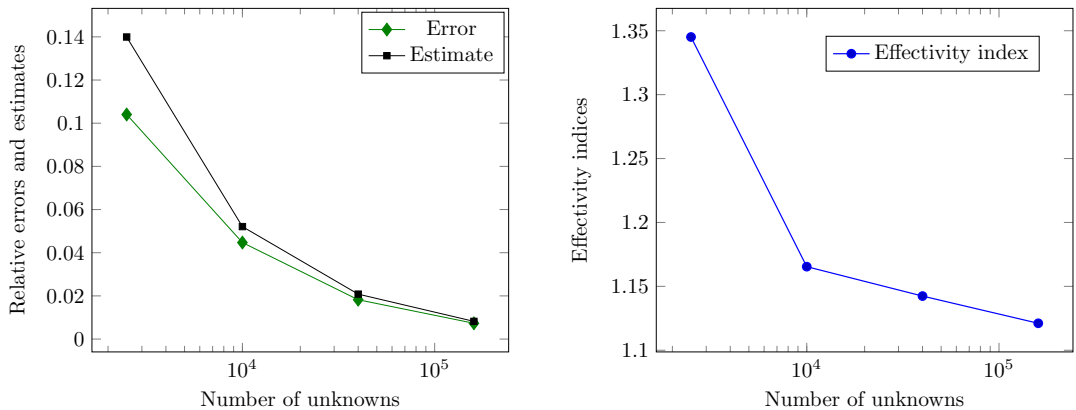

In the left part of Figure~\ref{fig_estim_reg_sol}, we plot the relative exact errors $\norm{\tu - \tu_h}{/\norm{\tu_h}}$ on each of the meshes $\T_0$--$\T_3$ together with their relative a posteriori error estimates $\eta/\norm{\tu_h}$ from Theorem~\ref{thm_est_Pois}, with respect to the number of mesh elements/number of unknowns. We indeed observe that the estimators give a guaranteed upper bound on the errors. In the right part of Figure~\ref{fig_estim_reg_sol} we then plot the effectivity indices defined as per~\eqref{eq_I_eff} by the ratio of the estimator to the error, 
\be \label{eq_eff_ind}
    I_{\mathrm{eff}} \eq \frac{\eta}{\norm{\tu - \tu_h}}.
\ee
We observe that they are remarkably close to one and tend to the optimal value of one with mesh refinement, showing asymptotic exactness as per Section~\ref{sec_ass_ex}.

\subsubsection{Singular solution}\label{sec_sing_sol}

\begin{figure}
   \centering
    \includegraphics[width=0.32\linewidth]{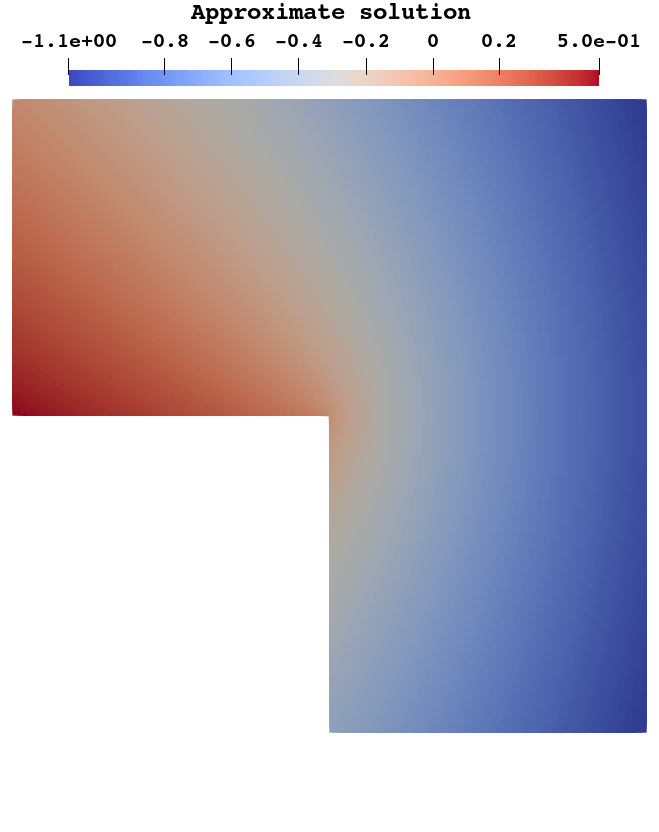}
\hspace{0.1cm}
    \includegraphics[width=0.311\linewidth]{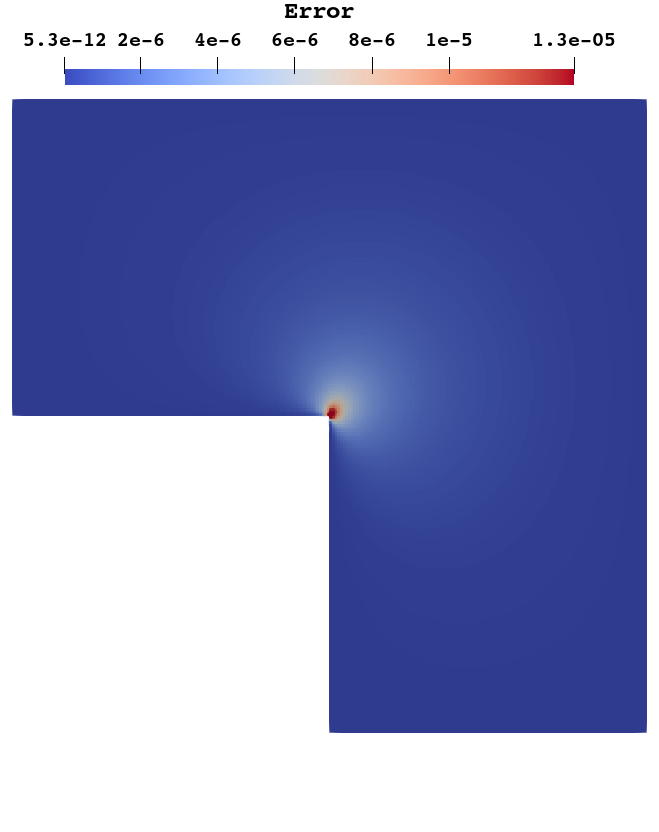}
\hspace{0.1cm}
    \includegraphics[width=0.311\linewidth]{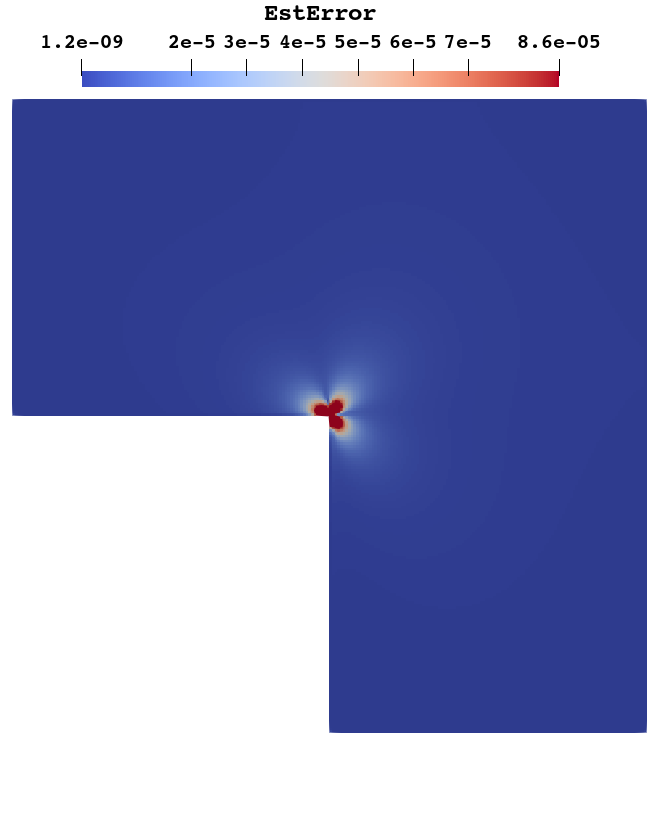}
    \caption{[Section~\ref{sec_sing_sol}, mesh $\T_2$] The approximate solution given by the values $p_\elm$, $\elm \in \Th$, ({\em left}), the exact errors $\norm{\tu - \tu_h}_\elm$ ({\em middle}), and the error estimators $\eta_\elm$ from Theorem~\ref{thm_est_Pois} ({\em right}).}
  \label{fig_sing_sol}
\end{figure}

We consider a benchmark test case defined on the L-shaped domain 
\[
    \Omega=(-1,1)\times (-1,1)\setminus (-1,0]\times (-1,0]
\]
with the exact solution 
\begin{align*}
p(r,\theta)=r^{\frac{2}{3}}\sin\left(\frac{2}{3}\theta + \frac{3}{2} \pi\right),
\end{align*}
where ($r,\theta$) are the polar coordinates. The Dirichlet boundary condition is partly inhomogeneous here and given by the value of the exact solution on $\pt \Om$. We illustrate in Figure~\ref{fig_sing_sol} the approximate solution given by the values $p_\elm$, $\elm \in \Th$, the elementwise errors $\norm{\tu - \tu_h}_\elm$, and the corresponding a posteriori error estimators $\eta_\elm$ from Theorem~\ref{thm_est_Pois}. We observe that the a posteriori error estimators $\eta_\elm$ detect perfectly the error related \cor{in the vicinity of} the singularity at $(0,0)$.

\begin{figure}
\centering
     \begin{tikzpicture}[scale=0.85]
      \begin{semilogxaxis}[
        xmin = 1e3,
        xmax = 1e5,
          xlabel = {Number of unknowns},
          ylabel = {\cor{Relative e}rrors and estimators},
          legend style ={at = { (0.99,0.99)}} 
        ]
        \addplot +[green!50!black,  mark options={solid}, mark=diamond*, mark size=3]table[x=N,y=relativeErr]{Figs/LShape.tex};
        \addplot +[mark=square*, mark size=1.5, mark options={solid}, black]table[x=N,y=relativeEst]{Figs/LShape.tex};
        \legend{Error, Estimate};
      \end{semilogxaxis}
    \end{tikzpicture}
\hspace{0.5cm}
     \begin{tikzpicture}[scale=0.85]
      \begin{semilogxaxis}[
        xmin = 1e3,
        xmax = 1e5,
        ymin = 1.25,
        ymax = 1.45,
          xlabel = {Number of unknowns},
          ylabel = {Effectivity indices},
          legend style ={at = { (0.92,0.92)}} 
        ]
        \addplot +[line width=0.2mm] table[x=N,y=effId]{Figs/LShape.tex};
        \legend{Effectivity index};
      \end{semilogxaxis}
    \end{tikzpicture}
  \caption{[Section~\ref{sec_sing_sol}, uniformly refined meshes $\T_0$--$\T_2$] Relative errors $\norm{\tu - \tu_h}/\norm{\tu_h}$ and relative estimates $\eta/\norm{\tu_h}$ from Theorem~\ref{thm_est_Pois} ({\em left}), effectivity indices $I_{\mathrm{eff}}$ from~\eqref{eq_eff_ind} ({\em right})}
\label{fig_sing_estim}
\end{figure}
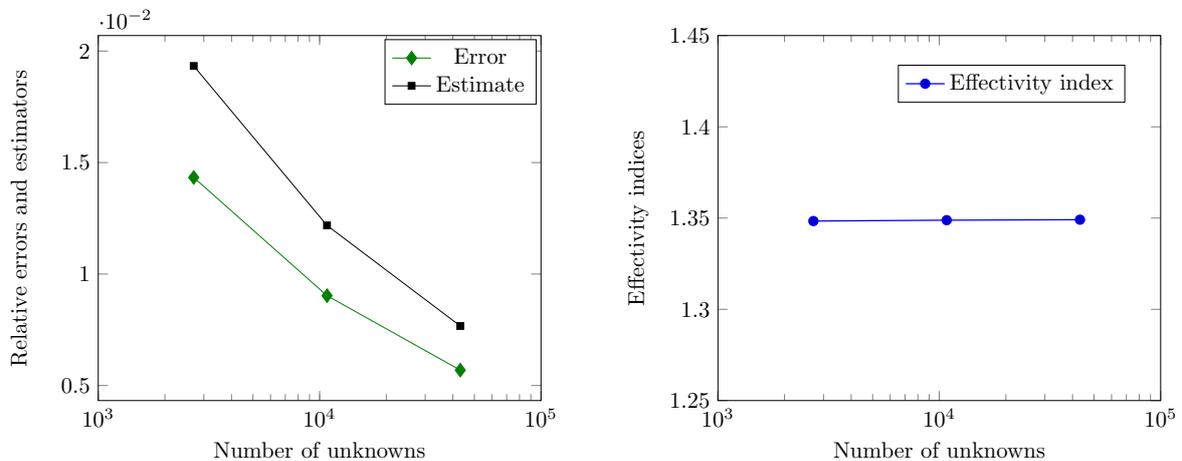

The convergence histories on the sequence of the uniformly-refined meshes $\T_0 - \T_2$ are depicted in the left part of Figure~\ref{fig_sing_estim}. The right part of Figure~\ref{fig_sing_estim} then displays the corresponding effectivity indices from~\eqref{eq_eff_ind}. We see that, for this test case, the effectivity indices are still close to the optimal value of one, even with singular solution, though no asymptotic exactness is observed here.

\subsection{Bibliographic resources}\label{sec_biblio_Pois}

One of the first a posteriori error analyses for finite volume methods was performed by Angermann~\cite{Anger_bal_apost_FV_CD_95}. Later, Achdou~\eal\ \cite{Ach_Ber_Coq_FV_Darcy_03}, Afif~\eal\ \cite{Af_Ber_Mghaz_Verf_a_post_FV_el_03}, Lazarov and Tomov~\cite{Laz_Tom_a_post_FVE_CRD_02}, Carstensen~\eal\ \cite{Cars_Laz_Tom_expl_aver_a_post_FV_05}, and Amaziane~\eal\ \cite{Am_Ber_El_Os_MGH_a_post_VC_FV_CD_09} also considered the so-called vertex-centered finite volume discretizations, closely related to lowest-order finite elements. Cell-centered finite volumes of the form~\eqref{eq_FV_scheme} were then addressed in Agouzal and Oudin~\cite{Ag_Ou_a_post_FV_simpl_00}, Achdou~\eal\ \cite{Ach_Ber_Coq_FV_Darcy_03}, Nicaise~\cite{Nic_a_post_FV_dif_05}, Vohral{\'{\i}}k~\cite{Voh_apost_FV_08}, and Erath and Praetorius~\cite{Er_Pra_a_post_FV_09}.

\section{Steady linear pure diffusion problems. Polytopal meshes and inexpensive implementation and evaluation} \label{sec_Darcy}

We describe here first a model steady linear problem, the Darcy flow problem. This extends the Poisson equation~\eqref{eq_Pois} from Section~\ref{sec_Pois} and forms the basis for complex multiphase compositional porous media flows as that of Section~\ref{sec_MP_MC}. We also focus on extending the finite volume discretizations on simplicial meshes from Section~\ref{sec_FV_Pois} to general polytopal meshes. Finally, we investigate how the guaranteed a posteriori error estimates of the form of Theorem~\ref{thm_est_Pois} can be implemented and evaluated in the most inexpensive way, and this on polytopal meshes. In particular, the flux and potential reconstructions of Sections~\ref{sec_fl_rec_Pois}--\ref{sec_pot_rec_Pois} will only be performed here virtually, avoiding their factual construction. Instead, the a posteriori error estimates will only need {\em one flux value per polytopal face} and {\em one potential value per vertex of a (virtual) simplicial submesh}. The former are directly available from the numerical scheme (these are the \cor{generalizations of the} face normal fluxes $U_{\elm,\sd}$ in~\eqref{eq_FVs}) and the latter will be trivially obtained by some postprocessing and averaging of the elementwise potential values \cor{generalizing} $p_\elm$ in~\eqref{eq_FVs}. 

\subsection{Singlephase steady linear Darcy flow}

Consider the problem of finding $p: \Om \ra \RR$ such that
\bse \label{eq_Darcy} \begin{empheq}[box=\widefbox]{align}
    \ds - \Dv(\Km \Gr p) & = f \qquad \mbox{ in } \, \Om, \label{eq_Darcy_eq} \\
    \ds p & = 0 \qquad \mbox{ on } \, \pt \Om. \label{eq_Darcy_BC}
\end{empheq} \ese
Here $f$ is a source term and $\Km$ is a symmetric, positive
definite, and bounded diffusion tensor with values in $\RR^{d \times
d}$; we suppose for simplicity here that both $f$ and $\Km$ are piecewise
constant with respect to the polytopal mesh $\T_H$. 

\subsection{Weak solution and its properties} \label{sec_prop_Darcy}

The Sobolev space $\Hoo$ from Definition~\ref{def_Hoo} still gives a proper mathematical setting to define a unique solution $p$ of problem~\eqref{eq_Darcy}, which is $p \in \Hoo$ such that
\be \label{eq_Darcy_WF}
    (\Km \Gr p, \Gr v) = (f, v) \qquad \forall v \in \Hoo.
\ee
From the pressure head $p$, we can define the Darcy velocity
\be \label{eq_velocity}
    \tu \eq - \Km \Gr p;
\ee
\cor{recall that we} refer to $p$ as {\em potential} and to $\tu$ as {\em flux}. As in Proposition~\ref{pr_prop_WS}, it follows from~\eqref{eq_Darcy_WF} that 
\be \label{eq_tu_pol}
    \tu \in \Hdv, \quad \Dv \tu = f. 
\ee
Thus, just like in Remark~\ref{rem_p_u}, the scalar-valued pressure head $p$ is trace continuous, \cf\ Figure~\ref{fig_Ho}, whereas the vector-valued Darcy velocity $\tu$ is normal-trace continuous, \cf\ Figure~\ref{fig_H_div}. Moreover, $\tu$ is locally mass conservative, in equilibrium with the load, satisfying $\Dv \tu = f$.

\subsection{Spatial discretization error}\label{sec_err_NL}

Applying the abstract reflections of Section~\ref{sec_sp_err}, we will evaluate the spatial discretization error as
\be \label{eq_en_norm_Darcy}
    \norm{\tu - \tu_h}_{\Km^{-\frac{1}{2}}}
\ee
for an approximate Darcy velocity $\tu_h$.
Here, the energy norm on $\om \subset \Om$, for ${\tv} \in \left[L^2(\om)\right]^d$, writes as
\be \label{eq_norm_Darcy}
\norm{\tv}_{\Km^{-\frac{1}{2}},\om} \eq \norm{\Km^\mft \tv}_{L^2(\om)} = \left\{\int_\om \left|\Km^\mft(\tx) {\tv}(\tx)\right|^2 \dx \right\}^{\frac12};
\ee
when $\om = \Om$, we \cor{as usual} drop the subscript $\om$.

\subsection{Generic discretizations on a polytopal mesh}

Let $\T_H$ be a polytopal mesh of $\Om$ in the sense of
Section~\ref{sec_meshes}. Define the load vector $\alg{F} \eq \{\alg{F}_\elm\}_{\elm \in \T_H} \in \RR^{|\T_H|}$
by
\be \label{eq_F}
    \alg{F}_\elm \eq (f,1)_\elm \qquad \forall \elm \in \T_H.
\ee
We consider any numerical discretization of~\eqref{eq_Darcy} that can be written under the following abstract lowest-order locally conservative form:

\bas[Locally conservative discretization on a polytopal mesh]
\label{as_polyt_disc} Find the algebraic vector $\alg{P} \eq \{\Pelm\}_{\elm \in \T_H} \in
\RR^{|\T_H|}$ such that the flux balance
\be \label{eq_flux_balance}
    \boxed{\sum_{\sd \in \FK} \Usd \tn_\elm \scp
    \tn_\sd = \alg{F}_\elm \quad \forall \elm \in \T_H}
\ee
is satisfied. Here $\alg{U} \eq \{\Usd\}_{\sd \in \F_H} \in \RR^{|\F_H|}$, and $\Usd \in \RR$ for each face $\sd \in \F_H$ approximates the normal flux $\<\tu \scp \tn_\sd, 1\>_\sd$ over the face $\sd$ and depends linearly on $\alg{P}$.
The unknowns $\Pelm \in \RR$ for each element $\elm \in \T_H$ approximate the potential
$p$ in the element $\elm$. \eas

Assumption~\ref{as_polyt_disc} is very generic (allows for any lowest-order locally conservative method). Note in particular that for the a posteriori error analysis performed below, we {\em do not need to know} how the flux unknowns $\alg{U}$ depend on the pressure head unknowns $\alg{P}$. Recall that on admissible simplicial meshes, the easiest example for this relation is the two-point formula~\eqref{eq_dif_fl_int}--\eqref{eq_dif_fl_Dir}. 

Many lowest-order locally conservative methods take a form comprised in Assumption~\ref{as_polyt_disc} but more specific. We consider them
separately under the following assumption, in the spirit of mimetic finite differences~\cite[Theorem~5.1]{Brez_Lip_Shash_MFD_cvg_05}, mixed finite elements~\cite[Theorems~7.2 and~7.3]{Voh_Wohl_MFE_1_unkn_el_rel_13}, and the unifying frameworks~\cite{Dro_Ey_Gal_Her_un_10, Dro_Eym_Her_grad_sch_16}:

\bas[Locally conservative saddle-point discretization on a polytopal mesh]
\label{as_polyt_disc_mod} Find the algebraic vectors $\alg{U} \eq \{\Usd\}_{\sd \in \F_H} \in
\RR^{|\F_H|}$ and $\alg{P} \eq \{\Pelm\}_{\elm \in \T_H} \in \RR^{|\T_H|}$
such that
\be \label{eq_polyt_disc}
 \boxed{
 \begin{pmatrix}
    \matr{A}  & \matr{B}^{\mathrm{t}} \\
    \matr{B}  & 0  \\
 \end{pmatrix}
 \begin{pmatrix}
    \alg{U} \\
    \alg{P} \\
 \end{pmatrix}
=
 \begin{pmatrix}
    \alg{0} \\
    - \alg{F} \\
 \end{pmatrix},}
\ee
where \ben[label=\arabic*)]
\item the unknown $\Usd \in \RR$ for each face $\sd \in \F_H$ approximates the normal
flux $\<\tu \scp \tn_\sd, 1\>_\sd$ over the face $\sd$ and the unknown $\Pelm \in \RR$ for each element $\elm \in \T_H$ approximates the potential $p$ in the element $\elm$;
\item the matrix $\matr{B} \in \RR^{|\T_H| \times |\F_H|}$ has a full rank; 
\item for each polytopal cell $\elm \in \T_H$ and face $\sd \in
\F_H$, $\matr{B}_{\elm, \sd} = - \tn_\elm \scp \tn_\sd$ if $\sd$ is a face of $\elm$, $\sd \in \FK$, and $\matr{B}_{\elm, \sd} = 0$ otherwise; 
\item the matrix $\matr{A} \in \RR^{|\F_H| \times |\F_H|}$ is invertible; 
\item for faces $\sd, \sdt \in \F_H$, $\matr{A}_{\sd, \sdt} = 0$ if $\sd$ and $\sdt$ are not faces of the same element $\elm \in \T_H$, $\{\sd, \sdt\} \not \subset \FK$ for some $\elm \in \T_H$, and 
\be \label{eq_A_struct}
    \matr{A}_{\sd, \sdt} = \sum_{\elm \in \T_H,
    \{\sd, \sdt\} \subset \FK} \tn_\elm \scp \tn_\sd \, \tn_\elm \scp \tn_\sdt (\widehat{\matr{A}}_\elm)_{\sd, \sdt}
\ee
otherwise, where $\widehat{\matr{A}}_\elm \in \RR^{|\FK|\times|\FK|}$ are symmetric and positive definite (the ``element matrices'' of the given
method for each $\elm \in \T_H$); 
\item on each polytopal cell $\elm \in \T_H$, there exists a lifting $\tilde \tu_h|_\elm \in \Hdvi{\elm}$ of the face normal fluxes $\alg{U}_\elm^{\mathrm {ext}} = \{\Usd\}_{\sd \in \FK}$ to the interior of $\elm$ such that
\bse \label{eq_lift_K} \ba 
    \<\tilde \tu_h|_\elm \scp \tn_\sd, 1\>_\sd & = \Usd \quad \forall \sd \in \FK, \label{eq_norm_fl_K}\\
    \Dv \tilde \tu_h|_\elm & = f|_\elm, \label{eq_div_K}\\
    \norm{\tilde \tu_h}_{\Km^{-\frac{1}{2}},\elm}^2 & = (\alg{U}_\elm^{\mathrm {ext}})^{\mathrm{t}} \widehat{\matr{A}}_{\elm} \alg{U}_\elm^{\mathrm
    {ext}}. \label{eq_en_K}
\ea \ese
\een
\eas

\subsection{Face normal fluxes} \label{sec_fl_faces}

In Theorem~\ref{thm_est_Darcy} below, we will need \cor{a vector of real values with one} value $\Usd$ per face $\sd \in \F_{H,h}$ (approximate face normal fluxes): 

\bd[Face normal fluxes] \label{def_face_fluxes}
Let a discretization scheme of the form of Assumption~\ref{as_polyt_disc} or~\ref{as_polyt_disc_mod} be given, leading to the face normal fluxes $\alg{U} \eq \{\Usd\}_{\sd \in \F_H} \in \RR^{|\F_H|}$. For each polytopal face $\sd \in \F_H$ and each simplicial subface $\sdt \in \F_{H,h}$ \cor{as per Section~\ref{sec_meshes}}, define the face normal flux
\be \label{eq_conv}
    \Usdt \eq \frac{\Usd}{|\sd|} |\sdt| \quad \forall \sdt \in \F_{H,h}, \, \sdt \subset \sd, \, \sd \in \F_H
\ee
and collect them in the element vectors $\alg{U}_\elm^{\mathrm {ext}} \in \RR^{|\FKhext|}$ 
\be \label{eq_frh_sd_vec}
    \alg{U}_\elm^{\mathrm {ext}} \eq \{\Usd\}_{\sd \in \FKhext}
\ee
for each polytopal cell $\elm \in \T_H$.
\ed

Remark that if the faces of the simplicial mesh $\Th$ do not subdivide the
faces of the polytopal mesh $\T_H$, then $\F_{H,h} = \F_H$ and consequently $\Usdt = \Usd$ for all $\sd \in \FK$ and $\alg{U}_\elm^{\mathrm {ext}} = \{\Usd\}_{\sd \in \FK}$. 
\cor{Figures~\ref{fig:notation_types} and~\ref{fig:notations} give an illustration in two space dimensions.}

\subsection{Fictitious flux reconstruction by lifting the face normal fluxes} \label{sec_fl_rec}

In Definition~\ref{def_fr_Pois}, \cor{on simplicial meshes,} we have reconstructed a discrete $\Hdv$-conforming flux $\tu_h$ by lifting the face normal fluxes of a finite volume scheme. We extend here this procedure to general polytopal meshes. In contrast to Theorem~\ref{thm_est_Pois}, however, Theorem~\ref{thm_est_Darcy} below will not need $\tu_h$ to evaluate the posteriori error estimate, so one does {\em not need to perform} physically/on the computer this reconstruction. Congruently, we may not have a closed formula for the reconstruction, which on top may not be discrete. This is the reason for our naming {\em fictitious}. 

\subsection*{A (fictitious) flux reconstruction under Assumption~\ref{as_polyt_disc}} 

We first treat the generic case of Assumption~\ref{as_polyt_disc}; the case of
Assumption~\ref{as_polyt_disc_mod} is treated below.
Following the concept of lifting operators used in, \eg, Eymard~\eal\
\cite[Section~1.2]{Eym_Gal_Her_grad_01}, Kuznetsov and
Repin~\cite{Kuzn_Rep_MFE_pol_03}, Brezzi~\eal\
\cite[Theorem~5.1]{Brez_Lip_Shash_MFD_cvg_05}, Beir{\~a}o da
Veiga~\cite[Section~2.1]{Bei_res_a_post_MFD_08}, Beir{\~a}o da Veiga and
Manzini~\cite[Section~2.4]{Bei_Manz_a_post_MFD_08},
Kuznetsov~\cite{Kuzn_MFE_pol_08},
Vohral\'ik~\cite[Section~3.2]{Voh_apost_FV_08}, or Sboui~\eal\
\cite{Sbou_Jaff_Rob_comp_MFE_09}, we will extend the fluxes $\Usd$, $\sd
\in \FK$, to the interior of each polytopal element $\elm \in \T_H$. The
obtained flux $\tu_h|_\elm$ will be the approximation of \cor{the exact flux} $\tu$ from~\eqref{eq_velocity} for which we will estimate the \cor{energy} error $\norm{\tu - \tu_h}_{\Km^{-\frac{1}{2}}}$. 
Let us first present a motivation. 

\bmo[Motivation for the fictitious flux reconstruction] \label{mo_fl_rec}
Following~\cite[equation~(3.14)]{Voh_apost_FV_08}, consider in each element $\elm \in \T_H$ the infinite-dimensional Darcy problem: for the fluxes $\{\Usd\}_{\sd \in \FK}$ and the potential $\Pelm$ from Assumption~\ref{as_polyt_disc}, find $\tilde p_\elm: \elm \ra \RR$ such that
\begin{subequations} \label{eq_loc_probl} \begin{align}
     - \Dv (\Km \Gr \tilde p_\elm) & = {f|_\elm}, \\
    \frac{(\tilde p_\elm,1)_\elm}{|\elm|} & = \Pelm, \\
    - \Km \Gr \tilde p_\elm \cdot \tn_{\sd} & = \frac{\Usd}{|\sd|}
        \quad \forall \sd \in
        \FK. \label{eq_loc_probl_3}
\end{align} \end{subequations}
We have actually already seen this problem, in Definition~\ref{def_postpr} of Section~\ref{sec_pot_rec_Pois}. There, on a simplicial mesh, the solution was \cor{easily} practically computable and discrete: the second-order polynomial denoted as $\pth|_\elm$ (indeed, \eqref{eq_postpr} gives a solution to~\eqref{eq_loc_probl}) on a simplicial mesh. Unfortunately, for a general polytopal mesh element $\elm \in \T_H$, \eqref{eq_loc_probl} amounts to a partial differential equation on $\elm$, with the solution $\tilde p_\elm$ living in the infinite-dimensional Sobolev space $\Hoi{\elm}$ of Definition~\ref{def_Ho} and the associated flux $\tilde \tu_\elm \eq - \Km \Gr \tilde p_\elm$ living in the infinite-dimensional Sobolev space $\Hdvi{\elm}$ of Definition~\ref{def_Hdv}.
\emo

Based on Motivation~\ref{mo_fl_rec}, we would like to produce a discrete version of~\eqref{eq_loc_probl}, which will in particular enable us to derive (easily) computable a posteriori error estimates. Following~\cite[Theorem~5.1]{Voh_apost_FV_08} and~\cite[Theorems~7.2 and~7.3]{Voh_Wohl_MFE_1_unkn_el_rel_13}, we consider the mixed finite element method on the simplicial submesh $\TK$ of $\elm$ to {\em approximate the solution} $\tilde p_\elm$ of~\eqref{eq_loc_probl}. 
Note that this is to be done separately for each $\elm \in \T_H$\cor{, so it is an intrinsically parallel procedure}. 
Recall Definition~\ref{def_face_fluxes} and consider the following sets and spaces of {\em piecewise polynomials} over the {\em simplicial submesh} $\TK$ of $\elm$:
\be \label{eq_spaces}
    \begin{array}{l}
        \tV_{h, \mathrm{N}}^\elm \eq \{\tv_h \in  \RT_0(\TK) \cap \Hdvi{\elm}; \, \<\tv_h
        \scp \tn_\sd, 1\>_\sd = \Usd \quad \forall \sd \in \FKhext\}, \\[1mm]
        \tV_{h, \mathrm{0}}^\elm \eq \{\tv_h \in \RT_0(\TK) \cap \Hdvi{\elm}; \, \<\tv_h
         \scp \tn_\sd, 1\>_\sd = 0 \quad \forall \sd \in \FKhext\}, \\[1mm]
        Q_{h, \mathrm{N}}^\elm \eq \{q_h \in \PP_0(\TK); \, \frac{(q_h, 1)_\elm }{|\elm|} = \Pelm \}, \\[1mm]
        Q_{h, 0}^\elm \eq \{q_h \in \PP_0(\TK); \, (q_h, 1)_\elm = 0 \}.
    \end{array}
\ee
Here, $\tV_{h, \mathrm{N}}^\elm$ and $\tV_{h, \mathrm{0}}^\elm$ are subsets/subspaces of the lowest-order Raviart--Thomas spaces as defined in Section~\ref{sec_RT}, whereas $Q_{h, \mathrm{N}}^\elm$ and $Q_{h, 0}^\elm$ are subsets/subspaces of the piecewise polynomial spaces as defined in Section~\ref{sec_pw_pol}. We then let:

\bd[Flux reconstruction on polytopal meshes under Assumption~\ref{as_polyt_disc}] \label{def_fr} \cor{Let $f$ be constant on each $\elm \in \T_H$ and} let a discretization scheme of the form of Assumption~\ref{as_polyt_disc} be given, leading to the normal face fluxes $\alg{U}_\elm^{\mathrm {ext}} = \{\Usd\}_{\sd \in \FK}$, $\elm \in \T_H$, and consequently to $\{\Usd\}_{\sd \in \FKhext}$ as per Definition~\ref{def_face_fluxes}. For all $\elm \in \T_H$, define $\tu_h|_\elm$ by the mixed finite element discretization of~\eqref{eq_loc_probl}, described by the constrained quadratic minimization problem
\be \label{eq_arg_min}
    \tu_h|_\elm \eq \arg\min_{\tv_h \in \tV_{h, \mathrm{N}}^\elm, \, \Dv \tv_h = \mathrm{constant} =f|_\elm} \norm{\tv_h}_{\Km^{-\frac{1}{2}},\elm}^{2}.
\ee
Equivalently, writing the Euler--Lagrange conditions of~\eqref{eq_arg_min} and imposing the divergence constraint with a Lagrange multiplier, $\tu_h|_\elm \in \tV_{h, \mathrm{N}}^\elm$ together with $p_h|_\elm \in
Q_{h, \mathrm{N}}^\elm$ are given by
\bse \label{eq_fl_equil} \bat{2}
    (\Km^{-1} \tu_h, \tv_h)_\elm - (p_h, \Dv \tv_h)_\elm & = 0
    & \qquad & \forall \tv_h \in \tV_{h, \mathrm{0}}^\elm, \label{eq_fl_equil_1}\\
    - (\Dv \tu_h,q_h)_\elm & = -(f,q_h)_\elm = 0 & \qquad &\forall q_h \in Q_{h, 0}^\elm. \label{eq_fl_equil_2}
\eat\ese
\ed

It is worth noting that the divergence constraint in~\eqref{eq_arg_min} can be written with either a generic constant, or more specifically with $f|_\elm$; indeed, the assumption that $f$ is constant on each $\elm \in \T_H$, the definition~\eqref{eq_spaces} of $\tV_{h, \mathrm{N}}^\elm$, and the flux
balance~\eqref{eq_flux_balance} give the equivalence. 
Congruently, \eqref{eq_fl_equil_2} can be equivalently written as
\[
    (\Dv \tu_h,q_h)_\elm = (f,q_h)_\elm \qquad  \forall q_h \in \PP_0(\TK),
\]
which is \cor{further} equivalent to $\Dv \tu_h|_\elm = f|_\elm$.
Note also that since the normal fluxes $\Usd$ are univalued on the faces from $\F_{H,h}$, the resulting flux reconstruction $\tu_h$ has the normal trace continuous over all interior mesh faces from $\F_{H,h}$ and thus belongs to the
Raviart--Thomas space $\RT_0(\Th) \cap \Hdv$ associated with the
simplicial mesh $\Th$ of the entire domain $\Om$. In consequence, \cor{we have:} 

\bl[\cor{Flux reconstruction on polytopal meshes under Assumption~\ref{as_polyt_disc}}] \label{lem_fr} \cor{Let $f$ be constant on each $\elm \in \T_H$. Then} the flux reconstruction of Definition~\ref{def_fr} satisfies 
\[
    \tu_h \in \RT_0(\Th) \cap \Hdv \quad \text{ with } \quad \Dv \tu_h = f.
\]
\el 

\begin{figure}
\centerline{\includegraphics[height=0.3\textwidth]{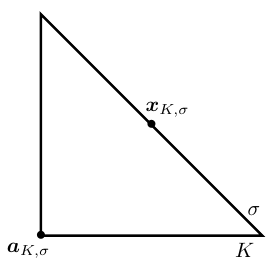} \qquad \includegraphics[height=0.35\textwidth]{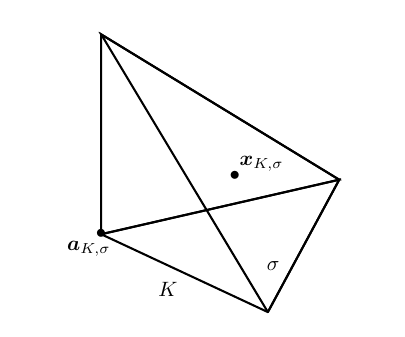}}
\caption{Face \cor{barycentres} $\tx_{\elm, \sd}$ and vertices $\ver_{\elm, \sd}$ opposite to a face $\sd \in \FK$. Simplex $\elm \in \Th$ in space dimensions $d=2$ (left) and $d=3$ (right)}\label{fig_face_bar_vert}
\end{figure}

\bco[Matrix form of problem~\eqref{eq_fl_equil}] \label{co_fr} Let $\tv_\sd$,
$\sd \in \FKhint$, be the basis functions of the space $\tV_{h,
\mathrm{0}}^\elm$, and $\tv_\sd$, $\sd \in \FKhext$, be the remaining basis
functions of the space $\RT_0(\TK) \cap \Hdvi{\elm}$\cor{, see~\eqref{eq_RT0_bas} and Figures~\ref{fig_RTN_bas} and~\ref{fig_face_bar_vert}}. Similarly, let $q_i$, $1 \leq i \leq
|Q_{h, 0}^\elm|$, be the basis functions of the space $Q_{h, 0}^\elm$ \cor{(piecewise constants with mean value $0$)}.
Consider these basis functions of $\tV_{h, \mathrm{0}}^\elm$ and $Q_{h,
0}^\elm$ as test functions $\tv_h$ and $q_h$ in~\eqref{eq_fl_equil} and develop 
\[
    \tu_h|_\elm = \underbrace{\sum_{\sd \in \FKhext} (\alg{U}_\elm^{\mathrm {ext}})_\sd \tv_\sd}_{\tu_{h,\elm}^{\mathrm {ext}}} + \underbrace{\sum_{\sd
\in \FKhint} (\alg{U}_\elm^{\mathrm {int}})_\sd \tv_\sd}_{\tu_{h, \elm}^{\mathrm {int}}} 
\]
and 
\[
    p_h|_\elm = \Pelm + \underbrace{\sum_{i=1}^{|Q_{h, 0}^\elm|} (\Peo)_i q_i}_{p_{h,\elm}^0},
\]
where $\{\Usd\}_{\sd \in \FK}$ and $\Pelm$ are the given data of the problem~\eqref{eq_loc_probl} and \cor{the data} $\alg{U}_\elm^{\mathrm {ext}} \eq \{\Usd\}_{\sd \in \FKhext}$ \cor{appearing in problem~\eqref{eq_arg_min}} are given by Definition~\ref{def_face_fluxes}.
Then~\eqref{eq_fl_equil} in matrix form corresponds to: find
$\alg{U}_\elm^{\mathrm {int}} \in \RR^{|\tV_{h, \mathrm{0}}^\elm|}$ and $\Peo
\in \RR^{|Q_{h, 0}^\elm|}$ such that
\begin{equation} \label{eq_mixed_eq_disc_loc}
 \left (
 \begin{array}{cc}
    \matr{A}_\elm^{\mathrm {int}, \mathrm {int}}  & (\matr{B}_\elm^{0,\mathrm{int}})^{\mathrm{t}} \\
    \matr{B}_\elm^{0,\mathrm{int}}  & 0  \\
 \end{array}
 \right )
 \left (
 \begin{array}{l}
    \alg{U}_\elm^{\mathrm {int}} \\
    \Peo \\
 \end{array}
 \right )
=
 \left (
 \begin{array}{l}
    - \matr{A}_\elm^{\mathrm {int}, \mathrm {ext}} \alg{U}_\elm^{\mathrm {ext}} \\
    - \matr{B}_\elm^{0,\mathrm{ext}} \alg{U}_\elm^{\mathrm {ext}} \\
 \end{array}
 \right ),
\end{equation}
where
\bse\label{eq_A}\bat{3}
    \matr{A}_\elm^{\mathrm {int}, \mathrm {int}} & \in \RR^{|\FKhint|\times|\FKhint|}, \qquad & (\matr{A}_\elm^{\mathrm {int}, \mathrm {int}})_{\sd, \sdt} & \eq (\Km^{-1} \tv_\sdt, \tv_\sd)_\elm \qquad & & \sd, \sdt \in \FKhint, \label{eq_A_1} \\
    \matr{A}_\elm^{\mathrm {int}, \mathrm {ext}} & \in \RR^{|\FKhint|\times|\FKhext|} & (\matr{A}_\elm^{\mathrm {int}, \mathrm {ext}})_{\sd, \sdt} & \eq (\Km^{-1} \tv_\sdt, \tv_\sd)_\elm
    &&\sd \in \FKhint, \, \sdt \in \FKhext, \label{eq_A_2}\\
    \matr{B}_\elm^{0,\mathrm{int}} & \in \RR^{|Q_{h, 0}^\elm| \times |\FKhint|} & (\matr{B}_\elm^{0,\mathrm{int}})_{i,\sdt} & \eq - (\Dv \tv_\sdt,q_i)_\elm
    && 1 \leq i \leq |Q_{h, 0}^\elm|, \, \sdt \in \FKhint, \label{eq_A_3}\\
    \matr{B}_\elm^{0,\mathrm{ext}} & \in \RR^{|Q_{h, 0}^\elm| \times |\FKhext|} & (\matr{B}_\elm^{0,\mathrm{ext}})_{i,\sdt} & \eq - (\Dv \tv_\sdt,q_i)_\elm
    && 1 \leq i \leq |Q_{h, 0}^\elm|, \, \sdt \in \FKhext. \label{eq_A_4}
\eat\ese
It will also be useful to define
\[
    (\matr{A}_\elm^{\mathrm {ext}, \mathrm {ext}})_{\sd, \sdt} \eq (\Km^{-1} \tv_\sdt, \tv_\sd)_\elm
    \qquad \sd, \sdt \in \FKhext.
\]
\eco

\cor{Definition~\ref{def_fr} is practical and we can solve problem~\eqref{eq_fl_equil} or~\eqref{eq_mixed_eq_disc_loc} to have $\tu_h$ available. We, however, want to avoid this (unless we want to compute the energy error $\norm{\tu - \tu_h}_{\Km^{-\frac{1}{2}}}$ in model cases with known exact solution $\tu$). In Theorem~\ref{thm_est_Darcy} below, we will principally only need the energy $\norm{\tu_h}_{\Km^{-\frac{1}{2}},\elm}$. In this respect, we observe that} 
Definition~\ref{def_fr} lifts the information from the boundary of each polytopal
element $\elm \in \T_H$ given by the fluxes $\alg{U}_\elm^{\mathrm {ext}}$ to
the interior of the element $\elm$. 
It is thus clear that $\norm{\tu_h}_{\Km^{-\frac{1}{2}},\elm}$ must only depend on
$\alg{U}_\elm^{\mathrm {ext}}$ \cor{and can be obtained directly therefrom}. It turns out that this link can be
expressed by a single element matrix constructed solely from the geometry of the simplicial submesh $\TK$ of $\elm$, using Construction~\ref{co_fr}. 
As in~\cite[proof of Theorem~7.3]{Voh_Wohl_MFE_1_unkn_el_rel_13}, we indeed
have:

\bl[Energy norm and MFE element matrix] \label{lem_en_norm_MFE} Let a discretization scheme of the form of Assumption~\ref{as_polyt_disc} be given, leading to the normal face fluxes $\alg{U}_\elm^{\mathrm {ext}} = \{\Usd\}_{\sd \in \FK}$, $\elm \in \T_H$, and consequently to $\{\Usd\}_{\sd \in \FKhext}$ as per Definition~\ref{def_face_fluxes}. For each
$\elm \in \T_H$, let $\tu_h|_\elm \in \tV_{h, \mathrm{N}}^\elm$ be given by
Definition~\ref{def_fr}. Relying on Construction~\ref{co_fr}, define the {\em element matrix} of the mixed finite element method
\be \label{eq_loc_matr_MFE}
    \widehat{\matr{A}}_{\mathrm{MFE},\elm} \eq \matr{A}_\elm^{\mathrm {ext}, \mathrm {ext}} - \left (
 \begin{array}{l}
    \matr{A}_\elm^{\mathrm {int}, \mathrm {ext}} \\
    \matr{B}_\elm^{0,\mathrm{ext}} \\
 \end{array}
 \right )^{\mathrm{t}}
 \left (
     \begin{array}{cc}
    \matr{A}_\elm^{\mathrm {int}, \mathrm {int}}  & (\matr{B}_\elm^{0,\mathrm{int}})^{\mathrm{t}} \\
    \matr{B}_\elm^{0,\mathrm{int}}  & 0  \\
 \end{array}
 \right )^{-1}  \left (
 \begin{array}{l}
    \matr{A}_\elm^{\mathrm {int}, \mathrm {ext}} \\
    \matr{B}_\elm^{0,\mathrm{ext}} \\
 \end{array}
 \right )\cor{,}
\ee
\cor{where the building matrices are given by~\eqref{eq_A}.}
Then
\be \label{eq_en}
    \norm{\tu_h}_{\Km^{-\frac{1}{2}},\elm}^2 = (\alg{U}_\elm^{\mathrm {ext}})^{\mathrm{t}} \widehat{\matr{A}}_{\mathrm{MFE},\elm} \alg{U}_\elm^{\mathrm {ext}}.
\ee
\el

\bp Let $\elm \in \T_H$ \cor{be a fixed polytopal mesh element}. Following Construction~\ref{co_fr}, decompose
$\tu_h|_\elm = \tu_{h,\elm}^{\mathrm {ext}} + \tu_{h, \elm}^{\mathrm {int}}$
and $p_h|_\elm = \Pelm + p_{h,\elm}^0$. Note that choosing $\tu_{h,
\elm}^{\mathrm {int}} \in \tV_{h, \mathrm{0}}^\elm$ as the test function
in~\eqref{eq_fl_equil_1} and $p_{h,\elm}^0 \in Q_{h, 0}^\elm$ as the test
function in~\eqref{eq_fl_equil_2}, one has
\[
    (\Km^{-1} \tu_h, \tu_{h, \elm}^{\mathrm {int}})_\elm \cor{\reff{eq_fl_equil_1}=(p_h, \Dv
    \tu_{h, \elm}^{\mathrm {int}})_\elm =} (p_{h,\elm}^0, \Dv
    \tu_{h, \elm}^{\mathrm {int}})_\elm \cor{\reff{eq_fl_equil_2}=} - (p_{h,\elm}^0, \Dv \tu_{h,
    \elm}^{\mathrm {ext}})_\elm,
\]
where we have \cor{also} used $(\Pelm, \Dv \tu_{h, \elm}^{\mathrm {int}})_\elm = 0$ \cor{in the middle equality}.
Consequently,
\ban
    \norm{\tu_h}_{\Km^{-\frac{1}{2}},\elm}^2 & = {(\Km^{-1} \tu_h, \tu_h)_\elm} \\
    & = (\Km^{-1} \tu_h, \tu_{h,\elm}^{\mathrm {int}})_\elm + (\Km^{-1} \tu_{h,\elm}^{\mathrm {int}}, \tu_{h,\elm}^{\mathrm {ext}})_\elm
        + (\Km^{-1} \tu_{h,\elm}^{\mathrm {ext}}, \tu_{h,\elm}^{\mathrm {ext}})_\elm \\
    & = - (p_{h,\elm}^0, \Dv \tu_{h, \elm}^{\mathrm {ext}})_\elm + (\Km^{-1} \tu_{h,\elm}^{\mathrm {int}}, \tu_{h,\elm}^{\mathrm {ext}})_\elm +
        (\Km^{-1} \tu_{h,\elm}^{\mathrm {ext}}, \tu_{h,\elm}^{\mathrm {ext}})_\elm\\
    & = (\alg{U}_\elm^{\mathrm {ext}})^{\mathrm{t}} \left (
 \begin{array}{l}
    \matr{A}_\elm^{\mathrm {int}, \mathrm {ext}} \\
    \matr{B}_\elm^{0,\mathrm{ext}} \\
 \end{array}
 \right )^{\mathrm{t}} \left (
 \begin{array}{l}
    \alg{U}_\elm^{\mathrm {int}} \\
    \Peo \\
 \end{array}
 \right ) + (\alg{U}_\elm^{\mathrm {ext}})^{\mathrm{t}} \matr{A}_\elm^{\mathrm {ext}, \mathrm {ext}} \alg{U}_\elm^{\mathrm
    {ext}}.
\ean
Combining this with the MFE matrix form~\eqref{eq_mixed_eq_disc_loc} and the
MFE element matrix definition~\eqref{eq_loc_matr_MFE} yields the claim~\eqref{eq_en}. \ep

\subsection*{A fictitious flux reconstruction under Assumption~\ref{as_polyt_disc_mod}} 

There exist other ways how to lift the face normal fluxes
$\alg{U}_\elm^{\mathrm {ext}}$ into the interior of each polytope $\elm \in
\T_H$ than via Definition~\ref{def_fr}. For instance, one could use a finite volume discretization of~\eqref{eq_loc_probl} (like~\eqref{eq_FV_scheme}) in place of the mixed finite element one~\eqref{eq_fl_equil}. Under Assumption~\ref{as_polyt_disc_mod}\cor{, which defines $\tilde \tu_h$}, a much more direct and simple procedure suggests itself: to use directly the lifting operator from~\eqref{eq_lift_K}. This is in particular the theoretical vehicle in mimetic finite differences according to, \eg, \cite[Theorem~5.1]{Brez_Lip_Shash_MFD_cvg_05}. In contrast to~\eqref{eq_arg_min}, \cor{where we can solve for $\tu_h$ if needed,} the reconstruction $\tilde \tu_h$ \cor{we now introduce} may be truly fictitious in the sense that it may not be practically possible to construct it. \cor{Thus, we will not be able to compute} the actual error $\norm{\tu - \tilde \tu_h}_{\Km^{-\frac{1}{2}}}$ in model cases with known exact flux $\tu$, but otherwise this not be an issue (we will still be able to compute the a posteriori error estimate of Theorem~\ref{thm_est_Darcy} below). In particular, thanks to~\eqref{eq_en_K}, we have the following trivial counterpart of Lemma~\ref{lem_en_norm_MFE}: 

\bl[Energy norm and given discretization scheme element matrix]
\label{lem_en_norm_MFD} Let a discretization scheme of the form of
Assumption~\ref{as_polyt_disc_mod} be given, leading to the element matrices
$\widehat{\matr{A}}_\elm$ and normal face fluxes $\alg{U}_\elm^{\mathrm
{ext}} = \{\Usd\}_{\sd \in \FK}$, $\elm \in \T_H$. Then the lifting $\tilde \tu_h$ of~\eqref{eq_lift_K} satisfies
\be \label{eq_en_mod}
    \norm{\tilde \tu_h}_{\Km^{-\frac{1}{2}},\elm}^2 = (\alg{U}_\elm^{\mathrm {ext}})^{\mathrm{t}} \widehat{\matr{A}}_{\elm} \alg{U}_\elm^{\mathrm
    {ext}}.
\ee
\el

Following~\cite[Theorem~7.3]{Voh_Wohl_MFE_1_unkn_el_rel_13}, the mixed finite element element matrix $\widehat{\matr{A}}_{\mathrm{MFE},\elm}$ defined by~\eqref{eq_loc_matr_MFE}
belongs to the same family as the element matrices $\widehat{\matr{A}}_\elm$
of mimetic finite differences, mixed finite volumes, and hybrid finite
volumes of Assumption~\ref{as_polyt_disc_mod} when the simplicial
faces of $\Th$ do not subdivide the faces of the polytopal mesh $\T_H$.
This leads to the following strong interconnection of Lemmas~\ref{lem_en_norm_MFE} and~\ref{lem_en_norm_MFD} \cor{(this is an illuminating observation not needed in practice/thory)}:

\bl[Approximate energy norm] \label{lem_en_norm_appr} Let the assumptions of Lemmas~\ref{lem_en_norm_MFE} and~\ref{lem_en_norm_MFD} and of~\cite[Theorem~5.1]{Brez_Lip_Shash_MFD_cvg_05} be satisfied. Let the simplicial
faces of $\Th$ do not subdivide the faces of the polytopal mesh $\T_H$.
Then there exist two constants $c,C>0$ independent of the matrices $\widehat{\matr{A}}_{\elm}$, $\widehat{\matr{A}}_{\mathrm{MFE},\elm}$, and the mesh size but possibly depending on the shape-regularity of the simplicial submesh $\T_K$ such that
\be \label{eq_en_approx_det}
    c \underbrace{\cor{\norm{\tilde \tu_h}_{\Km^{-\frac{1}{2}},\elm}^2}}_{(\alg{U}_\elm^{\mathrm {ext}})^{\mathrm{t}} \widehat{\matr{A}}_{\elm} \alg{U}_\elm^{\mathrm
    {ext}}} \leq \underbrace{\norm{\tu_h}_{\Km^{-\frac{1}{2}},\elm}^2}_{\cor{(\alg{U}_\elm^{\mathrm {ext}})^{\mathrm{t}} \widehat{\matr{A}}_{\mathrm{MFE},\elm} \alg{U}_\elm^{\mathrm {ext}}}} \leq C \underbrace{\cor{\norm{\tilde \tu_h}_{\Km^{-\frac{1}{2}},\elm}^2}}_{(\alg{U}_\elm^{\mathrm {ext}})^{\mathrm{t}}
    \widehat{\matr{A}}_{\elm} \alg{U}_\elm^{\mathrm {ext}}}
\ee
that we denote by
\be \label{eq_en_approx}
    \norm{\tu_h}_{\Km^{-\frac{1}{2}},\elm}^2 \approx
    (\alg{U}_\elm^{\mathrm {ext}})^{\mathrm{t}} \widehat{\matr{A}}_{\elm}
    \alg{U}_\elm^{\mathrm {ext}}.
\ee
\el

\bp Any scheme of the form of Assumption~\ref{as_polyt_disc_mod} with
the lifting operator according to~\cite[Theorem~5.1]{Brez_Lip_Shash_MFD_cvg_05} satisfies the equality~(5.7) from~\cite{Brez_Lip_Shash_MFD_cvg_05}, \ie, \eqref{eq_en_K}. As mixed finite elements
on polytopal meshes belong to the mimetic finite difference family by
\cite[Theorem~7.3]{Voh_Wohl_MFE_1_unkn_el_rel_13}, this implies~\eqref{eq_en_approx_det}. \ep

\subsection{Potential point values} \label{sec_pot_points}

\cor{In Theorem~\ref{thm_est_Darcy} below, we will need a vector of real values with one value $\Sver$ per vertex $\ver \in \Vh$ of the simplicial submesh $\Th$ of $\T_H$ and one real value $\Ssd$ per face $\sd \in \F_{H,h}$ of the simplicial submesh $\Th$ of $\T_H$ (approximate point values of the pressure heads). We identify them here for the two Assumptions~\ref{as_polyt_disc} and~\ref{as_polyt_disc_mod}.}

\subsection*{Potential point values under Assumption~\ref{as_polyt_disc}} 

Here, under the general Assumption~\ref{as_polyt_disc}, we define the vectors $\{\Sver\}_{\ver \in \Vh}$ and $\{\Ssd\}_{\sd \in \F_{H,h}}$ by average values of $\Pelm$ in the neighboring polytopal elements:

\bd[Potential point values on polytopal meshes under Assumption~\ref{as_polyt_disc}] \label{def_pr_values} Let a discretization scheme of the form of Assumption~\ref{as_polyt_disc} be given, leading to the elementwise pressure heads $\alg{P} \eq \{\Pelm\}_{\elm \in \T_H}$. For each vertex $\ver$ of the simplicial mesh $\Th$
lying on $\pt \elm$ for some $\elm \in \T_H$ but not on $\pt \Om$, let $\T_\ver$ denote the set of
polytopal elements $\elm \in \T_H$ sharing $\ver$. We set
\bse \label{eq_prh_el} \be \label{eq_prh_el_1}
    \Sver \eq \frac{1}{|\T_\ver|} \sum_{\elm \in \T_\ver} \Pelm\cor{,}
\ee
\cor{where, recall, $|\T_\ver|$ is the cardinality of the set $\T_\ver$, \ie, the number of 
polytopal elements $\elm \in \T_H$ sharing $\ver$.} 
We also set 
\ba \label{eq_prh_el_2}
    \Sver & \eq \Pelm \text{ for \cor{any} vertex $\ver$ of $\Th$ lying inside some $\elm \in \T_H$},\\
    \Sver & \eq 0 \text{ for any vertex $\ver$ lying simultaneously on $\pt \elm$ for some $\elm \in \T_H$ and on $\pt \Om$,}
\ea \ese
where the latter condition is set in accordance with the boundary condition~\eqref{eq_Darcy_BC}. 
We will also need the face values $\{\Ssd\}_{\sd \in \F_{H,h}}$
\be \label{eq_prh_sd}
    \Ssd \eq \frac 1 d \sum_{\ver \in \Ve} \Sver
\ee
for any simplicial face $\sd \in \F_{H,h}$, where, recall, $\Ve$ collects the vertices of the given face $\sd \in \F_{H,h}$ \cor{and $d$ is the space dimension}.
Then, for each polytopal cell $\elm \in \T_H$, we let the element vector $\Selm \in \RR^{|\VKh|}$ collect the values $\Sver$ for the vertices $\ver$ of the simplicial mesh $\TK$ lying
on the boundary $\pt \elm$ and the value\cor{(s)} inside $\elm$
\be \label{eq_prh_el_vec}
    \Selm \eq \{\Sver\}_{\ver \in \VKh}
\ee
and the element vector $\SeExt \in \RR^{|\FKhext|}$ collect the values associated
with faces $\sd \in \FKhext$ 
\be \label{eq_prh_sd_vec}
    \SeExt \eq \{\Ssd\}_{\sd \in \FKhext}.
\ee
\ed

Figure~\ref{fig:nodal} gives an illustration in two space dimensions. 

\def\sz{1}

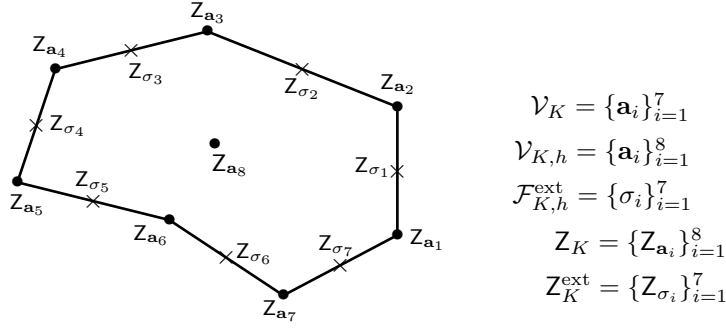
\begin{figure}
\centering
\begin{tikzpicture}
\draw (2.6,0.5)  node{$\bullet$} ; 
          \draw (2.8,0.5) node[below]{\small$\mathsf{Z}_{{\bf a}_{8}} $};
       \draw[line width=\sz pt] (5,-0.7) node{$\bullet$} 
                     -- (5,1) node[midway]{$\times$}  ;
        \draw  (5.4,-0.5) node[below]{\footnotesize$ \mathsf{Z}_{{\bf a}_1}$};
         \draw (5.05,0.2) node[left]{\footnotesize$  \mathsf{Z}_{\sigma_1} $};
        \draw[line width=\sz pt] (5,1)  node{$\bullet$} 
                     -- (2.5,2)  node[midway]{$\times$} ;
        \draw  (4.65,1.3) node[right]{\footnotesize$ \mathsf{Z}_{{\bf a}_2}$};
         \draw (3.75,1.5) node[below]{\footnotesize$  \mathsf{Z}_{\sigma_2} $};
        \draw[line width=\sz pt] (2.5,2)  node{$\bullet$} 
                      -- (0.5,1.5)  node[midway]{$\times$} ;
        \draw  (2.15,2.25) node[right]{\footnotesize$ \mathsf{Z}_{{\bf a}_3}$};
         \draw (1.7,1.7) node[below]{\footnotesize$  \mathsf{Z}_{\sigma_3} $};
        \draw[line width=\sz pt] (0.5,1.5)  node{$\bullet$} 
                     -- (0.,0)  node[midway]{$\times$} ;
        \draw  (0.8,1.8) node[left]{\footnotesize$ \mathsf{Z}_{{\bf a}_4}$};
         \draw (0.3,0.75) node[right]{\footnotesize$  \mathsf{Z}_{\sigma_4} $};
        \draw[line width=\sz pt] (0,0)  node{$\bullet$} 
                     -- (2,-0.5) node[midway]{$\times$} ;
        \draw  (0.5,-0.3) node[left]{\footnotesize$ \mathsf{Z}_{{\bf a}_5}$};
         \draw (1,-0.25) node[above]{\footnotesize$  \mathsf{Z}_{\sigma_5} $};
        \draw[line width=\sz pt] (2,-0.5)  node{$\bullet$} 
                     -- (3.5,-1.5) node[midway]{$\times$} ;
        \draw  (1.4,-0.7) node[right]{\footnotesize$ \mathsf{Z}_{{\bf a}_6}$};
         \draw (2.75,-0.95) node[right]{\footnotesize$  \mathsf{Z}_{\sigma_6} $};
        \draw[line width=\sz pt] (3.5,-1.5)  node{$\bullet$} 
                     -- (5,-0.7) node[midway]{$\times$} ;
        \draw  (3.1,-1.7) node[right]{\footnotesize$ \mathsf{Z}_{{\bf a}_7}$};
         \draw (4.2,-1.1) node[above]{\footnotesize$  \mathsf{Z}_{\sigma_7} $};
         \draw (7.8,1.) node{$\VK =  \{{\bf a}_i\}^7_{i=1} $} ;
         \draw (7.7,0.4) node{$ \VKh =  \{{\bf a}_i\}^8_{i=1}$} ;
         \draw (7.7,-0.2) node{$ \FKhext =  \{\sigma_i\}^7_{i=1}$} ;
         \draw (8.2,-0.8) node{$\mathsf{Z}_K = \{\mathsf{Z}_{{\bf a}_i}\}^8_{i=1}   $} ;
         \draw (8.15,-1.4) node{$ \mathsf{Z}^{\mathrm{ext}}_K = \{\mathsf{Z}_{\sigma_i}\}^7_{i=1}  $};
\end{tikzpicture}
\caption{Example of element vectors $\mathsf{Z}_K$ and $\mathsf{Z}^{\mathrm{ext}}_K$ of nodal {and facial potential point} values}
\label{fig:nodal}
\end{figure}

\subsection*{Potential point values under Assumption~\ref{as_polyt_disc_mod}} \label{sec_pot_rec_mod}

A more straightforward and precise potential reconstruction can be designed under Assumption~\ref{as_polyt_disc_mod}. The main idea is
that any scheme of the form~\eqref{eq_polyt_disc} can be hybridized,
following, \eg, Roberts and Thomas~\cite{Ro_Tho_91}, Brezzi and
Fortin~\cite{Brez_For_91}, or Droniou~\eal\ \cite{Dro_Ey_Gal_Her_un_10},
giving rise to one Lagrange multiplier $\Lambda_\sd$ per face $\sd \in \F_H$:

\bd[Lagrange multipliers under Assumption~\ref{as_polyt_disc_mod}] \label{def_Lag_mult} Let a discretization scheme of the form of Assumption~\ref{as_polyt_disc_mod} be given and fix a face $\sd \in \F_H$.
Consider the line associated with $\sd$ in the first block equation of~\eqref{eq_polyt_disc}. For an interior face $\sd$, let it be shared by two polytopes $\elm, \elmt \in \T_H$ such that $\tn_\sd$ points from $\elm$ to $\elmt$. From the structure of the matrices $\matr{A}$
and $\matr{B}$ supposed in Assumption~\ref{as_polyt_disc_mod}, it follows that

\[
    \sum_{\sdt \in \FK} \tn_\elm \scp \tn_\sdt (\widehat{\matr{A}}_\elm)_{\sd, \sdt} \Usdt - \Pelm
    = \sum_{\sdt \in \FKt} \tn_{\elmt, \sdt} \scp \tn_\sdt (\widehat{\matr{A}}_{\elmt})_{\sd, \sdt} \Usdt -
    \Pelmt,
\]
so that these expressions are univalued from both elements $\elm$ and $\elmt$. Such an expression is also clearly univalued on boundary faces. It allows to define the Lagrange multipliers
\be \label{eq_LM}
    \Lambda_\sd \eq \Pelm - \sum_{\sdt \in \FK} \tn_\elm \scp \tn_\sdt (\widehat{\matr{A}}_\elm)_{\sd, \sdt}
    \Usdt
\ee
for each element $\elm \in \T_H$ and each face $\sd \in \FK$.
\ed

\bd[Potential point values on polytopal meshes under Assumption~\ref{as_polyt_disc_mod}]
\label{def_hybr_pr} Let a discretization scheme of the form of Assumption~\ref{as_polyt_disc_mod} be given, leading to the normal face fluxes $\alg{U} = \{\Usd\}_{\sd \in \F_H}$ and the elementwise pressure heads $\alg{P} \eq \{\Pelm\}_{\elm \in \T_H}$. Let the Lagrange multipliers $\Lambda_\sd$ be given by Definition~\ref{def_Lag_mult}. 
For each vertex $\ver$ of the simplicial mesh $\Th$
lying on $\pt \elm$ for some $\elm \in \T_H$ but not on $\pt \Om$, let $\F_\ver$ denote the set of
polytopal faces $\sd \in \F_H$ sharing $\ver$ \cor{and $|\F_\ver|$ its cardinality}. We set
\bse \label{eq_prh_el_alt} \be \label{eq_prh_el_alt_1}
    \Sver \eq \frac{1}{|\F_\ver|} \sum_{\sd \in \F_\ver} \Lambda_\sd.
\ee
We also set 
\ba \label{eq_prh_el_alt_2}
    \Sver & \eq \Pelm \text{ for \cor{any} vertex $\ver$ of $\Th$ lying
inside some $\elm \in \T_H$},\\
    \Sver & \eq 0 \text{ for any vertex $\ver$ lying simultaneously on $\pt \elm$ for some $\elm \in \T_H$ and on $\pt \Om$,}
\ea \ese
where the latter condition is set in accordance with the boundary condition~\eqref{eq_Darcy_BC}. 
We then define the face values $\{\Ssd\}_{\sd \in \F_{H,h}}$ and the element vectors $\Selm \in \RR^{|\VKh|}$ and $\SeExt \in \RR^{|\FKhext|}$ by respectively~\eqref{eq_prh_sd}, \eqref{eq_prh_el_vec}, and~\eqref{eq_prh_sd_vec} as in Definition~\ref{def_pr_values}.
\ed

\subsection{Fictitious potential reconstruction by respecting the point values} \label{sec_pot_rec}

In Definition~\ref{def_pr_Pois}, we have reconstructed a discrete $\Hoo$-conforming potential $\prh$ from a finite volume approximation on a simplicial mesh. We extend here this procedure to general polytopal meshes. As in Section~\ref{sec_fl_rec}, in contrast to Theorem~\ref{thm_est_Pois}, however, Theorem~\ref{thm_est_Darcy} below will not need $\prh$ to evaluate the posteriori error estimate, so one again does {\em not need to perform} physically/on the computer this reconstruction. This is the reason for our naming {\em fictitious} (though we will give a closed formula for $\prh$ and $\prh$ will be discrete and possible to construct if desired).
Recall the notation from Section~\ref{sec_meshes}\cor{.} 

Given the point values $\Sver$, either prescribed by~\eqref{eq_prh_el} from Definition~\ref{def_pr_values} or by~\eqref{eq_prh_el_alt} of Definition~\ref{def_hybr_pr}, the potential reconstruction $\prh {\in \PP_1(\Th) \cap \Hoo}$ (piecewise
affine with respect to the simplicial submesh $\Th$ and $\Hoo$-conforming) is trivial:

\bd[Potential reconstruction on polytopal meshes] \label{def_pr} 
Let the potential point values $\Sver$ be given by either Definition~\ref{def_pr_values} or Definition~\ref{def_hybr_pr}.
Then the potential reconstruction $\prh {\in \PP_1(\Th) \cap \Hoo}$ is prescribed by
\[
    \prh(\ver) \eq \Sver \qquad \forall \ver \in \Vh.
\]
\ed

Note that by~\eqref{eq_prh_sd}, \cor{since $\prh$ is affine on each face,}
\be \label{eq_mean_val}
    (\SeExt)_\sd = \Ssd = \prh(\tx_\sd),
\ee
 \ie, the face values $\Ssd$ are the punctual values of the reconstruction $\prh$ in the face barycenters $\tx_\sd$, $\sd \in \FKhext$\cor{, \cf\ Figure~\ref{fig_face_bar_vert}}. 

Consider the usual hat basis functions $\psi_\ver$ of the space $\PP_1(\Th) \cap \Ho$ on the simplicial mesh $\Th$: $\psi_\ver$ is piecewise affine on $\Th$ and $\Ho$-conforming,
$\psi_\ver(\vertt) = 1$ if $\ver = \vertt$ and $0$ otherwise, where $\ver,
\vertt \in \Vh$ are the vertices of $\Th$. For any polytopal element $\elm \in \T_H$, we will need below the element {\em stiffness matrix} $\widehat{\matr{S}}_{\mathrm{FE},\elm} \in \RR^{|\VKh|
\times |\VKh|}$ defined by
\be \label{eq_loc_matr_FE}
    (\widehat{\matr{S}}_{\mathrm{FE},\elm})_{\ver, \vertt} \eq (\Km \Gr \psi_\vertt, \Gr
    \psi_\ver)_\elm \qquad \ver, \vertt \in \VKh
\ee
and the element {\em mass matrix} $\widehat{\matr{M}}_{\mathrm{FE},\elm} \in \RR^{|\VKh| \times |\VKh|}$ defined by
\be \label{eq_loc_matr_FE_mass}
    (\widehat{\matr{M}}_{\mathrm{FE},\elm})_{\ver, \vertt} \eq (\psi_\vertt, \psi_\ver)_\elm
    \qquad \ver, \vertt \in \VKh.
\ee
An immediate consequence is:

\bl[\cor{Energy norm and integration using the stiffness and mass matrices}] \label{lem_en_int_FE}
\cor{Let the potential reconstruction be given by Definition~\ref{def_pr} and the stiffness and mass matrices respectively by~\eqref{eq_loc_matr_FE} and~\eqref{eq_loc_matr_FE_mass}. Then there holds}
\be\label{eq_en_norm_pr}
    \norm{\Km \Gr \prh}_{\Km^{-\frac{1}{2}},\elm}^2 = \Selm^{\mathrm{t}} \widehat{\matr{S}}_{\mathrm{FE},\elm}
    {\Selm}
\ee
and
\be\label{eq_L2_norm_pr}
    (1, \prh)_\elm = \alg{1}^{\mathrm{t}}
    \widehat{\matr{M}}_{\mathrm{FE},\elm} {\Selm}.
\ee
\el

\br[Element matrices $\widehat{\matr{S}}_{\mathrm{FE},\elm}$ and
$\widehat{\matr{M}}_{\mathrm{FE},\elm}$] \label{rem_stiff_mass} To
obtain the finite element stiffness and mass matrices
$\widehat{\matr{S}}_{\mathrm{FE},\elm}$ of~\eqref{eq_loc_matr_FE} and
$\widehat{\matr{M}}_{\mathrm{FE},\elm}$ of~\eqref{eq_loc_matr_FE_mass}, one
could use a finite element assembly code linked to a simplicial mesh on the
polytopal cell $\elm \in \T_H$. Note, however, that there exist simple
analytical formulas for the basis functions $\psi_\ver$ of piecewise affine Lagrange
finite elements and for their gradients $\Gr \psi_\ver$ that merely
necessitate the position of the vertices of the polygonal element $\elm$ and
of the central point. In this sense, we say that neither $\Th$ nor $\TK$
need not be constructed in practice in order to obtain
$\widehat{\matr{S}}_{\mathrm{FE},\elm}$ and
$\widehat{\matr{M}}_{\mathrm{FE},\elm}$. \er

We finish this section by a remark on an alternative potential reconstruction in the spirit of Definitions~\ref{def_postpr}--\ref{def_pr_Pois}. This is still more precise, though slightly more computationally demanding; we will use it for comparison in numerical experiments below.

\br[Piecewise quadratic potential reconstruction] \label{rem_pq_quad_rec} Let
$(\tu_h|_\elm, p_h|_\elm)$ be given by Definition~\ref{def_fr} on all $\elm \in \T_H$.
We let $\pth$ be a piecewise
quadratic polynomial on the simplicial mesh $\Th$ given by
\[
    -\Km \Gr \pth|_\kappa = \tu_h|_\kappa, \qquad \frac{(\pth,1)_\kappa}{|\kappa|} = p_h|_\kappa
    \qquad \forall \kappa \in \Th.
\]
Then a piecewise quadratic potential reconstruction $\prh \in \PP_2(\Th) \cap \Hoo$ can
be obtained by averaging the values of $\pth$ in all Lagrange nodes of
freedom as per Definition~\ref{def_pr_Pois}. \er

\subsection{A guaranteed a posteriori error estimate with inexpensive implementation and evaluation}
\label{sec_a_post_darcy}

We now present our a posteriori error estimate for the model problem~\eqref{eq_Darcy}. Our main result is given for the error between the exact Darcy velocity $\tu \in \Hdv$ of~\eqref{eq_Darcy_WF}--\eqref{eq_velocity} and the reconstruction $\tu_h \in \RT_0(\Th) \cap \Hdv$ of Definition~\ref{def_fr} \cor{evaluated} as per~\eqref{eq_en_norm_Darcy}, under the general Assumption~\ref{as_polyt_disc}:

\begin{ctheorem}{A guaranteed a posteriori error estimate with inexpensive implementation and evaluation}{thm_est_Darcy}
\cor{Let $f$ be constant on each $\elm \in \T_H$ and} let the exact flux $\tu$ be given by~\eqref{eq_Darcy_WF}--\eqref{eq_velocity}. 
For any polytopal discretization satisfying Assumption~\ref{as_polyt_disc}, let, for each polytopal element $\elm \in \T_H$, the element vectors of face normal fluxes $\alg{U}_\elm^{\mathrm {ext}}$ be given by Definition~\ref{def_face_fluxes}.
Let the element vectors of potential point values $\Selm$ and $\SeExt$ be given by Definition~\ref{def_pr_values}.
Let finally the matrices $\widehat{\matr{A}}_{\mathrm{MFE},\elm}$,
$\widehat{\matr{S}}_{\mathrm{FE},\elm}$, and
$\widehat{\matr{M}}_{\mathrm{FE},\elm}$ be respectively defined
by~\eqref{eq_loc_matr_MFE}, \eqref{eq_loc_matr_FE},
and~\eqref{eq_loc_matr_FE_mass}. 
Then there holds
\begin{equation} \label{eq_est_Darcy}
    \norm{\tu - \tu_h}_{\Km^{-\frac{1}{2}}} \leq \Biggl\{\sum_{\elm \in \T_H} \eta_\elm^2\Biggr\}^\ft,
\end{equation}
where
\be \label{eq_est}
    \eta_\elm^2 \eq (\alg{U}_\elm^{\mathrm {ext}})^{\mathrm{t}} \widehat{\matr{A}}_{\mathrm{MFE},\elm} \alg{U}_\elm^{\mathrm {ext}} +
        \Selm^{\mathrm{t}} \widehat{\matr{S}}_{\mathrm{FE},\elm} \Selm
    + 2 (\alg{U}_\elm^{\mathrm {ext}})^{\mathrm{t}} \SeExt
    - 2 {\alg{F}_\elm |\elm|^{-1}} \alg{1}^{\mathrm{t}} \widehat{\matr{M}}_{\mathrm{FE},\elm} \Selm.
\ee
Here the flux reconstruction $\tu_h \in \RT_0(\Th) \cap \Hdv$ with $\Dv \tu_h = f$ is obtained following Definition~\ref{def_fr}. 
\end{ctheorem}

\bp
Recall \cor{from Lemma~\ref{lem_fr}} that $\tu_h \in \Hdv$ with $\Dv \tu_h = f$. Thus, \cor{considering the diffusion tensor $\Km$ in place of the identity matrix and proceeding} as for the Prager--Synge equality of Corollary~\ref{cor_Prag_Syng}, we have
\be \label{eq_Prag_Syng_2}
    \norm{\tu - \tu_h}_{\Km^{-\frac{1}{2}}} = \min_{v \in \Hoo} \norm{\tu_h + \Km \Gr v}_{\Km^{-\frac{1}{2}}}.
\ee
Consequently, for an arbitrary $\prh \in \Hoo$,
\[
    \norm{\tu - \tu_h}_{\Km^{-\frac{1}{2}}} \leq \norm{\tu_h + \Km \Gr \prh}_{\Km^{-\frac{1}{2}}}.
\]

We now choose for $\prh$ the fictitious potential reconstruction $\prh \in \PP_1(\Th) \cap \Hoo$ of Definition~\ref{def_pr}, continuous and piecewise affine with respect to
the simplicial submesh $\Th$ and given by the nodal values of the
vector $\Sel$ from Definition~\ref{def_pr_values}. \cor{We focus on the evaluation of $\norm{\tu_h + \Km \Gr \prh}_{\Km^{-\frac{1}{2}}}$.} Developing, we obtain, for each polytopal mesh element $\elm \in \T_H$,
\[
    \norm{\tu_h + \Km \Gr \prh}_{\Km^{-\frac{1}{2}},\elm}^2 = \norm{\tu_h}_{\Km^{-\frac{1}{2}},\elm}^2
    + 2(\tu_h, \Gr \prh)_\elm + \norm{\Km \Gr \prh}_{\Km^{-\frac{1}{2}},\elm}^2.
\]
We now use Lemma~\ref{lem_en_norm_MFE} for the first term
\cor{and~\eqref{eq_en_norm_pr} from Lemma~\ref{lem_en_int_FE}} for the last one. For the middle term, recall first
that the normal components of vector fields in $\RT_0(\Th) \cap \Hdv$ are constant on
each face. Thus the Green theorem together with~\eqref{eq_prh_sd} \cor{and its consequence~\eqref{eq_mean_val}}, \cor{$\Dv \tu_h = f$ and}~\eqref{eq_F}, and~\eqref{eq_L2_norm_pr} \cor{from Lemma~\ref{lem_en_int_FE}} give
\be \label{eq_Green}
    (\tu_h, \Gr \prh)_\elm = \<\tu_h \scp \tn, \prh\>_{\pt \elm} -
    (\Dv \tu_h, \prh)_\elm
    = (\alg{U}_\elm^{\mathrm {ext}})^{\mathrm{t}} \SeExt -
    {\alg{F}_\elm |\elm|^{-1}} \alg{1}^{\mathrm{t}} \widehat{\matr{M}}_{\mathrm{FE},\elm} \Selm.
\ee
Thus the proof is finished. \ep

Note that the a posteriori error estimators $\eta_\elm$ of Theorem~\ref{thm_est_Darcy} take a simple form of {\em matrix-vector multiplication} on each polytopal mesh element $\elm \in \T_H$, only need the element vectors $\alg{U}_\elm^{\mathrm {ext}}$, $\Selm$, and $\SeExt$ together with the element matrices $\widehat{\matr{A}}_{\mathrm{MFE},\elm}$,
$\widehat{\matr{S}}_{\mathrm{FE},\elm}$, and
$\widehat{\matr{M}}_{\mathrm{FE},\elm}$, but {\em not the reconstructions} $\tu_h$ and $\prh$ (which become fictitious), and yet~\eqref{eq_est_Darcy} delivers a {\em guaranteed upper bound} on the Darcy velocity error.

For schemes of the form of Assumption~\ref{as_polyt_disc_mod}, Definition~\ref{def_hybr_pr} can alternatively be used in place of Definition~\ref{def_pr_values} in Theorem~\ref{thm_est_Darcy}. In this case, however, it is probably still more attractive to proceed as:

\bc[A simple guaranteed estimate with the given element matrices
$\widehat{\matr{A}}_{\elm}$] \label{cor_eval_elm_matr} \cor{Let $f$ be constant on each $\elm \in \T_H$ and} let \cor{the exact flux} $\tu$ be given by~\eqref{eq_Darcy_WF}--\eqref{eq_velocity}. 
For any polytopal discretization satisfying Assumption~\ref{as_polyt_disc_mod}, let, for each polytopal element $\elm \in \T_H$, the element vectors of face normal fluxes $\alg{U}_\elm^{\mathrm {ext}}$ be given by Definition~\ref{def_face_fluxes}.
Let the element vectors of potential point values $\Selm$ and $\SeExt$ be given by Definition~\ref{def_hybr_pr}.
Let the matrices $\widehat{\matr{A}}_{\elm}$ be the element matrices from Assumption~\ref{as_polyt_disc_mod}. Let finally $\widehat{\matr{S}}_{\mathrm{FE},\elm}$ and
$\widehat{\matr{M}}_{\mathrm{FE},\elm}$ be respectively defined
by~\eqref{eq_loc_matr_FE} and~\eqref{eq_loc_matr_FE_mass}. 
Then there holds
\be \label{eq_est_Darcy_mod}
    \norm{\tu - \tilde \tu_h}_{\Km^{-\frac{1}{2}}} \leq \Biggl\{\sum_{\elm
    \in \T_H} \eta_\elm^2\Biggr\}^\ft,
\ee
where
\be \label{eq_est_mod}
    \eta_\elm^2 \eq (\alg{U}_\elm^{\mathrm {ext}})^{\mathrm{t}}
    \widehat{\matr{A}}_{\elm} \alg{U}_\elm^{\mathrm {ext}} +
        \Selm^{\mathrm{t}} \widehat{\matr{S}}_{\mathrm{FE},\elm} \Selm
    + 2 (\alg{U}_\elm^{\mathrm {ext}})^{\mathrm{t}} \SeExt
    - 2 {\alg{F}_\elm |\elm|^{-1}} \alg{1}^{\mathrm{t}} \widehat{\matr{M}}_{\mathrm{FE},\elm} \Selm,
\ee
\cor{\ie, the same form as~\eqref{eq_est} but with $\widehat{\matr{A}}_{\elm}$ in place of $\widehat{\matr{A}}_{\mathrm{MFE},\elm}$.} Here the fictitious flux reconstruction $\tilde \tu_h \in \RT_0(\Th) \cap \Hdv$ with $\Dv \tu_h = f$ is obtained following~\eqref{eq_lift_K}. 
\ec

\bp The proof is similar to that of Theorem~\ref{thm_est_Darcy}, using
the facts that $\tilde \tu_h \in \Hdv$, $\Dv \tilde \tu_h = f$, and
$\tilde \tu_h \scp \tn_\sd = \Usd / |\sd|$ for all $\sd \in \FK$
following~\eqref{eq_norm_fl_K}, and replacing
Lemma~\ref{lem_en_norm_MFE} by Lemma~\ref{lem_en_norm_MFD}. \ep

Finally, relying on Lemma~\ref{lem_en_norm_appr}, we also have:

\bc[A simple estimate with the given element matrices
$\widehat{\matr{A}}_{\elm}$]\label{cor_eval_elm_matr_appr} Let the
assumptions of Corollary~\ref{cor_eval_elm_matr} hold, let $\FKhext = \FK$,
\ie, the simplicial faces of $\Th$ do not subdivide the faces of the
polytopal mesh $\T_H$, and let the flux $\tu_h \in \RT_0(\Th) \cap \Hdv$ be
constructed by Definition~\ref{def_fr}. Then
\be \label{eq_est_Darcy_approx}
    \norm{\tu - \tu_h}_{\Km^{-\frac{1}{2}}} \lesssim \Biggl\{\sum_{\elm
    \in \T_H} \eta_\elm^2\Biggr\}^\ft,
\ee
where $\eta_\elm$ is given by~\eqref{eq_est_mod} and the ``approximately less
than or equal to'' denoted by $\lesssim$ stems from the
approximation~\eqref{eq_en_approx}. \ec

A few comments are in order:

\br[Comparison of Theorem~\ref{thm_est_Darcy} with
Corollaries~\ref{cor_eval_elm_matr} and~\ref{cor_eval_elm_matr_appr}]
Theorem~\ref{thm_est_Darcy} is applicable to any scheme of the form of
Assumption~\ref{as_polyt_disc} and relies on the the mixed finite element
matrices $\widehat{\matr{A}}_{\mathrm{MFE},\elm}$ that one needs to construct
by~\eqref{eq_loc_matr_MFE}. In contrast, Corollaries~\ref{cor_eval_elm_matr}
and~\ref{cor_eval_elm_matr_appr} are only applicable to schemes of the form
of Assumption~\ref{as_polyt_disc_mod}, but use directly the matrices
$\widehat{\matr{A}}_{\elm}$. The estimates~\eqref{eq_est_Darcy}
and~\eqref{eq_est_Darcy_approx} hold for the reconstruction $\tu_h$ that one
can compute in practice if needed by Definition~\ref{def_fr},
whereas~\eqref{eq_est_Darcy_mod} only holds for the generally unavailable
lifting $\tilde \tu_h$, as
in~\cite{Bei_res_a_post_MFD_08,Bei_Manz_a_post_MFD_08}. Finally,
both~\eqref{eq_est_Darcy} and~\eqref{eq_est_Darcy_mod} are guaranteed, with the $\leq$ inequality,
whereas~\eqref{eq_est_Darcy_approx} is not. In the numerical experiments in
Section~\ref{sec_num_Darcy} below, though, it is hard to distinguish the two
estimators~\eqref{eq_est} and~\eqref{eq_est_mod}. \er

\br[Simplicial submeshes and solution of local problems] We can avoid
the physical construction of the simplicial submeshes $\TK$ under the
conditions that 1) the given scheme takes the form of
Assumption~\ref{as_polyt_disc_mod}; 2) the element matrices
$\widehat{\matr{A}}_{\elm}$ are explicitly given; 3) the stiffness matrix
$\widehat{\matr{S}}_{\mathrm{FE},\elm}$ of~\eqref{eq_loc_matr_FE} and the
mass matrix $\widehat{\matr{M}}_{\mathrm{FE},\elm}$
of~\eqref{eq_loc_matr_FE_mass} of the conforming finite element method on
$\TK$ for each polytopal cell $\elm \in \T_H$ are constructed from the
geometry of $\TK$ only, in the sense of Remark~\ref{rem_stiff_mass}; 4) we use
the estimator $\eta_\elm$~\eqref{eq_est_mod} of
Corollary~\ref{cor_eval_elm_matr} or~\ref{cor_eval_elm_matr_appr}. Note,
however, that one may need a submesh of each polytopal cell $\elm$ already to
produce the scheme element matrix $\widehat{\matr{A}}_{\elm}$. Similarly, a
typical implementation of $\widehat{\matr{A}}_{\mathrm{MFE},\elm}$~\eqref{eq_loc_matr_MFE} will also need the submesh $\TK$, though it can also be constructed from the geometry of $\TK$ only. Remark finally that behind~\eqref{eq_loc_matr_MFE},
there is the local Neumann problem of Definition~\ref{def_fr}.\er

\br[Higher-order locally conservative methods] Higher-order locally conservative methods can be treated similarly to the above exposition but are not subject of this work. \er

\subsection{Numerical experiments} \label{sec_num_Darcy}

\begin{figure}
   \centering
    \includegraphics[width=0.4\linewidth]{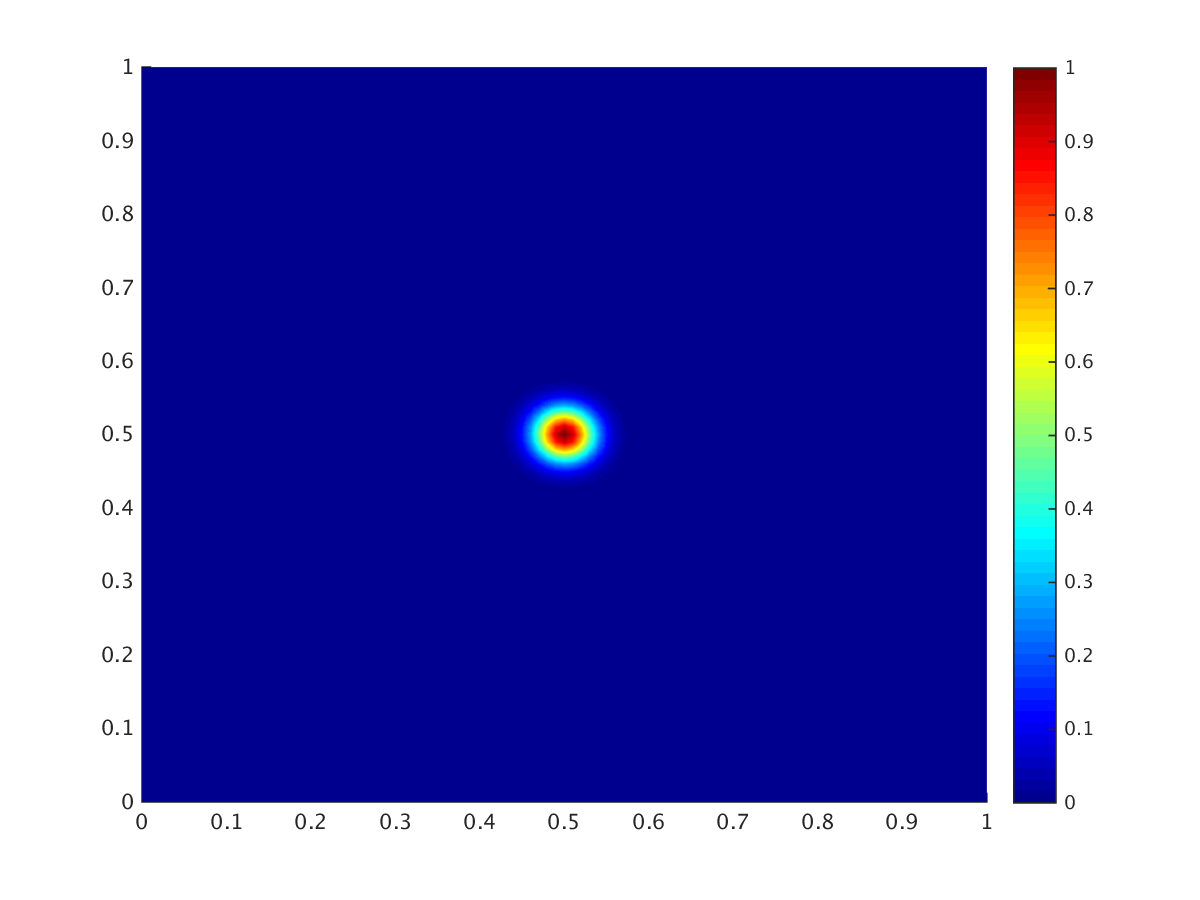}
\hspace{1.cm}
    \includegraphics[width=0.45\linewidth]{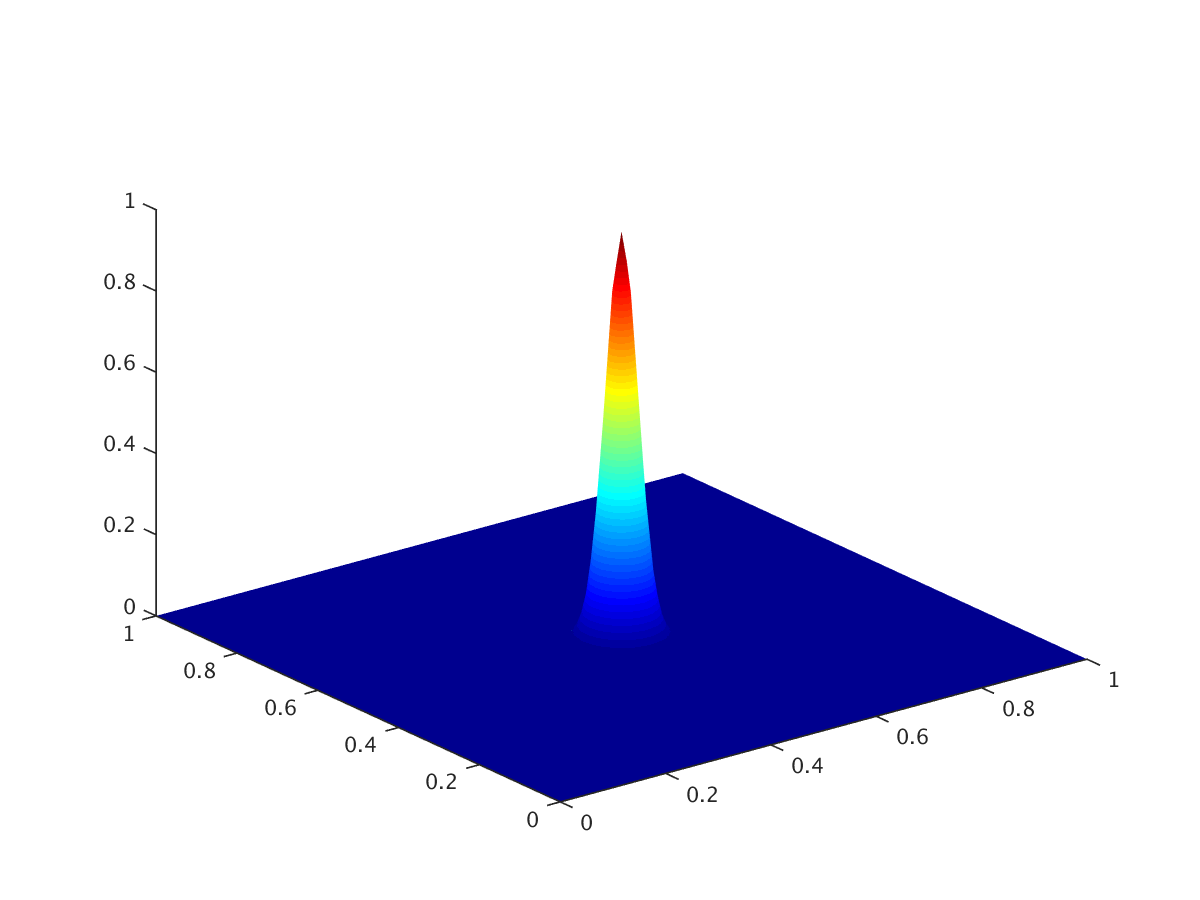}
    \caption{The \cor{exact solution~\eqref{eq_sol}}}
  \label{fig:analytic.sol}
\end{figure}

The purpose of this section is to numerically illustrate the performance of
the estimators of Theorem~\ref{thm_est_Darcy}, as well as the simplified
estimate~\eqref{eq_est_Darcy_approx} of
Corollary~\ref{cor_eval_elm_matr_appr} . We also compare this
methodology to a reference estimate on the simplicial submesh in particular
using the potential reconstruction of Remark~\ref{rem_pq_quad_rec}. The test
is taken from~\cite{Mitch_col_2D_adapt_13}, where a collection of
two-dimensional elliptic problems for testing adaptive grid refinement
algorithms is proposed. We approximate $-\Delta p = {f}$ on the space
domain $\Omega = (0,1)^2$ where the analytical solution is
\be \label{eq_sol}
    p(x,y) = 2^{4\alpha}x^{\alpha}(1-x)^{\alpha}y^{\alpha}(1-y)^{\alpha}, \quad \alpha = 200,
\ee
see Figure~\ref{fig:analytic.sol}. 
The source term $f$ prescribed correspondingly and we neglect here that it is not piecewise constant, \ie, the data oscillation term as the second term in~\eqref{eq_est_Pois_2}. We use homogeneous Dirichlet
boundary conditions.

The mesh $\T_H$ consists of general polygonal elements\cor{; together with} the fictitious
triangular submesh $\Th$, it is shown in Figure~\ref{fig:mesh}. We consider the
hybrid finite volume (HFV) discretization (Droniou~\eal\
\cite{Eym_Gal_Her_SUSHI_10}
or~\cite[Section~2.2]{Dro_Ey_Gal_Her_un_10}), taking the form~\eqref{eq_polyt_disc} \cor{of Assumption~\ref{as_polyt_disc_mod}} with the matrix $\matr{A}$ formed by local element
matrices $\widehat{\matr{A}}_\elm$. We compare three versions of a posteriori
error estimates: \cor{a) Theorem~\ref{thm_est_Darcy}}, with the
estimators $\eta_\elm$ evaluated via the matrices $\widehat{\matr{A}}_{\mathrm{MFE},\elm}$,
$\widehat{\matr{S}}_{\mathrm{FE},\elm}$, and $\widehat{\matr{M}}_{\mathrm{FE},\elm}$  (polygonal MFE estimate); \cor{b) Corollary~\ref{cor_eval_elm_matr_appr}}, using the element matrix $\widehat{\matr{A}}_{\elm}$ of the HFV scheme that is already available (polygonal HFV estimate); \cor{c) estimate~\eqref{eq_est_Darcy}} with the estimators $\eta_\elm$ replaced by the expression
$\norm{\Km^\mft\tu_h + \Km^\ft \Gr \prh}_\elm$, where $\tu_h|_\elm {\in
\tV_{h, \mathrm{N}}^\elm}$ is given in Definition~\ref{def_fr} and $\prh \in
\PP_2(\Th) \cap \Hoo$ is described in Remark~\ref{rem_pq_quad_rec} (triangular estimate\cor{, sharpest-possible, for comparison}); In the first two cases, \cor{we use Definition~\ref{def_pr_values} to obtain the vectors of the potential point values} $\Selm$ and $\SeExt$ by~\eqref{eq_prh_el} and~\eqref{eq_prh_sd}.

Note that if the mixed finite element method on polygonal meshes
of~\cite[Theorem~7.2]{Voh_Wohl_MFE_1_unkn_el_rel_13} was used instead of the
HFV discretization, the \cor{third} procedure would be,
following~\cite[Remark~7.3]{Voh_Wohl_MFE_1_unkn_el_rel_13}, fully equivalent
to solving the problem~\eqref{eq_Darcy} directly on the simplicial mesh $\Th$
by the lowest-order Raviart--Thomas mixed finite element method
and applying the a posteriori error estimates of Theorem~\ref{thm_est_Pois}.
This is why the triangular estimate serves here as a reference a posteriori error estimate. The first procedure only uses the
piecewise affine potential reconstruction $\prh$ of Definition\cor{s~\ref{def_pr_values} and}~\ref{def_pr},
which allows for the simple matrix form~\eqref{eq_est} of the estimators. The
second procedure is definitely the easiest choice in practice, where only the
already available element matrices $\widehat{\matr{A}}_\elm$ are used and
there is no need to construct the mixed finite element matrices
$\widehat{\matr{A}}_{\mathrm{MFE},\elm}$ via~\eqref{eq_loc_matr_MFE}.

\begin{figure}
   \centering
    \subfloat[Energy error]{
    \begin{tabular}{c}
    \includegraphics[width=0.22\linewidth]{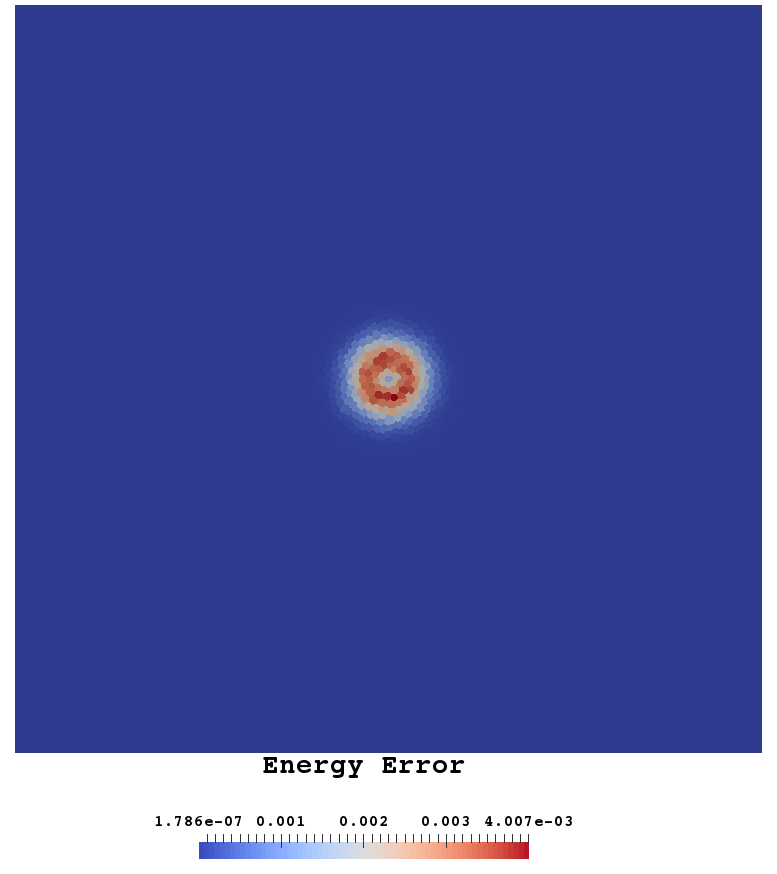}     \\
    \includegraphics[width=0.22\linewidth]{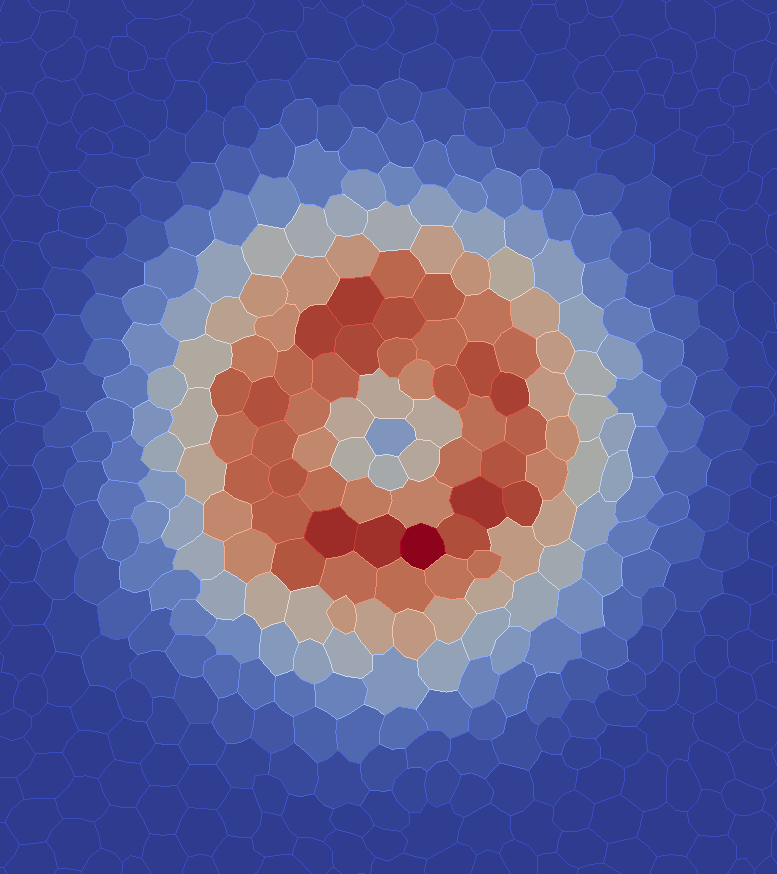}
    \end{tabular}
}
    \subfloat[Triangular estimate]{
    \begin{tabular}{c}
    \includegraphics[width=0.22\linewidth]{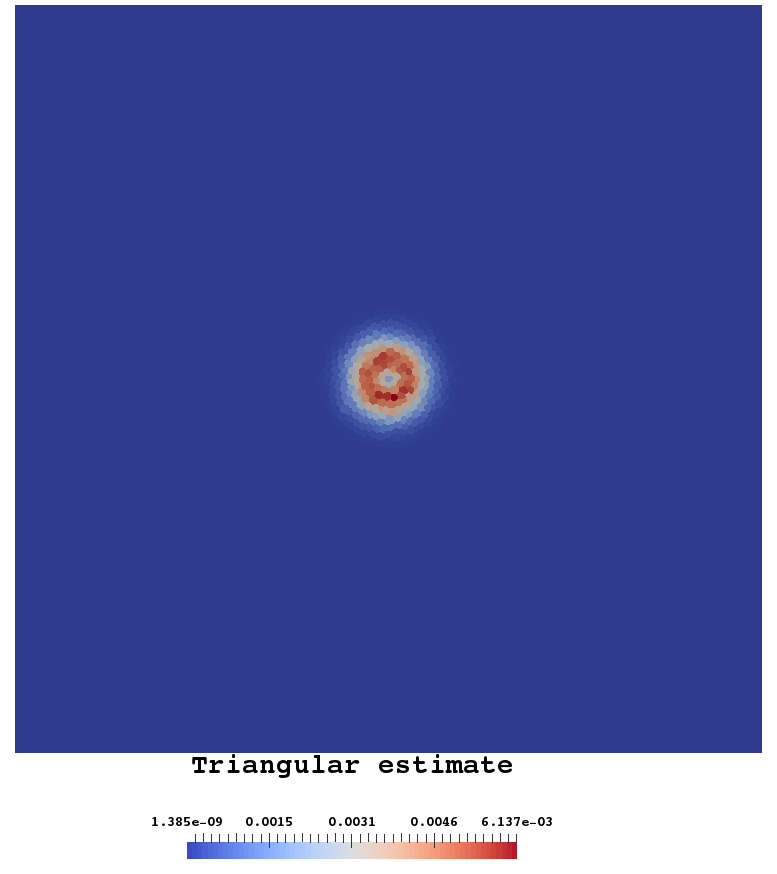}     \\
    \includegraphics[width=0.22\linewidth]{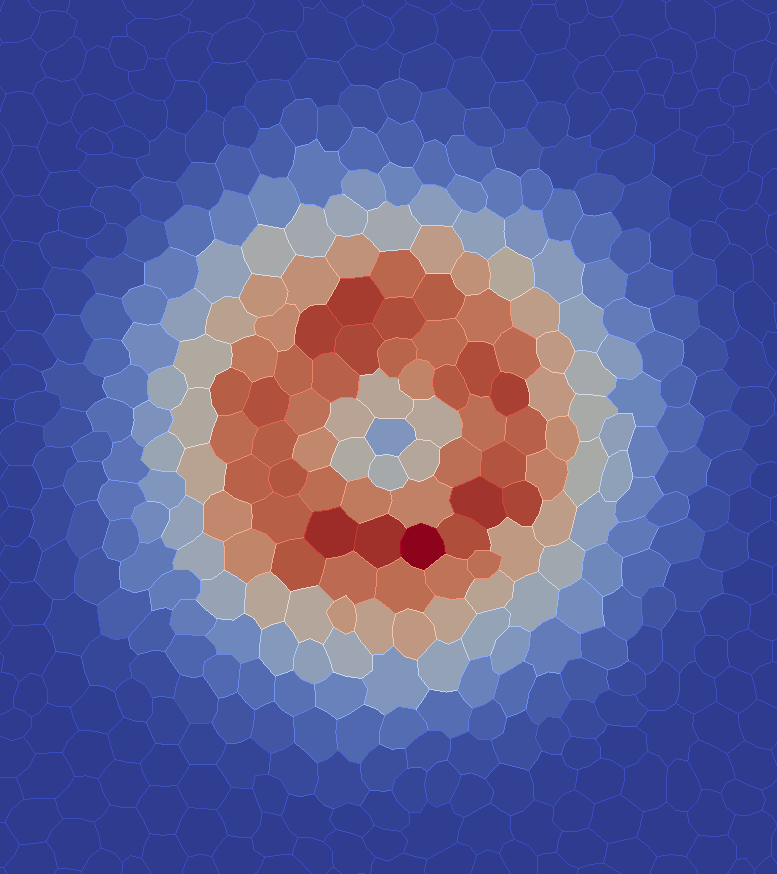}
    \end{tabular}
}
    \subfloat[Polygonal MFE estimate]{
    \begin{tabular}{c}
    \includegraphics[width=0.22\linewidth]{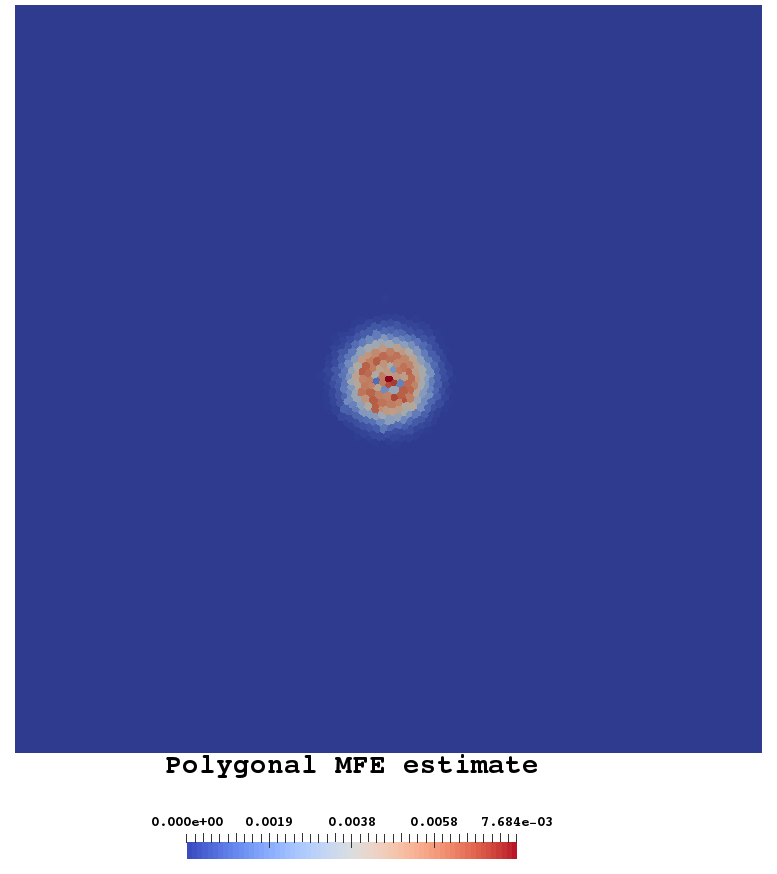}     \\
    \includegraphics[width=0.22\linewidth]{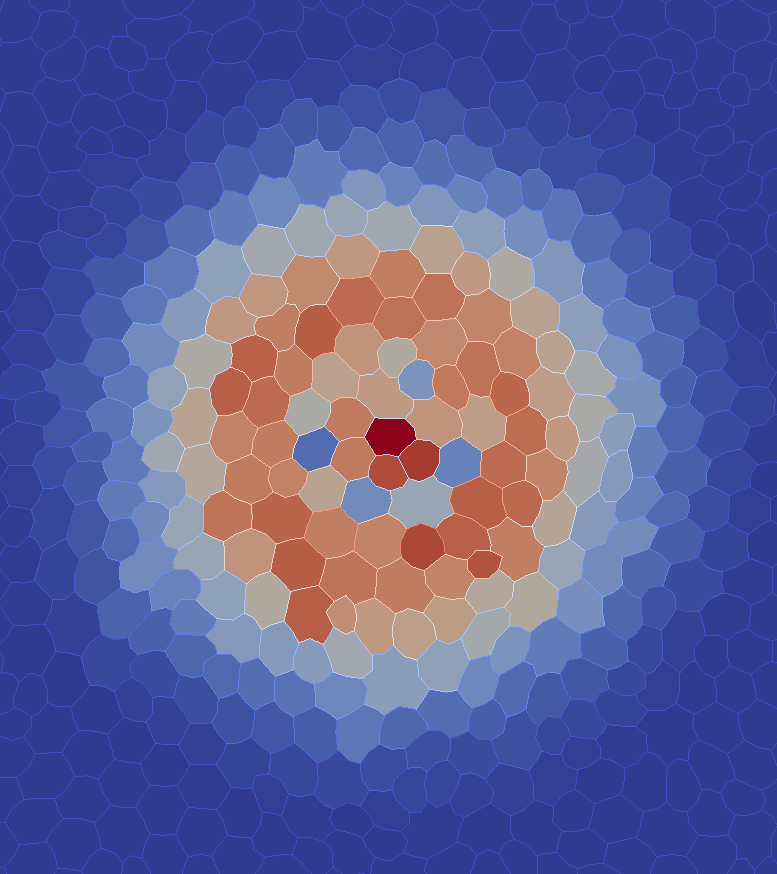}
    \end{tabular}
}
    \subfloat[Polygonal HFV estimate]{
    \begin{tabular}{c}
    \includegraphics[width=0.22\linewidth]{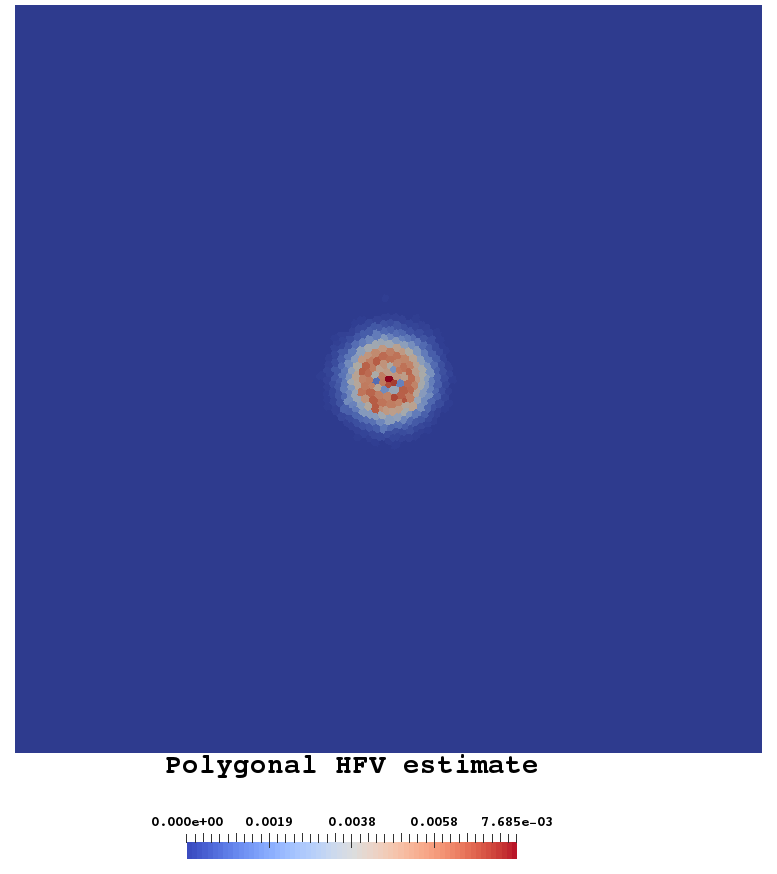}     \\
    \includegraphics[width=0.22\linewidth]{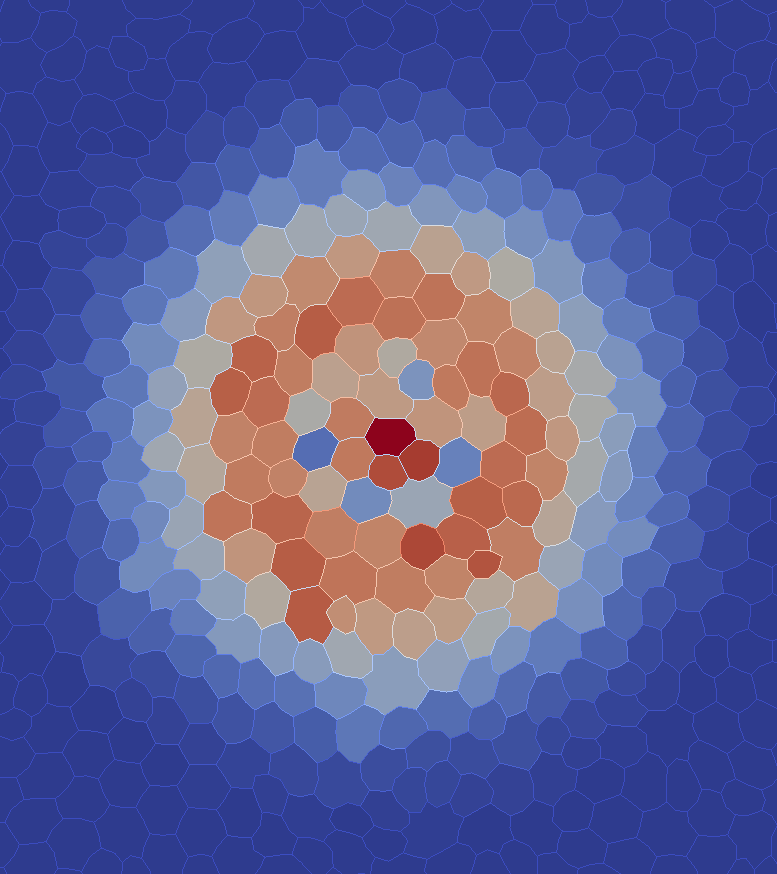}
    \end{tabular}
}
    \caption{Actual and estimated error distributions, entire domain ({\em top}) and center zoom ({\em bottom})}
  \label{fig:dist.estim}
\end{figure}

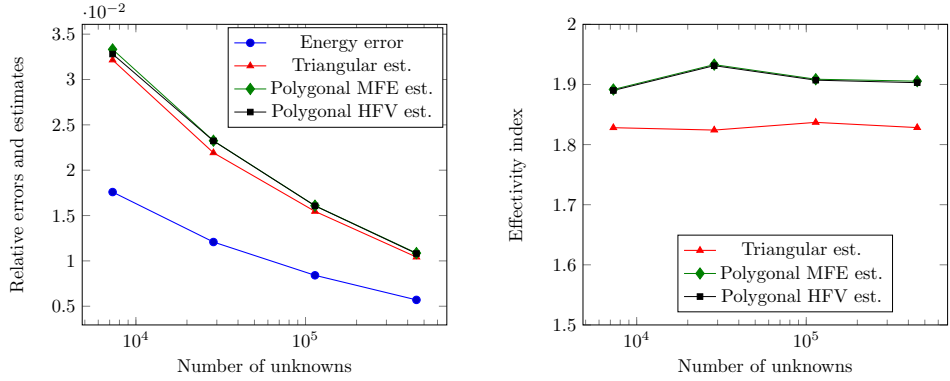
\begin{figure}
\centering
     \begin{tikzpicture}[scale=0.7]
      \begin{semilogxaxis}[
         max space between ticks=30,
          xlabel = {Number of unknowns},
          ylabel = {\cor{Relative errors and estimates}},
          legend style ={at = { (0.99,0.99)}} 
        ]
        \addplot +[line width=0.2mm] table[x=N,y=error]{Figs/relative.tex};
        \addplot +[line width=0.2mm,mark=triangle*] table[x=N,y=estimate]{Figs/relative.tex};
        \addplot +[green!50!black,  mark options={solid}, mark=diamond*, mark size=3]table[x=N,y=new]{Figs/relative.tex};
        \addplot +[mark=square*, mark size=1.5, mark options={solid}, green!0!black]table[x=N,y=modif]{Figs/relative.tex};
        \legend{Energy error, Triangular est., Polygonal MFE est., Polygonal HFV est.};
      \end{semilogxaxis}
    \end{tikzpicture}
\hspace{0.5cm}
     \begin{tikzpicture}[scale=0.7]
      \begin{semilogxaxis}[
        ymin = 1.5,
        ymax = 2,
         max space between ticks=30,
          xlabel = {Number of unknowns},
          ylabel = {Effectivity index},
          legend style ={at = { (0.85,0.3)}} 
        ]
        \addplot +[line width=0.2mm,mark=triangle*,  mark options={solid},red] table[x=N,y=estimate]{Figs/index.txt};
        \addplot +[mark=diamond*, mark size=3,  mark options={solid},green!50!black]table[x=N,y=new]{Figs/index.txt};
        \addplot +[mark=square*, mark size=1.5,  mark options={solid},green!0!black]table[x=N,y=modif]{Figs/index.txt};
        \legend{Triangular est., Polygonal MFE est., Polygonal HFV est.};
      \end{semilogxaxis}
    \end{tikzpicture}
  \caption{Relative error $\frac{\norm{\tu - \tu_h}}{\norm{\tu_h}}$ and relative estimators  and $\frac{\eta}{\norm{\tu_h}}$ (left) and effectivity indices (right), uniform mesh refinement}
\label{fig:estim.eff}
\end{figure}

In Figure~\ref{fig:dist.estim} we compare the actual and predicted error
distributions. The energy error and the \cor{sharp} triangular estimate distributions match perfectly. The polygonal MFE estimate and the
polygonal HFV estimates give similar results and match also well with the
energy error, \cor{though the error distribution around the peak is less well captured}. We depict in Figure~\ref{fig:estim.eff} the error and estimates
as a function of the total number of unknowns and the corresponding
effectivity indices for a uniform mesh refinement. This is in this test
performed as follows: we refine uniformly the triangular submesh $\Th$,
giving actually rise to a Delaunay triangulation on each step, and then we
merge triangles into polygons. All the three estimators behave in a similar
way, with a slight advantage for the triangular estimate. The graphs
confirm in particular that replacing the mixed finite element matrix
$\widehat{\matr{A}}_{\mathrm{MFE},\elm}$ in~\eqref{eq_est} by
$\widehat{\matr{A}}_\elm$ in~\eqref{eq_est_mod} has a very small
influence.

\begin{figure}
   \centering
     \begin{tikzpicture}[scale=0.7]
      \begin{axis}[
      xlabel = {Number of unknowns},
          ylabel = {\cor{Relative errors and estimates}},
          legend style ={at = { (0.99,0.99)}} 
        ]
        \addplot +[line width=0.2mm] table[x=N,y=error]{Figs/relativeAd.tex};
        \addplot +[line width=0.2mm,mark=triangle*] table[x=N,y=estimate]{Figs/relativeAd.tex};
        \addplot +[mark=diamond*,  mark options={solid}, mark size=3, green!50!black]table[x=N,y=new]{Figs/relativeAd.tex};
        \addplot +[mark=square*,  mark options={solid}, mark size=1.5, green!0!black]table[x=N,y=modif]{Figs/relativeAd.tex};
        \addplot +[line width=0.2mm] table[x=N,y=error]{Figs/relativeAd.tex};
        \legend{Energy error triang., Triangular est., Polygonal MFE est., Polygonal HFV est., Energy error polyg.};
      \end{axis}
    \end{tikzpicture}
\hspace{0.5cm}
     \begin{tikzpicture}[scale=0.7]
      \begin{axis}[
        ymin = 1,
        ymax = 2,
        xmin = 5e3,
          xlabel = {Number of unknowns},
          ylabel = {Effectivity index},
          legend style ={at = { (0.62,0.3)}} 
        ]
        \addplot +[line width=0.2mm,mark=triangle*,  mark options={solid},red] table[x=N,y=estimate]{Figs/indexAd.txt};
        \addplot +[mark=diamond*, mark size=3,  mark options={solid},green!50!black]table[x=N,y=new]{Figs/indexAdM.txt};
        \addplot +[mark=square*, mark size=1.5,  mark options={solid},green!0!black]table[x=N,y=modif]{Figs/indexAdM.txt};
        \legend{Triangular est., Polygonal MFE est., Polygonal HFV est.};
      \end{axis}
    \end{tikzpicture}
  \caption{Relative error $\frac{\norm{\tu - \tu_h}}{\norm{\tu_h}}$ and relative estimators  and $\frac{\eta}{\norm{\tu_h}}$ (left) and effectivity indices (right), adaptive mesh refinement}
\label{fig:ad.estim.eff}
\end{figure}
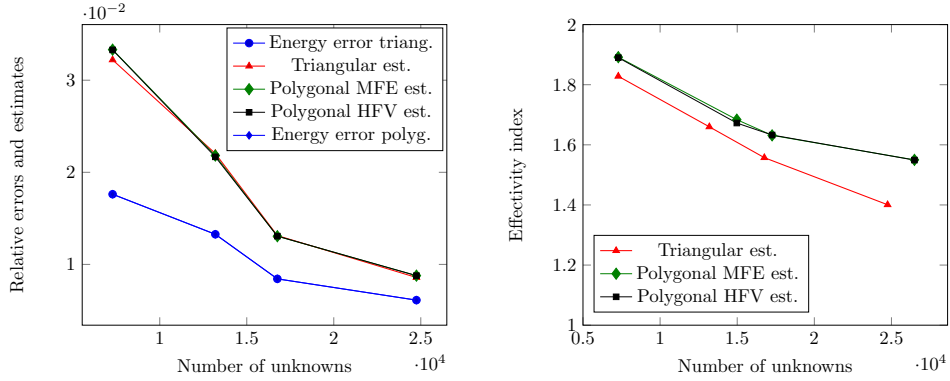

Figure~\ref{fig:ad.estim.eff} shows the results for adaptive mesh refinement,
achieved by using the local distribution of the predicted error as an
indicator to refine only the cells of the mesh where the error is important. More precisely, we refine the cells $\elm \in \T_H$ such that $\eta_\elm \geq 0.7 \max_{\elmt \in \T_H} \eta_\elmt$.
We observe quasi-identical values of the polygonal MFE and
polygonal HFV estimates; consequently, we obtain the same number of unknowns
at each step of adaptivity and almost identical final energy error for these
two estimators. The left part of Figure~\ref{fig:ad.estim.eff} displays the
results where we present the estimators and two energy errors, one resulting
with adaptivity based on the polygonal estimates and another with adaptivity
based on the triangular estimates. \cor{Almost undistinguishable} results are obtained. Finally, adaptive mesh
refinement leads to both smaller error and better effectivity indices in
comparison with uniform mesh refinement.

\subsection{Bibliographic resources}\label{sec_biblio_Darcy}

Literature on a posteriori analysis on polytopal meshes is much less plentiful than that discussed in Section~\ref{sec_biblio_Pois} for simplicial meshes. Beir{\~a}o da Veiga~\cite{Bei_res_a_post_MFD_08} and Beir{\~a}o da Veiga and Manzini~\cite{Bei_Manz_a_post_MFD_08} derive a posteriori error estimates for low-order mimetic finite difference methods; extensions are presented in Antonietti~\eal\
\cite{Ant_Bei_Lov_Ver_a_post_MFD_13}. Omnes~\eal\ \cite{Omn_Pen_Ros_a_post_ddfv_lap_09} considered the discrete duality finite volume method. 
\cor{Recent contributions on the subject include Munar~\eal\ \cite{Mun_Cang_Vel_a_post_mixed_VEM_24} or Li~\eal\ \cite{Li_Mu_Ye_a_post_WG_19}, see also the references therein.}
Here we have followed Vohral{\'{\i}}k and Yousef~\cite{Voh_Yous_polyt_18}, developing the idea of the local Neumann problems~\eqref{eq_loc_probl} from Vohral\'ik~\cite[Section~3.2]{Voh_apost_FV_08}.
Numerous \cor{additional} results were obtained for more complex model problems, which we discuss in Section~\ref{sec_biblio_MS_MC} below.

\section{Steady nonlinear pure diffusion problems. Iterative linearization and algebraic errors} \label{sec_Darcy_NL}

As a passage between the model steady linear Darcy flow of
Section~\ref{sec_Darcy} and the complex multiphase compositional Darcy flow
of Section~\ref{sec_MP_MC}, we consider here a steady nonlinear Darcy flow.
We showcase that our methodology applies here in a similar simple fashion as
in Section~\ref{sec_Darcy}. We in particular give all the details how inexact
linearization and linear solvers can be taken into account, following the
concept in Ern and Vohral{\'{\i}}k~\cite{Ern_Voh_adpt_IN_13}, giving rise to
a simple {\em adaptive inexact Newton method} on {\em polytopal meshes}. Its
particularity is that it gives a guaranteed upper bound on the total error in
the fluxes on {\em each linearization} and {\em each algebraic} solver {\em step},
distinguishes the {\em different error components}, and is obtained by {\em
simple multiplications} of the same element matrices as in the linear case by
vectors of face normal fluxes and potential values immediately available in
each lowest-order locally conservative polytopal scheme.

\subsection{Singlephase steady nonlinear Darcy flow}\label{sec_Darcy_NL_ass}

For simplicity of exposition, let us consider the following quasi-linear
version of~\eqref{eq_Darcy}: find $p: \Om \ra \RR$ such that
\bse \label{eq_Darcy_NL} \begin{empheq}[box=\widefbox]{align}
    \ds - \Dv(\Km(|\Gr p|) \Gr p) & = f \qquad \mbox{ in } \, \Om, \\
    \ds p & = 0 \qquad \mbox{ on } \, \pt \Om.
\end{empheq} \ese
We suppose that the nonlinearity can be inverted in the sense that
\be \label{eq_inv}
    \tv = - \Km(|\tw|) \tw \quad \Longleftrightarrow \quad \tw = - \Kmt (|\tv|) \tv
\ee
for all $\tv, \tw \in \RR^d$.
We suppose {\em strong monotonicity} and {\em Lipschitz-continuity}, \ie,
that there exist two positive constants $c_\Kmt, C_\Kmt$ so that for all
$\tv, \tw \in \RR^d$,
\bse \label{eq_mon_Lip} \ba
    c_\Kmt |\tv - \tw|^2 & \leq (\tv - \tw) \scp (\Kmt(|\tv|) \tv - \Kmt(|\tw|) \tw), \label{eq_mon} \\
    |\Kmt(|\tv|) \tv - \Kmt(|\tw|) \tw | & \leq C_\Kmt  |\tv - \tw|. \label{eq_cont}
\ea \ese
Moreover, we let $\Km$ and $\Kmt$ take symmetric values and, namely to
arrive at a simple matrix-vector multiplication form of the estimates, we
also suppose that, for all $\tv, \tw \in \RR^d$,
\be \label{eq_min_max_eig}
    c_\Kmt |\tv|^2 \leq \tv \scp \Kmt(|\tw|) \tv, \qquad |\Kmt(|\tw|) \tv| \leq C_\Kmt  |\tv|.
\ee
Note that for a linear problem, $\Kmt = \Km^{-1}$, $c_\Kmt$ is simply the
smallest eigenvalue of $\Kmt$, \ie, the reciprocal of the largest eigenvalue
of $\Km$, and $C_\Kmt$ is the largest eigenvalue of $\Kmt$, \ie, the
reciprocal of the smallest eigenvalue of $\Km$.

\begin{figure}
\centerline{\includegraphics[width=0.34\textwidth]{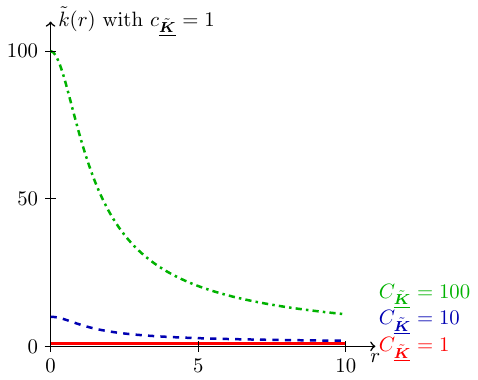} \hspace*{-0.3cm} \includegraphics[width=0.34\textwidth]{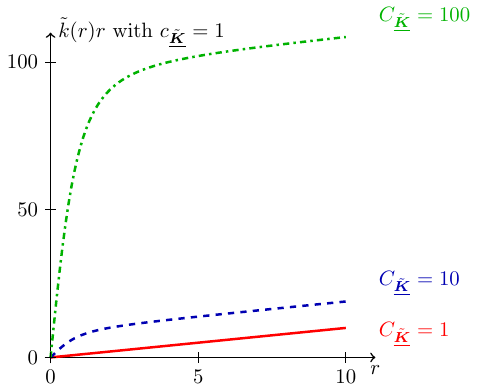} \hspace*{-0.3cm} \includegraphics[width=0.34\textwidth]{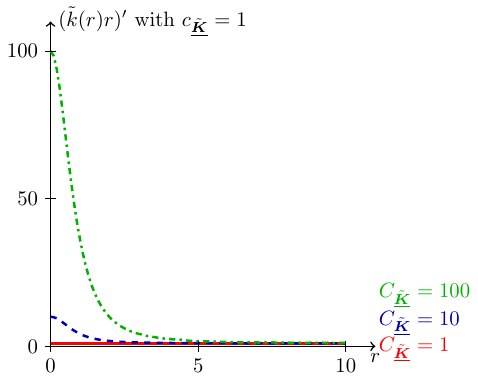}}
\caption{Prototypical monotonous and Lipschitz-continuous function $\tilde k(r) \eq c_\Kmt +\frac{C_\Kmt - c_\Kmt}{\sqrt{1+r^2}}$ for $c_\Kmt = 1$ and different values of $C_\Kmt$ (left), $\tilde k(r) r$ (middle), and $(\tilde k(r) r)'$ (right)}\label{fig_nonlin}
\end{figure}

A prototypical example satisfying~\eqref{eq_mon_Lip}--\eqref{eq_min_max_eig} is\cor{, for two positive constant $c_\Kmt$, $C_\Kmt$,}
\be\label{eq_nl_example}
    \tilde k(r) \eq c_\Kmt +\frac{C_\Kmt - c_\Kmt}{\sqrt{1+r^2}}, \qquad  \Kmt \eq \tilde k \Idd.
\ee
In this case, it can be easily verified that 
\[
    c_\Kmt \leq \tilde k(r) \leq C_\Kmt, \qquad c_\Kmt \le (\tilde k(r)r)' \leq C_\Kmt.
\]
An illustration is provided in Figure~\ref{fig_nonlin}.
A much wider class of nonlinearities is considered below in Section~\ref{sec_model}.

\subsection{Weak solution and its properties} \label{sec_prop_Darcy_NL}

In a generalization of Section~\ref{sec_prop_Darcy}, the weak formulation of~\eqref{eq_Darcy_NL} looks for $p \in \Hoo$ such that
\be \label{eq_Darcy_WF_NL}
    (\Km(|\Gr p|) \Gr p, \Gr v) = (f, v) \qquad \forall v \in \Hoo.
\ee
From the pressure head $p$, the Darcy velocity is given by
\be \label{eq_velocity_NL}
    \tu \eq - \Km(|\Gr p|) \Gr p,
\ee
and
\be \label{eq_inv_p}
    \Gr p = - \Kmt(|\tu|) \tu
\ee
from~\eqref{eq_inv}.
Still as in Proposition~\ref{pr_prop_WS}, it follows from~\eqref{eq_Darcy_WF_NL} that 
\be \label{eq_tu_Darcy_NL}
    \tu \in \Hdv, \quad \Dv \tu = f. 
\ee

\subsection{Total error and its spatial discretization, linearization, and algebraic components}

In the Poisson model problem in Section~\ref{sec_Pois}, we measure the error as $\norm{\tu - \tu_h}$, see~\eqref{eq_sp_dis_err}. In the Darcy problem in Section~\ref{sec_Darcy}, the intrinsic error measure becomes $\norm{\tu - \tu_h}_{\Km^{-\frac{1}{2}}}$, see~\eqref{eq_en_norm_Darcy}. From the mathematical setting in Section~\ref{sec_Darcy_NL_ass}, a natural extension is here to consider
\be \label{eq_en_norm_Darcy_NL}
    c_\Kmt^\ft \norm{\tu - \tu_h^{k,i}}
\ee
as the error, where the weight $c_\Kmt$ in front of the $\Lt$-norm is the monotonicity constant from~\eqref{eq_mon}. 

In~\eqref{eq_en_norm_Darcy_NL}, we have denoted the approximate solution by $\tu_h^{k,i}$, anticipating that our numerical approximation will here consist in {\em spatial discretization} (with characteristic mesh size $h$), {\em iterative linearization} (with iteration index $k$), and {\em iterative algebraic resolution} (with iteration index $i$). Our goal will be to {\em distinguish the error components} as in~\eqref{eq_err_comps} in Section~\ref{sec_err_comps} and to design an a posteriori error estimate of the form~\eqref{eq_est_rel_comps_st}, following the general concepts outlined in Section~\ref{sec_a_post_props}.

\subsection{Generic discretizations on a polytopal mesh, iterative linearization, and iterative algebraic system solution} \label{sec_model_lin.alg_NL}

Proceeding as generally as possible, we accommodate any discretization of the form of
Assumption~\ref{as_polyt_disc}. Let $\alg{F}_\elm \eq (f,1)_\elm$ for all $\elm \in \T_H$, as in~\eqref{eq_F}. For the model nonlinear problem~\eqref{eq_Darcy_NL}, the flux
balance~\eqref{eq_flux_balance} now takes the form: find the algebraic vector $\alg{P} \eq \{\Pelm\}_{\elm \in \T_H} \in \RR^{|\T_H|}$ such that
\be \label{eq_flux_balance_NL}
    \boxed{\sum_{\sd \in \FK} (\alg{U}(\alg {P}))_\sd \tn_\elm \scp \tn_\sd =
    \alg{F}_\elm \quad \forall \elm \in \T_H.}
\ee
Here $(\alg{U}(\alg {P}))_\sd \in \RR$ for each face $\sd \in \F_H$ again approximates the normal flux $\<\tu \scp \tn_\sd, 1\>_\sd$ over the face $\sd$, but it depends in a {\em nonlinear way} on the pressures $\alg{P}$, which is denoted by $\alg{U}(\alg P)$. It follows that
\eqref{eq_flux_balance_NL} is a system of $|\T_H|$ {\em nonlinear algebraic equations}
for the $|\T_H|$ unknowns $\alg {P}$.

The solution to~\eqref{eq_flux_balance_NL} is typically sought by an {\em iterative linearization} such as the Newton or the fixed point method. Given an initial guess $\alg {P}^0$, on each linearization step $k \geq 1$, this leads to a system of {\em linear algebraic equations} of the form: find the algebraic vector $\alg {P}^k \eq \{\alg {P}^k_\elm\}_{\elm \in \T_H} \in \RR^{|\T_H|}$ such that
\be \label{eq_flux_balance_NL_lin}
    \boxed{\sum_{\sd \in \FK} (\alg{U}^{k-1}(\alg {P}^k))_\sd \tn_\elm \scp \tn_\sd =
    \alg{F}_\elm \quad \forall \elm \in \T_H.}
\ee
Here, $\alg{U}^{k-1}(\alg{P}^k)$ are linearized face normal fluxes.
These are obtained from the previous iterate $\alg{P}^{k-1}$ and $\alg{U}(\alg {P^{k-1}})$ but depend in a linear (affine) way on $\alg{P}^k$. For the Newton linearization, in particular,
\be \label{eq_Newton_lin}
    (\alg{U}^{k-1}(\alg {P}^k))_\sd \eq \sum_{\elm' \in \T_H} \frac{\partial (\alg{U}(\alg {P}^{k-1}))_\sd}{\partial \alg{P}_{\elm'}}
  \cdot \big(\alg {P}^k_{\elm'}-\alg {P}^{k-1}_{\elm'}\big) + (\alg{U}(\alg {P^{k-1}}))_\sd.
\ee
%
%
In Section~\ref{sec_MP_MC} below, \eqref{eq:linear-sys} provides a general Newton linearization for a more involved problem.

The solution to~\eqref{eq_flux_balance_NL_lin} is typically sought by an {\em iterative algebraic solver}. Given an initial guess $\alg {P}^{k,0}$, \cor{typically the last algebraic solver iterate on the previous linearization step $k-1$ when $k > 1$,} on each algebraic solver step $i \geq 1$, the constant 
\be \label{eq_flux_balance_NL_lin_alg_def}
    (\alg{R})_\elm^{k,i} \eq \alg{F}_\elm - \sum_{\sd \in \FK} (\alg{U}^{k-1}(\alg {P}^{k,i}))_\sd \tn_\elm \scp \tn_\sd
\ee
is the {\em algebraic residual} on the element $\elm$,
the misfit of linearized fluxes $\alg{U}^{k-1}(\alg {P}^{k,i})$ not to
satisfy~\eqref{eq_flux_balance_NL_lin}. \cor{Consequently,} the flux balance \cor{on linearization step $k$ and algebraic step $i$} takes the form
\be \label{eq_flux_balance_NL_lin_alg}
    \boxed{\sum_{\sd \in \FK} (\alg{U}^{k-1}(\alg {P}^{k,i}))_\sd \tn_\elm \scp \tn_\sd =
    \alg{F}_\elm - (\alg{R})_\elm^{k,i} \quad \forall \elm \in \T_H\cor{,}}
\ee
\cor{\ie, \eqref{eq_flux_balance_NL_lin} up to the algebraic residual.}
In order to estimate the algebraic error, we
will, following~\cite{Gol_Strak_est_quad_94, Jir_Strak_Voh_a_post_it_solv_10,
Ern_Voh_adpt_IN_13}, also employ $j \geq 1$ additional algebraic solver steps
in~\eqref{eq_flux_balance_NL_lin_alg}, giving rise to
\be \label{eq_flux_balance_NL_lin_alg_j}
    \sum_{\sd \in \FK} (\alg{U}^{k-1}(\alg {P}^{k,i+j}))_\sd \tn_\elm \scp \tn_\sd =
    \alg{F}_\elm - (\alg{R})_\elm^{k,i+j} \quad \forall \elm \in \T_H.
\ee
\cor{The motivation is to take $j$ large enough to make the algebraic residual $(\alg{R})_\elm^{k,i+j}$ considerably smaller than $(\alg{R})_\elm^{k,i}$, see Remark~\ref{rem_mas_bal} below.}

\bex[Discretization of problem~\eqref{eq_Darcy_NL} by the finite volume scheme on a simplicial mesh with fixed-point linearization and an arbitrary iterative algebraic solver] \label{ex_FV_NL}
Consider the setting and notation of Section~\ref{sec_FV_Pois}. 
The cell-centered finite volume scheme for the nonlinear Darcy problem~\eqref{eq_Darcy_NL} on the simplicial mesh $\Th$ reads: find the real values $p_\elm$, $\elm \in \Th$, the approximations to the mean values of $p$ in the mesh elements $\elm$, such that
\bse\label{eq_FVs_NL} \begin{equation} \label{eq_FV_scheme_NL}
    \sum_{\sd \in \FK} \cor{\xi}_{\elm,\sd} U_{\elm,\sd} = (f,1)_\elm \qquad \forall \elm \in \Th,
\end{equation}
where the real number $\cor{\xi}_{\elm,\sd} U_{\elm,\sd}$ for each face $\sd \in \FK$ approximates the normal (out)flux $\<\tu \scp \tn_\elm, 1\>_\sd = \<(- \Km(|\Gr p|) \Gr p) \scp \tn_\elm, 1\>_\sd$ from $\elm$ over the face $\sd$. This is composed of the approximation of $\<- \Gr p \scp \tn_\elm, 1\>_\sd$ given, as in~\eqref{eq_dif_fl_int}--\eqref{eq_dif_fl_Dir}, by
\begin{alignat}{2}
    \label{eq_dif_fl_int_NL}
        U_{\elm,\sd} & \eq - \frac{|\sd_{\elm,\elmt}|}{d_{\elm,\elmt}}(p_\elmt - p_\elm) \qquad & & 
        \sd = \sd_{\elm,\elmt} \in \Fhint, \\
    \label{eq_dif_fl_Dir_NL}
        U_{\elm,\sd} & \eq - \frac{|\sd|}{d_{\elm,\sd}}(0-p_\elm) & & 
        \sd \in \FK \cap \Fhext
\end{alignat}
and of the approximation of $\Km(|\Gr p|)$ on $\sd$ given by
\begin{alignat}{2}
    \label{eq_coef_dif_fl_int_NL}
        \cor{\xi}_{\elm,\sd} & \eq \Km\left(\left|\frac{1}{2}(\cor{\bxi}_h|_\elm(\tx_{\elm, \sd}) + \cor{\bxi}_h|_\elmt(\tx_{\elm, \sd}))\right|\right)  \qquad & & 
        \sd = \sd_{\elm,\elmt} \in \Fhint, \\
    \label{eq_coef_dif_fl_Dir_NL}
        \cor{\xi}_{\elm,\sd} & \eq \Km(|\cor{\bxi}_h|_\elm(\tx_{\elm, \sd})|) & & 
        \sd \in \FK \cap \Fhext,
\end{alignat}
where the approximation \cor{of} $- \Gr p$ is, as in~\eqref{eq_fr_Pois_equiv}, given by
\be \label{eq_fr_Darcy_NL}
    (\cor{\bxi}_h|_\elm) (\tx) \eq \sum_{\sd \in \FK} U_{\elm,\sd} \frac{1}{d |\elm|}(\tx - \ver_{\elm, \sd}), \qquad \tx \in \elm, \, \elm \in \Th,
\ee\ese
see the illustrations in Figure~\ref{fig_face_bar_vert}.
Note that in contrast to~\eqref{eq_dif_fl_int}--\eqref{eq_dif_fl_Dir}, the face fluxes $\cor{\xi}_{\elm,\sd} U_{\elm,\sd}$ are now {\em nonlinear} functions of the pressure unknowns $p_\elm$ and $p_\elmt$ (unless $\Km({\cdot})$ is a constant function). Also note that~\eqref{eq_FVs_NL} writes in the form~\eqref{eq_flux_balance_NL} with $(\alg{U}(\alg {P}))_\sd = \cor{\xi}_{\elm,\sd} U_{\elm,\sd} \tn_\elm \scp \tn_\sd$ (recall that $\tn_\elm \scp \tn_\sd = \pm 1$ only determines the sign).

Iterative linearization of~\eqref{eq_FVs_NL} by the fixed point method then reads: given the initial values $p_\elm^0$, $\elm \in \Th$, on each linearization step $k \geq 1$, find the real values $p_\elm^k$, $\elm \in \Th$, such that
\bse\label{eq_FVs_NL_lin} \begin{equation} \label{eq_FV_scheme_NL_lin}
    \sum_{\sd \in \FK} \cor{\xi}_{\elm,\sd}^{k-1} U_{\elm,\sd}^k = (f,1)_\elm \qquad \forall \elm \in \Th,
\end{equation}
where 
\begin{alignat}{2}
    \label{eq_dif_fl_int_NL_lin}
        U_{\elm,\sd}^k & \eq - \frac{|\sd_{\elm,\elmt}|}{d_{\elm,\elmt}}(p_\elmt^k - p_\elm^k) \qquad \qquad \qquad \qquad \qquad \qquad & & 
        \sd = \sd_{\elm,\elmt} \in \Fhint, \\
    \label{eq_dif_fl_Dir_NL_lin}
        U_{\elm,\sd}^k & \eq - \frac{|\sd|}{d_{\elm,\sd}}(0-p_\elm^k) & & 
        \sd \in \FK \cap \Fhext,\\
    \label{eq_coef_dif_fl_int_NL_lin}
        \cor{\xi}_{\elm,\sd}^{k-1} & \eq \Km\left(\left|\frac{1}{2}(\cor{\bxi}_h^{k-1}|_\elm(\tx_{\elm, \sd}) + \cor{\bxi}_h^{k-1}|_\elmt(\tx_{\elm, \sd}))\right|\right)  & & 
        \sd = \sd_{\elm,\elmt} \in \Fhint, \\
    \label{eq_coef_dif_fl_Dir_NL_lin}
        \cor{\xi}_{\elm,\sd}^{k-1} & \eq \Km(|\cor{\bxi}_h^{k-1}|_\elm(\tx_{\elm, \sd})|) & & 
        \sd \in \FK \cap \Fhext, \\
\label{eq_fr_Darcy_NL_lin}
    (\cor{\bxi}_h^{k-1}|_\elm) (\tx) & \eq \sum_{\sd \in \FK} U_{\elm,\sd}^{k-1} \frac{1}{d |\elm|}(\tx - \ver_{\elm, \sd}), & & \tx \in \elm, \, \elm \in \Th.
\end{alignat}\ese
We see that~\eqref{eq_FVs_NL_lin} takes the form requested in~\eqref{eq_flux_balance_NL_lin}.

Finally, \eqref{eq_flux_balance_NL_lin_alg} for an arbitrary iterative algebraic solver follows by the definition~\eqref{eq_flux_balance_NL_lin_alg_def}.
\eex

\subsection{Face normal fluxes and potential point values} \label{sec_fl_faces_pot_NL}

The key for the guaranteed a posteriori error estimate with inexpensive implementation and evaluation for the linear Darcy problem~\eqref{eq_Darcy} in Theorem~\ref{thm_est_Darcy} were the face normal fluxes of Section~\ref{sec_fl_faces} and the potential point values of Section~\ref{sec_pot_points}.
We will now identify the face normal fluxes and potential point values for the present nonlinear Darcy problem~\eqref{eq_Darcy_NL}. In order to devise a guaranteed a posteriori error estimate distinguishing the error components, we will also identify {\em linearization error} and {\em algebraic error} face fluxes.

Let a polytopal cell $\elm \in \T_H$ with its face $\sd \in \FK$ be fixed.
For a given linearization step $k \geq 1$, algebraic step $i \geq 1$, and $j
\geq 1$ additional algebraic iterations, we set
\bse \label{eq_fl} \begin{empheq}[box=\widefbox]{align}
    (\alg{U}_{\elm}^{{k,i}})_\sd & \eq (\alg{U}(\alg {P}^{k,i}))_\sd, \label{eq_fl_disc}\\
    (\algUi{lin})_\sd & \eq (\alg{U}^{k-1}(\alg {P}^{k,i}))_\sd - (\alg{U}(\alg {P}^{k,i}))_\sd, \label{eq_fl_lin}\\
    (\algUi{alg})_\sd & \eq (\alg{U}^{k-1}(\alg {P}^{k,i+j}))_\sd - (\alg{U}^{k-1}(\alg
    {P}^{k,i}))_\sd. \label{eq_fl_alg}
\end{empheq} \ese

Observe that $\alg{U}^{k-1}(\alg {P}^{k,i+j})$ are directly available as employed in~\eqref{eq_flux_balance_NL_lin_alg_j} and that $\alg{U}^{k-1}(\alg {P}^{k,i})$ are directly available as employed in~\eqref{eq_flux_balance_NL_lin_alg}. In turn, $\alg{U}(\alg {P}^{k,i})$ are available upon plugging the current approximation $\alg {P}^{k,i}$ into the nonlinear flux function of~\eqref{eq_flux_balance_NL}. 
We call $\alg{U}_{\elm}^{{k,i}}$ {\em discretization face fluxes} since in limit $k,i \ra \infty$, they tend to $\alg{U}(\alg {P})$ of~\eqref{eq_flux_balance_NL}.
Similarly, we call $\algUi{lin}$ {\em linearization error face fluxes} since they vanish with linearization \cor{and algebraic solver} iterations $k\cor{,i} \ra \infty$. 
Finally, we call $\algUi{alg}$ the {\em algebraic error face fluxes} since they vanish with algebraic solver iterations $i \ra \infty$. 
Observe that
\[
    (\alg{U}_{\elm}^{{k,i}})_\sd + (\algUi{lin})_\sd + (\algUi{alg})_\sd = (\alg{U}^{k-1}(\alg {P}^{k,i+j}))_\sd,
\]
so that from~\eqref{eq_flux_balance_NL_lin_alg_j}
\be \label{eq_div_comps_disc}
    \sum_{\sd \in \FK} \left((\alg{U}_{\elm}^{{k,i}})_\sd + (\algUi{lin})_\sd + (\algUi{alg})_\sd\right) \tn_\elm \scp \tn_\sd =
    \alg{F}_\elm - (\alg{R})_\elm^{k,i+j} \quad \forall \elm \in \T_H.
\ee
The fluxes per simplex face are then obtained as in Definition~\ref{def_face_fluxes}.

The potential point values are obtained as in Definition~\ref{def_pr_values} or~\ref{def_hybr_pr} directly from the cell pressure heads $\alg {P}^{k,i}$. 
This gives the vectors $\Sel_\elm^{k,i}$ and $\Sel_\elm^{k,i,\mathrm{ext}}$ of {\em potential point values} in the vertices of the simplicial submesh $\TK$ of each polytopal element $\elm \in \T_H$. 

\subsection{Fictitious flux and potential reconstructions}
\label{sec_fl_pot_rec_NL}

We now extend the reconstructions of Sections~\ref{sec_fl_rec} and~\ref{sec_pot_rec} to the present nonlinear setting.

\cor{We have to face the fact that the diffusion tensor $\Km(|\Gr p|)$ is now unknown since $p$ is unknown}. We \cor{thus} modify Definition~\ref{def_fr} to
\be \label{eq_fl_NL}
    \tu_h^{k,i}\cor{|_\elm} \eq c_\Kmt^{-1} C_\Kmt^{2} \arg\min_{\tv_h \in \tV_{h, \mathrm{N}}^\elm, \, \Dv \tv_h = \cor{\mathrm{constant}} } \norm{\tv_h}_{\elm}^{2},
\ee
where the data in the spaces $\tV_{h, \mathrm{N}}^\elm$ are now given by $\alg{U}_{\elm}^{{k,i}}$ from~\eqref{eq_fl_disc} in place of $\Usd$ (we systematically use the convention~\eqref{eq_conv} (without change of the notation) to define the fluxes of the simplicial submesh faces in case they subdivide the original polytopal faces). \cor{Note that we cannot in general impose here $\Dv (\tu_h^{k,i})|_\elm = f|_\elm$.} This {\em discretization flux reconstruction} is again not needed to evaluate the a posteriori error estimators of Theorem~\ref{thm_estim_Darcy_NL} below (it would, however, be needed to evaluate the error $c_\Kmt^\ft \norm{\tu - \tu_h^{k,i}}$.) For appropriate scaling, we have replaced $\Km^{-1}$ by $c_\Kmt^{-1} C_\Kmt^2 \Idd$.
It is important to stress that just like in~\eqref{eq_arg_min}, problem~\eqref{eq_fl_NL}
is {\em linear} and can be written in a matrix form as
in~\eqref{eq_mixed_eq_disc_loc}, with the {\em same matrices}, just replacing $\Km^{-1}$ by $c_\Kmt^{-1} C_\Kmt^2 \Idd$. This in particular leads to the definition of the MFE element matrix $\widehat{\matr{A}}_{\mathrm{MFE},\elm}$ by~\eqref{eq_loc_matr_MFE} also in the present nonlinear setting\cor{, up to the replacement of $\Km^{-1}$ by $c_\Kmt^{-1} C_\Kmt^2 \Idd$}.

Similarly, we define the linearization error flux reconstruction $\tu_{\mathrm{lin},h}^{k,i}$ by the lifting~\eqref{eq_fl_NL}, employing $\algUi{lin}$ from~\eqref{eq_fl_lin} in the definition of the space $\tV_{h, \mathrm{N}}^\elm$, and the algebraic error flux reconstruction $\tu_{\mathrm{alg},h}^{k,i}$ by the lifting~\eqref{eq_fl_NL}, employing $\algUi{alg}$ from~\eqref{eq_fl_alg} in the definition of the space $\tV_{h, \mathrm{N}}^\elm$. 
The flux reconstructions $\tu_{\mathrm{lin},h}^{k,i}$ and $\tu_{\mathrm{alg},h}^{k,i}$ will only be used as a theoretical vehicle below.
From~\eqref{eq_div_comps_disc}, we observe by the Green theorem that
\be \label{eq_div_comps}
    \cor{\big(}\Dv\big(\tu_h^{k,i} + \tu_{\mathrm{lin},h}^{k,i} + \tu_{\mathrm{alg},h}^{k,i}\big)\cor{\big)|_\elm} =
    f|_\elm - |\elm|^{-1} (\alg{R})_\elm^{k,i+j} \quad \forall \elm \in \T_H.
\ee

Finally, the {\em potential reconstruction} $\prhki$ is simply obtained from the element vectors $\Sel_\elm^{k,i}$ and $\Sel_\elm^{k,i,\mathrm{ext}}$ assembled in Section~\ref{sec_fl_faces_pot_NL} as in Definition~\ref{def_pr}, \ie, setting 
\[
    \prhki(\ver) \eq \Sver^{k,i} \qquad \forall \ver \in \Vh.
\]
Congruently, for any $\elm \in \T_H$, we define the stiffness matrix $\widehat{\matr{S}}_{\mathrm{FE},\elm}$ as in~\eqref{eq_loc_matr_FE}, with $c_\Kmt^{-2} C_\Kmt \Idd$ in place of $\Km$, whereas the mass matrix $\widehat{\matr{M}}_{\mathrm{FE},\elm}$ is defined exactly as in~\eqref{eq_loc_matr_FE_mass}.

\subsection{A guaranteed a posteriori error estimate distinguishing the error
components with inexpensive implementation and evaluation} \label{sec_est_NL_comps}

Let $C_{\mathrm{F}}$ be the constant from the Friedrichs inequality
\be \label{eq_Fried}
    \norm{v} \leq C_{\mathrm{F}} h_\Om \norm{\Gr v} \qquad \forall v \in \Hoo;
\ee
there holds $1 / (\pi d) \leq C_{\mathrm{F}} \leq 1$, as $C_{\mathrm{F}}
h_\Om$ is the square root of the reciprocal of the smallest eigenvalue of the
Laplace operator on $\Om$ with homogenous Dirichlet boundary condition. The
main result of this section is:

\begin{ctheorem}{A guaranteed a posteriori error estimate distinguishing the error
components with inexpensive implementation and evaluation}{thm_estim_Darcy_NL} 
\cor{Let $f$ be constant on each $\elm \in \T_H$ and} let $\tu$ be given by~\eqref{eq_Darcy_WF_NL}--\eqref{eq_velocity_NL}. Consider any polytopal
discretization of the form~\eqref{eq_flux_balance_NL}, any
iterative linearization~\eqref{eq_flux_balance_NL_lin} on step $k \geq 1$, and any
iterative algebraic solver~\eqref{eq_flux_balance_NL_lin_alg} on step $i \geq 1$.
Consider $j \geq 1$ additional algebraic solver steps leading to~\eqref{eq_flux_balance_NL_lin_alg_j}. 
Let, for each polytopal element $\elm \in \T_H$, the element vectors of face normal fluxes $\alg{U}_{\elm}^{{k,i}}$, $\algUi{lin}$, and $\algUi{alg}$ be given by~\eqref{eq_fl}.
Let the element vectors of potential point values $\Sel_\elm^{k,i}$ and $\Sel_\elm^{k,i,\mathrm{ext}}$ be constructed from the values $\alg {P}^{k,i}$ following Section~\ref{sec_fl_faces_pot_NL}.
Let finally the element matrices $\widehat{\matr{A}}_{\mathrm{MFE},\elm}$ and
$\widehat{\matr{S}}_{\mathrm{FE},\elm}$ be defined by~\eqref{eq_loc_matr_MFE}
and~\eqref{eq_loc_matr_FE}, where one respectively takes the multiples of the
identity matrix $c_\Kmt^{-1} C_\Kmt^2 \Idd$ in place of $\Km^{-1}$ and
$c_\Kmt^{-2} C_\Kmt \Idd$ in place of $\Km$. Let also
$\widehat{\matr{M}}_{\mathrm{FE},\elm}$ be given
by~\eqref{eq_loc_matr_FE_mass}. Then there holds
\begin{equation}\label{eq:local.comps:a_NL}
    c_\Kmt^\ft \norm{\tu - \tu_h^{k,i}} \leq \Esti{sp}{k,i} + \Esti{lin}{k,i} + \Esti{alg}{k,i}
    + \Esti{rem}{k,i}
\end{equation}
with
\[
    \Esti{\bullet}{k,i} = \left\{\sum_{\elm \in \T_H}\left(\est{\bullet}{\elm}{k,i}\right)^2\right\}^\ft, \qquad \bullet = \{\mathrm{sp}, \, \mathrm{lin},
    \, \mathrm{alg}, \, \mathrm{rem} \},
\]
where the (spatial) {\em discretization estimators} are given by
\[
   \left( \est{sp}{\elm}{k,i}\right)^2 \! \eq \big(\alg{U}_{\elm}^{{k,i}}\big)^{\mathrm{t}} \widehat{\matr{A}}_{\mathrm{MFE},\elm} \alg{U}_{\elm}^{{k,i}} +
   (\Sel_\elm^{k,i})^{\mathrm{t}} \widehat{\matr{S}}_{\mathrm{FE},\elm} \Sel_\elm^{k,i}
    + 2 c_\Kmt^{-1} C_\Kmt \left[(\alg{U}_\elm^{k,i,\mathrm {ext}})^{\mathrm{t}} \Sel_\elm^{k,i,\mathrm{ext}} \!
    - {\alg{F}_\elm |\elm|^{-1}} \alg{1}^{\mathrm{t}} \widehat{\matr{M}}_{\mathrm{FE},\elm} \Sel_\elm^{k,i} \right]\!,
\]
the {\em linearization estimators} by
\[
    \left(\est{lin}{\elm}{k,i}\right)^2 \eq
    (\algUi{lin})^{\mathrm{t}} \widehat{\matr{A}}_{\mathrm{MFE},\elm} \algUi{lin},
\]
the {\em algebraic estimators} by,
\[
    \left(\est{alg}{\elm}{k,i}\right)^2 \eq (\algUi{alg})^{\mathrm{t}}
    \widehat{\matr{A}}_{\mathrm{MFE},\elm} \algUi{alg},
\]
and the {\em algebraic remainder estimators} by
\[
    \est{rem}{\elm}{k,i} \eq c_\Kmt^\mft C_\Kmt C_{\mathrm{F}} h_\Om |\elm|^\mft |(\alg{R})_\elm^{k,i+j}|.
\]
Here, the flux reconstruction $\tu_h^{k,i} \in \RT_0(\Th) \cap \Hdv$ is obtained following~\eqref{eq_fl_NL}. 
\end{ctheorem}

\br[\cor{Algebraic remainder estimators $\est{rem}{\elm}{k,i}$, choice of $j$ additional linear solver iterations, and mass balance}] \label{rem_mas_bal}
From~\eqref{eq_div_comps_disc} or~\eqref{eq_div_comps}, at ``convergence'' of the iterative algebraic solver (or for a direct solver)(recall that we neglect rounding errors), as the algebraic residuals $(\alg{R})_\elm^{k,i+j}$ vanish, we recover exact mass balance on each linearization step $k \geq 1$. More importantly, for a sufficiently high number $j$ of additional algebraic iterations, we ``almost'' have the mass balance on each linearization step $k \geq 1$. We typically choose $j$ adaptively, following~\cite{Gol_Strak_est_quad_94, Jir_Strak_Voh_a_post_it_solv_10,
Ern_Voh_adpt_IN_13}, so that this mass balance misfit related to \cor{the algebraic residuals} $(\alg{R})_\elm^{k,i+j}$ is negligible: we choose $j$ such that $\Esti{rem}{k,i}$ becomes negligible in comparison with the other estimators in~\eqref{eq:local.comps:a_NL}. 
The presence of the terms $(\alg{R})_\elm^{k,i+j}$ is the disadvantage of the $j$ additional iterations. In turn, the advantage of the $j$ of additional algebraic iterations is the simplicity of the formula~\eqref{eq_fl_alg}.
A more involved choice can be made following~\cite{Pap_Rud_Voh_Wohl_MG_20}. Then no $j$ additional iterations are needed, the term $|\elm|^{-1} (\alg{R})_\elm^{k,i+j}$ in~\eqref{eq_div_comps_disc} and~\eqref{eq_div_comps} and the term \cor{$\Esti{rem}{k,i}$} in~\eqref{eq:local.comps:a_NL} vanish, there is an exact mass balance on each linearization step $k \geq 1$ and each linear algebraic step $i \geq 1$, and the a posteriori bound~\eqref{eq:local.comps:a_NL} is typically slightly more precise. Remark finally that these two choices respectively correspond to nonzero and zero $\rho_h^{k,i}$ in the general discussion in Section~\ref{sec_mass_bal}.\er

\bp The proof follows in spirit of previous works by 
Kim~\cite{Kim_a_post_loc_cons_nonlin_07}, Jir{\'a}nek~\eal\
\cite{Jir_Strak_Voh_a_post_it_solv_10}, and Ern and
Vohral{\'{\i}}k~\cite{Ern_Voh_adpt_IN_13, Ern_Voh_p_rob_15}. Let $\pr \in \Hoo$
be the solution of
\be \label{eq_Darcy_WF_NL_proj}
    (\Km(|\Gr \pr|) \Gr \pr, \Gr v) = - (\tu_h^{k,i}, \Gr v) \qquad \forall v \in \Hoo
\ee
and set $\bzeta \eq - \Km(|\Gr \pr|) \Gr \pr$, similarly to~\eqref{eq_velocity_NL}. Remark that if there holds $\tu_h^{k,i} \in \Hdv$ with $\Dv \tu_h^{\cor{k,i}} = f$ (if~\eqref{eq_flux_balance_NL} is satisfied exactly, no iterative linearization and exact algebraic solve), then, by the Green theorem, $- (\tu_h^{k,i}, \Gr v) = (f, v)$, so that $\pr$ coincides with the solution $p$ of~\eqref{eq_Darcy_WF_NL}. 

By the triangle inequality,
\be \label{eq_triang}
    c_\Kmt^\ft \norm{\tu - \tu_h^{k,i}} \leq c_\Kmt^\ft \norm{\tu -\bzeta}
    + c_\Kmt^\ft \norm{\bzeta - \tu_h^{k,i}}.
\ee
For the first term above, the strong monotonicity implies
\ban
    c_\Kmt \norm{\tu - \bzeta}^2 \leq {} & (\Km(|\Gr p|) \Gr p - \Km(|\Gr \pr|) \Gr \pr, \Gr (p -
    \pr)) \tag*{\text{ by~\eqref{eq_mon} and~\eqref{eq_inv}}}\\
    = {} & (\Km(|\Gr p|) \Gr p + \tu_h^{k,i}, \Gr (p -
    \pr)) \tag*{\text{ by~\eqref{eq_Darcy_WF_NL_proj}}}\\
    = {} & (f,p - \pr) + (\tu_h^{k,i}, \Gr (p -
    \pr)) \tag*{\text{ by~\eqref{eq_Darcy_WF_NL}}}.
\ean
Note that by the Lipschitz-continuity~\eqref{eq_cont} together
with~\eqref{eq_inv},
\be \label{eq_bound}
    \norm{\Gr (p-\pr)} \leq C_\Kmt \norm{\tu - \bzeta}.
\ee
Thus,
\be \label{eq_dual_res} \bs
    c_\Kmt^\ft \norm{\tu - \bzeta} & = \frac{(f,p - \pr) + (\tu_h^{k,i}, \Gr (p -
    \pr))}{c_\Kmt^\ft \norm{\tu - \bzeta}} \\
    & \leq c_\Kmt^\mft C_\Kmt \max_{\vf \in \Hoo, \, \norm{\Gr \vf}=1} \{(f,\vf) + (\tu_h^{k,i}, \Gr \vf)\}.
\es \ee

To estimate~\eqref{eq_dual_res} by a computable quantity, we use~\eqref{eq_div_comps}.
Fix now $\vf \in \Hoo$ with
$\norm{\Gr \vf}=1$. Then adding and subtracting
$(\tu_{\mathrm{lin},h}^{k,i} + \tu_{\mathrm{alg},h}^{k,i}, \Gr \vf)$ and using
the Green theorem,
\ban
    (f,\vf) + (\tu_h^{k,i}, \Gr \vf) & = (f - \Dv( \tu_h^{k,i} + \tu_{\mathrm{lin},h}^{k,i} + \tu_{\mathrm{alg},h}^{k,i}),\vf)
        - (\tu_{\mathrm{lin},h}^{k,i} + \tu_{\mathrm{alg},h}^{k,i}, \Gr \vf) \\
    & = \sum_{\elm \in \T_H} \{|\elm|^{-1} ((\alg{R})_\elm^{k,i+j},\vf)_\elm - (\tu_{\mathrm{lin},h}^{k,i} + \tu_{\mathrm{alg},h}^{k,i}, \Gr
        \vf)_\elm\}.
\ean
\cor{The Cauchy--Schwarz inequality gives
\ban
    \sum_{\elm \in \T_H} -(\tu_{\mathrm{lin},h}^{k,i} + \tu_{\mathrm{alg},h}^{k,i}, \Gr \vf)_\elm & \leq \sum_{\elm \in \T_H} \norm{\tu_{\mathrm{lin},h}^{k,i} + \tu_{\mathrm{alg},h}^{k,i}}_\elm \norm{\Gr \vf}_\elm \leq \norm{\tu_{\mathrm{lin},h}^{k,i} + \tu_{\mathrm{alg},h}^{k,i}} \underbrace{\norm{\Gr \vf}}_{=1} \\
    & \leq \norm{\tu_{\mathrm{lin},h}^{k,i}} + \norm{\tu_{\mathrm{alg},h}^{k,i}}.
\ean
}
\cor{Similarly, again} the Cauchy--Schwarz inequality leads to
\be \label{eq_alg_rem_est}
    \sum_{\elm \in \T_H} |\elm|^{-1} |((\alg{R})_\elm^{k,i+j},\vf)_\elm | \leq \sum_{\elm \in \T_H} |\elm|^\mft |(\alg{R})_\elm^{k,i+j}| \norm{\vf}_\elm
    \leq \left\{\sum_{\elm \in \T_H} |\elm|^{-1} |(\alg{R})_\elm^{k,i+j}|^2\right\}^\ft \norm{\vf}_\Om.
\ee
Summing up these developments with the Friedrichs
inequality~\eqref{eq_Fried}, infer from~\eqref{eq_dual_res}
\[
    c_\Kmt^\ft \norm{\tu - \bzeta} \leq c_\Kmt^\mft C_\Kmt \left[C_{\mathrm{F}} h_\Om \left\{\sum_{\elm \in \T_H} |\elm|^{-1} |(\alg{R})_\elm^{k,i+j}|^2\right\}^\ft
    + \norm{\tu_{\mathrm{lin},h}^{k,i}} + \norm{\tu_{\mathrm{alg},h}^{k,i}} \right].
\]
Finally, to evaluate the \cor{norm $\norm{\tu_{\mathrm{lin},h}^{k,i}}$}, we use
\[
    \norm{\tu_{\mathrm{lin},h}^{k,i}}^2 = \sum_{\elm \in \T_H}
    \norm{\tu_{\mathrm{lin},h}^{k,i}}_{\elm}^2
\]
and employ the constant matrix $c_\Kmt^{-1} C_\Kmt^2 \Idd$ in place of
$\Km^{-1}$ in Lemma~\ref{lem_en_norm_MFE} to
infer
\be \label{eq_matr_NL}
    c_\Kmt^{-1} C_\Kmt^2 \norm{\tu_{\mathrm{lin},h}^{k,i}}_\elm^2 = \big(\algUi{lin}\big)^{\mathrm{t}}
    \widehat{\matr{A}}_{\mathrm{MFE},\elm} \algUi{lin}
\ee
for each polytopal element $\elm \in \T_H$. We proceed similarly for \cor{$\norm{\tu_{\mathrm{alg},h}^{k,i}}$}, which gives
\be \label{eq_matr_alg}
    c_\Kmt^{-1} C_\Kmt^2 \norm{\tu_{\mathrm{alg},h}^{k,i}}_\elm^2 = \big(\algUi{alg}\big)^{\mathrm{t}}
    \widehat{\matr{A}}_{\mathrm{MFE},\elm} \algUi{alg}.
\ee

We finally estimate the second term in~\eqref{eq_triang}. For $v \in \Hoo$
arbitrary, the strong monotonicity in particular leads to
\ban
    c_\Kmt \norm{\bzeta - \tu_h^{k,i}}^2 & \leq (\Km(|\Gr \pr|) \Gr \pr + \tu_h^{k,i}, \Gr \pr + \Kmt(|\tu_h^{k,i}|)
    \tu_h^{k,i}) \tag*{\text{ by~\eqref{eq_mon} and~\eqref{eq_inv}}} \\
    & = (\Km(|\Gr \pr|) \Gr \pr + \tu_h^{k,i}, \Gr v + \Kmt(|\tu_h^{k,i}|)
    \tu_h^{k,i}) \tag*{\text{ by~\eqref{eq_Darcy_WF_NL_proj}}} \\
    & \leq \norm{\bzeta - \tu_h^{k,i}} \norm {\Kmt(|\tu_h^{k,i}|) \tu_h^{k,i} + \Gr v} \tag*{\text{ by~Cauchy--Schwarz}}
\ean
so that,
\bse \label{eq_min_NC} \ba
    c_\Kmt^\ft \norm{\bzeta - \tu_h^{k,i}} & \leq c_\Kmt^\mft \min_{v \in \Hoo} \norm{\Kmt(|\tu_h^{k,i}|) \tu_h^{k,i} + \Gr
    v} \label{eq_min} \\
    & \leq c_\Kmt^\mft \norm{\Kmt(|\tu_h^{k,i}|) \tu_h^{k,i} + \Gr
    \prhki} \label{eq_NC}
\ea \ese
for an arbitrary $\prhki \in \Hoo$. We pick the potential reconstruction from Section~\ref{sec_fl_pot_rec_NL}. For the inexpensive implementation and evaluation, we proceed as in Section~\ref{sec_a_post_darcy}. Consequently, for each polytopal mesh
element $\elm \in \T_H$, employing here also the eigenvalue
hypothesis~\eqref{eq_min_max_eig},
\ban
    {} & \norm{\Kmt(|\tu_h^{k,i}|) \tu_h^{k,i} + \Gr \prhki}_{\elm}^2\\
    =  {} & \norm{[\Kmt(|\tu_h^{k,i}|)]^\ft \big([\Kmt(|\tu_h^{k,i}|)]^\ft
    \tu_h^{k,i} + [\Kmt(|\tu_h^{k,i}|)]^\mft \Gr \prhki\big)}_{\elm}^2\\
    \leq {} & C_\Kmt \norm{[\Kmt(|\tu_h^{k,i}|)]^\ft \tu_h^{k,i} + [\Kmt(|\tu_h^{k,i}|)]^\mft \Gr \prhki}_{\elm}^2 \tag*{\text{ by~\eqref{eq_min_max_eig}}}\\
    = {} & C_\Kmt \left((\Kmt(|\tu_h^{k,i}|)\tu_h^{k,i}, \tu_h^{k,i})_\elm + 2 (\tu_h^{k,i}, \Gr \prhki)_\elm + ([\Kmt(|\tu_h^{k,i}|)]^{-1} \Gr \prhki, \Gr
    \prhki)_\elm \right)\\
    \leq {} & C_\Kmt^2 \norm{\tu_h^{k,i}}_\elm^2 + 2 C_\Kmt (\tu_h^{k,i}, \Gr \prhki)_\elm
    + c_\Kmt^{-1} C_\Kmt \norm{\Gr \prhki}_\elm^2 \tag*{\text{ by~\eqref{eq_min_max_eig}}}.
\ean
For the first term above, we proceed as in~\eqref{eq_matr_NL}. For the second
one, the development~\eqref{eq_Green} remains unchanged. Finally, for the
third term, we use~\eqref{eq_en_norm_pr} where $\Km$ has been replaced by
$c_\Kmt^{-2} C_\Kmt \Idd$ to see that
\[
    c_\Kmt^{-2} C_\Kmt \norm{\Gr \prhki}_\elm^2 = (\Sel_\elm^{k,i})^{\mathrm{t}} \widehat{\matr{S}}_{\mathrm{FE},\elm} \Sel_\elm^{k,i}.
\]
The proof is finished by combining the two estimates on $c_\Kmt^\ft \norm{\tu
-\bzeta}$ and $c_\Kmt^\ft \norm{\bzeta - \tu_h^{k,i}}$. \ep

\br[Inexpensive implementation and evaluation]
Note that the estimate of Theorem~\ref{thm_estim_Darcy_NL} still takes the
same simple matrix-vector multiplication form of Theorem~\ref{thm_est_Darcy},
with in particular the same element matrices
$\widehat{\matr{A}}_{\mathrm{MFE},\elm}$,
$\widehat{\matr{S}}_{\mathrm{FE},\elm}$,
$\widehat{\matr{M}}_{\mathrm{FE},\elm}$ as in the liner case. Similarly to
Corollaries~\ref{cor_eval_elm_matr} and~\ref{cor_eval_elm_matr_appr}, a
further simplification is to replace $\widehat{\matr{A}}_{\mathrm{MFE},\elm}$
by the element matrices $\widehat{\matr{A}}_{\elm}$ of the given scheme whenever they are available, as per Assumption~\ref{as_polyt_disc_mod}.\er

\br[Distinction of error components] The distinction of the different error components is based on the face fluxes~\eqref{eq_fl}, directly available form any finite volume discretization~\eqref{eq_flux_balance_NL}--\eqref{eq_flux_balance_NL_lin_alg_j}. Beyond simplicity, a remarkable property is that all the error component estimators have the {\em same physical units} (that of energy flux error) and take the {\em same form}. This stands in contrast to the usual practice where the algebraic error is typically treated by the (relative) $L^2$-norm of the algebraic residual (here, from~\eqref{eq_err_NL}, we rather employ the \cor{energy-scaled} $\Hmo$-norm of the residual) and the linearization error estimator is often the $L^\infty$-norm of the difference of two consecutive iterates. This unified treatment of the different error components is a consequence of our holistic approach. In practice, the distinction of error components is to be used in an {\em adaptive inexact Newton algorithm} such as Algorithm~\ref{algo2} of the next section.\er

\br[Error structure] \label{rem_er_str} From~\eqref{eq_triang},
\eqref{eq_dual_res}, and~\eqref{eq_min}, we see that
\be \label{eq_err_NL} \bs
    c_\Kmt^\ft \norm{\tu - \tu_h^{k,i}} & \leq c_\Kmt^\mft
    C_\Kmt \max_{\vf \in \Hoo, \, \norm{\Gr \vf}=1} \{(f,\vf) + (\tu_h^{k,i}, \Gr \vf)\}
    + c_\Kmt^\mft \min_{v \in \Hoo} \norm{\Kmt(|\tu_h^{k,i}|) \tu_h^{k,i} + \Gr
    v}\\
    & \leq 2 c_\Kmt^\mft C_\Kmt \norm{\tu - \tu_h^{k,i}};
\es \ee
in the second inequality, we have used the Cauchy--Schwarz inequality
and~\eqref{eq_inv_p} together with~\eqref{eq_cont} for the second term. This means that the energy error
$c_\Kmt^\ft \norm{\tu - \tu_h^{k,i}}$ is equivalent to the sum of
the {\em dual norm of the residual} and {\em nonconformity} evaluated as the {\em
distance to the $\Hoo$ space}. This is an immediate extension of Theorem~\ref{thm_err_char} from the linear case. Compared to Corollary~\ref{cor_Prag_Syng} for the linear case with an exact algebraic solver giving $\Dv \tu_h = f$, this in particular enables to treat the
nonlinear case with {\em inexact linearization and \cor{algebraic} solvers}, so that $\Dv
\tu_h^{k,i} \neq f$. Our motivation here was to evaluate the error as the $\tLt$-norm, which gives rise to the weights $c_\Kmt$ and $C_\Kmt$ and the factor $2$ in~\eqref{eq_err_NL}. Alternatively, we could have stick to the {\em intrinsic error measure} of the form
\[
    \left\{ \max_{\vf \in \Hoo, \, \norm{\Gr \vf}=1} \{(f,\vf) + (\tu_h^{k,i}, \Gr \vf)\}^2
    + \min_{v \in \Hoo} \norm{\Kmt(|\tu_h^{k,i}|) \tu_h^{k,i} + \Gr v}^2\right\}^{\frac12},
\]
\cf\ \eqref{eq_err_char}\cor{, which would lead to constant-free estimates without the weights $c_\Kmt$ and $C_\Kmt$}. In this way, we define the error measure $\Nn{n,k,i}$ in~\eqref{eq:dual.meas} in the complex multiphase flow below. \er

\subsection{Numerical experiments} \label{sec_num_Darcy_NL}

In this section, we numerically illustrate the efficiency of our theoretical results of Section~\ref{sec_est_NL_comps} on two different examples. Our main goals are to asses the sharpness of the guaranteed bound~\eqref{eq:local.comps:a_NL} and to \cor{examine} the robustness of our estimates with respect to the ratio $C_\Kmt/c_\Kmt$ where we consider $C_\Kmt/c_\Kmt= 10^{i}, \, i \in \{1, ..., 6\}$. As in previous sections, the effectivity ind\cor{ex} is defined by the ratio of the estimator from Theorem~\ref{thm_estim_Darcy_NL} to the exact error, 
\be \label{eq_eff_ind_nl}
    I_{\mathrm{eff}} \eq \frac{\eta^{k,i}}{c_\Kmt^\ft\norm{\tu -  \tu_h^{k,i}}}.
\ee
Additionally, we study in the following examples a stopping criterion for the linearization algorithm based on balancing the error components of Theorem~\ref{thm_estim_Darcy_NL}. The criterion is that the linearization iteration is pursued until step $k$ such that
\be\label{eq_crit_lin}
 \Esti{lin}{k,i}  \leq \Gamma_{\rm lin} \cor{\Esti{sp}{k,i}},
\ee
with $\Gamma_{\rm lin} \in (0,1)$ \cor{(we focus on the linearization only, the algebraic solver is used with a very small tolerance yielding $\Esti{alg}{k,i}, \Esti{rem}{k,i} \approx 0$)}. We compare this criterion with a classic one based on a fixed threshold on the relative linearization residual
\begin{equation} \label{eq:st.crit.lin.clas}
    e_{\rm lin }^{k} \le 10^{-8},
\end{equation}
where the relative linearization residual, for a resolution of a system of nonlinear algebraic equations
$F(X)=0$ by the Newton method, is given by 
\[
 e_{\rm lin }^{k} \eq \frac{\norm{F(X^{k})}}{\norm{F(X^{0})}}.
\]
Note that the criterion~\eqref{eq_crit_lin} can be used in more general adaptive algorithm balancing all the error components \cor{such} as in Algorithm~\ref{algo2}.

\subsubsection{\cor{Regular} solution}\label{sec_smooth_sol}

\begin{figure}
   \centering
    \includegraphics[width=0.3\linewidth]{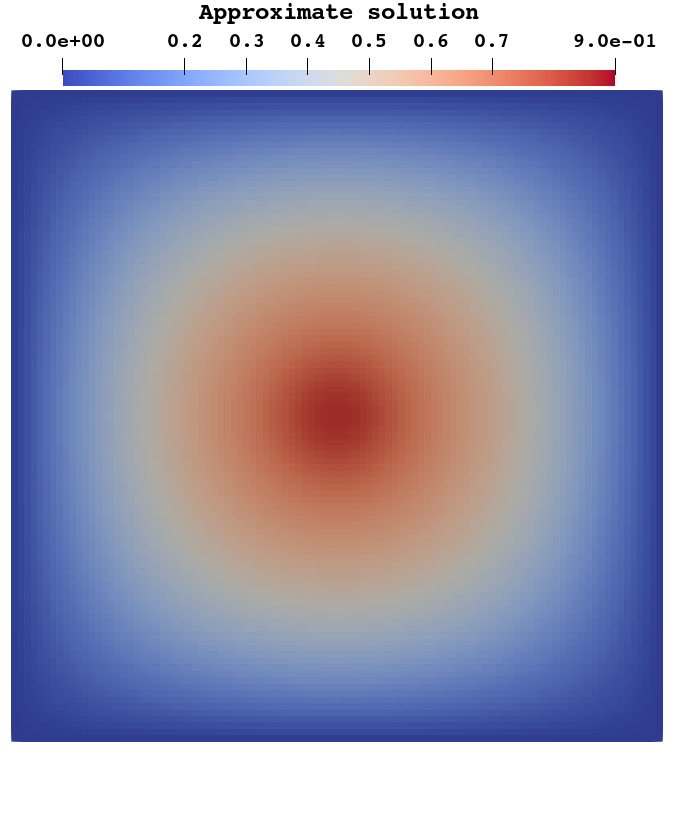}
\hspace{0.4cm}
    \includegraphics[width=0.3\linewidth]{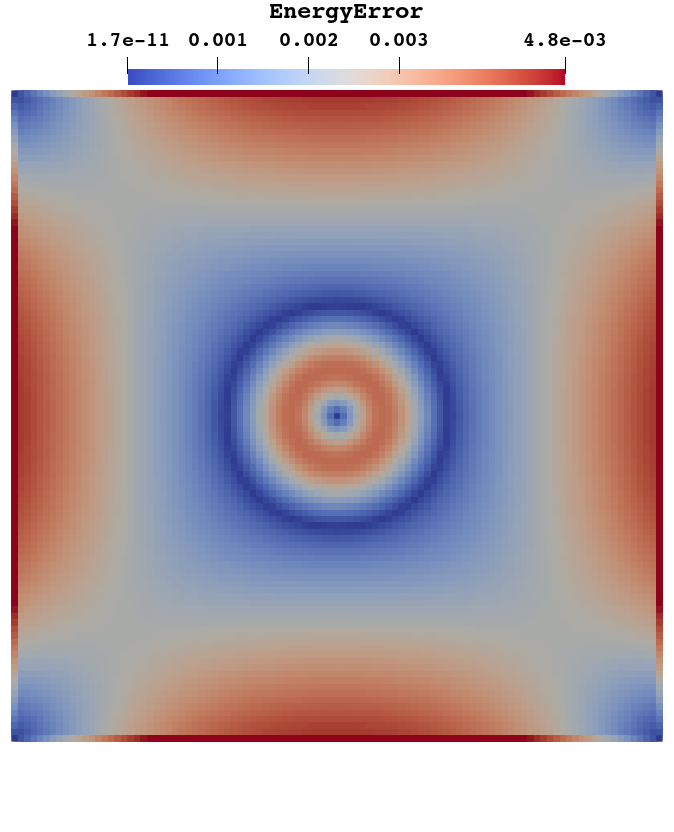}
\hspace{0.4cm}
    \includegraphics[width=0.3\linewidth]{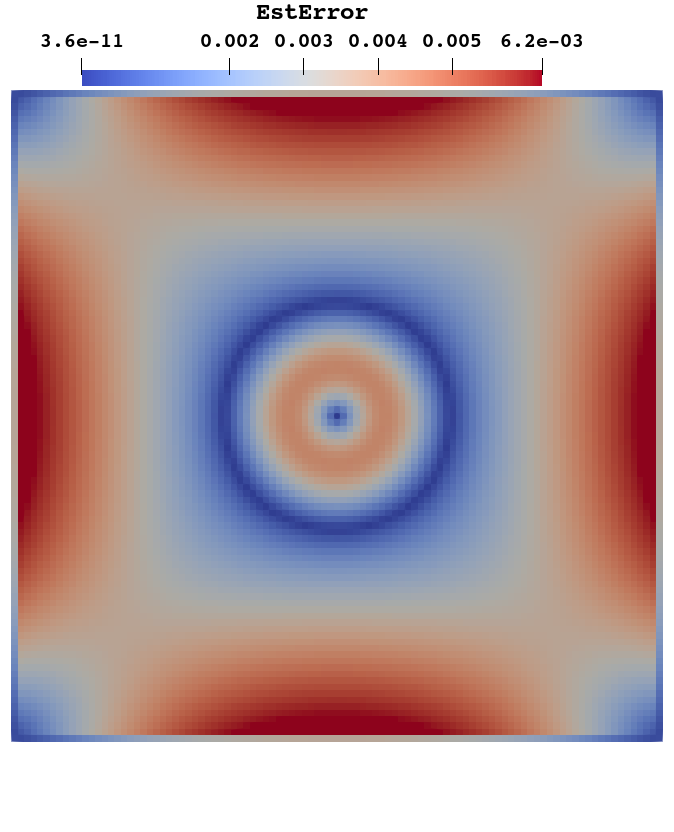}
    \caption{[Section~\ref{sec_smooth_sol}] The approximate solution given by the values $p_\elm^{\cor{k,i}}$, $\elm \in \Th$, ({\em left}) \cor{the energy errors $c_\Kmt^\ft \norm{\tu - \tu_h^{k,i}}_\elm$ ({\em middle}),} and the error estimators \cor{$\est{\mathrm{sp}}{\elm}{k,i} + \est{\mathrm{lin}}{\elm}{k,i} + \est{\mathrm{alg}}{\elm}{k,i} + \est{\mathrm{rem}}{\elm}{k,i}$} from Theorem~\ref{thm_estim_Darcy_NL} ({\em right}).}
  \label{fig:smooth.sol}
\end{figure}

In this test, we consider on the domain $\Omega = (0,1) \times (0,1)$ the exact solution
\be\label{eq_smooth_sol}
p(x,y)=16 x (1-x) y (1-y).
\ee
The source term $f$ is then constructed through~\eqref{eq_Darcy_NL}. We perform the numerical tests with the nonlinearity example given by \eqref{eq_nl_example} on a uniformly refined rectangular mesh of $10^4$ elements. Figure~\ref{fig:smooth.sol} shows the approximate solution given by the values $p^{k,i}_\elm$, $\elm \in \Th$, the elementwise errors $c_\Kmt^\ft \norm{\tu - \tu_h^{k,i}}_\elm$, and the corresponding a posteriori error estimators $\eta^{k,i}_\elm$ from Theorem~\ref{thm_estim_Darcy_NL}. 

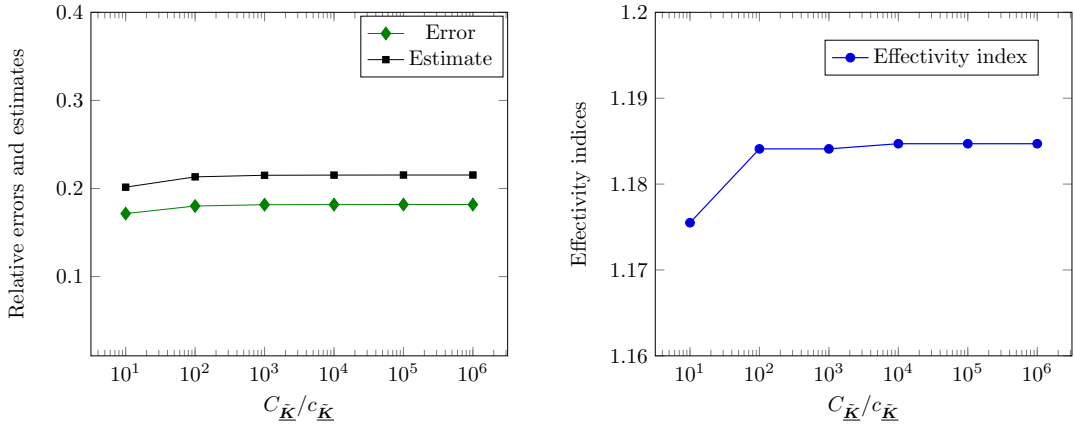
\begin{figure}
\centering
     \begin{tikzpicture}[scale=0.8]
      \begin{semilogxaxis}[
        ymin = 1e-2,
        ymax = 0.4,
         max space between ticks=30,
         yticklabel style={/pgf/number format/fixed},
          xlabel = {$C_\Kmt/c_\Kmt$},
          ylabel = {\cor{Relative e}rrors and estimates},
          legend style ={at = { (0.99,0.99)}}
        ]
        \addplot +[green!50!black,  mark options={solid}, mark=diamond*, mark size=3]
        table[x=fact,y=error]{Figs/smooth_sol_sharpness.tex};
        \addplot +[mark=square*, mark size=1.5, mark options={solid}, black]table[x=fact,y=estimate]{Figs/smooth_sol_sharpness.tex};
        \legend{Error, Estimate};
      \end{semilogxaxis}
    \end{tikzpicture}
\hspace{0.5cm}
     \begin{tikzpicture}[scale=0.8]
      \begin{semilogxaxis}[
        ymin = 1.16,
        ymax = 1.20,
         max space between ticks=30,
          xlabel = {$C_\Kmt/c_\Kmt$},
          ylabel = {Effectivity indices},
          legend style ={at = { (0.92,0.92)}} 
        ]
        \addplot +[line width=0.2mm] table[x=fact,y=index]{Figs/smooth_sol_sharpness.tex};
        \legend{Effectivity index};
      \end{semilogxaxis}
    \end{tikzpicture}
  \caption{[Section~\ref{sec_smooth_sol}], Relative errors $c_\Kmt^\ft\norm{\tu -  \tu_h^{k,i}}{/\norm{ \tu_h^{k,i}}}$ and relative estimates $\eta^{k,i}{/\norm{ \tu_h^{k,i}}}$ from Theorem~\ref{thm_estim_Darcy_NL} ({\em left}), effectivity indices $I_{\mathrm{eff}}$ from~\eqref{eq_eff_ind_nl} ({\em right})}
\label{fig_estim_smooth_sol}
\end{figure}

The robustness of our estimates with respect to the ratio $C_\Kmt/c_\Kmt$ is \cor{examined} in Figure~\ref{fig_estim_smooth_sol}, where we fix the value of $c_\Kmt=1$ and plot, for different values of $C_\Kmt \in \{10, 100, ..., 10^6\}$, the values of the relative energy error $c_\Kmt^\ft\norm{\tu -  \tu_h^{k,i}}{/\norm{ \tu_h^{k,i}}}$ and the corresponding relative total error estimate $\eta^{k,i}{/\norm{ \tu_h^{k,i}}}$ from Theorem~\ref{thm_estim_Darcy_NL}. The effectivity index plotted in the right part of Figure~\ref{fig_estim_smooth_sol} reveal that, in this case and for the choice of nonlinearity example~\eqref{eq_nl_example}, the estimate is robust and not affected by the values of the ratio $C_\Kmt/c_\Kmt$.  

\begin{figure}
   \centering
     \begin{tikzpicture}[scale=0.92]
      \begin{axis}[
          xlabel = {Newton iteration},
          ylabel = {Error component \cor{estimates}},
          ymode=log,
          legend style ={at = { (0.5,0.41)}} 
        ]
        \addplot +[line width=0.01mm] table[x=iter,y=estimate]{Figs/smooth_sol_iterations.tex};
        \addplot +[line width=0.01mm,mark=triangle*] table[x=iter,y=linEst]{Figs/smooth_sol_iterations.tex};
        \addplot +[line width=0.01mm,mark=diamond*] table[x=iter,y=linErr]{Figs/smooth_sol_iterations.tex};
        \node (AD) [ red!60, text=red, draw] at (155, -4.) {\small adaptive stopping criterion};
      \draw (AD) edge [->, shorten >=1pt, thick, red, bend right=1.5]
      (65, -4.);
      \node (UN) [ blue!60, text=blue, draw] at  (70, -19) {\small standard stopping criterion};
      \draw (UN) edge [->, shorten >=1pt, thick, blue, bend right=1.5]
      (205, -18.8);
        \legend{\small total estimator, \small lin. estimator, \small rel. lin. residual};
      \end{axis}
    \end{tikzpicture}
\hspace{0.3cm}
     \begin{tikzpicture}[scale=0.92]
      \begin{semilogxaxis}[
          ymin = 1,
          ymax = 40,
          xlabel = {$C_\Kmt/c_\Kmt$},
          ylabel = {Number of Newton iterations},
          legend style ={at = { (0.97,0.95)}} 
        ]
        \addplot +[line width=0.05mm] table[x=fact,y=classic]{Figs/smooth_iter.tex};
        \addplot +[line width=0.05mm,mark=triangle*] table[x=fact,y=adaptive]{Figs/smooth_iter.tex};
        \legend{Standard resolution, Adaptive resolution};
      \end{semilogxaxis}
    \end{tikzpicture}
\caption{[Section~\ref{sec_smooth_sol}], Standard linearization vs. adaptive linearization}
\label{fig:smooth.iter}
\end{figure}
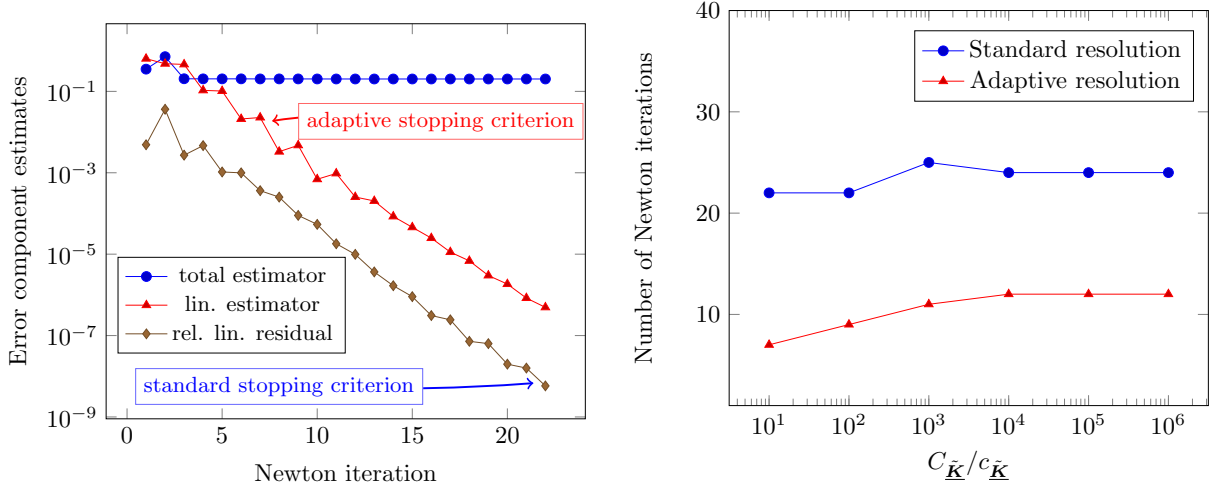

Figure~\ref{fig:smooth.iter} is dedicated to study the stopping criterion~\cor{\eqref{eq_crit_lin}} \cor{and its comparison with the usual stopping criterion~\eqref{eq:st.crit.lin.clas}}. The left part of Figure~\ref{fig:smooth.iter} depicts the evolution of the total estimator $\eta^{k,i}$, the linearization estimator $\eta_{\rm lin}^{k,i}$, and the relative linearization residual
$e_{\rm lin }^{k}$ as a function of the number of Newton iterations for a uniform fixed mesh of $10^4$ elements and with the ratio $C_\Kmt/c_\Kmt = 10$. We observe that the linearization estimator and the relative linearization residual steadily decrease, while the total estimator
stagnate starting from the third iteration. The stopping criterion~\eqref{eq_crit_lin} with $\Gamma_{\rm 1in}=0.1$, indicates that the resolution can be stopped at the sixth iteration, which already yields a sufficiently accurate approximate solution and helps avoid unnecessary additional iterations. The right part of Figure~\ref{fig:smooth.iter} shows a comparison of the required number of Newton iteration in order to satisfy the standard criterion $e_{\rm lin }^{k}\leq 10^{-8}$ and the adaptive one $ \Esti{lin}{k,i}  \leq 0.1 \cor{\Esti{sp}{k,i}}$ for different values of the ratio $C_\Kmt/c_\Kmt$. We see that we can save at least the half of Newton iterations on each resolution without altering the precision.

\subsubsection{Singular solution}\label{sec_nl_lshape}
In this test, we consider the L-shaped domain 
\[
    \Omega=(-1,1)\times (-1,1)\setminus (-1,0]\times (-1,0]
\]
with \cor{the exact solution 
\begin{align*}
p(r,\theta)=r^{\frac{2}{3}}\sin\left(\frac{2}{3}\theta + \frac{3}{2} \pi\right),
\end{align*}
where ($r,\theta$) are the polar coordinates.}
The Dirichlet boundary condition is partly inhomogeneous and given by the value of the exact solution on $\pt \Om$. The following results are carried out with the choice of the nonlinearity example given by~\eqref{eq_nl_example} and on a uniformly refined rectangular mesh of 10800 elements.

\begin{figure}
   \centering
    \includegraphics[width=0.3\linewidth]{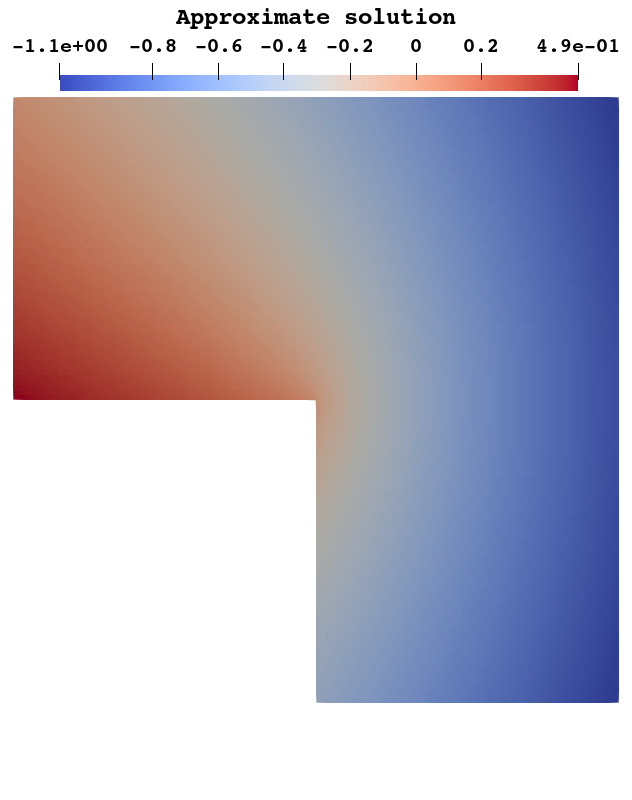}
\hspace{0.4cm}
    \includegraphics[width=0.3\linewidth]{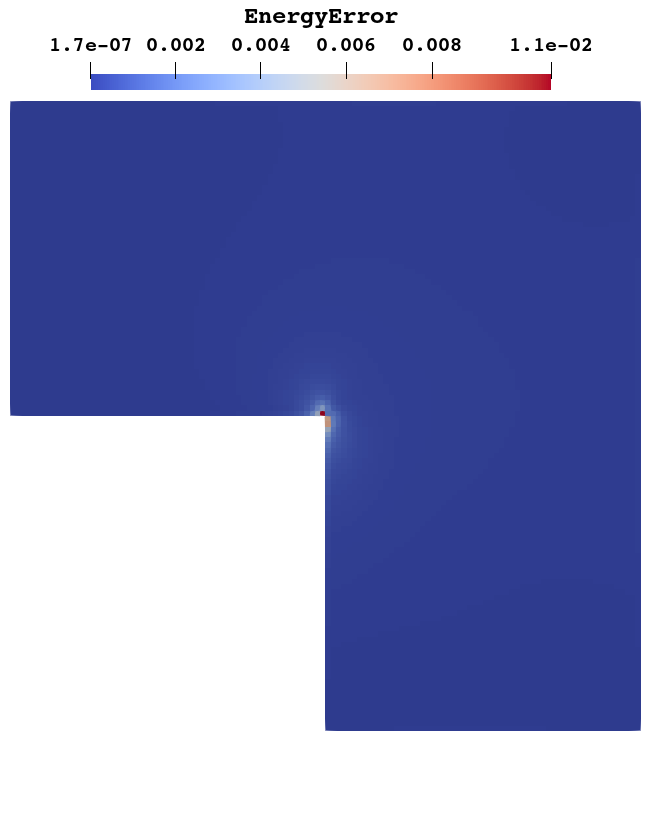}
\hspace{0.4cm}
\includegraphics[width=0.3\linewidth]{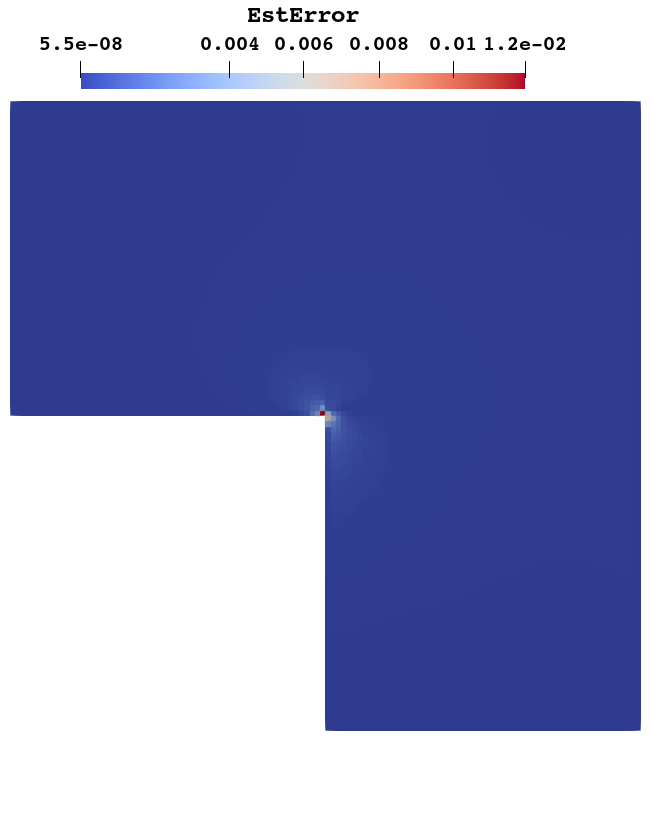}

    \caption{[Section~\ref{sec_nl_lshape}] The approximate solution given by the values $p_\elm^{\cor{k,i}}$, $\elm \in \Th$, ({\em left}), the energy errors $c_\Kmt^\ft\norm{\tu -  \tu_h^{k,i}}_\elm$ ({\em middle}), and the error estimators \cor{$\est{\mathrm{sp}}{\elm}{k,i} + \est{\mathrm{lin}}{\elm}{k,i} + \est{\mathrm{alg}}{\elm}{k,i} + \est{\mathrm{rem}}{\elm}{k,i}$} from Theorem~\ref{thm_estim_Darcy_NL} ({\em right}).}
  \label{fig:lshape.ln.sol}
\end{figure}

In Figure~\ref{fig:lshape.ln.sol}, we show the distribution of the approximate solution and we compare the distributions of the errors $c_\Kmt^\ft \norm{\tu - \tu_h^{k,i}}_\elm$ and the corresponding a posteriori error estimators $\eta^{k,i}_\elm$ from Theorem~\ref{thm_estim_Darcy_NL}. We can observe that the estimators identify perfectly the error in the zone around the singularity at (0,0).  

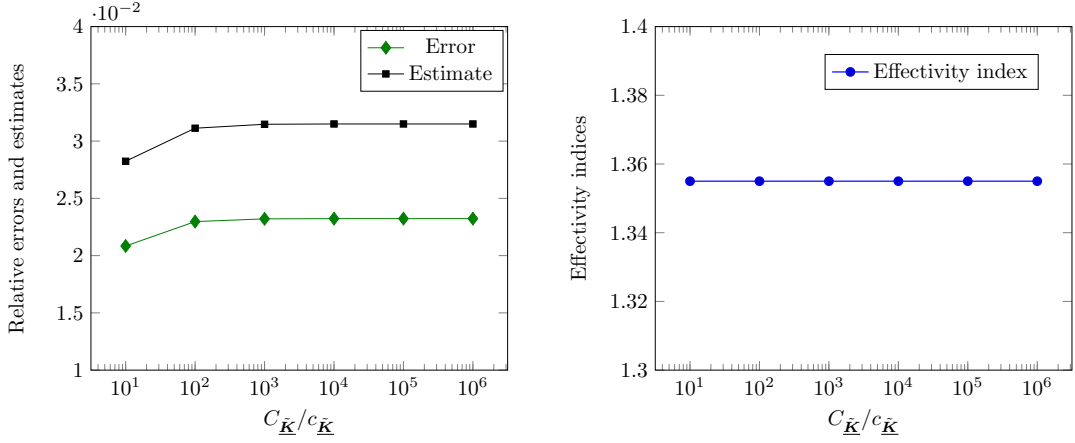
\begin{figure}
\centering
     \begin{tikzpicture}[scale=0.8]
      \begin{semilogxaxis}[
        ymin = 1e-2,
        ymax = 4e-2,        
         max space between ticks=30,
         yticklabel style={/pgf/number format/fixed},
          xlabel = {$C_\Kmt/c_\Kmt$},
          ylabel = {\cor{Relative e}rrors and estimates},
          legend style ={at = { (0.99,0.99)}}
        ]
        \addplot +[green!50!black,  mark options={solid}, mark=diamond*, mark size=3]
        table[x=fact,y=error]{Figs/LShapeNL.tex};
        \addplot +[mark=square*, mark size=1.5, mark options={solid}, black]table[x=fact,y=estimate]{Figs/LShapeNL.tex};
        \legend{Error, Estimate};
      \end{semilogxaxis}
    \end{tikzpicture}
\hspace{0.5cm}
     \begin{tikzpicture}[scale=0.8]
      \begin{semilogxaxis}[
        ymin = 1.3,
        ymax = 1.4,
         max space between ticks=30,
          xlabel = {$C_\Kmt/c_\Kmt$},
          ylabel = {Effectivity indices},
          legend style ={at = { (0.92,0.92)}} 
        ]
        \addplot +[line width=0.2mm] table[x=fact,y=index]{Figs/LShapeNL.tex};
        \legend{Effectivity index};
      \end{semilogxaxis}
    \end{tikzpicture}
  \caption{[Section~\ref{sec_nl_lshape}], Relative errors $c_\Kmt^\ft\norm{\tu -  \tu_h^{k,i}}{/\norm{ \tu_h^{k,i}}}$ and relative estimates $\eta^{k,i}{/\norm{ \tu_h^{k,i}}}$ from Theorem~\ref{thm_estim_Darcy_NL} ({\em left}), effectivity indices $I_{\mathrm{eff}}$ from~\eqref{eq_eff_ind_nl} ({\em right})}
\label{fig_lshape_ln_smooth_sol}
\end{figure}

In Figure~\ref{fig_lshape_ln_smooth_sol}, we study the effect of the ratio $C_\Kmt/c_\Kmt$ on our estimators. The left part of Figure~\ref{fig_lshape_ln_smooth_sol} shows the relative exact errors $c_\Kmt^\ft\norm{\tu -  \tu_h^{k,i}}{/\norm{ \tu_h^{k,i}}}$ together with their relative a posteriori error estimates $\eta^{k,i}{/\norm{ \tu_h^{k,i}}}$ from Theorem~\ref{thm_estim_Darcy_NL}. 
The right part of Figure~\ref{fig_lshape_ln_smooth_sol} displays the corresponding effectivity indices, given by the ratio of the estimate over the error~\eqref{eq_eff_ind_nl}. We see that the effectivity indices are not affected by the values of the
ratio $C_\Kmt/c_\Kmt$, which \cor{numerically gives} the robustness of the estimate. 

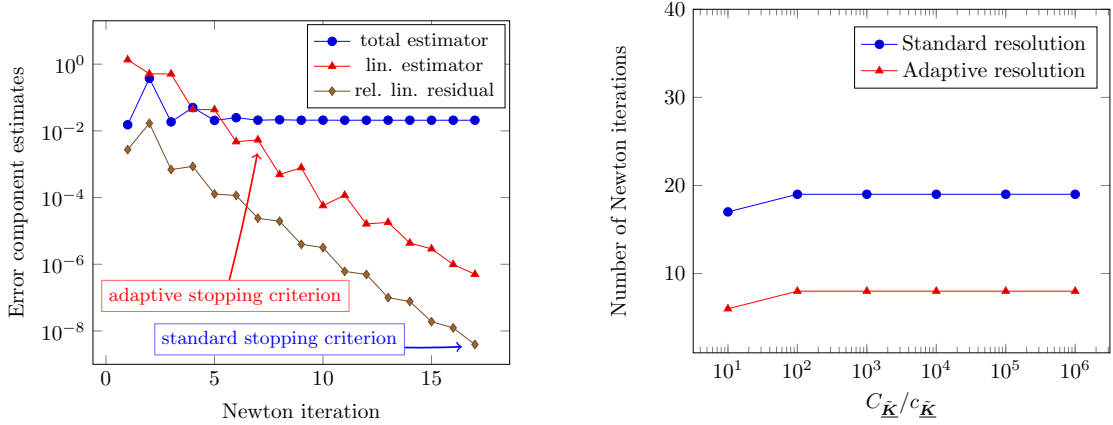
\begin{figure}
   \centering
     \begin{tikzpicture}[scale=0.8]
      \begin{axis}[
        ymin = 1e-9,
        ymax = 20,
          xlabel = {Newton iteration},
          ylabel = {Error component \cor{estimates}},
          ymode=log,
          legend style ={at = { (0.99,0.99)}} 
        ]
        \addplot +[line width=0.01mm] table[x=iter,y=estimate]{Figs/LShape_iter.tex};
        \addplot +[line width=0.01mm,mark=triangle*] table[x=iter,y=linEst]{Figs/LShape_iter.tex};
        \addplot +[line width=0.01mm,mark=diamond*] table[x=iter,y=linErr]{Figs/LShape_iter.tex};
        \node (AD) [ red!60, text=red, draw] at (45, -16) {\small adaptive stopping criterion};
      \draw (AD) edge [->, shorten >=1pt, thick, red, bend right=1.5]
      (60, -6);
      \node (UN) [ blue!60, text=blue, draw] at  (70, -19) {\small standard stopping criterion};
      \draw (UN) edge [->, shorten >=1pt, thick, blue, bend right=1.5]
      (155, -19.5);
        \legend{\small total estimator, \small lin. estimator, \small rel. lin. residual};
      \end{axis}
    \end{tikzpicture}
\hspace{1.cm}
     \begin{tikzpicture}[scale=0.8]
      \begin{semilogxaxis}[
          ymin = 1,
          ymax = 40,
          xlabel = {$C_\Kmt/c_\Kmt$},
          ylabel = {Number of Newton iterations},
          legend style ={at = { (0.97,0.95)}} 
        ]
        \addplot +[line width=0.05mm] table[x=fact,y=classic]{Figs/LShape_nb_iter.tex};
        \addplot +[line width=0.05mm,mark=triangle*] table[x=fact,y=adaptive]{Figs/LShape_nb_iter.tex};
        \legend{Standard resolution, Adaptive resolution};
      \end{semilogxaxis}
    \end{tikzpicture}
\caption{[Section~\ref{sec_nl_lshape}], Standard linearization vs. adaptive linearization}
\label{fig:lshape.ln.iter}
\end{figure}

In order to evaluate the potential of a stopping criterion for the linearization algorithm, we illustrate in the left part of Figure~\ref{fig:lshape.ln.iter} the evolution of the total estimator $\eta^{k,i}$, the linearization estimator $\eta_{\rm lin}^{k,i}$, and the relative linearization residual
$e_{\rm lin }^{k}$ as a function of the number of Newton iterations for a uniform fixed mesh of $10800$ elements and with the ratio $C_\Kmt/c_\Kmt = 10$. These results are given with a standard resolution stopped when $e_{\rm lin }^{k}\leq 10^{-8}$. The linearization estimator and the relative linearization residual decrease as expected while the total estimator start to stagnate from the fifth iteration. Considering the stopping criterion~\eqref{eq_crit_lin} with $\Gamma_{\rm lin}=0.1$ suggests that the resolution process can be halted after the seventh iteration, as it already provides a sufficiently accurate approximate solution while avoiding a significant number of unnecessary linearization iterations. In the right part of Figure~\ref{fig:lshape.ln.iter}, computational savings in
terms of linearization iterations is observed for different values of the ratio $C_\Kmt/c_\Kmt$. These results make it
apparent that using criteria based on balancing the error components gives a powerful tool for adaptivity algorithms. This is shown in the next section by considering Algorithm~\ref{algo2} applied \cor{to a} realistic problem and achieving a significant gain in terms of the computation effort.

\subsection{Bibliographic resources}\label{sec_biblio_Darcy_NL}

One of the first a posteriori error estimates for finite volume methods for nonlinear problems is presented in Bergam~\eal\ \cite{Ber_Mgh_Verf_A_post_FV_NL_03}, though it concerns vertex-centered finite volumes, yielding $\Hoo$-conforming approximations as in the finite element case. Kim~\cite{Kim_a_post_loc_cons_nonlin_07} was one of the first to consider locally conservative methods in a Prager--Synge-type framework \cor{of the form~\eqref{eq_min_NC}}.

The distinction of different error components (and the associated adaptivity) has been initiated in 
Han~\cite{Han_a_post_NL_94}, Becker~\eal\ \cite{Beck_John_Ran_95}, Kaltenbacher~\cite{Kalt_a_post_Newt_98}, Axelsson and Kaporin~\cite{Axel_Kap_err_est_stop_crit_CG_01}, Patera and R{\o}nquist~\cite{Pat_Ronq_gen_out_it_err_01}, Arioli~\cite{Ari_st_crit_FE_04}, Bakushinsky and Smirnova~\cite{Bak_Smi_stop_crit_FP_06}, Chaillou and Suri~\cite{Chai_Sur_comp_a_post_NL_06, Chai_Sur_a_post_lin_err_mon_NL_07}, Picasso~\cite{Pic_stop_crit_09}, Ahusborde~\eal\ \cite{Ahu_Aza_Bel_Ber_mod_red_Darc_13}, Bernardi~\eal\ \cite{Ber_Reb_Ver_Cor_DD_09, Ber_Dak_Mans_Say_a_post_it_NL_15}, Jir{\'a}nek~\eal\ \cite{Jir_Strak_Voh_a_post_it_solv_10}, El Alaoui~\eal\ \cite{El_Al_Ern_Voh_a_post_NL_11}, Nordbotten~\eal\ \cite{Nordb_Keile_Sand_mass_cons_DD_12}, Arioli~\eal\ \cite{Ari_Lies_Mied_Strak_13, Ar_Geor_Log_st_crit_cvg_13}, Rannacher and Vihharev~\cite{Ran_Vir_adapt_NL_bal_13}, Ern and Vohral{\'{\i}}k~\cite{Ern_Voh_adpt_IN_13}, Keilegavlen and Nordbotten~\cite{Keil_Nord_inex_lin_alg_14, Keil_Nord_inex_FV_15}, Rey~\eal\ \cite{Rey_Rey_Gos_a_post_DD_14, Rey_Rey_Gos_a_post_DD_16, Rey_Rey_Gos_a_post_DD_GO_15}, and Becker~\eal\ \cite{Beck_Cap_Luce_st_crit_flux_rec_15}, see also the references therein. 
It has been reconsidered more recently, namely in the framework of conforming finite elements, in Heid and Wihler~\cite{Heid_Wih_it_lin_NL_20}, Heid~\eal\ \cite{Heid_Wihl_Praet_en_contr_cvg_21}, Gantner~\eal\ \cite{Gant_Hab_Praet_Schi_opt_cost_NL_21}, Haberl~\eal\ \cite{Hab_Praet_Schim_Voh_inex_Newt_opt_cost_21}, Kumar~\eal\ \cite{Kum_Kay_Nord_Rep_a_post_Biot_21}, Chaudhry~\eal\ \cite{Chaud_Estep_Taven_a_post_DD_21}, Dolgov and Vejchodsk\'{y}~\cite{Dolg_Vejch_a_post_low_rank_21}, Liu and Keyes~\cite{Liu_Keyes_a_post_nonl_prec_21}, F{\'e}votte~\eal\ \cite{Fev_Rapp_Voh_adapt_reg_24}, Bringmann~\eal\ \cite{Bring_Feisch_Mir_Praet_Strei_lin_cvg_inex_25}\cor{, Harnist~\eal\ \cite{Har_Mitr_Rapp_Voh_en_difs_rob_24}, Mitra and Vohral{\'{\i}}k~\cite{Mitr_Voh_NL_rob_it_norm_23}, Stokke~\eal\ \cite{Stok_Mitr_Stor_Both_Radu_adapt_lin_23}, and Pape\v{z}~\cite{Pap_alg_err_24}.} The classical books of Kelley~\cite{Kelley_it_methods_95}, Saad~\cite{Saad_03}, or Olshanskii and Tyrtyshnikov~\cite{Olsh_Tyrt_it_methods_14} do not discuss \cor{total} a posteriori error estimates during iterative linearization and iterative linear algebraic solvers \cor{iteration} but rather \cor{only} treat the nonlinear algebra level (for iterative linearization) and linear algebra level (for iterative linear algebraic solvers), not addressing the discretization/partial differential equation level.

\section{Unsteady nonlinear coupled degenerate advection--diffusion--reaction problems. Temporal discretization error} \label{sec_MP_MC}

In this section, the methodology derived in Sections~\ref{sec_Pois}--\ref{sec_Darcy_NL} is applied to polytopal discretizations of a model real-life {\em multiphase compositional Darcy flow} in porous media. In the following, we describe briefly the multiphase model,
discuss its discretization by an implicit finite volume scheme, the
linearization by the Newton \cor{or combined Newton and fixed point} method, and the algebraic solution of the arising
linear system. We explain that the required fictitious flux and potential reconstructions
stay unchanged with respect to Section~\ref{sec_Darcy_NL} and finally develop our a posteriori error estimate with inexpensive implementation and evaluation distinguishing the different error components. 

By decomposing the estimators into the space, time, linearization, and algebraic error components, we can design \cor{an algorithm involving adaptivity of both solvers of meshes. We namely} formulate criteria for stopping the iterative algebraic solver and the iterative
linearization solver when the corresponding error components do not affect
significantly the total error. Moreover, the spatial and temporal error
components are balanced respectively by time step and space mesh adaptation, materializing the abstract considerations of Section~\ref{sec_adapt}.
Denote by $\Est{sp}{n,k,i}, \Est{tm}{n,k,i}, \Est{lin}{n,k,i}$, and $\Est{alg}{n,k,i}$
respectively the {total} spatial, temporal, linearization, and algebraic
estimators, on time step $n$, linearization step $k$, and algebraic solver
step $i$. Let $\param{lin},
\param{alg} \in (0,1)$ and $\Param{tm}> \param{tm}>0$ be user-given fixed balancing
parameters and fix the fractions \cor{(percentages)} of cells to refine
$\cor{\delta}_{\mathrm{ref}}$ and to derefine $\cor{\delta}_{\mathrm{deref}}${, $0 <
\cor{\delta}_{\mathrm{deref}} < \cor{\delta}_{\mathrm{ref}} < 1$}.
Following~\cite{Ern_Voh_adpt_IN_13, Voh_Whee_a_post_2P_13,
Canc_Pop_Voh_a_post_2P_14, Di_Pi_Voh_Yous_a_post_comp_14,
Di_Pi_Voh_Yous_a_post_therm_14} and the references therein, we in particular derive in this
section the following {\em fully adaptive algorithm}:

\begin{algo}[Adaptive stopping criteria and adaptive time and space mesh refinement]\label{algo2} $ $
{\em
  \begin{algorithmic}
    \STATE  Set $n \eq 0$.
    \WHILE[time loop] {$t^n\le t_{\rm{F}}$}
     \STATE Set $n \eq n+1$.
     \LOOP[spatial and temporal errors balancing loop]
     \STATE  Set $k \eq 0$.
     \LOOP[linearization loop]
     \STATE Set $k \eq k+1$.
     \STATE  Set up the linear system.
     \STATE  Set $i \eq 0$.
     \LOOP[algebraic solver loop]
     \STATE Perform a step of the iterative algebraic solver {and} set $i \eq i+1$.
     \STATE Evaluate the different a posteriori error estimators.
     \STATE {\bf Terminate if:}
$
          \Est{alg}{n,k,i} \le
          \param{alg}
          \Est{sp}{n,k,i}.
$
     \ENDLOOP \, \{algebraic solver loop\}
     \STATE {\bf Terminate if:}
$          \Est{lin}{n,k,i} \le
          \param{lin} \Est{sp}{n,k,i}.
$
     \ENDLOOP \, \{linearization loop\}
     \STATE {\bf Terminate if}
        \begin{align*}
        {\Est{sp}{\elm,}{n{,k,i}}} & {\geq \cor{\delta}_{\mathrm{ref}}
    \max_{\elmt \in \T_H^n}\big\{ \Est{sp}{\elmt,}{n{,k,i}}\big\} \qquad \forall \elm \in \T_H^n} \\
          \param{tm} \Est{sp}{n,k,i} & \le
          \Est{tm}{n,k,i} \le
          \Param{tm} \Est{sp}{n,k,i};
        \end{align*} \vspace{-0.5cm}
    \STATE{\bf else}
    \STATE Refine the cells $\elm \in \T_H^n$ such that $\Est{sp}{\elm,}{n{,k,i}} \geq \cor{\delta}_{\mathrm{ref}}
    \max_{\elmt \in \T_H^n}\big\{ \Est{sp}{\elmt,}{n{,k,i}}\big\}$.

    \STATE Derefine the cells $\elm \in \T_H^n$ such that $\Est{sp}{\elm,}{n{,k,i}}
 \leq
        \cor{\delta}_{\mathrm{deref}} \max_{\elmt \in \T_H^n}\big\{ \Est{sp}{\elmt,}{n{,k,i}}\big\}$.
    \STATE {Refine the time step if $\Est{tm}{n,k,i} > \Param{tm} \Est{sp}{{\rm
    t}}{n,k,i}$, derefine the time step if $\param{tm}\Est{sp}{n,k,i} > \Est{tm}{n,k,i}$.}
       \ENDLOOP \, \{spatial and temporal errors balancing loop\}
      \STATE Update data.
    \ENDWHILE \, \{time loop\}
  \end{algorithmic}}
\end{algo}

\subsection{Multiphase compositional unsteady Darcy flow}
\label{sec_model}

We consider the multiphase compositional Darcy flow model in porous media~\cite{Acs_Dol_Fark_comp_mod_85, Aziz_Set_petr_res_sim_79, Coats_eq_state_comp_80, Chen_res_sim_07, Young_Steph_comp_res_sim_83}. 
We use the generalization of Coats' formulation~\cite{Coats_impl_composit_89} to an arbitrary number of phases by~\cite{Eym_Guich_Her_Mas_MC_12}.

Let $\Pp=\{p\}$ be the set of {\em phases}, $\C=\{c\}$ the set
of {\em components}, and, for a given phase $p \in \Pp$, let $\C_p\subset\C$
be the set of its components. $S_p$ then denotes the {\em saturation} of the
phase $p$ and $C_{p,c}$ the {\em molar fraction} {of the component $c$ in the
phase $p$}. For a given component $c \in \C$, denote by $\Pp_c$ the set of
the phases which contain $c$.

We denote by $P$ the {\em reference pressure} such that the {\em phase
pressures} $P_p$, $p \in \Pp$, are expressed as
\begin{equation}\label{pressure}
  P_p \eq P +  P_{\mathrm{c}_p}(S_p),
\end{equation}
where $P_{\mathrm{c}_p} :\RR \ra \RR$ is a given generalized capillary pressure function.
We collect the unknowns of the model in the vector 
\be \label{eq_unkns}
    \sol {\eq} (P,(S_p)_{p\in \Pp},(C_{p,c})_{p \in \Pp, c \in \C_p}).
\ee
For a phase $p \in \Pp$, let $k_{{\rm r},p}:\RR \ra \RR$ be the relative permeability function of the unknown saturation $S_p$. We denote by $\phi$ the porosity of the medium, by $\perm$ the permeability tensor, by $\mu_p$ the dynamic viscosity, by $\cor{\varsigma}_p$ the molar density, and by $\rho_p$ the mass density. 
For simplicity, we suppose that all these porous media parameters are given as constants on each mesh element and on each time step\cor{; $\phi$ is actually supposed to be constant in time}. \cor{We also define} the mobility by $\nu_p\eq \cor{\varsigma}_p\frac{k_{{\rm r},p}}{\mu_p}$. \cor{Finally}, let $q_c$ be the {\em source or sink} term, again piecewise constant on the space-time mesh. 

Let $t_{\mathrm F}>0$ be the final simulation time.
The system of governing equations is posed in $\Omega \times (0,t_{\mathrm
F})$ by
\begin{equation}\label{composit}
  \boxed{\pt_t \Ll_c + \DIV\VEC{\theta}_c = q_c,
  \qquad \forall c \in \C,}
\end{equation}
where $\Ll_c$ is the {\em amount} (in moles) of component $c$ per unit volume,
\begin{equation}\label{eq:lc}
  \Ll_c \eq \phi \sum_{p \in \Pp_c} \cor{\varsigma}_p S_p C_{p,c}.
\end{equation}
In~\eqref{composit}, for each component $c\in\C$, the {\em component flux}
$\VEC{\theta}_c$ has the following expression:
\begin{equation}\label{eq:comp.flux}
  \VEC{\theta}_c \eq \sum_{p \in \Pp_c} \VEC{\theta}_{p,c},\qquad
  \VEC{\theta}_{p,c}
  \eq \VEC{\theta}_{p,c}(\sol)
  = \nu_pC_{p,c}\tv_p({P_p}),
\end{equation}
where for all phases $p \in \Pp$, $\tv_p({P_p})$ represents the average {\em phase
velocity} given by Darcy's law,
\begin{equation}\label{eq:darcy.vel}
  \tv_p({P_p})
  \eq -\perm\left(\GRAD P_p + \rho_p g\GRAD z\right).
\end{equation}
Here\cor{, $z$ is the vertical space coordinate and} $g$ is the gravitation acceleration constant. Note at this
occasion that only the first part of the Darcy velocity,
\begin{equation}\label{eq:darcy.vel_pot}
  \tu_p({P_p})
  \eq -\perm \GRAD P_p,
\end{equation}
is a gradient of the phase pressure. We will use this observation in our
error measure below. We assume that no-flow boundary conditions are
prescribed for all the component fluxes,
\begin{equation}\label{eq:bc}
  \cor{\VEC{\theta}_c} \SCAL\normall[\Omega] = 0 \qquad
  \text{on $\partial \Omega \times (0,t_{\mathrm F})$} \qquad \forall c\in\C,
\end{equation}
where $\partial\Omega$ denotes the boundary of $\Omega$ and $\normall[\Omega]$
its unit outward normal.
At $t=0$, we prescribe the {\em initial amount} of each
component,
\begin{equation} \label{eq:init}
\Ll_c(\cdot, 0)=\Ll_c^0 \qquad \forall c \in \C.
\end{equation}
The previous PDEs is supplemented by a system of algebraic equations imposing
the volume conservation
\begin{equation}\label{compositA}
\sum_{p \in \Pp } S_p = 1,
\end{equation}
the conservation of the quantity of matter
\begin{equation}\label{compositB}
  \sum_{c \in \C_p} C_{p,c} = 1 \qquad \forall p \in \Pp,
\end{equation}
and local
thermodynamic equilibrium expressed by
\begin{equation}\label{eq:fug}
\sum_{c\in\C}(|\Pp_c| - 1) = \sum_{p\in\Pp}|\C_p| - |\C|
\end{equation}
equalities of fugacities, where we refer to~\cite{Di_Pi_Voh_Yous_a_post_comp_14} for details. 

\subsection{Space-time mesh}

We suppose that the time simulation interval $[0, t_{\mathrm F}]$ is partitioned into time intervals $I^n \eq [t^{n-1}, t^n]$, $1 \leq n \leq N$, given by the discrete times $0 = t^0 < \ldots < t^n < \ldots < t^N = t_{\mathrm F}$. The discrete time steps are then $\tau^n \eq t^n- t^{n-1}$. Then, we let $(\T_H^n)_{0 \leq n \leq N}$ be a family of polytopal meshes of the space domain $\Om$ defined in the sense of Section~\ref{sec_meshes}; we suppose that
$\T_H^n$ has been obtained from $\T_H^{n-1}$ by refinement of some elements and
by coarsening of some other. For all $0\leq n \leq N$, we define $\cor{\mathcal{F}}^n_H$,
$\Fh^n$, $\cor{\mathcal{F}}^n_{H,h}$, $\FKhextn$, and $\FKhintn$ as in Section~\ref{sec_meshes}.

\subsection{Weak solution and its properties} \label{sec_prop_Darcy_MP}

At this stage, we need to characterize a weak solution for the multiphase compositional
model~\eqref{pressure}--\eqref{eq:fug}.
We define
\begin{subequations}
  \begin{align}
  \label{eq:Y.X}
  X & \eq L^2((0,t_{\mathrm F}); H^1(\Omega)),\\
  Y & \eq H^1((0,t_{\mathrm F}); L^2(\Omega)).
  \end{align}
\end{subequations}
We equip the space $X$ with the norm
\begin{equation} \label{eq:X.norm}
    \norm{\varphi}_X := \Biggl\{ \sum_{n=1}^N \norm{\varphi}_{X^n}^2 
    \Biggr\}^{\frac12}, \qquad     
    \norm{\varphi}_{X^n}^2 \eq \int_{I^n} \sum_{\elm \in \T_H^n} \big\{ h_\elm^{-2} \norm{\varphi}_{\elm}^2 +
    \norm{\Km^{\frac{1}{2}}\GRAD \varphi}_{\elm}^2 \big\} \dt
    \qquad \varphi \in X,
\end{equation}
where, we recall, $h_\elm$ is the diameter of the cell $\elm \in \T_H^n$. Below, it will be convenient to use $X^n \eq L^2(I^n;\Ho)$. 
We suppose that there exists a unique weak solution of~\eqref{pressure}--\eqref{eq:fug} in the following sense:

\begin{assumption}[Weak solution] \label{ass:sol.reg}
There exists a weak solution $\sol$ of~\eqref{pressure}--\eqref{eq:fug} which
can be characterized as follows:
\begin{subequations}\label{eq:weak}
\begin{align}
    & l_c \in Y \qquad \forall c \in \C, \label{eq:weak.lc} \\
    & P_p \in X \qquad \forall p \in \Pp, \label{eq:weak.Pp} \\
    & \VEC{\theta}_c \in [L^2((0,t_{\mathrm F}); L^2(\Omega))]^d \qquad \forall c \in \C, \label{eq:weak.Phic} \\
    & \int_0^{t_{\mathrm F}}  \left\{ ( \partial_t l_c, \varphi )(t)
      - (\VEC{\theta}_c, \GRAD\varphi)(t)
      \right\} \d t = \int_0^{t_{\mathrm F}} (q_c,\varphi)(t) \d t \qquad \forall \varphi \in X, \, \forall c \in \C, \label{eq:weak:eqt} \\
    & \text{the initial condition~\eqref{eq:init} holds},\\
    & \text{and the algebraic closure equations~\eqref{compositA}--\eqref{eq:fug} hold},
\end{align}
\end{subequations}
where $P_p$, $l_c$, and $\VEC{\theta}_c$ are defined, respectively,
by~\eqref{pressure}, \eqref{eq:lc}, and~\eqref{eq:comp.flux}.
\end{assumption}

We have formulated the existence and uniqueness of a weak solution as an assumption, since, in contrast to the model problems of Sections~\ref{sec_Pois}--\ref{sec_Darcy_NL}, the problem~\eqref{pressure}--\eqref{eq:fug} is too abstract and complicated to prove it.
In simplified settings, with possibly only two phases present and each
phase composed of a single component, such results can be found
in~\cite{Kro_Luck_2P_84, Chav_Jaff_res_sim_86, Ant_Kaz_Mon_BVP_90,
Chen_deg_2P_I_01, Canc_Gal_Por_TP_cap_het_09, Khal_Saad_2P_compr_mis_11, Am_Jur_Zgal_2P_ex_11} and the references therein.

\begin{remark}[Component fluxes] \label{rem_fluxes} Since we assume $q_c$ piecewise constant and thus  $q_c\in L^2((0,t_{\mathrm F});$ $L^2(\Omega))$, it follows from~\eqref{eq:weak} that actually
\begin{subequations} \label{eq_reg_fluxes} \begin{alignat}{2}
    \VEC{\theta}_c & \in L^2((0,t_{\mathrm F}); {\mathbf H}(\mathrm{div},\Omega)) \qquad & & \forall c \in
        \C, \label{eq_fluxes_Hdv} \\
    \DIV \VEC{\theta}_c & = q_c - \partial_t l_c & & \forall c \in
        \C, \label{eq_fluxes_div} \\
    \VEC{\theta}_c \SCAL\tn_\Om & = 0 \qquad
        \text{on $\partial \Omega \times (0,t_{\mathrm F})$} & & \forall
        c\in\C \label{eq_fluxes_BC},
\end{alignat} \end{subequations}
so that the component fluxes $\VEC{\theta}_c$ have the normal trace
continuous in a proper sense, the governing equation~\eqref{composit}
is satisfied with a weak divergence, and the boundary
conditions~\eqref{eq:bc} hold in the normal trace sense. This extends Proposition~\ref{pr_prop_WS}, \eqref{eq_tu_pol}, and~\eqref{eq_tu_Darcy_NL} to the multiphase compositional model~\eqref{pressure}--\eqref{eq:fug}.
\end{remark}

\subsection{Total error and its spatial discretization, temporal discretization, linearization, and algebraic components}\label{sec:comp.measure}

Suppose that Assumption~\ref{ass:sol.reg} holds, so that we dispose of a weak solution. Then, we describe the distance of the current approximate solution to the exact solution by the following intrinsic error measure:
\bse\label{eq:int_err}\begin{equation}\label{eq:dual.meas}
  \Nn{n,k,i} \eq \Bigg\{
  \sum_{c \in \C}\big(\Nn{n,k,i}_c\big)^2
  \Bigg\}^{\frac12} + \Bigg\{
  \sum_{p \in \Pp}\big(\Nn{n,k,i}_p\big)^2\Bigg\}^{\frac12},
\end{equation}
where
\begin{equation}
  \label{eq:Nc}
  \Nn{n,k,i}_c \eq
  \sup_{\varphi \in X^n, \norm{\varphi}_{X^n}=1} \int_{I^n} \big\{
  (q_c,\varphi)
  -(\pt_t \Ll^{n,k,i}_{c,h\tau}, \varphi)
  + \big(\Phicht[n,k,i], \GRAD \varphi \big)
  \big\} \dt
\end{equation}
is the {\em dual norm of the residual} of the weak formulation for each component flux $\Phicht[n,k,i]$ on the time interval $I^n$ and where 
\begin{equation}
  \label{eq:Np}
  \Nn{n,k,i}_p \eq
  \inf_{\cor{\varsigma}_p\in X^n} \Bigg \{ \sum_{c \in\C_p} \int_{I^n}
  \Bigg\{{\sum_{\elm \in \T_H^n}} \big(\nu_{p,\elm}^{n,k,i}  C^{n,k,i}_{p,c,\elm} \big)^2
  \norm{\tu_{p,h\tau}^{n,k,i} + \Km \Gr \cor{\varsigma}_p}_{\Km^{-\frac{1}{2}},\elm}^2 \Bigg\}\dt
    \Bigg\}^{\frac12}
\end{equation}
\ese
evaluates the {\em nonconformity} (the distance of $\tu_{p,h\tau}^{n,k,i}$ to $\Km \Gr \Hoo$) for each Darcy phase flux $\tu_{p,h\tau}^{n,k,i}$. 
Indeed, from~\cor{\eqref{eq_fluxes_div}}, the first term should be zero and from~\eqref{eq:darcy.vel_pot}, $\tu_p({P_p})$ should be minus a
gradient of a scalar-valued function, which is not necessarily the case for
the discrete $\tu_{p,h\tau}^{n,k,i}$ given at time $t^n$.

\br[Error measure $\Nn{n,k,i}$] The error measure $\Nn{n,k,i}$ is a
natural extension of the energy error from the simpler cases. Indeed, component by component and phase by phase, it extends Theorem~\ref{thm_err_char} from the Poisson case as well as Remark~\ref{rem_er_str} from the steady nonlinear case. The sum $\sum_{n=1}^N (\Nn{n,k,i})^2$ (if the initial condition is satisfied exactly) then corresponds to the square of the energy error for the heat equation, see~\cite[Theorem~2.1 and equation~(2.7)]{Ern_Sme_Voh_heat_HO_Y_17}.
\er

In~\eqref{eq:int_err}, we have denoted the approximate component fluxes by $\Phicht[n,k,i]$ and the approximate Darcy phase fluxes by $\tu_{p,h\tau}^{n,k,i}$, anticipating that our numerical approximation will here consist in {\em temporal discretization} (with characteristic time step $\tau$ and discrete times $t^n$), {\em spatial discretization} (with characteristic mesh size $h$), {\em iterative linearization} (with iteration index $k$), and {\em iterative algebraic resolution} (with iteration index $i$). Our goal will be to {\em distinguish the error components} as in~\eqref{eq_err_comps} in Section~\ref{sec_err_comps}, to design an a posteriori error estimate of the form~\eqref{eq_est_rel_comps_unst}, and to use them in the fully adaptive Algorithm~\ref{algo2}.

\subsection{Backward Euler time \cor{discretization} and locally conservative space discretization on polytopal meshes}
\label{sec:disc}

To discretize problem~\eqref{pressure}--\eqref{eq:fug}, we consider a fully implicit backward Euler time stepping and an abstract locally conservative space discretization on a polytopal mesh of the form of Assumption~\ref{as_polyt_disc}, using phase-upwinding.
For all $1\le n\le N$, we let $\solcM \eq (\solM)_{\elm\in\T_H^n}$, with, following~\eqref{eq_unkns},
\be \label{eq_unkns_disc}
    \solM \eq ( P^n_{\elm}, (S^n_{p,\elm})_{ p\in \Pp}, (C^n_{p,c,\elm})_{ p\in
\Pp,c\in\C_p})
\ee
be the algebraic vector of discrete unknowns \cor{on an element $\elm$}.
System~\eqref{pressure}--\eqref{eq:fug} is then discretized as follows: for all time
steps $1\le n\le N$, all polytopal cells $\elm \in \T_H^n$, and each
component $c\in\C$, we require
 \begin{equation}
  \label{eq:discret}
  \boxed{\frac{|\elm|}{\tau^n}\big(\Ll_{c,\elm} (\solM) - \Ll_{c,\elm} (\solMp)\big)
  + \sum_{\sd \in {\FK \cap \F_H^{\rm int}}} \FcMs(\solcM) = |\elm| q_{c,\elm}^n.}
\end{equation}
This equation expresses the mass balance for the element $\elm$. Here
$q_{c,\elm}^n$ is the value of the source or sink $q_{c}$ on the element $\elm$ and time interval $I^n$ \cor{-- recall that we consider $q$ piecewise constant}. Recall also that $|\elm|$ stands for the \cor{volume} of the element $\elm$. 

In a discrete version of~\eqref{eq:lc}, \cor{we set}
\be \label{eq:lc:disc}
    \Ll_{c,\elm} (\solM) \eq \phi_\elm \sum_{p \in \Pp_c} \cor{\varsigma}^n_{p,\elm} S^n_{p,\elm} C^n_{p,c,\elm}.
\ee
Similarly, for each component $c\in\C$, following~\eqref{eq:comp.flux},
the total flux across the face $\sd$ results from the sum of the
corresponding fluxes for each phase $p\in\Pp_c$, \ie, for all elements $\elm\in\T_H^n$
and all faces $\sd \in {\FK \cap \F_H^{\rm int}}$ with $\sd \cor{\subset} \partial
\elm\cap\partial \elmt$,
\begin{equation}
  \label{eq:FcMs}
 \FcMs(\solcM) \eq \sum_{p\in\Pp_c}   \nu_{p,\elm^\uparrow_p}^{n}
   C^n_{p,c,\elm^{\uparrow}_p}  \FpMs(\solcM),\qquad
  \elm^{\uparrow}_p \eq \begin{cases}
    \elm & \text{if $ \FpMs(\solcM)\ge 0$}, \\
    \elmt & \text{otherwise},
  \end{cases}
\end{equation}
with $C_{p,c,\elm^\uparrow_p}^n$ and $\nu_{p,\elm^\uparrow_p}^{n}$ denoting,
respectively, the upstream molar fraction and the upstream mobility.
Let us also introduce
\begin{equation}
  \label{eq:FcMs:NL}
 \FcMsNL(\solcM) \eq \sum_{p\in\Pp_c}   \nu_{p,\elm}^{n}
   C^n_{p,c,\elm}  \FpMs(\solcM),
\end{equation}
\cor{where no upwinding is used.}
In~\eqref{eq:FcMs}--\eqref{eq:FcMs:NL}, we have employed a finite volume approximation of the normal component of the average phase velocity on face
$\sd$ given by
\begin{equation}
  \label{eq:FpMs}
  \FpMs(\solcM) \eq \FpMs(\{\solMK\}_{\elm' \in \Set_\sd})\eq
  \sum_{\elm' \in \Set_\sd} \tau_{\elm'}^\sd ( P^n_{p,\elm'} + \rho_{p,\sd}^{n}g z_{\elm'}),
\end{equation}
where, following~\eqref{pressure},
\begin{equation}\label{pressure_disc}
    P^n_{p,\elm} \eq P_\elm^n +  P_{\mathrm{c}_p}(S_{p,\elm}^n),
\end{equation}
and, for all $\sd \in {\FK \cap \F_H^{\rm int}}$, $\Set_\sd$ is the
flux stencil collecting the elements in $\T_H^{n}$ with nonzero flux contribution ($\Set_\sd = \{\elm,\elmt\}$ for $\sd$ such that $\sd = \sd_{\elm,\elmt} = \pt \elm \cap \pt \elmt \in \F_H^{\rm int}$) in the two-point scheme~\eqref{eq_dif_fl_int}. Moreover, for all $\elmt \in \Set_\sd$, $\tau_{\elmt}^\sd \in \RR$ is the
transmissibility coefficient of the face $\sd$ (if $\Km = \Idd$, $\tau_{\elm}^\sd = \frac{|\sd_{\elm,\elmt}|}{d_{\elm,\elmt}}$ and $\tau_{\elmt}^\sd = - \frac{|\sd_{\elm,\elmt}|}{d_{\elm,\elmt}}$ for the two-point scheme\cor{~\eqref{eq_dif_fl_int}}). 
Please remark that the boundary fluxes are set to zero to account for the homogeneous
\cor{no-flow} boundary conditions~\eqref{eq:bc}.

Finally, the initial condition comes from~\eqref{eq:init}, and the algebraic constraints~\eqref{compositA}--\eqref{eq:fug} are trivially discretized separately on all $\elm\in\T_H^n$, namely 
\begin{align}
\sum_{p \in \Pp } S^n_{p,\elm} & = 1,\label{compositA_disc}\\
\sum_{c \in \C_p} C^n_{p,c,\elm} & = 1 \qquad \forall p \in \Pp.\label{compositB_disc}
\end{align}

\subsection{Iterative linearization}
\label{sec:lin}

At this stage, we need to solve, at each time step, the system of nonlinear
algebraic equations resulting from the discretization~\eqref{eq:discret}--\eqref{compositB_disc}.
To this purpose, for all times $1 \leq n \leq N$, we apply an iterative linearization, in extension of~\eqref{eq_flux_balance_NL_lin}. We consider the {\em Newton
linearization}, as in~\eqref{eq_Newton_lin}. This generates, for an initial guess $\solcM[n,0]$ (typically given by the previous time step/initial condition), a sequence
$(\solcM[n,k])_{k\ge 1}$ with $\solcM[n,k]$ solution to the following {\em
system of linear algebraic equations}: for all components $c \in \C$ and all
mesh elements $\elm \in \T_H^n$,
\begin{equation}\label{eq:linear-sys}
  \sum_{\elm' \in \T_H^n} \frac{\partial \Rl_{c,\elm} }{\partial \sol^n_{\elm'}}\big(\solcM[n,k-1]\big)
  \cdot \big(\sol^{n,k}_{\elm'}-\sol^{n,k-1}_{\elm'}\big)
  + \Rl_{c,\elm}\big(\solcM[n,k-1]\big) = 0,
\end{equation}
with, for all $c \in \C$ and all $\elm\in\T_H^n$,
\[
   \Rl_{c,\elm}\big(\solcM\big) \eq
  \frac{|\elm|}{\tau^n}{\big(\Ll_{c,\elm} (\solM) - \Ll_{c,\elm} (\solMp)\big)}
  + \sum_{\sd \in {\FK \cap \F_H^{\rm int}}} \FcMs\big(\solcM\big)
  - |\elm|q_{c,\elm}^n.
\]
\cor{Developing the derivatives in~\eqref{eq:linear-sys}, we can write
\begin{equation}\label{eq:RcM.lin}
  \boxed{\frac{|\elm|}{\tau^n}\left( \Ll_{c,\elm}^{n,k} - \Ll_{c,\elm} (\solMp)
  \right)
  + \sum_{\sd \in {\FK \cap \F_H^{\rm int}}}
  \FcMs^{n,k}
  - |\elm| q_{c,\elm}^n = 0,}
 \end{equation}
where the linearization of the amount of component $\Ll_{c,\elm}
(\sol^{n,k}_{\elm}) $ of~\eqref{eq:lc:disc} is given by
\[
    \Ll_{c,\elm}^{n,k} \eq \Ll_{c,\elm} (\solMpp) + \frac{\partial \Ll_{c,\elm} }{\partial \sol^n_{\elm}}\big(\solMpp)
  \cdot \big(\sol^{n,k}_{\elm}-\sol^{n,k-1}_{\elm}\big),
\]
and where $\FcMs^{n,k}$ are appropriate linearizations of the component fluxes $\FcMs(\solcM[n,k])$ of~\eqref{eq:FcMs}.

\br[Blending Newton and fixed-point linearizations] \label{rem_blending_N_FP} In numerical experiments below in Section~\ref{sec_num_exp_MP_MC}, we will actually define $\FcMs^{n,k}$ by evaluation of the mobility $\nu_{p,\elm^\uparrow_p}^{n}$ and the molar fraction $C^n_{p,c,\elm^{\uparrow}_p}$ at the previous linearization iterate $k-1$ as
\begin{eqnarray*}
  \FcMs^{n,k} \eq\sum_{p\in\Pp_c} \biggl\{\nu_{p,\elm^\uparrow_p}^{n{,k-1}} C^{n{,k-1}}_{p,c,\elm^{\uparrow}_p} \biggl[ \FpMs({\solcM[n,k-1]})
  + \sum_{\elm' \in \T_H^n} \frac{\partial \FpMs}{\partial
    \sol^n_{\elm'}}\big(\solcM[n,k-1]\big)\cdot\big(
  \sol^{n,k}_{\elm'}-\sol^{n,k-1}_{\elm'} \big) \biggr] \biggr\},
\end{eqnarray*}
blending thus the Newton and the fixed-point linearizations.
\er}

\subsection{Iterative algebraic system solution}
\label{sec:alg}

Now, for a given time step $1\le n\le N$ and a given \cor{linearization} iteration $k\ge 1$, \cor{\eqref{eq:RcM.lin}} is a system of linear algebraic equations. To find its approximate solution, we use an {\em
iterative algebraic solver}. Let $\solcM[n,k,0]$ be an initial guess. Typically, we take $\solcM[n,k,0] = \solcM[n,k-1]$. This gives a sequence
$(\solcM[n,k,i])_{i\ge 1}$ solving~\eqref{eq:RcM.lin} up to the residuals \cor{$\Rl^{n,k,i}_{c,\elm}$} for all $c \in \C$ and all $\elm \in \T_H^n$, \ie, 
\begin{equation}\label{eq:RcM.lin.alg}
  \boxed{\frac{|\elm|}{\tau^n}\left( \Ll_{c,\elm}^{n,k,i} - \Ll_{c,\elm} (\solMp)
  \right)
  + \sum_{\sd \in {\FK \cap \F_H^{\rm int}}}
  \FcMs^{n,k,i}
  - |\elm| q_{c,\elm}^n \qe \Rl^{n,k,i}_{c,\elm},}
 \end{equation}
\cf~\eqref{eq_flux_balance_NL_lin_alg_def}--\eqref{eq_flux_balance_NL_lin_alg}. Here 
\be \label{eq:l.lin.alg}
  \Ll_{c,\elm}^{n,k,i} \eq \Ll_{c,\elm} (\solMpp) + \frac{\partial \Ll_{c,\elm} }{\partial \sol^n_{\elm}}\big(\solMpp)
  \cdot \big(\sol^{n,k,i}_{\elm}-\sol^{n,k-1}_{\elm}\big),
\ee
and, in the case of Remark~\ref{rem_blending_N_FP},
\begin{eqnarray}\label{eq:FcMs.lin}
  \FcMs^{n,k,i} \eq\sum_{p\in\Pp_c} \biggl\{\nu_{p,\elm^\uparrow_p}^{n{,k-1}} C^{n{,k-1}}_{p,c,\elm^{\uparrow}_p} \biggl[ \FpMs({\solcM[n,k-1]})
  + \sum_{\elm' \in \T_H^n} \frac{\partial \FpMs}{\partial
    \sol^n_{\elm'}}\big(\solcM[n,k-1]\big)\cdot\big(
  \sol^{n,k,i}_{\elm'}-\sol^{n,k-1}_{\elm'} \big) \biggr] \biggr\}.
\end{eqnarray}
We will also consider $j \geq 1$ additional algebraic solver steps
in~\eqref{eq:RcM.lin.alg}, similarly as we have done it
in~\eqref{eq_flux_balance_NL_lin_alg_j}.

\subsection{Face normal fluxes and potential point values} \label{sec_fl_faces_pot_MP_MC}

We proceed here as in Section~\ref{sec_fl_faces_pot_NL}. In particular, in extension of~\eqref{eq_fl}, we rely on the finite volume {\em face normal fluxes}. These employ the approximate phase normal velocities~\eqref{eq:FpMs}, \cor{the total normal flux nonlinear expressions~\eqref{eq:FcMs:NL} and their upwindings~\eqref{eq:FcMs}, as well as the linearizations~\eqref{eq:FcMs.lin} of the latter.} 
For all phases $p\in\Pp$, components $c\in\C$, time levels $n \geq 1$, mesh elements $\elm \in \T_H^n$ and their faces not lying on the boundary of $\Om$\cor{, $\sd \in {\FK \cap \F_H^{\rm int}}$,} linearization steps $k \geq 1$, and algebraic iterations $i \geq 1$, at times $t \in I^n$ where necessary, we define
  \begin{subequations}
    \label{eq:local}
    \begin{align}
(\alg{U}_{\elm,p}^{t,n,k,i})_\sd & \eq {\frac{t - t^{n-1}}{\tau^n}}
 \sum_{\elm' \in \Set_\sd} \tau_{\elm'}^\sd  P^{n,k,i}_{p,\elm'} {+ \frac{t^n - t}{\tau^n}
 \sum_{\elm' \in \Set_\sd} \tau_{\elm'}^\sd  P^{n-1}_{p,\elm'}}, \label{eq:local.fl}\\
(\algPt{tm})_\sd &  \eq  {\frac{t^n - t}{\tau^n}} \big(\FcMsNL(\solcM[n,k,i]) - \FcMsNL(\solcM[n-1])\big),\label{eq:local.comp.fl}\\
(\algP{{upw}})_\sd  & \eq \FcMs(\solcM[n,k,i]) - \FcMsNL(\solcM[n,k,i]), \label{eq:local.upw}\\
(\algP{lin})_\sd & \eq \FcMs^{n,k,i}  - \FcMs(\solcM[n,k,i]),  \label{eq:local.lin}\\
(\algP{alg})_\sd & \eq \FcMs^{n,k,i+j}  - \FcMs^{n,k,i} \label{eq:local.alg}.
\end{align}
  \end{subequations}
All these \cor{face} normal fluxes are set to zero on faces $\sd$ located on the boundary
of $\Om$, in accordance with the \cor{no-flow} boundary condition~\eqref{eq:bc}, with the exception of $(\alg{U}_{\elm,p}^{t,n,k,i})_\sd$ that is set to $\<\perm \rho_p g\GRAD z
\SCAL \normall[\Omega],1\>_\sd$\cor{, in accordance with~\eqref{eq:darcy.vel}--\eqref{eq:darcy.vel_pot}}. Note that~\eqref{eq:local.fl} is an affine\cor{-in-time} interpolation of the previous time step $P^{n-1}_{p,\elm}$ and the current $P^{n,k,i}_{p,\elm}$ values.

We also need to define the phase pressure point values $\Sel_{p,\elm}^{t,n,k,i}$ and $\Sel^{{\rm ext},{t,n,k,i}}_{p,\elm}$. We again proceed as in Section~\ref{sec_fl_faces_pot_NL}, using directly the cells pressure values in Definition~\ref{def_pr_values} or~\ref{def_hybr_pr}. More
precisely, the starting values are the phase pressure values obtained from
the finite volume scheme, the current $P^{n,k,i}_{p,\elm}$ and the previous
time step $P^{n-1}_{p,\elm}$, that we combine for each time $t \in I^n$ to get an affine\cor{-in-time} interpolation
\be \label{eq_p_pres}
    \frac{t - t^{n-1}}{\tau^n} P^{n,k,i}_{p,\elm} + \frac{t^n - t}{\tau^n}
    P^{n-1}_{p,\elm}.
\ee
Note that all quantities of~\eqref{eq:local} and~\eqref{eq_p_pres} are {\em
readily available} from the {\em finite volume discretization} of
Section~\ref{sec:disc}, more precisely on each linearization step $k\geq1$
and each algebraic solver step $i\geq1$ as described in
Sections~\ref{sec:lin}--\ref{sec:alg}.

\subsection{Fictitious flux and potential reconstructions} \label{sec:fl.rec}

We now identify the fictitious flux and potential reconstructions, following Section~\ref{sec_fl_pot_rec_NL}. 

Let a phase $p\in\Pp$, a component $c\in\C$, a time step $1 \leq n \leq N$, a mesh element $\elm \in \T_H^n$, a linearization iteration $k\ge 1$, and an algebraic solver iteration $i\geq 1$
be fixed. For the fictitious flux reconstructions, we \cor{lift the face normal fluxes from~\eqref{eq:local}, using them to define the spaces $\tV_{h, \mathrm{N}}^\elm$ in~\eqref{eq_arg_min} of Definition~\ref{def_fr}.} We define:

\bi

\item The {\em Darcy phase fluxes reconstructions} $\tu_{p,h}^{n,k,i}|_\elm\in\RT_0(\TK) \cap \Hdvi{\elm}$, as a discrete counterpart
of~\eqref{eq:darcy.vel_pot} at the discrete time $t^n$, using the face fluxes $\alg{U}_{\elm,p}^{t,n,k,i}$ of~\eqref{eq:local.fl} at $t^n$.

\item The {\em component fluxes reconstructions} $\Phich[n,k,i]|_\elm\in\RT_0(\TK) \cap \Hdvi{\elm}$, as a discrete counterpart of~\eqref{eq:comp.flux} at the discrete time $t^n$, prescribed by $\FcMsNL(\solcM[n,k,i])$ appearing in~\cor{\eqref{eq:local.comp.fl}}.

\item The {\em upwinding \cor{error} component fluxes reconstructions} $\ffd^{n,k,i}|_\elm\in\RT_0(\TK) \cap \Hdvi{\elm}$ prescribed by \cor{$\algP{{upw}}$ of}~\eqref{eq:local.upw}.

\item The {\em linearization error flux reconstruction}
$\ffl^{n,k,i}|_\elm\in\RT_0(\TK) \cap \Hdvi{\elm}$ prescribed by $\algP{lin}$ of~\eqref{eq:local.lin}.

\item The {\em algebraic error flux reconstruction}
$\ffr^{n,k,i}|_\elm\in\RT_0(\TK) \cap \Hdvi{\elm}$ prescribed by $\algP{alg}$ of~\eqref{eq:local.alg}. 

\ei

The error measure $\Nn{n,k,i}$ of~\eqref{eq:int_err} works with approximate Darcy phase fluxes $\tu_{p,h\tau}^{n,k,i}$ together with approximate component fluxes $\Phicht[n,k,i]$. 
We define them as affine-in-time on the time interval $I^n$, given respectively by $\tu_{p,h}^{n,k,i}$ and $\Phich[n,k,i]$ at the current time $t^n$, and by $\tu_{p,h}^{n-1}$ and $\Phich[n-1]$ obtained in the same way but on time previous time $t^{n-1}$\cor{, once the linearization and algebraic iterations have been stopped}. \cor{All $\tu_{p,h\tau}^{n,k,i}$, $\Phicht[n,k,i]$, and} the other fluxes $\ffd^{n,k,i}$, $\ffl^{n,k,i}$, and $\ffr^{n,k,i}$ are only used for the proof of Theorem~\ref{thm_estim_MP_MC} below \cor{and {\em need not be constructed in practice} (whence, recall, the naming fictitious)}.
It namely follows from~\eqref{eq:RcM.lin.alg}, on algebraic step $i+j$,
the above flux reconstructions, and the Green theorem that
there holds
\be \label{eq_div_comps_MP_MC}
    q_{c,\elm}^n - \frac{\Ll_{c,\elm}^{n,k,i+j} - \Ll_{c,\elm} (\solMp)}
    {\tau^n} - \cor{\big(} \DIV(\cor{\Phich[n,k,i] +} \ffd^{n,k,i}+\ffl^{n,k,i}+\ffr^{n,k,i})\cor{\big)|_\elm} =  - |\elm|^{-1} \Rl^{n,k,i+j}_{c,\elm} \quad\forall \elm\in\T_H^n,
\ee
similarly as in~\eqref{eq_div_comps}. 

The fictitious reconstruction of the phase pressures is again achieved as in Definition~\ref{def_pr}, \ie, setting 
\[
    \prht^{n,k,i}(\ver) \eq \alg{Z}_{p,\ver}^{t,n,k,i} \qquad \forall \ver \in \Vh
\]
from the potential point values identified in Section~\ref{sec_fl_faces_pot_MP_MC}. It is piecewise
affine with respect to the simplicial submeshes $\Th^n$ and $\Ho$-conforming and piecewise affine in time. In the error measure $\Nn{n,k,i}$ of~\eqref{eq:int_err}, we also need $\Ll^{n,k,i}_{c,h\tau}$. It is defined as affine-in-time on each time interval $I^n$, given
respectively by the values $\Ll_{c,\elm} (\solMp)$ and $\Ll_{c,\elm}
(\solMppp)$ at times $t^{n-1}$ and $t^n$, for all $\elm \in \T_H^n$. This in particular gives
\be \label{eq_lc_der}
    \pt_t \Ll^{n,k,i}_{c,h\tau}|_\elm = \frac{\Ll_{c,\elm}
(\solMppp) - \Ll_{c,\elm} (\solMp)}{\tau^n}
\ee
on each \cor{time interval} $I^n$.

\subsection{A guaranteed a posteriori error estimate distinguishing the error
components with inexpensive implementation and evaluation} \label{sec:apost}

Proceeding as for the steady linear problem in Section~\ref{sec_a_post_darcy} and
the steady nonlinear problem in Section~\ref{sec_est_NL_comps} (under
Assumption~\ref{as_polyt_disc}), we obtain a posteriori error estimates for the problem~\eqref{pressure}--\eqref{eq:fug}. \cor{We first treat the simpler case where there is no coarsening of the space meshes between the time steps and postpone the general case to Remark~\ref{rem_coars} below.}

\begin{ctheorem}{A guaranteed a posteriori error estimate distinguishing the error
components with inexpensive implementation and evaluation}{thm_estim_MP_MC}
For the multiphase compositional Darcy flow~\eqref{pressure}--\eqref{eq:fug}, \cor{let the weak solution $\sol$ satisfy~\eqref{ass:sol.reg}.} Consider any polytopal discretization of the form~\eqref{eq:discret}, \cor{an iterative} linearization~\eqref{eq:linear-sys} on step $k \geq 1$, and any
iterative algebraic solver~\eqref{eq:RcM.lin.alg} on step $i \geq 1$, for a given time step $1\le n\le N$.
Consider $j \geq 1$ additional algebraic solver steps.
Let, for each polytopal element $\elm \in \T_H^n$, the element vectors of the face normal fluxes $\alg{\Theta}_{\bullet, \elm, c}^{{n,k,i}}$, $\bullet = {\rm upw}, {\rm lin}, {\rm alg}$, together with $\algPt{tm}$ and $\alg{U}_{\elm,p}^{t,n,k,i}$ be given by~\eqref{eq:local}.
Let the element vectors of potential point values $\Sel_{p,\elm}^{t,n,k,i}$ and $\Sel^{{\rm ext},{t,n,k,i}}_{p,\elm}$ be constructed from the time-interpolated values~\eqref{eq_p_pres} of $P^{n,k,i}_{p,\elm}$ and $P^{n-1}_{p,\elm}$ following Section~\ref{sec_fl_faces_pot_MP_MC}. \cor{Consider the approximate solution $\tu_{p,h\tau}^{n,k,i}, \Phicht[n,k,i]$ as defined in Section~\ref{sec:fl.rec} and the intrinsic error measure $\Nn{n,k,i}$ of~\eqref{eq:int_err}.}
Let \cor{finally} the element matrices
$\widehat{\matr{A}}_{\mathrm{MFE},\elm}$,
$\widehat{\matr{S}}_{\mathrm{FE},\elm}$, and
$\widehat{\matr{M}}_{\mathrm{FE},\elm}$ be respectively defined
by~\eqref{eq_loc_matr_MFE}, \eqref{eq_loc_matr_FE},
and~\eqref{eq_loc_matr_FE_mass}. Let $\T_H^n \subset \T_H^{n-1}$, \ie,
there is no coarsening. Then

\begin{equation}\label{eq:local.comps:a}
    \Nnn[{n,k,i}] \leq \Bigg\{\sum_{c \in \C} \big(\est{sp}{c}{n,k,i} + \est{tm}{c}{n,k,i}
    + \est{lin}{c}{n,k,i} + \est{alg}{c}{n,k,i} + \est{rem}{c}{n,k,i}\big)^2 \Bigg\}^{\frac12}
\end{equation}
with
\begin{equation}\label{eq:local.comps:b}
   \est{\bullet}{c}{n,k,i} \eq  \Bigg\{ {\delta_{\bullet}} \int_{I^n} \sum_{\elm \in \T_H^n}
    \big(\est{\bullet}{\elm,c}{n,k,i} \big)^2 \dt\Bigg\}^{\frac12},
\end{equation}
where $\bullet = {\rm sp}, {\rm tm}, {\rm lin}, {\rm alg}, {\rm rem}$
and $\delta_{\bullet} \eq 2$, except for $\delta_{\rm sp} \eq 4$. Here, for
$c \in \C$, we prescribe the elementwise {\em spatial estimators} \bse
\label{eq:simple.all.est} \be \label{eq:simple.sp.est}
    \est{sp}{\elm,c}{n,k,i} \eq \est{{upw}}{\elm,c}{n,k,i}+ \Bigg \{ \sum_{p\in\Pp_c} \left(  \est{NC}{\elm,c,p}{t,n,k,i}\right)^2\Bigg\}^{\frac12},
\ee
with
\be \label{eq:simple.NC.est}\bs
   \left( \est{NC}{\elm,c,p}{t,n,k,i}\right)^2  \eq {} & \left(\nu_{p,\elm}^{n,k,i}  C^{n,k,i}_{p,c,\elm}
   \right)^2 \biggl[
    \big(\alg{U}_{\elm,p}^{t,n,k,i}\big)^{\mathrm{t}} \widehat{\matr{A}}_{\mathrm{MFE},\elm} \alg{U}_{\elm,p}^{t,n,k,i} +
     \big(\Sel_{p,\elm}^{t,n,k,i}\big)^{\mathrm{t}} \widehat{\matr{S}}_{\mathrm{FE},\elm} \Sel_{p,\elm}^{t,n,k,i} \\
    {} &   +2 \big(\alg{U}_{\elm,p}^{t,n,k,i}\big)^{\mathrm{t}} \Sel^{{\rm ext},{t,n,k,i}}_{p,\elm}
    - 2 {\sum_{\sd \in \FK} (\alg{U}_{\elm,p}^{t,n,k,i}
    )_\sd |\elm|^{-1}}  {\alg{1}}^{\mathrm{t}} \widehat{\matr{M}}_{\mathrm{FE},\elm} \Sel_{p,\elm}^{t,n,k,i}
    \biggr],
\es \ee
and the {\em upwinding estimators}
\be \left(\est{{upw}}{\elm,c}{n,k,i}\right)^2 \eq (\algP{{upw}})^{\mathrm{t}}
\widehat{\matr{A}}_{\mathrm{MFE},\elm} (\algP{{upw}}) \ee
together with the {\em temporal estimators}
\be \label{eq:simple.temp.est}
    \left(\est{tm}{\elm,c}{n,k,i}\right)^2 \eq
    (\algPt{tm})^{\mathrm{t}} \widehat{\matr{A}}_{\mathrm{MFE},\elm} \algPt{tm} ,
\ee
the {\em linearization
estimators} \be \label{eq:simple.lin.est}
    \left(\est{lin}{\elm,c}{n,k,i}\right)^2 \eq
    (\algP{lin})^{\mathrm{t}} \widehat{\matr{A}}_{\mathrm{MFE},\elm}  \algP{lin}
    + h_\elm (\tau^n)^{-1} \norm{\Ll_{c,\elm}^{n,k,i} - \Ll_{c,\elm} (\solMppp)}_{\elm},
\ee the {\em algebraic estimators}
\be \label{eq:simple.alg.est} \left(\est{alg}{\elm,c}{n,k,i}\right)^2 \eq
(\algP{alg})^{\mathrm{t}} \widehat{\matr{A}}_{\mathrm{MFE},\elm} \algP{alg}
+ h_\elm (\tau^n)^{-1} \norm{\Ll_{c,\elm}^{n,k,i+j} - \Ll_{c,\elm}^{n,k,i}}_{\elm}, \ee
and the {\em algebraic remainder estimators} by
\be \label{eq:simple.rem.est}
    \est{rem}{\elm,c}{n,k,i} \eq h_\elm |\elm|^\mft |\Rl^{n,k,i+j}_{c,\elm}|.
\ee
\ese
\end{ctheorem}

\bp Let the time step $1\le n\le N$, the iterative linearization step $k \geq 1$, and the iterative algebraic solver $i \geq 1$ be fixed. We bound the two building terms $\Nn{n,k,i}_c$ and $\Nn{n,k,i}_p$ of the error measure $\Nn{n,k,i}$ from~\eqref{eq:int_err} separately, for each component $c \in \C$ and each phase $p \in \Pp$. We rely on the fictitious flux and potential reconstructions of Section~\ref{sec:fl.rec} and namely on~\eqref{eq_div_comps_MP_MC}.

We first bound $\Nn{n,k,i}_c$. Let $\varphi \in X^n$ with $\norm{\varphi}_{X^n}=1$ be fixed. We have, using~\eqref{eq_lc_der}, \eqref{eq_div_comps_MP_MC}, and the Green theorem
\ban
{} & (q_c,\varphi) -(\pt_t \Ll^{n,k,i}_{c,h\tau}, \varphi) + \big(\Phicht[n,k,i], \GRAD \varphi \big) \\
= {} & \sum_{\elm \in \T_H^n}\left(q_c - \frac{\Ll_{c,\elm}^{n,k,i+j} - \Ll_{c,\elm} (\solMp)}
    {\tau^n} -\DIV(\cor{\Phich[n,k,i] +} \ffd^{n,k,i}+\ffl^{n,k,i}+\ffr^{n,k,i}),\varphi \right)_\elm \\
{} & + \sum_{\elm \in \T_H^n}\left(\frac{\Ll_{c,\elm}^{n,k,i+j} - \Ll_{c,\elm} (\solMp)}
    {\tau^n} - \frac{\Ll_{c,\elm} (\solMppp) - \Ll_{c,\elm} (\solMp)}{\tau^n},\varphi\right)_\elm \\
{} & + (\DIV(\cor{\Phich[n,k,i] +} \ffd^{n,k,i}+\ffl^{n,k,i}+\ffr^{n,k,i}),\varphi) + \big(\Phicht[n,k,i], \GRAD \varphi \big) \\
= {} & - \sum_{\elm \in \T_H^n}(|\elm|^{-1} \Rl^{n,k,i+j}_{c,\elm} ,\varphi )_\elm + \sum_{\elm \in \T_H^n}\left(\frac{\Ll_{c,\elm}^{n,k,i+j} - \Ll_{c,\elm} (\solMppp)}
    {\tau^n},\varphi\right)_\elm \\
{} & + \big(\Phicht[n,k,i] - \cor{\Phich[n,k,i] -} \ffd^{n,k,i} - \ffl^{n,k,i} - \ffr^{n,k,i}, \GRAD \varphi \big)
\ean
We now treat the terms herein separately. First, the triangle inequality gives
\ban
    {} & \big(\Phicht[n,k,i] - \cor{\Phich[n,k,i] -} \ffd^{n,k,i} - \ffl^{n,k,i} - \ffr^{n,k,i}, \GRAD \varphi \big) \\
\leq {} &
    \sum_{\elm \in \T_H^n}\left\{\norm{\Phicht[n,k,i] - \Phich[n,k,i]}_{\Km^{-\frac{1}{2}},\elm} + 
    \norm{\ffd^{n,k,i}}_{\Km^{-\frac{1}{2}},\elm} \right. \\
    {} & \left.+ \norm{\ffl^{n,k,i}}_{\Km^{-\frac{1}{2}},\elm} + \norm{\ffr^{n,k,i}}_{\Km^{-\frac{1}{2}},\elm}\right\} \norm{\Km^{\frac{1}{2}}\GRAD \varphi}_\elm.
\ean
Second, the Cauchy--Schwarz inequality gives
\[
    \sum_{\elm \in \T_H^n}\left(\frac{\Ll_{c,\elm}^{n,k,i+j} - \Ll_{c,\elm} (\solMppp)}
    {\tau^n},\varphi\right)_\elm \leq \sum_{\elm \in \T_H^n} h_\elm (\tau^n)^{-1} \norm{\Ll_{c,\elm}^{n,k,i+j} - \Ll_{c,\elm} (\solMppp)}_\elm h_\elm^{-1} \norm{\varphi}_\elm.
\]
And third, 
\[
    - \sum_{\elm \in \T_H^n}(|\elm|^{-1} \Rl^{n,k,i+j}_{c,\elm} ,\varphi )_\elm \leq \sum_{\elm \in \T_H^n} h_\elm |\elm|^{-1/2} |\Rl^{n,k,i+j}_{c,\elm}| h_\elm^{-1} \norm{\varphi}_\elm.
\]
Thus, regrouping these estimates and using the Cauchy--Schwarz inequality together with~\eqref{eq:X.norm}, we have 
\be\label{eq_est_Nc}\bs
    \Nn{n,k,i}_c \leq {} & \sup_{\varphi \in X^n, \norm{\varphi}_{X^n}=1} \int_{I^n} \sum_{\elm \in \T_H^n} \big(\est{tm}{\elm,c}{n,k,i} + \est{{upw}}{\elm,c}{n,k,i} + \est{lin}{\elm,c}{{n,k,i}} + \est{alg}{\elm,c}{{n,k,i}} + \est{rem}{\elm,c}{n,k,i}\big) \\
    {} & \big\{ h_\elm^{-2} \norm{\varphi}_{\elm}^2 +
    \norm{\Km^{\frac{1}{2}}\GRAD \varphi}_{\elm}^2 \big\}^{\frac12} \dt \\
    \leq {} & \sup_{\varphi \in X^n, \norm{\varphi}_{X^n}=1} \left\{\int_{I^n} \sum_{\elm \in \T_H^n} \big(\est{tm}{\elm,c}{n,k,i} + \est{{upw}}{\elm,c}{n,k,i} + \est{lin}{\elm,c}{{n,k,i}} + \est{alg}{\elm,c}{{n,k,i}} + \est{rem}{\elm,c}{n,k,i}\big)^2 \dt\right\}^{\frac12} \underbrace{\norm{\varphi}_{X^n}}_{=1},
\es\ee
where, for all components $c\in\C$ and all $\elm \in \T_H^n$, 
\begin{subequations} \label{eq:est.comps.M} \begin{align}
\label{eq:est.tm.M}    \est{tm}{\elm,c}{n,k,i} & \eq
    \norm{\Phicht[n,k,i] - \Phich[n,k,i]}_{\Km^{-\frac{1}{2}},\elm},\\
\label{eq:est.upw.M} \est{{upw}}{\elm,c}{n,k,i} & \eq
\norm{\ffd^{n,k,i}}_{\Km^{-\frac{1}{2}},\elm},\\
\label{eq:est.lin.M}
      \est{lin}{\elm,c}{{n,k,i}} & \eq
      \norm{\ffl^{n,k,i}}_{\Km^{-\frac{1}{2}},\elm}
      {+ h_\elm (\tau^n)^{-1} \norm{\Ll_{c,\elm}^{n,k,i} - \Ll_{c,\elm} (\solMppp)}_\elm},\\
\label{eq:est.alg.M}    \est{alg}{\elm,c}{{n,k,i}} & \eq
\norm{\ffr^{n,k,i}}_{\Km^{-\frac{1}{2}},\elm} + h_\elm (\tau^n)^{-1} \norm{\Ll_{c,\elm}^{n,k,i+j} - \Ll_{c,\elm}^{n,k,i}}_\elm,\\
\label{eq:est.rem.M} \est{rem}{\elm,c}{n,k,i} & \eq
h_\elm |\elm|^\mft |\Rl^{n,k,i+j}_{c,\elm}|.
\end{align}

\cor{We now bound $\Nn{n,k,i}_p$.} Let
\begin{align}
    \label{eq:est.NC.m}
    \est{NC}{\elm,c,p}{t,n,k,i} & \eq
    \nu_{p,\elm}^{n,k,i}  C^{n,k,i}_{p,c,\elm}
    \norm{\tu_{p,h\tau}^{n,k,i} + \Km \Gr
    {\prht^{n,k,i}}}_{\Km^{-\frac{1}{2}},\elm}.
\end{align}\end{subequations}
We have
\be\label{eq_est_Np}
  \Nn{n,k,i}_p \leq
  \Bigg \{ \sum_{c \in\C_p} \int_{I^n}
  \Bigg\{{\sum_{\elm \in \T_H^n}} \big(\est{NC}{\elm,c,p}{t,n,k,i}\big)^2
  \Bigg\}\dt
  \Bigg\}^{\frac12}.
\ee

We now remark that the estimators in~\eqref{eq:est.tm.M}--\eqref{eq:est.NC.m} and~\eqref{eq:simple.NC.est}--\eqref{eq:simple.rem.est} are indeed the same, proceeding as in Theorem~\ref{thm_est_Darcy} to evaluate the nonconformity estimator
given by~\eqref{eq:est.NC.m} and by applying Lemma~\ref{lem_en_norm_MFE} to
evaluate the upwinding, temporal, linearization, and algebraic estimators.
Regrouping the spatial discretization estimators
\begin{align}
    \label{eq:est.sp.M}
    \est{sp}{\elm,c}{n,k,i} & \eq \est{{upw}}{\elm,c}{n,k,i}+ \Bigg \{ \sum_{p\in\Pp_c} \left(
    \est{NC}{\elm,c,p}{t,n,k,i}\right)^2\Bigg\}^{\frac12}
\end{align}
and all the estimators per component $c \in \C$, using the Cauchy--Schwarz inequality, and recalling~\eqref{eq:int_err}, we obtain~\eqref{eq:local.comps:a}.\ep

\br[Sharper estimates]Sharper estimates, namely avoiding the factors $\delta_{\bullet}$, (but not identifying/regrouping the error component estimators) immediately follow from~\eqref{eq_est_Nc} and~\eqref{eq_est_Np}.\er

\br[Coarsening] \label{rem_coars} If a mesh element $\elm \in \T_H^n$ has
been obtained by coarsening of elements of the previous mesh $\T_H^{n-1}$
collected in the set $\T^{n-1}_{H,\elm}$, then the estimators
$\est{NC}{\elm,c,p}{t,n,k,i}$ of~\eqref{eq:simple.NC.est} and
$\est{tm}{\elm,c}{n,k,i}$ of~\eqref{eq:simple.temp.est} need to be evaluated
using the matrices $\widehat{\matr{A}}_{\mathrm{MFE},\elm}$,
$\widehat{\matr{S}}_{\mathrm{FE},\elm}$, and
$\widehat{\matr{M}}_{\mathrm{FE},\elm}$ stemming from a common refinement of
the element $\elm$ by the simplicial meshes $\Th^{n-1}|_{\T^{n-1}_{H,\elm}}$
and $\TK$ and not $\TK$ solely. In this sense, Theorem~\ref{thm_estim_MP_MC}
also holds when $\T_H^n \not \subset \T_H^{n-1}$.\er

\subsection{Numerical experiments}\label{sec_num_exp_MP_MC}

The proposed a
posteriori error estimate framework summarized in Theorem~\ref{thm_estim_MP_MC} and
Algorithm~\ref{algo2} 
has been implemented in a reservoir
prototype simulator~\cite{Ric_ARCEOR_11}, a thermal multi-purpose simulator written
in C++, which is a part of the next generation IFPEn research simulators
based on the {\em Arcane} framework \cite{Grosp_Lelan_ARCANE_09}. The
execution platform is a public computer Intel Core i7, 8 cores, 3.7Ghz with
16GB of memory.

\subsubsection{Setup} \label{sec_num_MP_MC}

We study two different test cases. The first one
is taken from~\cite{Christ_Blunt_SPE_10_01}, relying on the tenth SPE comparative
solution project model. It is an incompressible water-oil two-phase flow
problem built on a Cartesian regular geometry. This test corresponds to the
layer 85 of SPE10. We choose here first the Cartesian regular mesh so as to
compare our new approach with approaches already validated on this type of
meshes. The second case is a simulation of a black-oil model. For this test,
the problem is built on a three-dimensional corner-point geometry (distorted
grids), well-known and most often used in reservoir simulation due to the
flexibility allowing a good representation of reservoir description,
see~\cite{Ding_Lemo_corner_point_95} and the references therein. 

In our computational experiments presented below, we perform a
nonconforming-type mesh refinement and allow one hanging node per cell face
for the space mesh adaptation. \cor{The resulting mesh is a polytopal mesh in the sense of Section~\ref{sec_meshes}.} Coarsening is then obtained by
retrogressive cell agglomeration. During the coarsening and refinement phases, we use some upscaling and interpolating mechanisms
to ensure the mass balance. Geological data (rock properties like porosity and
permeability) are stored at the finest mesh level as data on a cloud
of space points. These data are used during mesh adaptation by standard
upscaling algorithms to compute the values of the rock model properties on
the newly created cells following~\cite{brevet-amc}.

\subsubsection{Two-phase Darcy flow}\label{sec_2P}

\label{sec:two-phase} 
\begin{figure}
   \centering
    \includegraphics[width=0.4\linewidth]{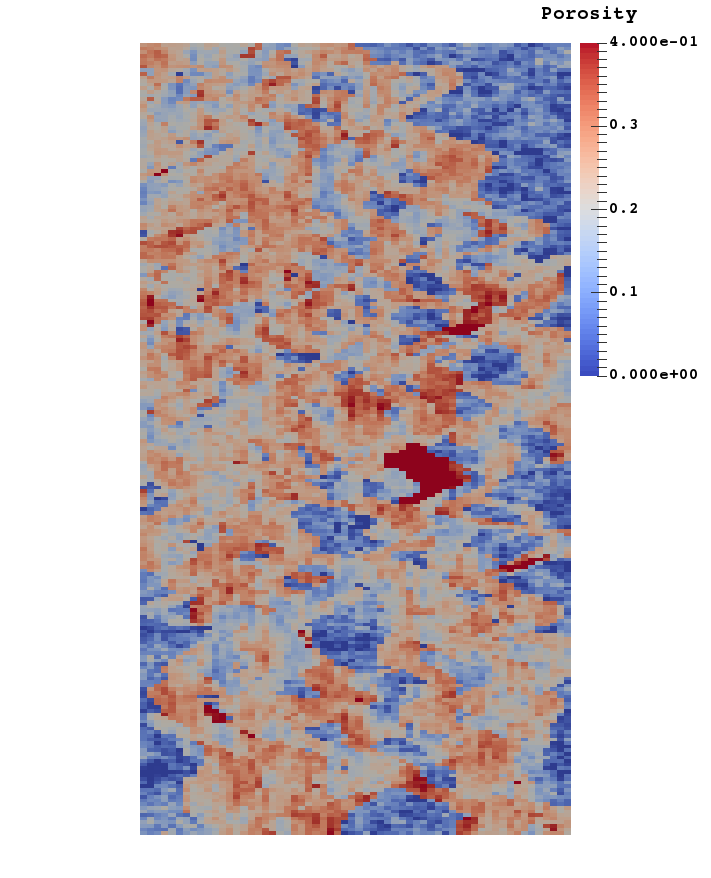}
    \hspace{1cm}
    \includegraphics[width=0.4\linewidth]{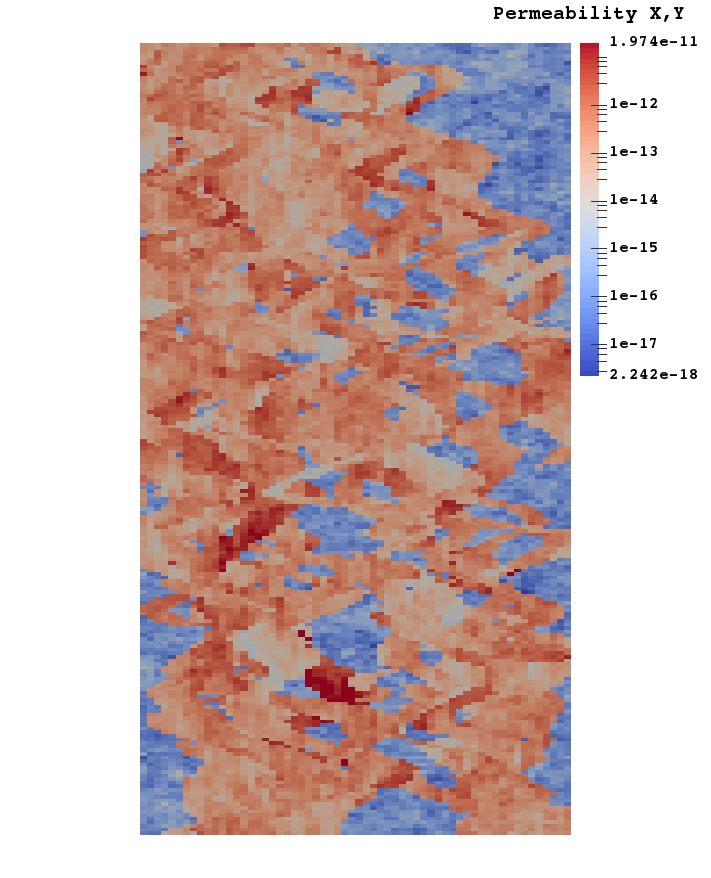}
    \caption{\cor{[Section~\ref{sec_2P}]} Porosity ({\em left}) and permeability ({\em right}), layer 85 of the 10th SPE case}
  \label{fig:spe10.perm}
\end{figure}

We consider a two-dimensional spatial domain
discretized by a grid of $60 \times 220$ rectangular cells of size $6.096$m
in the $x$ direction and $3.048$m in the $y$ direction. We choose the
initial time step as $\tau^0=4.32 \times 10^4$s, which equals to $0.5$ days,
and the process is simulated to $t_F=2000$ days. The reservoir is initially
saturated with hydrocarbons and we consider the injection of water by a well
located at the center of the grid. Four production wells are placed at the
four corners of the domain. Therefore, we have a water component ${\rm W}$
and an oil component ${\rm O}$ collected in {the set of components} $\C=\{\rm
W,O\}$ and two phases $\Pp=\{\rm w,o\}$ corresponding to water and oil. {The
model is actually simplified in that the components can be identified with
the phases, so that $\Pp_{\rm W} = \{{\rm w}\}$ and $\Pp_{\rm O} = \{{\rm
o}\}$, $C_{{\rm w},{\rm W}}=C_{{\rm o},{\rm O}}=1$, $C_{{\rm o},{\rm
W}}=C_{{\rm w},{\rm O}}=0$, and the vector of unknowns $\sol$ reduces to
$(P,(S_p)_{ p\in \Pp})$. Note that then the accumulation term $\Ll_c$
becomes linear in the only unknown $S_p$, and one particular consequence is
that the second term in the definition~\eqref{eq:simple.lin.est} of
$\est{lin}{\elm,c}{n,k,i}$ vanishes}. The porosity $\phi$ and the
permeability field $\perm$ (scalar coefficient times an identity matrix) are
shown in Figure~\ref{fig:spe10.perm}. The other parameters of
Section~\ref{sec_model} are given by (see~\cite{Christ_Blunt_SPE_10_01} for
more details):

\bi

\item $\mu_{\rm {o}} = 10^{-3}$ Pa$\cdot$s and  $\mu_{\rm {w}} = 0.3 \cdot 10^{-3}$ Pa$\cdot$s,

\item $\cor{\varsigma}_{\rm o} = \cor{\varsigma}_{\rm w} =1 $ mole$\cdot$m$^{-3}$,

\item $S_{\rm wr} = S_{\rm or} = 0.2$,

\item
\[
k_{{{\rm r}, {\rm w}}}(S_{{\rm w}}) = \left(\frac{ S_{{\rm w}} - S_{\rm wr}}{1 - S_{\rm wr} - S_{\rm or}}\right)^2 \qquad
\text{and} \qquad
k_{{{\rm r}, {\rm o}}}(S_{{\rm o}}) = \left(\frac{S_{{\rm o}} - S_{\rm or}}{1 - S_{\rm wr} - S_{\rm or}}\right)^2,
\]

\item $P_{\mathrm{c}_p}(S_p) = 0$,

\item there is no gravitational force, $z=0$, so that the mass densities
$\rho_p$ need not be specified.

\ei

We consider a discretization by the implicit multi-point finite volume scheme
of Section~\ref{sec:disc} with the linearization detailed in
Section~\ref{sec:lin}, using the \cor{blending of Newton and fixed-point} of Remark~\ref{rem_blending_N_FP}. For the linear solver needed in Section~\ref{sec:alg}, we use the Bi-Conjugated
Gradient Stabilized (BiCGStab)~\cite{vd_Vorst_BiCGStab_92} with an ILU$\{0\}$
preconditioner.

\begin{figure}
   \centering
    \subfloat[Water saturation]{
    \begin{tabular}{c}
    \includegraphics[width=0.26\linewidth]{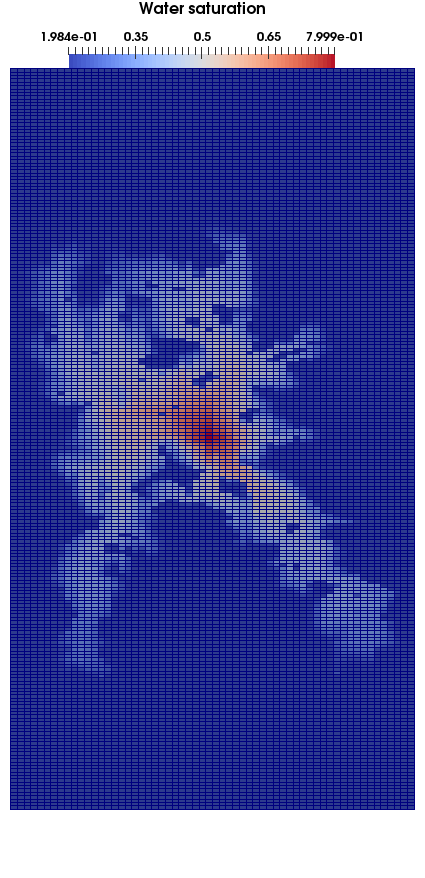}
    \end{tabular} \hspace{0.2cm}
}
    \subfloat[Cartesian estimate of~\cite{Yous_PhD_13}]{
    \hspace{0.7cm}\begin{tabular}{c}
    \includegraphics[width=0.26\linewidth]{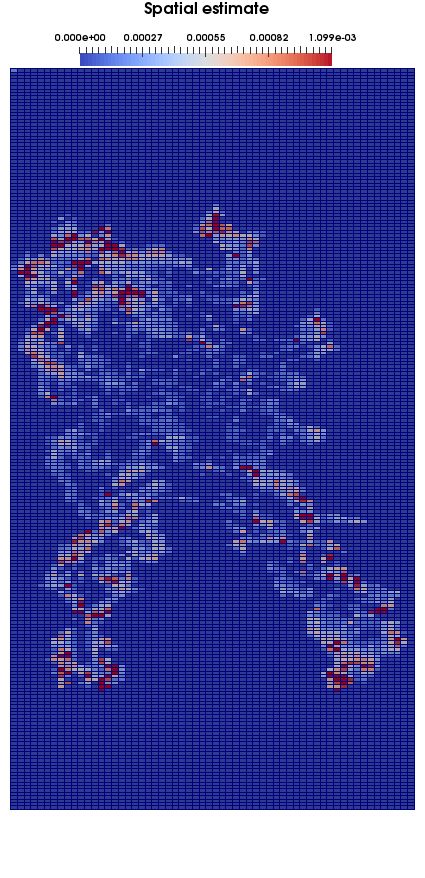}
    \end{tabular} \hspace{0.2cm}
}
    \subfloat[Polygonal simplified estimate~\eqref{eq:simple.sp.est}]{
    \hspace{0.7cm}\begin{tabular}{c}
    \includegraphics[width=0.26\linewidth]{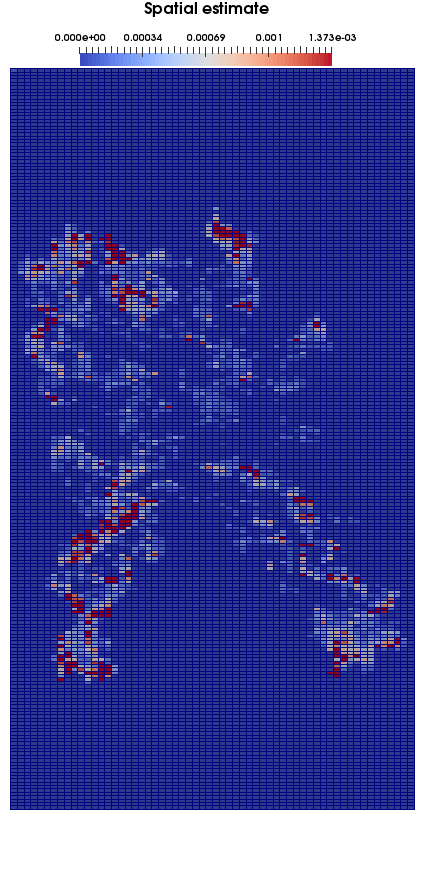}
    \end{tabular}
}

    \caption{\cor{[Section~\ref{sec_2P}]} Results on the fine mesh at $500$ days}
  \label{fig:spe10.estim}
\end{figure}

Figure~\ref{fig:spe10.estim} shows, at $500$ days of the simulation, the
evolution of the approximate water saturation, the spatial estimator computed
by the formula on Cartesian grid proposed
in~\cite[Chapter~4]{Yous_PhD_13}, and the spatial estimator given
by~\eqref{eq:simple.sp.est}. Similar behavior is observed for
the two a posteriori estimators. Furthermore, we see that they both detect
well the error following the saturation front despite the strong
heterogeneity of the domain.

\begin{figure}
   \centering
    \includegraphics[width=0.33\linewidth]{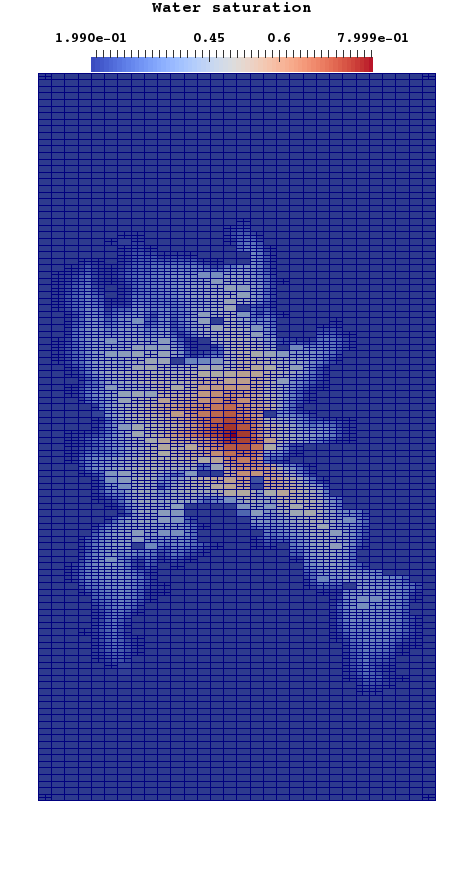}
    \hspace{1cm}
    \includegraphics[width=0.33\linewidth]{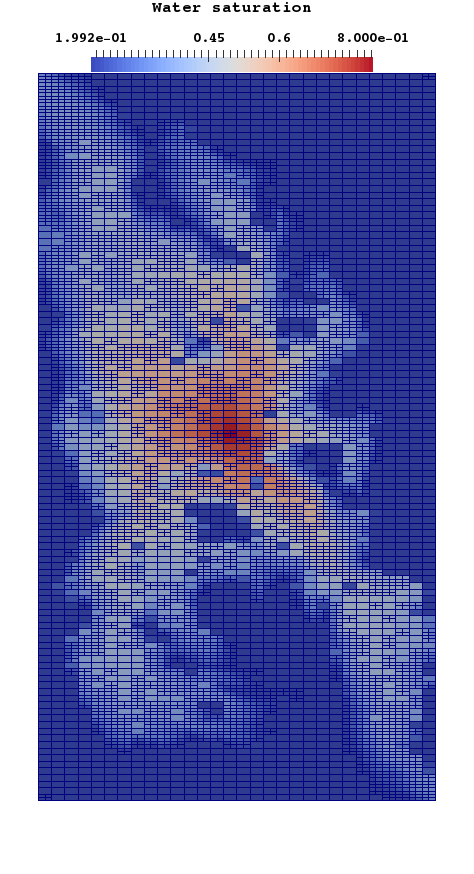}
    \caption{\cor{[Section~\ref{sec_2P}]} Results on the adaptive mesh at $400$ days and $1100$ days }
  \label{fig:spe10.front}
\end{figure}

\begin{figure}
   \centering
     \begin{tikzpicture}[scale=0.7]
      \begin{axis}[
          xlabel = {Time (seconds)},
          ylabel = {Cumulated oil rate ($m^3$)},
          legend style ={at = { (0.99,0.35)}} 
        ]
        \addplot +[line width=0.06mm] table[x=time,y=rate]{Figs/Wellprod_COIL_Fine.txt};
        \addplot +[line width=0.1mm,mark=triangle*] table[x=time,y=rate]{Figs/Wellprod_COIL_Cart.txt};
        \addplot +[mark=diamond*, mark size=1, green!50!black]table[x=time,y=rate]{Figs/Wellprod_COIL_CPG.txt};
        \addplot +[mark=square*, mark size=1., green!0!black]table[x=time,y=rate]{Figs/Wellprod_COIL_Coarse.txt};
        \legend{Fine mesh, AMR Cartesian est., AMR polygonal est., Coarse mesh};
      \end{axis}
    \end{tikzpicture}
\hspace{1.cm}
     \begin{tikzpicture}[scale=0.7]
      \begin{axis}[
          xlabel = {Time (seconds)},
          ylabel = {Water cut ($m^3$)},
          legend style ={at = { (0.65,0.95)}} 
        ]
        \addplot +[line width=0.05mm] table[x=time,y=rate]{Figs/Wellprod_CWAT_Fine.txt};
        \addplot +[line width=0.1mm,mark=triangle*] table[x=time,y=rate]{Figs/Wellprod_CWAT_Cart.txt};
        \addplot +[mark=diamond*, mark size=1, green!50!black]table[x=time,y=rate]{Figs/Wellprod_CWAT_CPG.txt};
        \addplot +[mark=square*, mark size=1., green!0!black]table[x=time,y=rate]{Figs/Wellprod_CWAT_Coarse.txt};
        \legend{Fine mesh, AMR Cartesian est., AMR polygonal est., Coarse mesh};
      \end{axis}
    \end{tikzpicture}
\caption{\cor{[Section~\ref{sec_2P}]} Fine mesh, coarse mesh, and adaptive meshes: cumulated oil rate ({\em left}) and water cut ({\em right})}
\label{fig:spe10.coil}
\end{figure}
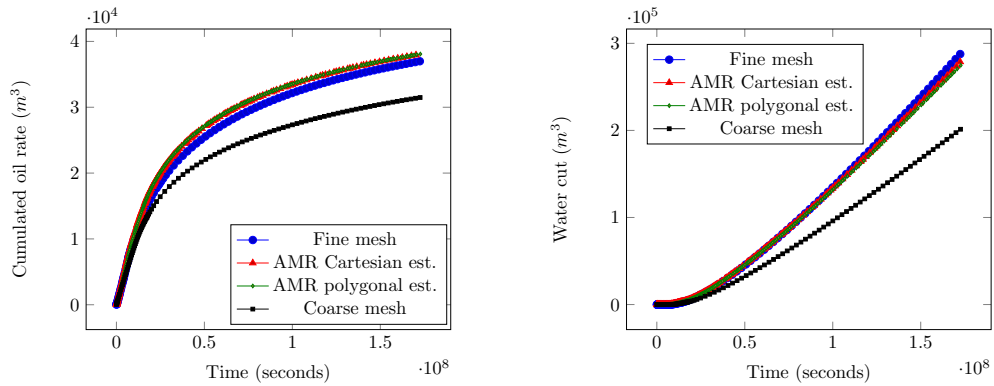

The previous result motivates an adaptive mesh refinement/coarsening (AMR)
strategy based on these estimators. We now verify that the simple a
posteriori estimators on polyhedral mesh of Theorem~\ref{thm_estim_MP_MC}
give the expected results, while comparing it with the already validated case
in~\cite[Chapter~4]{Yous_PhD_13}. We focus on the mesh space
adaptivity and we do not act on the time step. Actually, in our simulator,
the time step size is controlled by the behavior of the Newton\cor{/fixed-point} algorithm
(increased systematically in general, by multiplying by 2, and divided by 2
if the Newton\cor{/fixed-point} algorithm diverges), so the space-mesh adaptation only will
enable an easy comparison with the non-adaptive code. We first apply
Algorithm~\ref{algo2} with
$\cor{\delta}_{\mathrm{ref}}=0.7,\cor{\delta}_{\mathrm{deref}}=0.2$, and ``exact''
(\ie, negligible-error) algebraic and linearization solvers. On the
coarse scale, the domain is discretized by a grid of  $30  \times 110$ cells
and we allow one refinement level. Figure~\ref{fig:spe10.front} shows the
evolution of the approximate water saturation and of the meshes at two
different simulation times. We remark that the refinement follows the
saturation front as times evolves. Additionally, \cor{since} we have a
model with highly heterogeneous permeability, we are lead to perform a slow
derefinement process in the zone abandoned by the front of water saturation.
Figure~\ref{fig:spe10.coil} depicts the cumulated oil rate (left) and the
water-cut\footnote{The ratio of water produced compared to the volume of
total liquids produced.} (right) during the simulation. We compare there the
results on the fine grid, the results on the coarse grid, and the
results of the two AMR strategies. We remark that the accuracy of the results
on the fine grid is almost recovered by the AMR strategy, and appears much
better than the coarse-grid result.

The details of the efficiency of the AMR strategy based on the a posteriori
error estimator can finally be appreciated in Table~\ref{tab:spe10.perf}. We
compare the global CPU time of the simulation for the different strategies.
We detail the CPU time spent on the evaluation of the estimators and the mesh
adaptation. We remark a low cost of the estimators evaluation compared with
the total computation CPU time, thanks to the use of the simple and
fast-to-evaluate form of a posteriori error estimates on polygonal meshes of
Theorem~\ref{thm_estim_MP_MC}. In Table~\ref{tab:spe10.perf}, when applying
the mesh adaptation, the overall CPU time is the sum of the resolution
time, AMR time, and estimators evaluation time. We remark that applying the
AMR strategy on this two-dimensional test case leads to a gain factor in the
overall CPU time at around 2, for both the polygonal estimate proposed
in this paper and the Cartesian estimate already validated on this type of
meshes. Note that, however, that the simulation is slightly faster when
adapting using the Cartesian estimate, due to the fact that we can here
directly compute the $\RT$ basis functions on the rectangular cells.

\begin{table}
\centering
\begin{tabular}{|c|c|c|c|c|}
\hline
\rowcolor[gray]{0.7} \bf  - & \bf Resolution & \bf AMR & \bf Estimators evaluation & \bf Gain factor\\
\hline
Fine &603s&-&-&-\\
\hline
AMR Cartesian est.&229s&39s&19s&2.1\\
\hline
AMR polygonal est.&242s&46s&27s&1.9\\
\hline
\end{tabular}
\caption{\cor{[Section~\ref{sec_2P}]} Fine grid vs. adaptive mesh refinements}
\label{tab:spe10.perf}
\end{table}

\subsubsection{Three-phases, three-components Darcy flow}\label{sec_BO}

\label{sec:BO} In this section, we present a simulation of a black-oil model.
Here, we have three phases constituted by water, oil, and gas, represented by
lowercase letters w, o, g as indices, respectively. The oil phase contains
two types of components: nonvolatile oil and volatile oil, which we call here
oil component and gas component, respectively. This is due to the fact that
in this model, the hydrocarbon components are divided into light and heavy
components. The light component can dissolve into the liquid oil phase or
volatilize in the gas phase according to the pressure and temperature. The
{gas phase only contains the gas components and the} water phase only
contains the water component. The components are represented by uppercase
letters W for the water component, O for the oil component, and G for the gas
component. Therefore, we have a problem with three phases $\Pp=\{\rm w,o,g\}$
and three components $\C=\{\rm W,O,G\}$.

\begin{figure}
   \centering
    \includegraphics[width=0.45\linewidth]{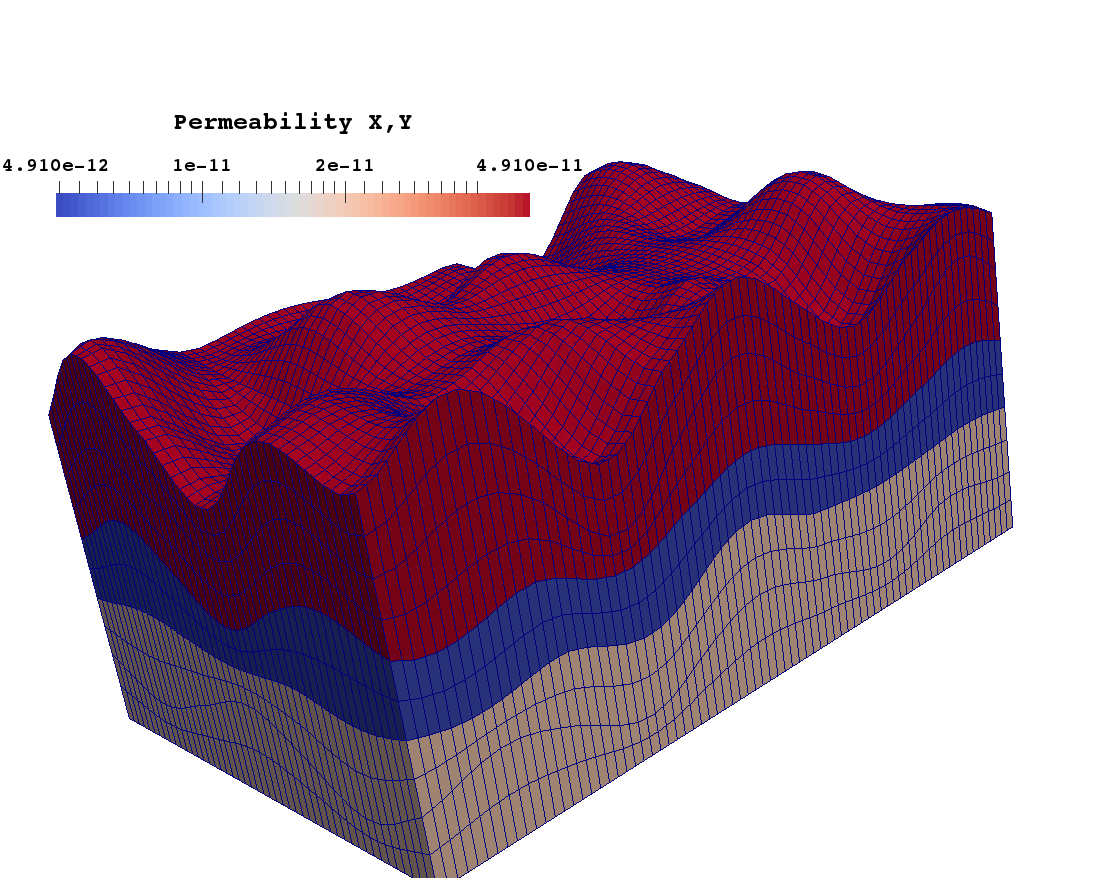}
    \hspace{0.7cm}
    \includegraphics[width=0.45\linewidth]{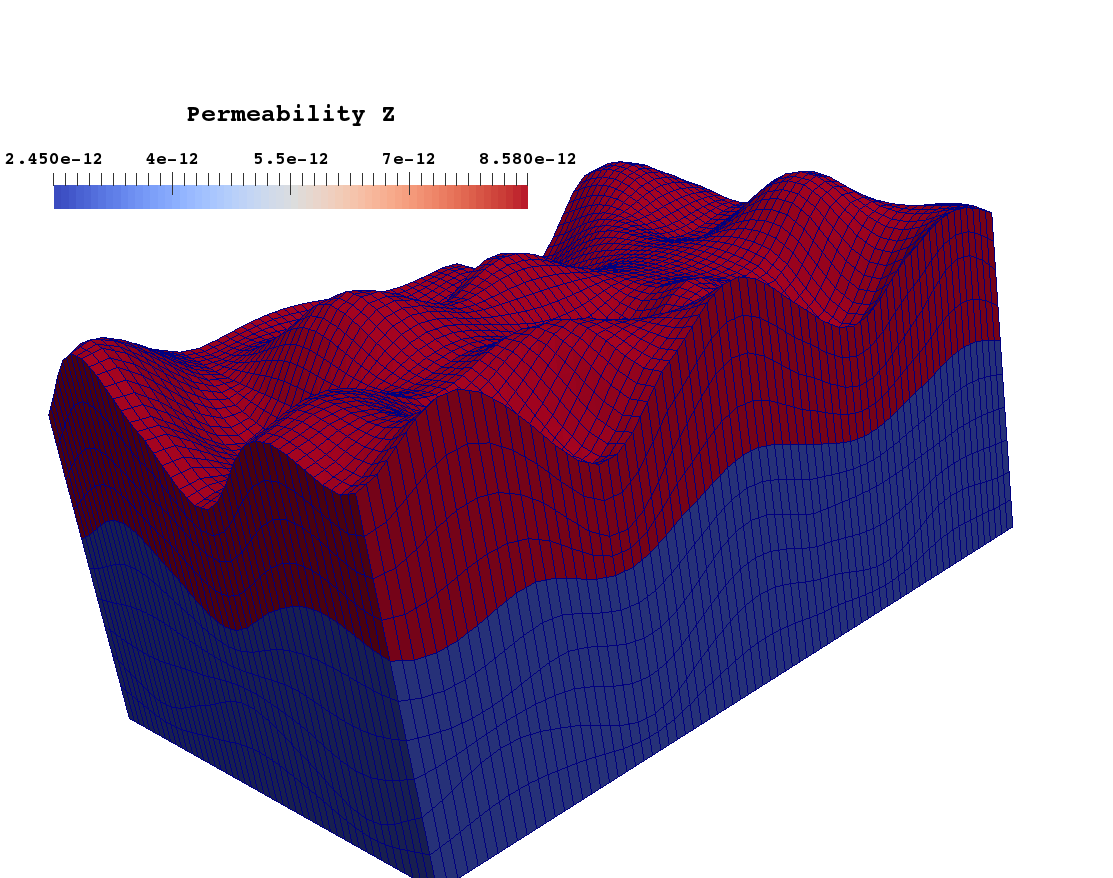}
    \caption{\cor{[Section~\ref{sec_BO}]} Permeability}
  \label{fig:BO.perm}
\end{figure}

\subsubsection*{Model setting}

\begin{figure}
   \centering
    \includegraphics[width=0.43\linewidth]{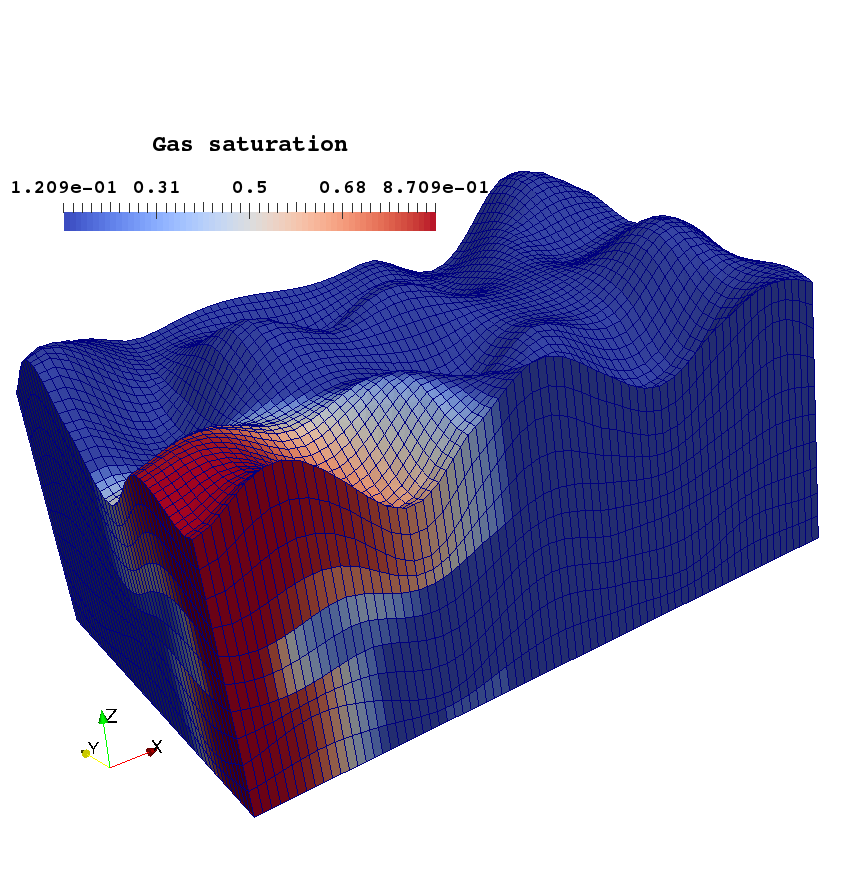}
    \hspace{0.7cm}
    \includegraphics[width=0.43\linewidth]{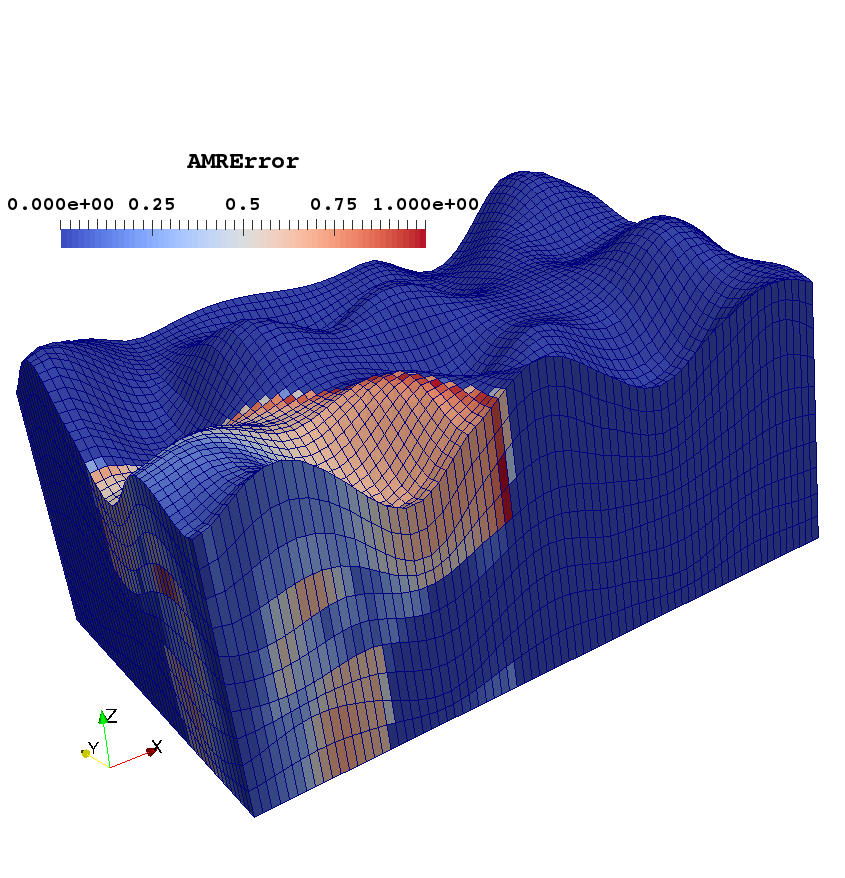}
    \caption{\cor{[Section~\ref{sec_BO}]} Results on the fine mesh at $1000$ days, gas saturation ({\em left}) and normalized polygonal estimate ({\em right})}
  \label{fig:BO.estim}
\end{figure}

The reservoir considered in this test case is a $3$-dimensional domain $\Om
\eq 4750\text{m} \times 3000\text{m} \times 114\text{m}$ discretized by a
corner-point geometry grid. We consider a heterogeneous anisotropic reservoir
with porosity $0.3$ and permeability $\perm$ in a form of a diagonal
matrix with the $x$ and $y$ components identical and forming three horizontal
layers, and the $z$ component forming two horizontal layers, see
Figure~\ref{fig:BO.perm}. We consider a gas injection in a reservoir
initially unsaturated. A vertical gas injection well perforates a corner of
the reservoir in the $z$ direction and a production well is located in
the opposite corner. On the fine scale, the domain is discretized by a grid
of $76 \times 48 \times 10$ {elements} and on the coarse scale by a grid of
$38 \times 24 \times 5 $ cells, leading to one refinement level. The process
is simulated to $t_F = 2000$ days with initial time step $\tau^0=4.32 \times
10^4$s, which equals to $0.5$ days. Data, constraints, and pressure-volume
properties are adapted from the first SPE comparative solution project model
(SPE1) designed to simulate a three-dimensional black-oil reservoir, given in
\cite[Tables 1,2, and 3]{Odeh_comp_81}. We consider a discretization by the
implicit multi-point finite volume scheme of Section~\ref{sec:disc} with the
\cor{blended Newton and fixed-point} linearization detailed in Section~\ref{sec:lin}. For the linear
solver of Section~\ref{sec:alg}, we again use the BiCGStab~\cite{vd_Vorst_BiCGStab_92} with an
ILU$\{0\}$ preconditioner. Figure~\ref{fig:BO.estim} shows the evolution of
the gas saturation and of the spatial estimator at $1000$ days. Note that, for
the spatial estimator, the data are normalized by max value in order to have
a $[0,1]$ range. We observe that the spatial estimator follows the saturation
front though the heterogeneous anisotropic medium with time evolution.

\subsubsection*{Adaptive {space mesh refinement and stopping criteria for the linear solver}}

In the standard resolution of our reservoir prototype simulator, the chosen
grid is the fine-scale one and the initial time step, as mentioned before, is
chosen as $\tau^0=4.32 \times 10^4$s. As in the test of Section
\ref{sec:two-phase}, the time step is increased systematically, by
multiplying by $2$, but also controlled by the convergence of the 
linearization loop in such a way that we divide it by $2$ if a divergence of
the \cor{linearization} algorithm occurs. We thus stick to this setting and only
focus on the stopping criteria for the linear solver (and not for the \cor{linearization} 
one) and on mesh adaptivity in space (and not in time). For the adaptive
resolution, we start on the coarse-level grid allowing to one refinement
level and we fix $\cor{\delta}_{\mathrm{ref}}=0.7,\cor{\delta}_{\mathrm{deref}}=0.2$ in
Algorithm~\ref{algo2}.

\begin{figure}
   \centering
    \includegraphics[width=0.43\linewidth]{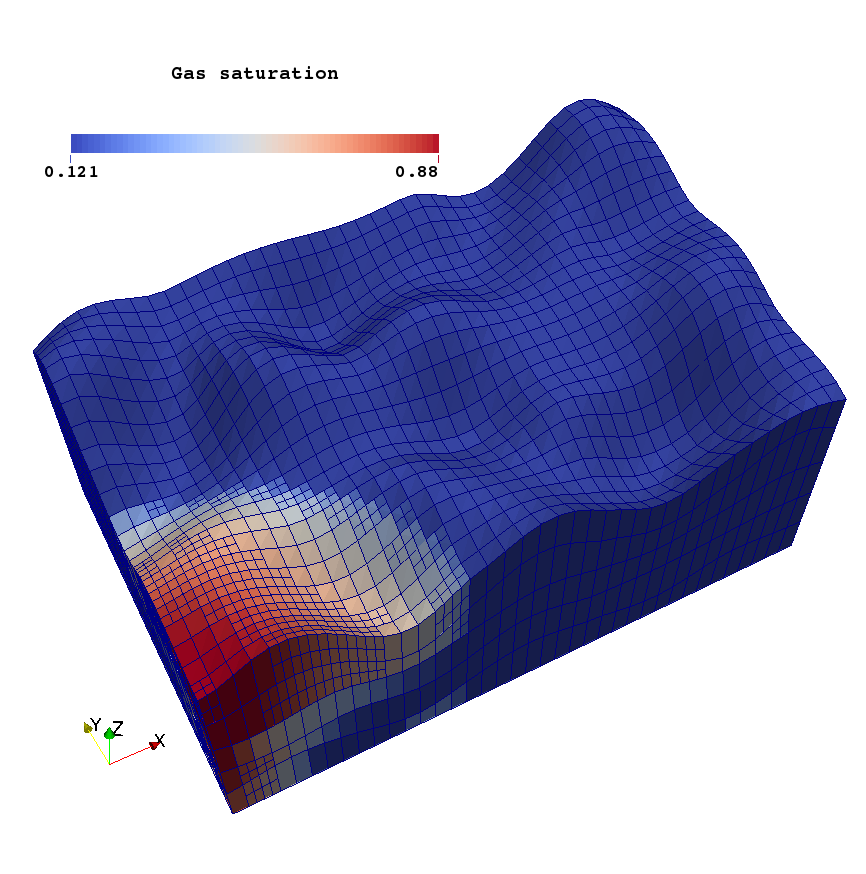}
    \hspace{0.7cm}
    \includegraphics[width=0.43\linewidth]{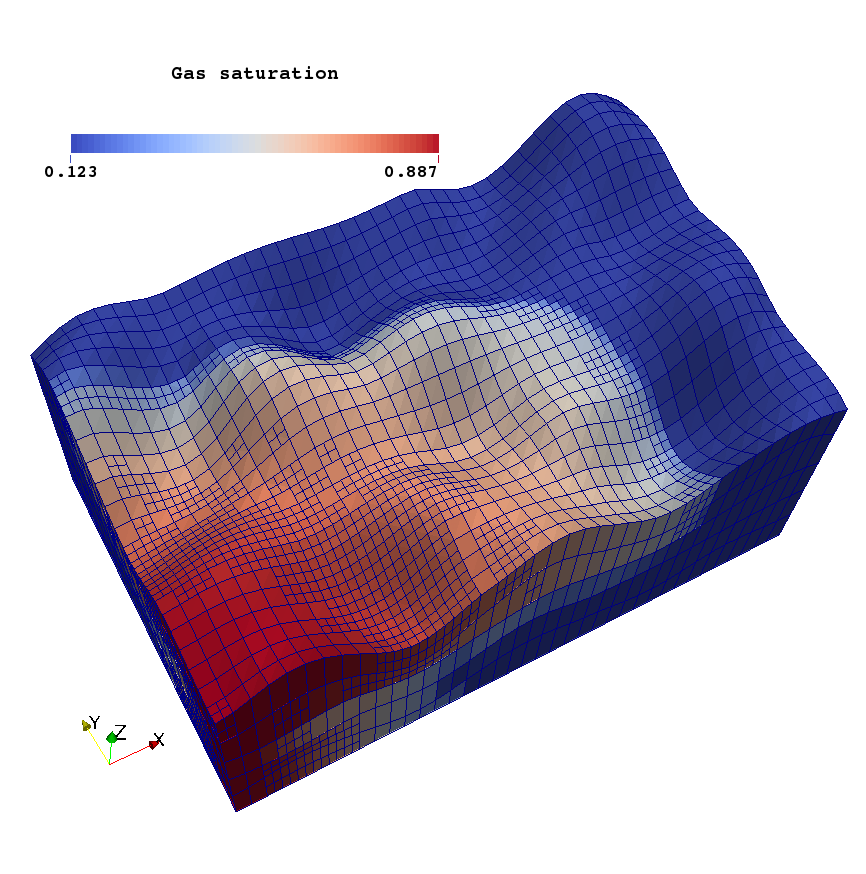}
    \caption{\cor{[Section~\ref{sec_BO}]} Gas saturation: results on the adaptive mesh at $500$ days ({\em left}) and $1500$ days ({\em right})}
  \label{fig:BO.amr}
\end{figure}

Figure~\ref{fig:BO.amr} illustrates the evolution of the approximate gas
saturation at two different time steps. We remark that the refinement follows
the saturation front as times evolves. Additionally, the fact that the light
component G (gas component) can dissolve into the liquid oil phase or
volatilize in the gas phase, we are lead to perform some localized refinement
in zones abandoned by the front of gas saturation.

\begin{figure}
   \centering
     \begin{tikzpicture}[scale=0.7]
      \begin{axis}[
          xlabel = {BiCGStab iteration},
          ylabel = {Error component \cor{estimates}},
          ymode=log,
          legend style ={at = { (0.72,0.27)}} 
        ]
        \addplot +[line width=0.01mm] table[x=iter,y=space]{Figs/estimatorsN.txt};
        \addplot +[line width=0.01mm,mark=triangle*] table[x=iter,y=alg]{Figs/estimatorsN.txt};
        \addplot +[line width=0.01mm,mark=diamond*] table[x=iter,y=residu]{Figs/estimatorsN.txt};
        \node (AD) [ red!60, text=red, draw] at (150, -9.9) {\small adaptive stopping criterion};
      \draw (AD) edge [->, shorten >=1pt, thick, red, bend right=1.5]
      (190, -6.);
      \node (UN) [ blue!60, text=blue, draw] at  (400, -3) {\small standard stopping criterion};
      \draw (UN) edge [->, shorten >=1pt, thick, blue, bend right=1.5]
      (570, -13.8);
        \legend{total estimator, algebraic estimator, relative algebraic residual};
      \end{axis}
    \end{tikzpicture}
\hspace{1.cm}
     \begin{tikzpicture}[scale=0.7]
      \begin{axis}[
          xlabel = {Time (seconds)},
          ylabel = {Number of BiCGStab iterations},
          legend style ={at = { (0.97,0.95)}} 
        ]
        \addplot +[line width=0.05mm] table[x=time,y=iter]{Figs/ref_iter_per_time.txt};
        \addplot +[line width=0.05mm,mark=triangle*] table[x=time,y=iter]{Figs/ca_iter_per_time.txt};
        \legend{Standard resolution, Adaptive resolution};
      \end{axis}
    \end{tikzpicture}
\caption{\cor{[Section~\ref{sec_BO}]} Standard resolution vs. adaptive resolution: total estimator and its algebraic component
({\em left}) and number of BiCGStab iterations per time step ({\em right})}
\label{fig:bo.iter}
\end{figure}
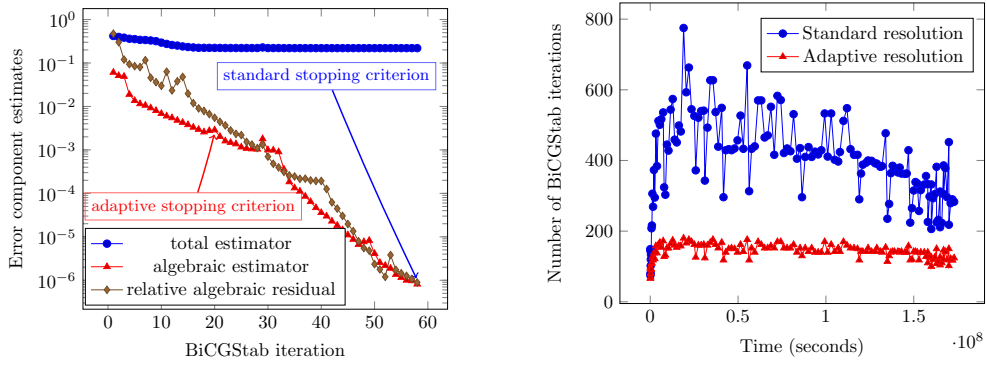

We show in the left part of Figure~\ref{fig:bo.iter}, at a fixed time of
$500$ days and for the first \cor{linearization} iteration, the evolution of the total
estimator $\left(\est{sp}{c}{n,k,i}+\est{tm}{c}{n,k,i}+\est{lin}{c}{n,k,i} +
\est{alg}{c}{n,k,i} + \est{rem}{c}{n,k,i}\right)$ with $c={\rm G}$
the gas component and the estimators given
in~\eqref{eq:local.comps:b}--\eqref{eq:simple.all.est}, the algebraic
estimator $\est{alg}{c}{n,k,i}$ of~\eqref{eq:local.comps:b},
\eqref{eq:simple.alg.est} (with $j=20$ additional steps), and the relative algebraic residual
given by
\[
\err{alg}{}{n,k,i} \eq \frac{\norm{\matr{A}^{n,k-1} \alg{X}^{n,k,i} - \alg{B}^{n,k-1}}}{\norm{\alg{B}^{n,k-1}}}
\]
with $\matr{A}^{n,k-1} \alg{X}^{n,k} = \alg{B}^{n,k-1}$ being
the linear system resulting from the $k$-th iteration of the \cor{blended Newton/fixed-point linearization} at
time step $t^n$, see Section~\ref{sec:lin}. For the standard resolution,
we stop the algebraic iteration using a fixed threshold $\err{alg}{}{n,k,i}
\leq 10^{-6}$, following our usual practice. In the adaptive
resolution based on Algorithm~\ref{algo2}, we fix $\param{alg}=10^{-2}$.

We remark, in the left part of Figure~\ref{fig:bo.iter}, that the algebraic
estimator steadily decreases, while the total estimator almost stagnates
after about a third of total number of iterations necessary to converge
using the standard stopping criterion. In the right part of
Figure~\ref{fig:bo.iter}, we depict the cumulated number of BiCGStab
iterations at each time step (the sum of the necessary number of BiCGStab
iterations at each \cor{linearization} iteration of the time step). We observe a
significant gain with the adaptive stopping criterion.

\begin{figure}
   \centering
     \begin{tikzpicture}[scale=0.7]
      \begin{axis}[
          xlabel = {Time (seconds)},
          ylabel = {Cumulated oil rate ($m^3$)},
          legend style ={at = { (0.95,0.35)}} 
        ]
        \addplot +[line width=0.01mm] table[x=time,y=rate]{Figs/BO_Well_WPROD_COIL_ref.txt};
        \addplot +[line width=0.01mm,mark=triangle*] table[x=time,y=rate]{Figs/BO_Well_WPROD_COIL.txt};
        \legend{Standard resolution, Adaptive resolution};
      \end{axis}
    \end{tikzpicture}
\hspace{1.cm}
     \begin{tikzpicture}[scale=0.7]
      \begin{axis}[
          xlabel = {Time (seconds)},
          ylabel = {Number of cells},
          legend style ={at = { (0.65,0.65)}} 
        ]
        \addplot +[line width=0.05mm] table[x=time,y=rate]{Figs/cell_nb_ref.txt};
        \addplot +[line width=0.05mm,mark=triangle*] table[x=time,y=rate]{Figs/cell_nb_bo.txt};
        \legend{Standard resolution, Adaptive resolution};
      \end{axis}
    \end{tikzpicture}
\caption{\cor{[Section~\ref{sec_BO}]} Standard resolution vs. adaptive resolution: cumulated oil rate ({\em left}) and number of cells ({\em right})}
\label{fig:bo.coil}
\end{figure}
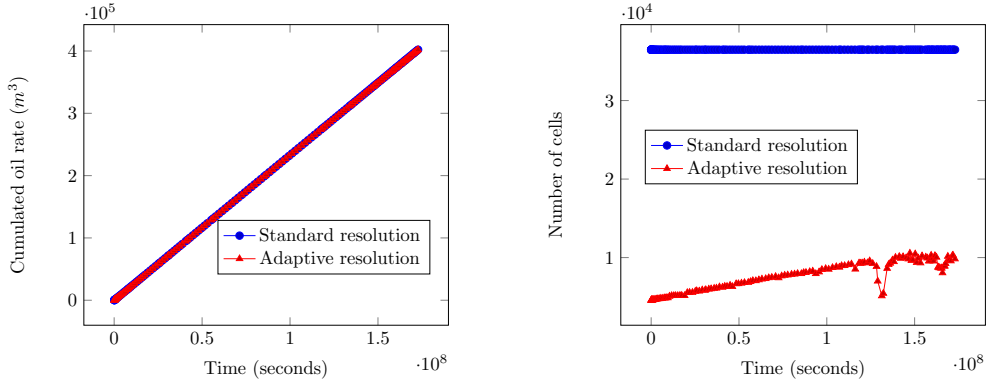

We compare in Figure~\ref{fig:bo.coil} the number of cells and the cumulated
rate of oil production resulting from both the standard and the adaptive
resolution. We can observe in the left part of Figure~\ref{fig:bo.coil} that
the adaptive algorithm does not have any significant influence on the
accuracy of production. The right part of Figure~\ref{fig:bo.coil} shows an
important reduction in the number of cells via the adaptive resolution
compared with the standard one.

The necessary linear solver steps and the different CPU times of the adaptive
and standard resolutions are collected in Table~\ref{tab:BO.perf}. It
shows that the total number of linear solver steps is reduced by $70 \%$,
which is an important gain. An important gain is also observed in the CPU
time. We in particular remark that the time spent for the adaptation (AMR and
estimators evaluation) only represents a small part of the actual resolution
time. Globally, a reduction factor of 3.8 of the overall CPU time is obtained
by comparing the adaptive resolution with the standard one. Important
increase of this reduction factor is still to be expected when also the
\cor{linearization} solver stopping criteria and time step adaptation are used.

\begin{table}
\centering
\begin{tabular}{| c | c | c | c | c | c |}
\hline
\rowcolor[gray]{0.7}     & {\bf  Linear solver  } & {\bf Resolution  }&{ \bf AMR  } & {\bf Estimators }&{ \bf Gain  }\\
\rowcolor[gray]{0.7}      & {\bf steps } & {\bf time } &  {\bf time } &{ \bf  evaluation} & {\bf factor} \\
\hline
Standard resolution &66386&1023s&-&-&-\\
\hline
Adaptive  resolution&20184&201s&42s&26s&3.8\\
\hline
\end{tabular}
\caption{\cor{[Section~\ref{sec_BO}]} Comparison between standard and adaptive resolutions}
\label{tab:BO.perf}
\end{table}

\subsection{Bibliographic resources}\label{sec_biblio_MS_MC}

Problem~\eqref{pressure}--\eqref{eq:fug} is highly involved. For this reason, the existence and uniqueness of the weak solution is supposed in Assumption~\ref{ass:sol.reg}. For simplified problems, existence and uniqueness of the weak solution together with existence, uniqueness, and a priori error estimates for numerical discretizations have been derived in Alt and Luckhaus~\cite{Alt_Luck_quas_83}, Alt~\eal\ \cite{Alt_Luck_Vis_PM_84}, Kr{\"o}ner and Luckhaus~\cite{Kro_Luck_2P_84}, Otto~\cite{Otto_L1_contr_96, Otto_L1_contr_unsat_97}, Knabner and Otto~\cite{Knab_Otto_un_00}, 
Chen~\cite{Chen_deg_2P_I_01}--\cite{Chen_deg_2P_II_02} and Chen and Ewing~\cite{Chen_Ew_deg_01}, Sma\"{\i}~\cite{Smai_2P_exist_09}, Caro~\eal\ \cite{Caro_Saad_Saad_ex_2P_14}, Cao and Pop~\cite{Cao_Pop_2P_dyn_cap_15}, Ruiz-Baier and Lunati~\cite{RuizB_Lunat_MFE_FV_MP_16}, Feireisl~\eal\ \cite{Fei_Hilh_Petz_Tak_var_dens_fl_16}, 
Kou~\eal\ \cite{Kou_Sun_Wang_dec_en_st_MP_18}, Murphy and Walkington~\cite{Murp_Walk_convexity_2P_19}, Bene\v{s}~\cite{Ben_ex_2P_therm_22}, and Brenner~\eal\ \cite{Bren_Chorf_Mass_VAG_2P_disc_cap_22, Bren_Mass_Quenj_Dron_2P_disc_cap_22}, see also the references therein. Some classical books on the subject are those of Bear and Bachmat~\cite{Bear_Bach_90}, Chavent and Jaffr{\'e}~\cite{Chav_Jaff_res_sim_86}, and Chen~\eal\ \cite{Chen_Huan_Ma_comput_MP_06}.

For model (unsteady nonlinear) advection--diffusion--reaction problems, a posteriori analyses have been initiated in Eriksson and Johnson~\cite{Er_John_adpt_SD_93}, S\"{u}li~\cite{Sul_a_post_hyp_FV_96, Sul_a_post_hyp_FE_99}, Angermann~\eal\ \cite{An_Knab_Th_err_est_FV_DD_98}, Verf{\"u}rth~\cite{Verf_RD_rob_a_post_98, Verf_CD_a_post_98, Verf_rob_a_post_CD_05}, Kr{\"o}ner and Ohlberger~\cite{Kron_Ohl_a_post_FV_CL_00}, Ohlberger~\cite{Ohl_a_post_FV_vert_CRD_01, Ohl_a_post_FV_cell_CRD_01}, Ohlberger and Rohde~\cite{Ohl_Rohde_a_post_FV_CD_02}, Kr{\"o}ner~\eal\ \cite{Kro_Kuth_Ohl_Rohde_a_post_FV_hyp_ADR_03}, Nicaise~\cite{Nic_a_post_FV_CDR_06}, and Ern~\eal\ \cite{Ern_Steph_Voh_apost_DG_10}, see also the references therein. 
A detailed attention has been paid to time-evolution problems of parabolic type in order to derive a posteriori error estimates robust with respect to the final time. This has been achieved in Verf{\"u}rth~\cite{Ver_a_post_heat_03, Verf_rob_a_post_CD_nonstat_05, Tob_Verf_a_post_uns_stab_15}, see also Bergam~\eal\ \cite{Ber_Ber_Mgh_a_post_par_04}, 
Ern and Vohral{\'{\i}}k~\cite{Ern_Voh_a_post_par_10}, Hilhorst and Vohral{\'{\i}}k~\cite{Hilh_Voh_apost_FV_FE_11}, Georgoulis~\eal\ \cite{Geor_Lak_a_post_DG_11, Georg_Lak_Wih_a_post_hp_DG_21, Georg_Makr_a_post_par_23}. Lately, locally space--time efficient a posteriori error estimates robust with respect to the final time were obtained in Ern~\eal\ \cite{Ern_Sme_Voh_heat_HO_Y_17}, see also Ern~\eal\ \cite{Ern_Sme_Voh_heat_HO_X_19} and Mitra and Vohral{\'{\i}}k~\cite{Mitra_Voh_Richards_24}.

Finally, for results on a posteriori error estimates and adaptivity for numerical discretizations of porous media flows and namely multiphase compositional flows, we refer to Saad and Zhang~\cite{Saad_Zhan_adpt_2P_97}, Chen and Ewing~\cite{Chen_Ew_deg_adapt_03}, Jenny~\eal\ \cite{Jen_Lee_Tchel_adpt_MS_FV_MP_04}, Chen and Liu~\cite{Chen_Liu_a_post_MFE_misc_08}, Klieber and Rivi{\`e}re~\cite{Klieb_Riv_adpt_2P_DG_06}, Chueh~\eal\ \cite{Chue_Sec_Bang_Dji_ML_adpt_2P_10, Chue_Dji_Bang_split_adpt_2P_13}, Vohral{\'{\i}}k and Wheeler~\cite{Voh_Whee_a_post_2P_13}, El Hassouni and Mghazli~\cite{Has_Mghaz_adapt_2P_13}, Bernardi~\eal\ \cite{Ber_El_Al_Mghaz_a_post_Rich_14}, Canc{\`e}s~\eal\ \cite{Canc_Pop_Voh_a_post_2P_14}, Di~Pietro~\eal\ \cite{Di_Pi_Voh_Yous_a_post_comp_14, Di_Pi_Voh_Yous_a_post_therm_14}, Faigle~\eal\ \cite{Fai_Hel_Aav_Flem_MP_adpt_14, Fai_Elf_Hel_Beck_Flem_Gei_MP_adpt_15}, Henning~\eal\ \cite{Hen_Ohl_Schw_adpt_het_MS_2P_15}, Buhr~\eal\ \cite{Buhr_Engw_Ohl_Rave_ArbiLoMod_17}, Ganis~\eal\ \cite{Gan_Pench_Whee_adpt_fully_coupled_19}, Beaude~\eal\ \cite{Beau_Brenn_Lop_Mass_Smai_form_MP_19}, Ahmed~\eal\ \cite{Ahm_Rad_Nord_adpat_st_crit_Biot_19, Ahm_Nord_Rad_adpat_st_crit_poromech_20}, and Li~\eal\ \cite{Li_Leun_Whee_adapt_2P_20}, Li and Wheeler~\cite{Li_Whee_dyn_coupl_22}, Ben Gharbia~\eal\ \cite{BenGha_Dab_Mar_Voh_a_post_compl_MP_20, BenGharb_Ferz_Voh_Yous_sem_sm_Newt_reg_23}, Deucher~\eal\ \cite{Deuch_Franc_Moy_Tchel_adpt_MS_FV_MP_25}, see also the references therein.

\bibliographystyle{acm_mod}
\bibliography{biblio}

\end{document}